\documentclass[11pt]{article}
\usepackage{ebezier}
\usepackage{amssymb,amsmath,mathrsfs,color} % font used for R in Real numbers
 \allowdisplaybreaks

\catcode`!=11
\let\!int\int \def\int{\displaystyle\!int}
\let\!lim\lim \def\lim{\displaystyle\!lim}
\let\!sum\sum \def\sum{\displaystyle\!sum}
\let\!sup\sup \def\sup{\displaystyle\!sup}
\let\!inf\inf \def\inf{\displaystyle\!inf}
\let\!cap\cap \def\cap{\displaystyle\!cap}
\let\!max\max \def\max{\displaystyle\!max}
\let\!min\min \def\min{\displaystyle\!min}
\let\!frac\frac \def\frac{\displaystyle\!frac}
\catcode`!=12

\newtheorem{theorem}{\bf Theorem}[section]
\newtheorem{lemma}{\bf Lemma}[section]
\newtheorem{definition}{\bf Definition}[section]
\newtheorem{proposition}{\bf Proposition}[section]
\newtheorem{remark}{\bf Remark}[section]

\def\Proof{\it{Proof.}\rm\quad}

\def\pd#1#2{\frac{\partial#1}{\partial#2}}

\catcode`@=11
\@addtoreset{equation}{section}
\catcode`@=12

\begin{document}

\title{Smooth Transonic Flows in De Laval Nozzles}
\author{Chunpeng Wang\footnote{The research is supported by a grant from the Chinese National Sciences Foundation.}
\\
\small School of Mathematics, Jilin University, Changchun
130012, P. R. China
\\
\small The Institute of Mathematical Sciences, The Chinese
University of Hong Kong, Shatin, NT, Hong Kong
\\
\small (email: wangcp@jlu.edu.cn)
\\[2mm]
Zhouping Xin\footnote{The research is supported by Zheng Ge Ru Foundation, Hong Kong RGC Earmarked Research 
Grants CUHK-4011/11P and CUHK-4042/08P, a grant from the Croucher Foundation, and a Focus Area Grant from the 
Chinese University of Hong kong.}
\\
\small The Institute of Mathematical Sciences and Department of
Mathematics,
\\
\small The Chinese University of Hong Kong, Shatin, NT, Hong Kong
\\
\small (email: zpxin@ims.cuhk.edu.hk)
}
\date{}
\maketitle

\begin{abstract}
This paper concerns smooth transonic flows of Meyer type in finite de Laval nozzles,
which are governed by an equation of mixed type with degeneracy and singularity at the sonic state.
First we study the properties of sonic curves.
For any $C^2$ transonic flow of Meyer type,
the set of exceptional points is shown to be a closed line segment (may be empty or only one point).
Furthermore, it is proved that a flow with nonexceptional points is unstable
for a $C^1$ small perturbation in the shape of the nozzle.
Then we seek smooth transonic flows of Meyer type which satisfy physical boundary conditions
and whose sonic points are exceptional.
For such a flow, its sonic curve must be located at the throat of the nozzle
and it is strongly singular in the sense that
the sonic curve is a characteristic degenerate boundary in the subsonic-sonic region,
while in the sonic-supersonic region all characteristics from sonic points coincide,
which are the sonic curve and never approach the supersonic region.
It is proved that there exists uniquely such a smooth transonic flow near the throat of the nozzle,
whose acceleration is Lipschitz continuous,
if the wall of the nozzle is sufficiently flat.
\\[6pt]
\\
{\sl 2000 MR Subject Classification:} 76H05 35M12 76N10
\end{abstract}

\newpage

\tableofcontents

\newpage

\section{Introduction}

In this paper we study the compressible Euler
system of steady isentropic and
irrotational smooth transonic flows in a class of two-dimensional de Laval nozzles. Such problems
arise naturally in physical experiments and engineering
designs, and there exist large literatures on experiments and numerical simulations.
However, a rigorous mathematical theory remains to be completed.

A flow is called transonic if it contains both subsonic and supersonic regions,
which are usually connected by shocks and sonic curves.
Roughly speaking, transonic flows are governed by partial differential equations
of mixed type; furthermore, shocks and sonic curves by which
subsonic and supersonic regions are connected are free interfaces in general.
Two kinds of transonic flow patterns have received much attention.
One is transonic flows past a profile.
Bers \cite{Bers1} showed that for two
dimensional flows past an arbitrarily given profile,
the whole flow field will be subsonic outside the profile if the freestream Mach number $M_0$
within a certain interval $0<M_0<\hat M<1$;
furthermore, the maximum flow speed will tend to
the sound speed as $M_0\to\hat M$. A natural question arises whether
there also exists a smooth flow past the profile for $\hat M<M_0<1$.
By constructing explicit solutions in the hodograph plane, one can get some
smooth transonic flows for some profiles and for some values of $M_0$ (\cite{Bers}).
However, smooth transonic flows past a profile do not exist in general
and are unstable even they exist according to Morawetz \cite{MorzI,MorzII,MorzIII}.
There are few studies on transonic flows with shocks past a profile
and almost no rigorous results are known.
The other kind is transonic flows in a nozzle.
Transonic flows in nozzles were investigated first by Taylor \cite{Taylor}
and Meyer \cite{Meyer},
where some special solutions were shown by using power series expansions.
One usually considers, as we shall do in this paper,
a de Laval nozzle which is symmetric and whose cross section
decreases first and then increases. There are two types of smooth transonic flows in such a
nozzle, named Taylor and Meyer types.
In a transonic flow of Taylor type, there are two supersonic enclosures of the type encountered
in a transonic flow past a profile, and smooth transonic flows of such type do not exist in general
and are unstable with respect to small changes in the shape of the nozzle even they exist (\cite{Bers}).
In a transonic flow of Meyer type, the sonic curve extends from one wall of the nozzle to the other,
which is believed to be located near the throat of the nozzle (where the cross section is smallest).
Stable transonic flows of Meyer type should be subsonic upstream and supersonic downstream,
while a flow with the reversing direction is unstable with respect to
small changes in the shape of the nozzle.
Most known smooth transonic flows in nozzles are solved as solutions to the governing equations
and no boundary conditions are discussed (\cite{Bers}). So, the nozzles cannot be given in advance.
It is noticed that Kuz'min \cite{Kuz} solved the
perturbation problems of accelerating smooth transonic flows with some structural assumptions
in a class of nozzles
by using the principle of contracting mappings.
However, the physical meaning of the boundary conditions for these smooth transonic flows is not clear.
There also exists a stable discontinuous transonic flow pattern in a nozzle called transonic shocks,
which are supersonic upstream and turn to subsonic across shocks.
We refer to \cite{CKL,CFM,CFM1,CFM3,LXY,LXY1,LXY2,XY} for the existence and stability
of such transonic flows in flat and curved nozzles.
Particularly, \cite{LXY,LXY1,LXY2} solved the transonic shock pattern
described by Courant and Friedrichs:
Given the appropriately large receiver pressure, if the
upstream flow is still supersonic behind the throat of the nozzle, then at a certain place
in the divergent part of the nozzle a shock front intervenes and the gas is compressed
and slowed down to subsonic speed (\cite{CF}).
For these transonic shocks, the shocks are free and the flows are away from the sonic state.

The main goal of the present paper is to understand the behavior of transonic flows near the sonic state.
Thus we investigate smooth transonic flows in a class of two-dimensional solid de Laval nozzles.
As usual, it is assumed that the throat of the nozzle lies in the $y$-axis
and the nozzle is symmetric with respect to the $x$-axis.
For convenience, we consider only the upper part of the nozzle due to its symmetry
and the upper wall is given by
$$
\Gamma_{\text{upw}}:y=f(x),\quad l_-\le x\le l_+,
$$
where $l_-<0<l_+$ and $f\in C^1([l_-,l_+])$ satisfies
$$
f'(x)\left\{
\begin{aligned}
&<0,\quad &&l_-\le x<0,
\\
&=0,\quad &&x=0,
\\
&>0,\quad &&0<x\le l_+.
\end{aligned}
\right.
$$
For the given inlet
$$
\Gamma_{\text{\rm in}}:x=g(y),\quad 0\le y\le f(l_-),
$$
we seek a smooth transonic flow of Meyer type in the nozzle
whose velocity vector is along the normal direction at the inlet.
Since the flow is supersonic downstream,
we choose the outlet
$$
\Gamma_{\text{\rm out}}:x=t(y),\quad 0\le y\le f(l_+)
$$
as a free boundary where the velocity vector is along
the normal of the outlet.
Let $\Omega$ be the domain bounded by $\Gamma_{\text{upw}}$, the $x$-axis, $\Gamma_{\text{\rm in}}$
and $\Gamma_{\text{\rm out}}$. The transonic flow in $\Omega$ satisfies
the steady isentropic compressible Euler system:
\begin{align}
\label{euler-1}
&\pd{}{x}(\rho u)+\pd{}{y}(\rho v)=0,\quad&&(x,y)\in\Omega,
\\
\label{euler-2}
&\pd{}{x}(P+\rho u^2)+\pd{}{y}(\rho uv)=0,\quad&&(x,y)\in\Omega,
\\
\label{euler-3}
&\pd{}{x}(\rho uv)+\pd{}{y}(P+\rho v^2)=0,\quad&&(x,y)\in\Omega,
\end{align}
where $(u,v)$, $P$ and $\rho$ represent the velocity, pressure and
density of the flow, respectively. The flow is assumed to be
isentropic so that $P=P(\rho)$ is a smooth function. In
particular, for a polytropic gas with  the adiabatic exponent
$\gamma>1$,
\begin{align}
\label{euler-4}
P(\rho)=\frac1\gamma\rho^\gamma
\end{align}
is the normalized pressure. Assume further that the flow is
irrotational, i.e.
\begin{align}
\label{euler-5}
\pd{u}{y}=\pd{v}{x},\quad(x,y)\in\Omega.
\end{align}
Then the density $\rho$ is expressed in terms of the speed
$q$ according to the Bernoulli law (\cite{CF})
\begin{align}
\label{Bernoulli}
\rho(q^2)=\Big(1-\frac{\gamma-1}{2}q^2\Big)^{1/(\gamma-1)},\quad
q=\sqrt{u^2+v^2},\quad0<q<q_{max}=\Big(\frac2{\gamma-1}\Big)^{1/2}.
\end{align}
Summing up, the flow is governed by the system
\eqref{euler-1}--\eqref{Bernoulli}.
It is well known that the system
\eqref{euler-1}--\eqref{Bernoulli} can be transformed into the full potential equation
\begin{align}
\label{potential}
\mbox{div}(\rho(|\nabla\varphi|^2)\nabla\varphi)=0,\quad(x,y)\in\Omega,
\end{align}
where $\varphi$ is the velocity potential.

\vskip10mm
\hskip8cm
\setlength{\unitlength}{0.6mm}
\begin{picture}(250,65)
\put(-70,0){\vector(1,0){140}}
\put(0,40){\vector(0,1){17}}
%\put(0,-5){$O$}
\put(67,-4){$x$} \put(-4,53){$y$}
\put(0,10){$\Omega$}
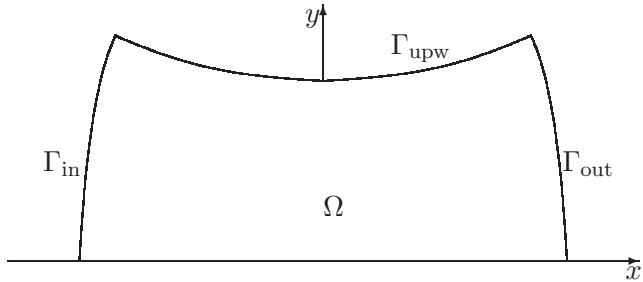
\put(-36,-12){Figure: the de Laval nozzle}

\put(-60,48){\cbezier(14,2)(30,-5)(40,-7)(60,-8)}
\put(60,48){\cbezier(-14,2)(-30,-5)(-40,-7)(-60,-8)}

%\put(-60,48){\qbezier(14,2)(30,-3)(40,-5)}
%\put(-60,48){\qbezier(40,-5)(50,-7)(60,-8)}
\put(-40,48){\qbezier(-6,2)(-12,-10)(-14,-48)}

%\put(60,48){\qbezier(-14,2)(-30,-3)(-40,-5)}
%\put(60,48){\qbezier(-40,-5)(-50,-7)(-60,-8)}
\put(40,48){\qbezier(6,2)(12,-10)(14,-48)}

\put(15,46){$\Gamma_{\text{upw}}$}
\put(-62,20){$\Gamma_{\text{in}}$}
\put(53,20){$\Gamma_{\text{out}}$}
\end{picture}
\vskip15mm

For a smooth transonic flow of Meyer type,
the governing equation \eqref{potential}
is elliptic in the subsonic region and hyperbolic in the supersonic region.
It should be noted that such a flow is singular near the sonic curve in the sense
that \eqref{potential} is degenerate elliptic in the subsonic-sonic region
while non-strictly hyperbolic in the sonic-supersonic region.
These factors cause the essential difficulties for the mathematical analysis.
Another difficulty is that the sonic curve is usually believed to be free.
So, understanding the behavior of sonic curves is not only
important for physical applications and engineering designs
but also one of the keys to study smooth transonic flows.
Bers \cite{Bers2,Bers} studied the continuation of a flow across a sonic curve, where
a subsonic-sonic flow was assumed to be given ahead.
He found out that the existence and the structure of exceptional points
play an important role in this continuation.
The definition of exceptional points is referred to \cite{Bers} (also $\S$ 2.2).
For a $C^2$ transonic flow, a sonic point is exceptional
is equivalent to that the velocity vector is orthogonal to the sonic curve at this point.
That is to say, exceptional points are characteristic degenerate in the subsonic-sonic region.
If there is not any exceptional point, then, the flow can be continued in a unique way across
the sonic curve as a supersonic flow without discontinuity.
Whereas, if there is a unique exceptional point,
the flow will be uniquely determined only in a strict subset of the same domain; furthermore,
if the flow can be continued at all into this excluded region, this continuation will not be unique.
Bers also mentioned in \cite{Bers} that it would
be interesting to know whether exceptional points are always isolated.
Here, the supersonic flow is solved as a Cauchy problem from the sonic curve
and there is no prescription of boundary conditions on the walls.
For a smooth transonic flow of Meyer type in a nozzle,
the supersonic flow satisfies not a Cauchy problem but an initial-boundary value problem.
By some precise analysis on sonic curves,
it is shown in the paper that for any $C^2$ transonic flow,
the set of exceptional points is a closed line segment
(may be empty or only one point).
Furthermore, exceptional points are strongly singular in the sonic-supersonic region
in the sense that there are two different characteristics from each nonexceptional point in the nozzle,
while all characteristics from interior exceptional points coincide
and they never approach the supersonic region locally.
Then, it is proved that a flow with nonexceptional points is unstable
for a $C^1$ small perturbation on the wall (even if the wall is still smooth).
This instability is weak since it is unknown whether the flow is unstable
if the wall is perturbed in $C^2$ or other smooth spaces.
So we seek a smooth transonic flow of Meyer type whose sonic points are exceptional.
For such a flow, its sonic curve must be located at the throat of the nozzle
and the potential on its sonic curve equals identically to a constant.
Another motivation arises from our early paper \cite{WX2}, where
the structural stability problem of a symmetric continuous subsonic-sonic flow
in a convergent nozzle with straight solid walls was proved and
the sonic curve of the continuous subsonic-sonic flow
is shown to be a free boundary where the potential equals identically
to a constant. In \cite{WX2} the continuous subsonic-sonic flow
is singular in the sense that while the speed is
continuous yet the acceleration blows up at the sonic state.
This singularity arises from the geometry of the nozzle
and it seems to vanish if the sonic curve is located at the throat of a de Laval nozzle.
It is noted that for a smooth transonic flow of Meyer type
whose sonic curve is located at the throat of the nozzle,
the potential at the sonic curve, which equals identically a constant, is free
although the location of the sonic curve in the nozzle is known.

For smooth transonic flows, the governing equation is elliptic in the subsonic region
and hyperbolic in the supersonic region;
furthermore, it is degenerate and singular at the sonic state.
So, to study a smooth transonic flow problem, one must determine and control precisely
the speed of the flow near the sonic state,
which plays an essential role in the mathematical analysis.
In \eqref{potential}, the speed of the flow is the absolute value of the gradient
of a solution. So, generally speaking,
it is very hard to estimate precisely the speed of the flow near the sonic state
in the physical plane even if the location of the sonic curve is known.
It is more convenient to study smooth transonic flows in the potential plane,
where the speed of the flow is a solution to the governing equation.
After confirming the sonic curve to be a free interface where the potential equals identically
to a constant, we decompose the smooth transonic flow problem into
a smooth subsonic-sonic flow problem with free boundary
and a smooth sonic-supersonic flow problem with fixed boundary in the potential plane,
which can be solved separately.
For the smooth subsonic-sonic flow problem, we encounter the two main difficulties in
the study on the continuous subsonic-sonic flow problem in \cite{WX2}.
One is that the problem is a characteristic degenerate free boundary problem
and the degeneracy occurs just at the free boundary.
The other is that at the inlet we should prescribe a Neumann boundary condition
instead of a Robin one with which the problem may be ill-posed.
However, besides these difficulties, there is an additional disadvantage in this study.
Different from the problem in \cite{WX2},
there are no background solutions,
which play an important role to determine and control the rate of the flow
from the subsonic region to the sonic state.
Furthermore, the boundary conditions at the inlet and the wall are nonlinear, nonlocal and implicit.
In the paper, we try to seek a subsonic-sonic flow
which tends uniformly to the sonic state along different streamlines
although there are no background solutions.
To this end, we assume that
the convergent part of the nozzle slopes so gently that
\begin{align}
\label{fo1}
f''(x)=o(x^2),\quad x\to0^-.
\end{align}
Indeed, \eqref{fo1} is necessary to
guarantee that the change of the speed along the stream direction is infinitesimal of higher order
than the one along the potential direction near the sonic state.
By a fixed point argument with many very precise elliptic estimates,
we obtain a subsonic-sonic flow which tends uniformly to the sonic state along different streamlines.
For the smooth sonic-supersonic flow problem,
the governing equation is a non-strictly hyperbolic equation with strong singularity at the sonic curve,
where its two eigenvalues coincide and the eigenvector space reduces
to a one-dimensional space.
Furthermore, the singularity at the sonic curve is so strong
that all characteristics from sonic points coincide,
which are the sonic curve and never approach the supersonic region.
And the boundary condition on the wall is nonlinear, nonlocal and implicit.
As mentioned above,
a crucial step to solve the problem is how to determine and control the rate of the flow
from the sonic state to supersonic region.
We seek a sonic-supersonic flow
which moves uniformly away from the sonic state along different streamlines in the paper.
For the sonic-supersonic flow in the divergent part of the de Laval nozzle,
\begin{align}
\label{fo2}
f''(x)=o(x^2),\quad x\to0^+
\end{align}
is necessary to
guarantee that the change of the speed along the stream direction is infinitesimal of higher order
than the one along the potential direction near the sonic state.
We solve this sonic-supersonic flow problem under the assumption \eqref{fo2}
by a fixed point argument and the method of characteristics.
Due to the strong singularity, the computations of the flow near the sonic curve are quite complicated.
Through some precise calculations and optimal estimates,
we are able to get a desired sonic-supersonic flow.

In the present paper, by solving a smooth subsonic-sonic flow problem with free boundary
and a smooth sonic-supersonic flow problem with fixed boundary separately,
we get a smooth transonic flow of Meyer type in the de Laval nozzle $\Omega$,
whose velocity vector is along the normal direction at the inlet
and which satisfies the slip condition on the wall,
under the assumptions \eqref{fo1}, \eqref{fo2} and $|l_\pm|$ being sufficiently small.
Different from the examples of smooth transonic flows of Meyer type by using power series expansions,
where the exceptional point is isolated,
the sonic curve of the smooth transonic flow in this paper is located at the throat and
each sonic point is exceptional.
Thus, this transonic flow pattern is strongly singular in the sense that
the sonic curve is a characteristic degenerate boundary in the subsonic-sonic region,
while in the sonic-supersonic region all characteristics from sonic points coincide,
which are the sonic curve and never approach the supersonic region.
It is surprising that there is a smooth transonic flow for this pattern with
so strong singularity.
Indeed, we get a smooth transonic flow
in the sense that the acceleration is Lipschitz continuous.
Furthermore, the transonic flow of this pattern is also shown to be unique.
It is noted that \eqref{fo1} and \eqref{fo2} are almost necessary for a $C^2$ transonic flow
whose sonic curve is located at the throat of the nozzle.
Moreover, the assumption $|l_\pm|$ being sufficiently small can be relaxed.
In the convergent part and the divergent part of the nozzle,
we get a smooth subsonic-sonic flow
and a smooth sonic-supersonic flow, respectively.
Using the same methods, one can get a smooth transonic flow of Meyer type
in a long de Laval nozzle if it possesses the same properties near the throat
(\eqref{fo1} and \eqref{fo2}, which are almost necessary)
and its walls slope gently away from the throat.
More generally, we believe that there also exist smooth transonic flows
for general finite de Laval nozzles satisfying \eqref{fo1} and \eqref{fo2}.
In the upstream, by more complicated elliptic estimates than the ones in the present paper,
we can get a subsonic-sonic flow for a very general convergent nozzle.
In the downstream, the flow can be extended for a class of long divergent nozzles
by solving an initial-boundary value problem of a strictly hyperbolic equation.
These questions will be dealt with in our forthcoming studies.

The paper is arranged as follows.
In $\S\, 2$ we analyze the structure and the location of the sonic curve
for $C^2$ transonic flows and show the instability
of transonic flows with nonexceptional points,
which motives us to seek a smooth transonic flow of Meyer type whose sonic points are exceptional.
We formulate the problem of smooth transonic flows and state the main results
(existence and uniqueness) of the paper in $\S\, 3$.
The smooth transonic flow problem can be decomposed into a smooth subsonic-sonic flow problem
with free boundary and a smooth sonic-supersonic flow problem with fixed boundary
in the potential plane,
which are investigated in $\S\, 4$ and $\S\, 5$, respectively.

\section{Sonic curves and instability of transonic flows with nonexceptional points}

In this section, we analyze the structure and the location of sonic curves
for $C^2$ transonic flows.
And we always assume that the sonic curve belongs to both the boundary of the subsonic region
and the boundary of the supersonic region.

\subsection{Governing equations}

Let us first rewrite the system \eqref{euler-1}--\eqref{Bernoulli} in the physical plane
as the Chaplygin equations in the potential plane.

Since the conservation of momentum,
\eqref{euler-2} and \eqref{euler-3}, can be derived from
\eqref{euler-1}, \eqref{euler-4}--\eqref{Bernoulli},
the flow is just governed by the conservation of mass and
the condition of irrotationality:
\begin{align}
\label{con-irro}
\pd{}{x}(\rho u)+\pd{}{y}(\rho v)=0,
\quad \pd{u}{y}-\pd{v}{x}=0,
\end{align}
where $\rho$ is related to the velocity $(u,v)$ by the Bernoulli
law \eqref{Bernoulli}. The sound speed $c$ is defined as
$$
c^2=P'(\rho)=\rho^{\gamma-1}=1-\frac{\gamma-1}{2}q^2.
$$
At the sonic state, the speed is
$$
c_*=\Big(\frac{2}{\gamma+1}\Big)^{1/2},
$$
which is called the critical speed in the sense that
the flow is subsonic ($q<c$) when $0<q<c_*$, sonic ($q=c$) when $q=c_*$
and supersonic ($q>c$) when $c_*<q<q_{max}$.

Define a velocity potential $\varphi$ and a stream
function $\psi$, respectively, by
$$
\pd\varphi x=u,\quad\pd\varphi y=v, \quad\pd\psi x=-\rho
v,\quad\pd\psi y=\rho u,
$$
which are
$$
\pd\varphi x=q\cos\theta,\quad\pd\varphi y=q\sin\theta,\quad
\pd\psi x=-\rho q\sin\theta,\quad\pd\psi y=\rho q\cos\theta
$$
in terms of polar coordinators in the velocity space,
where $\theta$, which is called a flow angle,
is the angle of the velocity inclination to
the $x$-axis. Direct calculations show that the system
\eqref{con-irro} can be reduced to the Chaplygin equations:
\begin{align}
\label{Chaplygin}
\pd\theta\psi+\frac{\rho(q^2)+2q^2\rho'(q^2)}{q\rho^2(q^2)}
\pd{q}\varphi=0,
\quad\frac{1}{q}\pd{q}\psi-\frac{1}{\rho(q^2)}\pd\theta\varphi=0
\end{align}
in the potential-stream coordinates $(\varphi,\psi)$. Note that
$$
\pd{(\varphi,\psi)}{(x,y)}=\pd\varphi x\pd\psi y-
\pd\varphi y\pd\psi x=\rho q^2.
$$
So the coordinates transformation
between the two coordinate systems
is valid at least in the absence of stagnation points. Eliminating
$\theta$ from \eqref{Chaplygin}
yields the following second-order quasilinear equation
\begin{align}
\label{eq2or}
\pd{^2A(q)}{\varphi^2}+\pd{^2B(q)}{\psi^2}=0,
\end{align}
where
\begin{align*}
A(q)=\int_{c_*}^q\frac{\rho(s^2)+2s^2\rho'(s^2)}{s\rho^2(s^2)}ds,\quad
B(q)=\int_{c_*}^q\frac{\rho(s^2)}{s}ds, \qquad 0<q<q_{max}.
\end{align*}
Here, $B(\cdot)$ is strictly increasing in $(0,q_{max})$,
while $A(\cdot)$ is strictly increasing in $(0,c_*]$ and
strictly decreasing in $[c_*,q_{max})$.
The inverse function of $B(\cdot)$ is denoted by $B^{-1}(\cdot)$,
while the inverse functions of $A(\cdot)$ lying in $(0,c_*]$ and
$[c_*,q_{max})$ are denoted by $A_{-}^{-1}(\cdot)$
and $A_{+}^{-1}(\cdot)$, respectively.

It can be checked easily that both \eqref{potential} and
\eqref{eq2or} are elliptic in the subsonic region $(q<c_*)$
and hyperbolic in the supersonic region $(q>c_*)$,
while degenerate and singular at the sonic state $(q=c_*)$. Therefore, the
governing equations \eqref{potential} and \eqref{eq2or} for
transonic flows are both second-order mixed type quasilinear equations
and are both degenerate and singular at sonic states.

\subsection{Exceptional points in the physical plane}

We begin with a description of sonic curves in \cite{Bers}. Let $\mathcal S$ be a sonic curve
of a $C^2$ transonic flow. The positive direction on $\mathcal S$ is defined by requiring that, if
one moves along $\mathcal S$ in this direction, the subsonic region
is located on the left.
Then, $\theta_s$, the derivative of $\theta$ with respect to the arc length on $\mathcal S$, satisfies
\begin{align}
\label{sonicc-1}
\theta_s=-\frac{\sin^2\vartheta}{c_*}\pd{q}{\nu},
\end{align}
where $\nu$ is the unit normal of $\mathcal S$ pointing into the supersonic region
and $\vartheta$ is the angle between the velocity vector and $\nu$.

\hskip85mm
\setlength{\unitlength}{0.6mm}
\begin{picture}(250,50)

\put(0,0){\cbezier(-30,-30)(-12,-15)(-28,35)(-40,40)}

\put(-25,15){\vector(2,1){15}}
\put(-25,15){\vector(2,-1){25}}

\put(-25,15){\qbezier(2.2,-1.42)(3.6,-0.5)(2.3,1)}

\put(-8,20){$\nu$}
\put(1,0){$(u,v)$}
\put(-20,13){$\vartheta$}
\put(-23,-20){$\mathcal S$}

\put(-50,-10){subsonic}
\put(-18,-10){supersonic}

\put(-22.8,-15){\vector(1,4){1}}
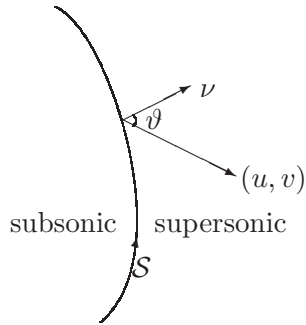
\put(-76,-42){Figure: the orientation of the sonic curve}

\end{picture}
\vskip30mm

It follows from \eqref{sonicc-1} that

\begin{lemma}
\label{lemmasoniccuver}
$\theta$ is nonincreasing along $\mathcal S$.
\end{lemma}

As in \cite{Bers}, points, where $\theta_s=0$, are called exceptional.
According to \eqref{sonicc-1}, a point is exceptional if and only if 
at this point either $\pd{q}{\nu}=0$
or the velocity vector is orthogonal to the sonic curve. 
Indeed, at a point on $\mathcal S$, if $\pd{q}{\nu}=0$, then 
the velocity vector is orthogonal to the sonic curve
owing to the following lemma by Gilbarg and Shiffman \cite{GS}.

\hskip77mm
\setlength{\unitlength}{0.6mm}
\begin{picture}(250,65)
\put(-30,18){\vector(1,0){60}}
\put(0,0){\vector(0,1){45}}
\put(0,18){\circle{30}}
\put(10.7,23.3){\circle*{1.5}}
\put(0,18){\vector(2,1){15}}
\put(10.7,23.3){\vector(1,2){3}}

\put(10,32){$(u,v)$}
\put(12,19){$P$}
\put(-5,46){$y$}
\put(30,13){$x$}
\put(-12,13.5){subsonic}
\put(-28,-14){Figure of Lemma \ref{GSlemma}}

\end{picture}
\vskip15mm

\begin{lemma}
\label{GSlemma}
Assume that $q=c_*$ at a point $P$ on the circumference of a circle in whose interior $q<c_*$.
If the velocity vector is not along the direction from the center of the circle to $P$, then
at $P$ the derivative of $q$ along this direction is positive.
\end{lemma}

\Proof
The case when the circle is the interior of a subsonic-sonic region was considered in \cite{GS}.
The authors derived a linear elliptic equation which is degenerate
at the sonic state from \eqref{potential}.
A comparison principle holds for this degenerate equation, which can be proved by the same
argument for uniformly elliptic equations.
Then, the lemma follows from the same proof of the Hopf lemma
and the auxiliary function has been constructed in \cite{GS}.
Here, it is used that the velocity vector is not along the direction from the center of the circle to $P$,
which shows that this direction
is not characteristic for the degenerate equation.
The proof is standard and thus omitted.
$\hfill\Box$\vskip 4mm

Therefore, one has

\begin{proposition}
\label{soniccurve1}
For a $C^2$ transonic flow, a point at the sonic curve is exceptional
if and only if the velocity vector
is orthogonal to the sonic curve at this point.
\end{proposition}

This proposition shows

\begin{remark}
For a $C^2$ transonic flow, exceptional points are characteristic degenerate in the subsonic-sonic region.
\end{remark}

Exceptional points are regarded to be singular in the following sense (\cite{Bers}). 
Assume that a subsonic flow is defined
in some domain whose boundary contains a smooth sonic curve $\mathcal S$.
If there is no exceptional point on $\mathcal S$, then, the flow can be continued in a unique way across
$\mathcal S$ as a supersonic flow without discontinuity. The flow will be determined in some neighborhood
of $\mathcal S$ contained between the two Mach lines from the two endpoints of $\mathcal S$.
However, if there is a unique exceptional point on $\mathcal S$, the flow will be uniquely determined only
in a subset of the same neighborhood, 
where the points between the two characteristics from the exceptional point
are excluded; furthermore, if the flow can be continued at all into this excluded region, 
this continuation will not be unique.
In this section, we will show that 
the set consisting of exceptional points for any $C^2$ transonic flow
of Meyer type is a closed line segment (may be empty or only one point),
while there is no exceptional point for any $C^2$ transonic flow of Taylor type.
Furthermore, the singularity of exceptional points in the sonic-supersonic region is so strong
that there are two different characteristics from each nonexceptional point in the nozzle,
while all characteristics from interior exceptional points coincide
and they never approach the supersonic region locally.
In the present paper, we get a smooth transonic flow of Meyer type
whose sonic points are all exceptional.
This transonic flow pattern is strongly singular in the sense that
the sonic curve is a characteristic degenerate boundary in the subsonic-sonic region,
while in the sonic-supersonic region all characteristics from sonic points coincide,
which are the sonic curve and never approach the supersonic region.
It is surprising that there is a smooth transonic flow for this pattern with
so strong singularity.
Indeed, we get a smooth transonic flow
in the sense that the acceleration is Lipschitz continuous.
This transonic flow pattern also answers the problem
mentioned in \cite{Bers} that it would
be interesting to know whether exceptional points are always isolated.

\subsection{Exceptional points in the potential plane}

It turns out to be more convenient to analyze exceptional points in the potential plane.
Transforming Proposition \ref{soniccurve1} from the Cartesian coordinates
to the potential-stream coordinates, one can get

\begin{proposition}
\label{soniccurve3}
For a $C^2$ transonic flow in the potential plane,
a sonic point is exceptional if and only if
$\pd{q}{\psi}=0$ at this point.
\end{proposition}

\Proof
On the sonic curve, it holds
\begin{align*}
v\pd{q}{x}-u\pd{q}{y}
=v\Big(u\pd{q}{\varphi}-\rho(c_*^2)v\pd{q}{\psi}\Big)
-u\Big(v\pd{q}{\varphi}+\rho(c_*^2)u\pd{q}{\psi}\Big)
=-\rho(c_*^2)c_*^2\pd{q}{\psi}.
\end{align*}
Therefore, at a sonic point the velocity vector
is orthogonal to the sonic curve in the physical plane
if and only if $\pd{q}{\psi}=0$ at this point.
Then, the proposition follows from Proposition \ref{soniccurve1}.
$\hfill\Box$\vskip 4mm

\hskip11mm
\setlength{\unitlength}{0.6mm}
\begin{picture}(250,65)

\put(80,18){\vector(1,0){60}}
\put(110,0){\vector(0,1){45}}
\put(110,18){\circle{30}}
\put(120.7,23.3){\circle*{1.5}}
\put(110,18){\vector(2,1){15}}

\put(122,19){$P$}
\put(105,46){$\psi$}
\put(140,13){$\varphi$}
\put(98,13.5){subsonic}
\put(80,-14){Figure of Lemma \ref{WWlemma}}

\end{picture}
\vskip15mm

Similar to the proof of Lemma \ref{GSlemma}, it can be shown that

\begin{lemma}
\label{WWlemma}
Consider a $C^2$ flow in the potential plane.
Assume that $q=c_*$ at a point $P$ on the circumference
of a circle in whose interior $q<c_*$.
If the direction from the center of the circle to $P$
is not parallel to the $\varphi$-axis,
then $\pd{q}{\mu}>0$ at $P$, where $\mu$ is any direction which forms an acute angle with
the direction from the center of the circle to $P$.
\end{lemma}

\Proof
Without loss of generality, we assume that $P=(\varphi_0,\psi_0)$ is the unique sonic point
on the circumference of the circle.
Otherwise, we can consider the flow in a smaller circle
which intersects the original circle at $P$.
For convenience, it is assumed that the circle is centered at the origin with radius $r$.
Since the direction from the center of the circle to $P$
is not parallel to the $\varphi$-axis, $\psi_0\not=0$.
Set
$$
\overline q(\varphi,\psi)=c_*-\varepsilon
\Big(\mbox{e}^{-M(\varphi^2+\psi^2)}-\mbox{e}^{-Mr^2}\Big),
\quad(\varphi,\psi)\in\overline B_1\cap\overline B_2
$$
with positive constants $\varepsilon\,(\le c_*/2)$ and $M$ to be determined, where
$$
B_1=\big\{(\varphi,\psi):\varphi^2+\psi^2<r^2\big\},
\quad B_2=\Big\{(\varphi,\psi):(\varphi-\varphi_0)^2+(\psi-\psi_0)^2<\frac14\psi_0^2\Big\}.
$$
Note that
\begin{align*}
&\frac{c_*}{2}<\overline q(\varphi,\psi)<c_*,
\quad(\varphi,\psi)\in B_1\cap B_2
\end{align*}
and
\begin{align*}
&\pd{^2A(\overline q)}{\varphi^2}+\pd{^2B(\overline q)}{\psi^2}
\\
\le&
A'(\overline q(\varphi,\psi))\varepsilon(2M-4M^2\varphi^2)\mbox{e}^{-M(\varphi^2+\psi^2)}
+B'(\overline q(\varphi,\psi))\varepsilon(2M-4M^2\psi^2)\mbox{e}^{-M(\varphi^2+\psi^2)}
\\
\le&
-M\varepsilon\big(B'(\overline q(\varphi,\psi))(M\psi_0^2-2)-2A'(\overline q(\varphi,\psi))\big)
\mbox{e}^{-M(\varphi^2+\psi^2)}
\\
\le&
-M\varepsilon\Big((M\psi_0^2-2)\min_{c_*/2\le s\le c_*}B'(s)-2\max_{c_*/2\le s\le c_*}A'(s)\Big)
\mbox{e}^{-M(\varphi^2+\psi^2)},
\quad(\varphi,\psi)\in B_1\cap B_2.
\end{align*}
Therefore, there exists a sufficiently large $M>0$ such that
$\overline q$ is a supersolution to \eqref{eq2or} in $B_1\cap B_2$ for
each $\varepsilon\in(0,c_*/2]$. Choose $\varepsilon\in(0,c_*/2]$ so small that
$$
q(\varphi,\psi)\le\overline  q(\varphi,\psi),
\quad(\varphi,\psi)\in B_1\cap\partial B_2.
$$
Then, it follows from a comparison principle that
$$
q(\varphi,\psi)\le\overline q(\varphi,\psi),
\quad(\varphi,\psi)\in B_1\cap B_2,
$$
which leads to
$$
\pd{q}{\mu}(\varphi_0,\psi_0)
\ge\pd{\overline q}{\mu}(\varphi_0,\psi_0)>0.
$$
Here, the comparison principle
is not the classical one since \eqref{eq2or} is degenerate at the sonic state
and it can be proved similar to Proposition 3.1 in \cite{WX2}.
$\hfill\Box$\vskip 4mm

As a corollary of Proposition \ref{soniccurve3} and Lemma \ref{WWlemma}, it holds that

\begin{lemma}
\label{soniccurveform}
For a $C^2$ transonic flow in the potential plane,
if an interior point in the nozzle is exceptional, then the sonic curve near this point is
a graph of function with respect to $\psi$.
\end{lemma}

\Proof
If the lemma is false, then $\pd{q}{\varphi}=0$ at this exceptional point.
Note that Proposition \ref{soniccurve3} implies $\pd{q}{\psi}=0$ at this exceptional point.
So the derivative of $q$ at this exceptional point along any direction is zero.
Since this exceptional point is an interior point in the nozzle,
Lemma \ref{WWlemma} shows that the normal of the sonic curve at this exceptional point
is parallel to the $\varphi$-axis, which yields that
the sonic curve near this point is
a graph of function with respect to $\psi$.
$\hfill\Box$\vskip 4mm

If an exceptional point is on the wall of the nozzle, one can get that

\begin{lemma}
\label{soniccurveform2}
Consider a $C^2$ transonic flow in the potential plane with the sonic curve given by
$$
S: \varphi=S_1(s),\,\psi=S_2(s),\,(S'_1(s))^2+(S'_2(s))^2>0,\quad 0\le s\le1,\quad S_1,S_2\in C^2([0,1]).
$$
Assume that $S_2(0)=\inf_{(0,1)}S_2<\sup_{(0,1)}S_2=S_2(1)$
and the subsonic region is located on the left of $S$.

{\rm(i)} If $S'_1(0)\ge0$ and $(S_1(0),S_2(0))$ is exceptional, then $S$
near $(S_1(0),S_2(0))$ is a graph of function with respect to $\psi$.

{\rm(ii)} If $S'_1(1)\le0$ and $(S_1(1),S_2(1))$ is exceptional, then $S$
near $(S_1(1),S_2(1))$ is a graph of function with respect to $\psi$.
\end{lemma}

\Proof
We will prove (ii) only since (i) can be dealt with similarly.
If $S'_1(1)=0$, then $S'_2(1)\not=0$ and the conclusion of the lemma is obvious.
Therefore, if (ii) is false, then $S'_1(1)<0$ and $S$ near $(S_1(1),S_2(1))$ is a graph of function with respect to $\varphi$,
which is denoted by
$$
\psi=h(\varphi),\quad S_1(1)\le\varphi\le S_1(1)+\tau\qquad(\tau>0)
$$
with
\begin{align}
\label{soniccurveform2-1}
h'(S_1(1))=0.
\end{align}
It follows from $S_2(1)=\sup_{(0,1)}S_2$ that
$$
h(\varphi)\le h(S_1(1)),\quad S_1(1)\le\varphi\le S_1(1)+\tau,
$$
which, together with \eqref{soniccurveform2-1}, leads to
$$
-\infty\le h''(S_1(1))\le0.
$$
If $-\infty<h''(S_1(1))\le0$, then there exists a circle whose interior is located at the subsonic region
and which intersects with $S$ at $(S_1(1),S_2(1))$.
Thus Lemma \ref{WWlemma} shows $\pd{q}{\psi}(S_1(1),S_2(1))>0$, which contradicts $(S_1(1),S_2(1))$ is exceptional
by Proposition \ref{soniccurve3}.
Therefore,
\begin{align}
\label{soniccurveform2-2}
h''(S_1(1))=-\infty.
\end{align}
Due to \eqref{soniccurveform2-1} and \eqref{soniccurveform2-2},
there exists a number $\varepsilon\in(0,\tau]$ such that
$$
h'(\varphi)<0,\quad S_1(1)<\varphi<S_1(1)+\varepsilon,
$$
which shows that $S$ near $(S_1(1),S_2(1))$ is a graph of function with respect to $\psi$.
$\hfill\Box$\vskip 4mm

\begin{remark}
As will be shown later (see Remark \ref{structure2}),
neither $S'_1(0)>0$ in {\rm (i)} nor $S'_1(1)<0$ in {\rm (ii)} can occur.
\end{remark}

\begin{remark}
In Lemma \ref{soniccurveform2}, if either $S'_1(0)<0$ in {\rm (i)}
or $S'_1(1)>0$ in {\rm (ii)}, then the similar proof is invalid
and it is unknown whether the conclusion of the lemma holds.
At this moment these two cases cannot be excluded for transonic flows in a nozzle.
\end{remark}

\subsection{Behavior of sonic curves near exceptional points}

Consider a $C^2$ transonic flow with the sonic curve $S$ in the potential plane.
To study the behavior of the sonic curve near exceptional points,
we start with a geometric property of the sonic curve.

\hskip44mm
\setlength{\unitlength}{0.6mm}
\begin{picture}(250,65)
\put(-35,0){\vector(1,0){70}}
\put(0,33){\vector(0,1){15}}

\put(0,33){\qbezier(-35,0)(0,0)(35,0)}

\put(-25,15){subsonic}
\put(5,15){supersonic}

\put(0,5){\cbezier(0,0)(2,4)(5,16)(1,23)}

\put(-7,46){$\psi$}
\put(30,-5){$\varphi$}
\put(-32,-14){Figure of Lemma \ref{lemma-aaa1} (i)}

\put(75,0){\vector(1,0){70}}
\put(110,33){\vector(0,1){15}}

\put(110,33){\qbezier(-35,0)(0,0)(35,0)}

\put(80,15){supersonic}
\put(115,15){subsonic}

\put(115,5){\cbezier(0,0)(-2,4)(-5,16)(-1,23)}

\put(103,46){$\psi$}
\put(140,-5){$\varphi$}
\put(78,-14){Figure of Lemma \ref{lemma-aaa1} (ii)}

\end{picture}
\vskip15mm

\begin{lemma}
\label{lemma-aaa1}
Assume that $S_0:\varphi=h_0(\psi)\,(\psi_1\le\psi\le\psi_2,\,\psi_1<\psi_2)$ is a portion of $S$.

{\rm (i)} If the flow is subsonic
in $\{(\varphi,\psi):\mu_1\le\varphi<h_0(\psi),\psi_1\le\psi\le\psi_2\}$
with a number $\mu_1<\min_{[\psi_1,\psi_2]}h_0$, then
$$
h_0(\psi)\ge\min\{h_0(\psi_1),h_0(\psi_2)\},\quad \psi_1\le\psi\le\psi_2.
$$

{\rm (ii)} If the flow is subsonic
in $\{(\varphi,\psi):h_0(\psi)<\varphi\le\mu_2,\psi_1\le\psi\le\psi_2\}$
with a number $\mu_2>\max_{[\psi_1,\psi_2]}h_0$, then
$$
h_0(\psi)\le\max\{h_0(\psi_1),h_0(\psi_2)\},\quad \psi_1\le\psi\le\psi_2.
$$
\end{lemma}

\Proof
We prove (i) only and the proof of (ii) is similar. Otherwise,
there exists $\psi_0\in(\psi_1,\psi_2)$ such that
$$
h_0(\psi_0)=\min_{[\psi_1,\psi_2]}h<\min\{h_0(\psi_1),h_0(\psi_2)\}.
$$
Set
$$
\bar q(\varphi,\psi)=A^{-1}_{-}(\delta(\varphi-h_0(\psi_0))),\quad
\mu_1\le\varphi\le h_0(\psi_0),\,\psi_1\le\psi\le\psi_2,
$$
where $\delta$ is a positive number so small that
\begin{align*}
\sup_{(\psi_1,\psi_2)}A(q(\mu_1,\cdot))\le-\delta(h_0(\psi_0)-\mu_1),\quad
\sup_{(\mu_1,h_0(\psi_0))}
\big\{A(q(\cdot,\psi_1)),A(q(\cdot,\psi_2))\big\}\le-\delta(h_0(\psi_0)-\mu_1).
\end{align*}
Then, it follows from a comparison principle that
\begin{align*}
q(\varphi,\psi)\le\bar q(\varphi,\psi)=A^{-1}_{-}(\delta(\varphi-h_0(\psi_0))),\quad
\mu_1\le\varphi\le h_0(\psi_0),\,\psi_1\le\psi\le \psi_2,
\end{align*}
which leads to $\pd{q}{\varphi}(h_0(\psi_0),\psi_0)=+\infty$.
Here, the comparison principle
is not the classical one since \eqref{eq2or} is degenerate at the sonic state
and it can be proved similar to Proposition 3.1 in \cite{WX2}.
$\hfill\Box$\vskip 4mm

\hskip75mm
\setlength{\unitlength}{0.6mm}
\begin{picture}(250,65)
\put(-35,5){\vector(1,0){70}}
\put(0,38){\vector(0,1){15}}

\put(0,35){\qbezier(-35,0)(0,0)(35,0)}

\put(10,23){\qbezier(5,-10)(0,0)(4,8)}
\put(12,23){\circle*{1.5}}

\put(12,23){\qbezier[18](0,-10)(0,-1)(0,8)}
\put(12,31){\qbezier[30](-30,0)(-15,0)(0,0)}
\put(12,13){\qbezier[30](-30,0)(-15,0)(0,0)}
\put(-18,23){\qbezier[18](0,-10)(0,-1)(0,8)}

\put(-14,20){subsonic}
\put(13,20){$(h_0(\psi_0),\psi_0)$}

\put(-7,51){$\psi$}
\put(30,0){$\varphi$}
\put(-45,-9){Figure of the proof of Lemma \ref{lemma-aaa1}}

\end{picture}
\vskip10mm

The following lemma and proposition describe the behavior of the sonic curve near an exceptional point.

\hskip70mm
\setlength{\unitlength}{0.6mm}
\begin{picture}(250,65)
\put(-35,0){\vector(1,0){70}}
\put(0,33){\vector(0,1){15}}

\put(0,33){\qbezier(-35,0)(0,0)(35,0)}

\put(10,30){\qbezier(-2,-10)(-0.2,-3)(0,0)}
\put(10,30){\circle*{1.5}}
\put(10,30){\qbezier[12](0,0)(0,-6)(0,-12)}
\put(-16,21){subsonic}
\put(13,21){supersonic}
\put(11.5,27.5){$(h_0(\psi_2),\psi_2)$}

\put(-5,3){\qbezier(-2,10)(-0.2,3)(0,0)}
\put(-5,3){\circle*{1.5}}
\put(-5,3){\qbezier[12](0,12)(0,6)(0,0)}
\put(-32,10){subsonic}
\put(-2,10){supersonic}
\put(-4,2){$(h_0(\psi_1),\psi_1)$}

\put(-7,46){$\psi$}
\put(30,-5){$\varphi$}
\put(-32,-14){Figure of Lemma \ref{lemma-aaa2}}

\end{picture}
\vskip15mm

\begin{lemma}
\label{lemma-aaa2}
Assume that $S_0:\varphi=h_0(\psi)\,(\psi_1\le\psi\le\psi_2,\,\psi_1<\psi_2)$
is a portion of $S$ and the flow is subsonic
on the left of $S_0$.

{\rm (i)} If $\pd{q}{\psi}(h_0(\psi_1),\psi_1)=0$ and $h'_0(\psi_1)=0$,
then there exists a positive constant $\varepsilon\,(\le\psi_2-\psi_1)$ such that
$$
h_0(\psi)\le h_0(\psi_1),\quad \psi_1\le\psi\le\psi_1+\varepsilon.
$$

{\rm (ii)} If $\pd{q}{\psi}(h_0(\psi_2),\psi_2)=0$ and $h'_0(\psi_2)=0$,
then  there exists a positive constant $\varepsilon\,(\le\psi_2-\psi_1)$ such that
$$
h_0(\psi)\le h_0(\psi_2),\quad \psi_2-\varepsilon\le\psi\le\psi_2.
$$
\end{lemma}

\Proof
We prove (ii) only and (i) can be proved similarly.
Otherwise, Lemma \ref{lemma-aaa1} shows that there exists $\psi_0\in[\psi_1,\psi_2)$ such that
\begin{align}
\label{lemma-aaa2-1}
h_0(\psi)>h_0(\psi_2),\quad h'_0(\psi)\le 0,\quad\psi_0<\psi<\psi_2.
\end{align}
Furthermore, \eqref{lemma-aaa2-1}, Proposition \ref{soniccurve3} and Lemma \ref{WWlemma} lead to
\begin{align}
\label{lemma-sc2-2}
\pd{q}{\psi}(h_0(\psi),\psi)\ge0,\quad\psi_0\le\psi\le\psi_2
\end{align}
and
\begin{align}
\label{lemma-sc2-2-0}
\Big\{\psi\in(\psi_0,\psi_2):\pd{q}{\psi}(h_0(\psi),\psi)>0\Big\}
\mbox{ is dense in any left neighberhood of }\psi_2.
\end{align}
Owing to \eqref{eq2or}, one gets that
\begin{align}
\label{lemma-sc2-14}
B'(c_*)\pd{^2q}{\psi^2}(h_0(\psi),\psi)=
-A''(c_*)\Big(\pd{q}{\varphi}(h_0(\psi),\psi)\Big)^2
-B''(c_*)\Big(\pd{q}{\psi}(h_0(\psi),\psi)\Big)^2\ge0,
\quad\psi_0\le\psi\le\psi_2.
\end{align}
Therefore,
\begin{align}
\label{lemma-sc2-15}
\pd{q}{\varphi}(h_0(\psi_2),\psi_2)=0.
\end{align}
Otherwise, \eqref{lemma-sc2-14} yields
$\pd{^2q}{\psi^2}(h_0(\psi_2),\psi_2)>0$ and hence
$$
\frac{d}{d\psi}\Big(\pd{q}{\psi}(h_0(\psi),\psi)\Big)\Big|_{\psi=\psi_2}
=\pd{^2q}{\varphi\partial\psi}(h_0(\psi_2),\psi_2)h'_0(\psi_2)
+\pd{^2q}{\psi^2}(h_0(\psi_2),\psi_2)>0,
$$
which contradicts $\pd{q}{\psi}(h_0(\psi_2),\psi_2)=0$ and \eqref{lemma-sc2-2}.

\hskip70mm
\setlength{\unitlength}{0.6mm}
\begin{picture}(250,65)
\put(-35,0){\vector(1,0){70}}
\put(0,33){\vector(0,1){15}}

\put(0,33){\qbezier(-35,0)(0,0)(35,0)}

\put(7,30){\qbezier(7,-20)(2,-16)(0,0)}
\put(-8,30){\qbezier[15](0,0)(7.5,0)(15.1,0)}
\put(5,10){\qbezier[7](0,0)(3.5,0)(7.2,0)}

\put(10,30){\cbezier(3,-5.6)(2,-5.5)(1,-5)(0,0)}
\put(10,30){\cbezier(3,-5.6)(5,-6)(5.5,-8)(6,-12)}
\put(10,30){\cbezier(9,-15.5)(8,-15.3)(7,-15)(6,-12)}
\put(10,30){\cbezier(9,-15.5)(10,-16)(12,-17)(13,-20)}

\put(-8,30){\cbezier(3,-5.6)(2,-5.5)(1,-5)(0,0)}
\put(-8,30){\cbezier(3,-5.6)(5,-6)(5.5,-8)(6,-12)}
\put(-8,30){\cbezier(9,-15.5)(8,-15.3)(7,-15)(6,-12)}
\put(-8,30){\cbezier(9,-15.5)(10,-16)(12,-17)(13,-20)}

\put(11,5){$h$}
\put(23,12){$h_0$}
\put(-19,12){$h_0-\mu$}
\put(2,16){$G$}

\put(10,30){\circle*{1.5}}
\put(-30,18){subsonic}
\put(25,18){supersonic}
\put(11.5,27.5){$(h_0(\psi_2),\psi_2)$}

\put(-7,46){$\psi$}
\put(30,-5){$\varphi$}
\put(-44,-14){Figure of the proof of Lemma \ref{lemma-aaa2}}

\end{picture}
\vskip15mm

Assume that $\mu>0$, $0<\tau\le\psi_2-\psi_0$ and $h\in C^2([\psi_2-\tau,\psi_2])$ satisfy
\begin{gather}
\label{lemma-sc2-3}
h'(\psi_2)=0,\quad h_0(\psi)-\mu<h(\psi)<h_0(\psi),\quad h''(\psi)\ge0,\quad\psi_2-\tau\le\psi\le\psi_2,
\\
\label{lemma-sc2-4}
q(\varphi,\psi)<c_*,\quad h_0(\psi)-\mu\le\varphi\le h(\psi),\,\psi_2-\tau\le\psi\le\psi_2,
\\
\label{lemma-sc2-5}
\pd{q}{\psi}(\varphi,\psi_2-\tau)>0,\quad h_0(\psi_2-\tau)-\mu\le\varphi\le h_0(\psi_2-\tau).
\end{gather}
Set $G=\big\{(\varphi,\psi):h_0(\psi)-\mu<\varphi<h(\psi),\,
\psi_2-\tau<\psi<\psi_2\big\}$ and
$$
\bar q(\varphi,\psi)=c_*+\delta(\psi-\psi_2)(h(\psi)-\varphi),\quad
(\varphi,\psi)\in\overline G,
$$
where $\delta$, which depends only on $\mu$, $\tau$,
$\sup_{(\psi_2-\tau,\psi_2)}q(h_0(\cdot)-\mu,\cdot)$,
$\inf_{(h_0(\psi_2-\tau)-\mu,h_0(\psi_2-\tau))}\pd{q}{\psi}(\cdot,\psi_2-\tau)$
and $\|h'\|_{L^\infty((\psi_2-\tau,\psi_2))}$, is a positive number so small that
\begin{align}
\label{lemma-sc2-6}
q(h_0(\psi)-\mu,\psi)\le\bar q(h_0(\psi)-\mu,\psi),\quad\psi_2-\tau<\psi<\psi_2
\end{align}
and
\begin{align}
\label{lemma-sc2-7}
\pd{q}{\psi}(\varphi,\psi_2-\tau)\ge\pd{\bar q}{\psi}(\varphi,\psi_2-\tau),
\quad h_0(\psi_2-\tau)-\mu<\varphi<h(\psi_2-\tau).
\end{align}
Noting
$$
\pd{^2\bar q}{\varphi^2}(\varphi,\psi)=0,\quad
\pd{^2\bar q}{\psi^2}(\varphi,\psi)=\delta(2h'(\psi)+(\psi-\psi_2)h''(\psi))\le0,
\quad(\varphi,\psi)\in G
$$
due to \eqref{lemma-sc2-3},
one gets that $\bar q$ is a supersolution to \eqref{eq2or} in $G$.
Since $q$ is strictly subsonic in $G$ from \eqref{lemma-sc2-4},
the classical comparison principle holds.
Then, it follows from \eqref{lemma-sc2-6} and \eqref{lemma-sc2-7} that
\begin{align}
\label{lemma-sc2-8}
q(\varphi,\psi)\le c_*+\delta(\psi-\psi_2)(h(\psi)-\varphi),\quad
h_0(\psi)-\mu\le\varphi\le h(\psi),\,\psi_2-\tau\le\psi\le\psi_2.
\end{align}

We now derive a contradiction.
Owing to \eqref{lemma-aaa2-1}--\eqref{lemma-sc2-2-0},
there exist $\mu>0$, $0<\tau\le\psi_2-\psi_0$, $\{h_k\}_{k=1}^\infty\subset C^2([\psi_2-\tau,\psi_2])$
and $\{\tau_n\}_{n=1}^\infty\subset(0,\tau)$ satisfying
\begin{gather}
\label{lemma-sc2-9}
h_0(\psi)-\mu<h_k(\psi)\le h_{k+1}(\psi)<h_0(\psi),\quad\psi_2-\tau\le\psi\le\psi_2,
\quad k=1,2,\cdots,
\\
\label{lemma-sc2-10}
h'_k(\psi_2)=0,\quad h''_k(\psi)\ge0,\quad\psi_2-\tau\le\psi\le\psi_2,
\quad k=1,2,\cdots,
\\
\label{lemma-sc2-10-1}
q(\varphi,\psi)<c_*,\quad h_0(\psi)-\mu\le\varphi<h_0(\psi),\,
\psi_2-\tau\le\psi\le\psi_2,
\\
\label{lemma-sc2-11}
\pd{q}{\psi}(\varphi,\psi_2-\tau)>0,\quad h_0(\psi_2-\tau)-\mu\le\varphi\le h_0(\psi_2-\tau),
\\
\label{lemma-sc2-12}
\lim_{n\to\infty}\tau_n=0\quad\mbox{ and }
\quad\lim_{k\to\infty}h_{k}(\psi_2-\tau_n)=h_0(\psi_2-\tau_n)\mbox{ for } n=1,2,\cdots.
\end{gather}
Note that \eqref{lemma-sc2-9}--\eqref{lemma-sc2-11} imply that
$h_k$ satisfies \eqref{lemma-sc2-3}--\eqref{lemma-sc2-5} for each positive integer $k$.
It follows from \eqref{lemma-sc2-8} that
\begin{align*}
q(\varphi,\psi)\le c_*+\delta(\psi-\psi_2)(h_k(\psi)-\varphi),\quad
h_0(\psi)-\mu\le\varphi\le h_k(\psi),\,\psi_2-\tau\le\psi\le\psi_2,
\quad k=1,2,\cdots,
\end{align*}
where $\delta>0$ is independent of $k$. This,
together with \eqref{lemma-sc2-12} and \eqref{lemma-sc2-15}, leads to
\begin{align*}
\pd{q}{\varphi}(h_0(\psi_2-\tau_n),\psi_2-\tau_n)\ge\delta\tau_n,\quad
n=1,2,\cdots
\end{align*}
and thus
\begin{align}
\label{lemma-sc2-13}
\pd{^2q}{\varphi\partial\psi}(h_0(\psi_2),\psi_2)\le-\delta.
\end{align}
Due to \eqref{lemma-sc2-14} and \eqref{lemma-sc2-13},
there exists a number $\varepsilon\in(0,\tau)$ such that
$$
\frac{d}{d\psi}\Big(\pd{q}{\psi}(h_0(\psi),\psi)\Big)
=\pd{^2q}{\varphi\partial\psi}(h_0(\psi),\psi)h'_0(\psi)
+\pd{^2q}{\psi^2}(h_0(\psi),\psi)
\ge-\frac{\delta}{2}h'_0(\psi),\quad
\psi_2-\varepsilon<\psi<\psi_2,
$$
which leads to
$$
\pd{q}{\psi}(h_0(\psi),\psi)-\pd{q}{\psi}(h_0(\psi_2),\psi_2)
\le-\frac{\delta}{2}(h_0(\psi)-h_0(\psi_2)),\quad
\psi_2-\varepsilon\le\psi\le\psi_2.
$$
This, together with $\pd{q}{\psi}(h_0(\psi_2),\psi_2)=0$,
\eqref{lemma-aaa2-1} and \eqref{lemma-sc2-2},
yields that
$$
\pd{q}{\psi}(h_0(\psi),\psi)=0,\quad
\psi_2-\varepsilon\le\psi\le\psi_2,
$$
which contradicts \eqref{lemma-sc2-2-0}.
$\hfill\Box$\vskip 4mm

Let us remove the restrictions $h'_0(\psi_1)=0$ and $h'_0(\psi_2)=0$
in Lemma \ref{lemma-aaa2}.

\hskip3mm
\setlength{\unitlength}{0.6mm}
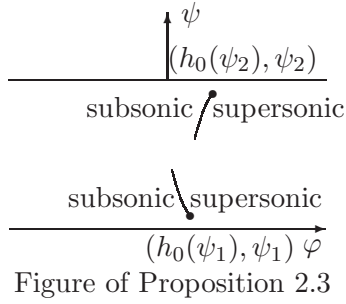
\begin{picture}(250,65)

\put(75,0){\vector(1,0){70}}
\put(110,33){\vector(0,1){15}}

\put(110,33){\qbezier(-35,0)(0,0)(35,0)}

\put(120,30){\qbezier(-4,-10)(-2,-2)(0,0)}
\put(120,30){\circle*{1.5}}
\put(92,25){subsonic}
\put(120,25){supersonic}
\put(110,36){$(h_0(\psi_2),\psi_2)$}

\put(115,3){\qbezier(-4,10)(-2,2)(0,0)}
\put(115,3){\circle*{1.5}}
\put(88,5){subsonic}
\put(115,5){supersonic}
\put(105,-6){$(h_0(\psi_1),\psi_1)$}

\put(113,46){$\psi$}
\put(140,-5){$\varphi$}
\put(76,-14){Figure of Proposition \ref{lemma-aaa3}}

\end{picture}
\vskip15mm

\begin{proposition}
\label{lemma-aaa3}
Assume that $S_0:\varphi=h_0(\psi)\,(\psi_1\le\psi\le\psi_2,\,\psi_1<\psi_2)$
is a portion of $S$ and the flow is subsonic
on the left of $S_0$.

{\rm (i)} If $\pd{q}{\psi}(h_0(\psi_1),\psi_1)=0$,
then there exists a positive constant $\varepsilon\,(\le\psi_2-\psi_1)$ such that
$$
h_0(\psi)\le h_0(\psi_1),\quad \psi_1\le\psi\le\psi_1+\varepsilon.
$$

{\rm (ii)} If $\pd{q}{\psi}(h_0(\psi_2),\psi_2)=0$,
then there exists a positive constant $\varepsilon\,(\le\psi_2-\psi_1)$ such that
$$
h_0(\psi)\le h_0(\psi_2),\quad \psi_2-\varepsilon\le\psi\le\psi_2.
$$
\end{proposition}

\Proof
We prove (ii) only and the proof of (i) is similar.
If the lemma is false, then Lemma \ref{lemma-aaa2} implies that
$-\infty\le h'_0(\psi_2)<0$ and thus
\begin{align}
\label{lemma-aaa3-0-1}
\pd{q}{\varphi}(h_0(\psi_2),\psi_2)=0.
\end{align}
Since the flow is subsonic on the left of $(h_0(\psi_2),\psi_2)$, \eqref{lemma-aaa3-0-1}
implies
\begin{align}
\label{lemma-aaa3-0-2}
\pd{^2q}{\varphi^2}(h_0(\psi_2),\psi_2)\le0.
\end{align}
Choose a positive number $\tau\,(<\psi_2-\psi_1)$ such that
\begin{align}
\label{lemma-aaa3-0}
h'_0(\psi)<0,\quad \psi_2-\tau\le\psi<\psi_2.
\end{align}
It follows from Lemma \ref{WWlemma} and \eqref{lemma-aaa3-0} that
$\pd{q}{\psi}(h_0(\psi_2-\tau),\psi_2-\tau)>0$.
Let $\mu$ be a positive number so small that
$\big\{(\varphi,\psi):h_0(\psi)-\mu\le\varphi<h_0(\psi),\psi_2-\tau<\psi<\psi_2\big\}$ belongs to the subsonic region and
$$
\pd{q}{\psi}(\varphi,\psi_2-\tau)\ge\frac12\pd{q}{\psi}(h_0(\psi_2-\tau),\psi_2-\tau),\quad
h_0(\psi_2-\tau)-\mu<\varphi<h_0(\psi_2-\tau).
$$

We will use a similar method as in Lemma \ref{lemma-aaa2} to complete the proof. Set
\begin{align*}
{\mathscr D}=&\Big\{h\in C^2([\psi_2-\tau,\psi_2]):h_0(\psi)-\mu<h(\psi)<h_0(\psi),\,h'(\psi)<0
\\
&\qquad\qquad\mbox{ and }(\psi_2-\psi)h''(\psi)\ge2 h'(\psi)\mbox{ for each }\psi_2-\tau\le\psi\le\psi_2\Big\}.
\end{align*}
For any $0<\varepsilon<\tau$, if $h\in C^2([\psi_2-\varepsilon,\psi_2))$ satisfies
$$
(\psi_2-\psi)h''(\psi)\le2 h'(\psi)<0,\quad \psi_2-\varepsilon\le\psi<\psi_2,
$$
then $h$ must be unbounded.
Hence, there exist
$\{h_k\}_{k=1}^\infty\subset{\mathscr D}$ and $\{\tau_n\}_{n=1}^\infty\subset(0,\tau)$ such that
\begin{gather}
\label{lemma-scaaa2-12}
\lim_{n\to\infty}\tau_n=0\quad\mbox{ and }
\quad\lim_{k\to\infty}h_{k}(\psi_2-\tau_n)=h_0(\psi_2-\tau_n)\mbox{ for } n=1,2,\cdots.
\end{gather}
For each $k\ge1$, set
$$
\bar q_k(\varphi,\psi)=c_*+\delta(\psi-\psi_2)(h_k(\psi)-\varphi),\quad
h_0(\psi)-\mu\le\varphi\le h_k(\psi),\,\psi_2-\tau\le\psi\le\psi_2,
$$
where $\delta$, which is independent of $k$, is a positive number so small that
\begin{align}
\label{lemma-scaaa2-6}
q(h_0(\psi)-\mu,\psi)\le\bar q_k(h_0(\psi)-\mu,\psi),\quad\psi_2-\tau<\psi<\psi_2
\end{align}
and
\begin{align}
\label{lemma-scaaa2-7}
\pd{q}{\psi}(\varphi,\psi_2-\tau)\ge\pd{\bar q_k}{\psi}(\varphi,\psi_2-\tau),
\quad h_0(\psi_2-\tau)-\mu<\varphi<h_k(\psi_2-\tau).
\end{align}
Note that
$$
\pd{^2\bar q_k}{\psi^2}(\varphi,\psi)=\delta(2h'_k(\psi)+(\psi-\psi_2)h''_k(\psi))\le0,
\quad h_0(\psi)-\mu<\varphi<h_k(\psi),\,\psi_2-\tau<\psi<\psi_2
$$
due to the definition of ${\mathscr D}$.
Thus $\bar q_k$ is a supersolution to \eqref{eq2or}.
Then one can get from the classical comparison principle,
together with \eqref{lemma-scaaa2-6} and \eqref{lemma-scaaa2-7}, that
$$
q(\varphi,\psi)\le c_*+\delta(\psi-\psi_2)(h_k(\psi)-\varphi),\quad
h_0(\psi)-\mu\le\varphi\le h_k(\psi),\,\psi_2-\tau\le\psi\le\psi_2,\quad k=1,2,\cdots,
$$
which, together with \eqref{lemma-scaaa2-12}, yields
\begin{align}
\label{lemma-aaa3-1-0}
\pd{q}{\varphi}(h_0(\psi_2-\tau_n),\psi_2-\tau_n)\ge\delta\tau_n,\quad n=1,2,\cdots.
\end{align}
Then \eqref{lemma-aaa3-1-0}, \eqref{lemma-aaa3-0-1}, \eqref{lemma-aaa3-0-2} and
$-\infty\le h'_0(\psi_2)<0$ imply that
\begin{align}
-\delta\ge\frac{d}{d\psi}\Big(\pd{q}{\varphi}(h_0(\psi),\psi)\Big)\Big|_{\psi=\psi_2}
=&\pd{^2q}{\varphi^2}(h_0(\psi_2),\psi_2)h'_0(\psi_2)
+\pd{^2q}{\varphi\partial\psi}(h_0(\psi_2),\psi_2)
\nonumber
\\
\ge&\pd{^2q}{\varphi\partial\psi}(h_0(\psi_2),\psi_2).
\label{lemma-aaa3-1}
\end{align}
Additionally, \eqref{eq2or} yields
\begin{align}
\label{lemma-aaa3-2}
B'(c_*)\pd{^2q}{\psi^2}(h_0(\psi_2),\psi_2)=
-A''(c_*)\Big(\pd{q}{\varphi}(h_0(\psi_2),\psi_2)\Big)^2
-B''(c_*)\Big(\pd{q}{\psi}(h_0(\psi_2),\psi_2)\Big)^2\ge0.
\end{align}
It follows from \eqref{lemma-aaa3-1}, \eqref{lemma-aaa3-2} and $-\infty\le h'_0(\psi_2)<0$ that
\begin{align}
\label{lemma-aaa3-3}
\frac{d}{d\psi}\Big(\pd{q}{\psi}(h_0(\psi),\psi)\Big)\Big|_{\psi=\psi_2}
=\pd{^2q}{\varphi\partial\psi}(h_0(\psi_2),\psi_2)h'_0(\psi_2)
+\pd{^2q}{\psi^2}(h_0(\psi_2),\psi_2)\in(0,+\infty].
\end{align}
Due to \eqref{lemma-aaa3-3} and $\pd{q}{\psi}(h_0(\psi_2),\psi_2)=0$,
one can get that for sufficiently small $\varepsilon\in(0,\tau)$,
$$
\pd{q}{\psi}(h_0(\psi),\psi)<0,\quad\psi_2-\varepsilon<\psi<\psi_2,
$$
which contradicts \eqref{lemma-aaa3-0} and $q$ is subsonic
on the left of $S_0$.
$\hfill\Box$\vskip 4mm

According to Lemmas \ref{soniccurveform}, \ref{soniccurveform2} and Proposition \ref{lemma-aaa3},
one can get the structure of exceptional points.

\begin{remark}
\label{structure1}
Consider a $C^2$ transonic flow in the potential plane.
If an interior point in the nozzle is exceptional,
then near the exceptional point, the sonic curve can be represented as $\varphi=h(\psi)$
which achieves a local maximum (if the left of this point is subsonic)
or minimum (if the right of this point is subsonic) at this point.
\end{remark}

\begin{remark}
\label{structure2}
Under the assumptions of Lemma \ref{soniccurveform2},
one further gets $S'_1(0)=0$ in {\rm (i)} and $S'_1(1)=0$ in {\rm (ii)}.
\end{remark}

\hskip43mm
\setlength{\unitlength}{0.6mm}
\begin{picture}(250,65)
\put(-35,0){\vector(1,0){70}}
\put(0,33){\vector(0,1){15}}

\put(0,33){\qbezier(-35,0)(0,0)(35,0)}

\put(0,16.5){\qbezier(-2,-7)(-0.2,-3)(0,0)}
\put(0,16.5){\qbezier(-2,7)(-0.2,3)(0,0)}
\put(0,16.5){\qbezier[14](0,-7)(0,0)(0,7)}

\put(12.8,41.5){\vector(-1,-2){12}}

\put(0,16.5){\circle*{1.5}}
\put(2.8,41.5){exceptional point}
\put(-7,46){$\psi$}
\put(30,-5){$\varphi$}
\put(-30,15){subsonic}
\put(8,15){supersonic}
\put(-32,-14){Figure of Remark \ref{structure1}}

\put(75,0){\vector(1,0){70}}
\put(110,33){\vector(0,1){15}}
\put(110,33){\qbezier(-35,0)(0,0)(35,0)}

\put(135,33){\qbezier(-20,0)(-21,-7)(-24,-10)}
\put(135,33){\qbezier[10](-20,0)(-20,-5)(-20,-10)}

\put(122,0){\qbezier(-20,0)(-21,7)(-24,10)}
\put(122,0){\qbezier[10](-20,0)(-20,5)(-20,10)}
\put(115,33){\circle*{1.5}}
\put(102,0){\circle*{1.5}}
\put(110,36){$(S_1(1),S_2(1))$}
\put(92,-6){$(S_1(0),S_2(0))$}

\put(128,20){\vector(-1,1){12}}
\put(120.2,13.6){\vector(-4,-3){17}}

\put(100,14){exceptional point}
\put(103,46){$\psi$}
\put(140,-5){$\varphi$}
\put(75,3){subsonic}
\put(115,3){supersonic}
\put(88,28){subsonic}
\put(128,28){supersonic}
\put(78,-14){Figure of Remark \ref{structure2}}

\end{picture}
\vskip15mm

\begin{remark}
We never use supersonic regions in the discussion in $\S\,2.2,\,2.3,\,2.4$ 
except the orientation of sonic curves. So, the results in these three subsections 
(except Lemma \ref{lemmasoniccuver})
also hold for any $C^2$ subsonic-sonic flow or subsonic-sonic-subsonic flow 
whose sonic points is a curve.
\end{remark}

\subsection{Structure of sonic curves}

Consider a nozzle in the potential plane of the form
\begin{align*}
%\label{soniccurveg}
G=\big\{(\varphi,\psi):g_1(\psi)<\varphi<g_2(\psi),\psi_1<\psi<\psi_2\big\},
\end{align*}
where $\psi_1<\psi_2$ and
$$
g_1(\psi)<g_2(\psi),\quad\psi_1\le\psi\le\psi_2.
$$
The structure of the sonic curve of a $C^2$ transonic flow in $G$ is as follows.

\vskip10mm
\hskip40mm
\setlength{\unitlength}{0.6mm}
\begin{picture}(250,65)
\put(-35,0){\vector(1,0){70}}
\put(0,50){\vector(0,1){15}}
\put(0,50){\qbezier(-35,0)(0,0)(35,0)}

\put(-10,50){\qbezier(4,-8)(-6,-7)(-10,0)}
\put(-10,50){\cbezier(9,-4)(7,-4.2)(6,-8.2)(4,-8)}
\put(-10,50){\qbezier(11,-8)(10,-4.2)(9,-4)}
\put(-10,50){\qbezier(11,-8)(12.4,-13)(12.5,-14)}

\put(-10,50){\qbezier(12.5,-37.5)(12.5,-16)(12.5,-14)}

\put(-10,0){\cbezier(-7,0)(-6,16)(0,9.8)(1,9)}
\put(-10,0){\qbezier(7,5)(5,4)(1,9)}
\put(-10,0){\cbezier(7,5)(8,5.5)(11,7)(12.5,12)}

\put(2.5,36){\circle*{1.5}}
\put(2.5,12.5){\circle*{1.5}}

\put(-25,17){subsonic}
\put(8,17){supersonic}

\put(-4,26){$S_e$}
\put(-15,38){$S_+$}
\put(-24,6){$S_-$}

\put(2.45,23){\vector(0,1){1}}
\put(-1.5,6.1){\vector(1,1){1}}
\put(0.9,43){\vector(-1,2){1}}

\put(-7,63){$\psi$}
\put(30,-5){$\varphi$}
\put(-32,-14){Figure of Theorem \ref{soniclocation}}

\put(85,0){\vector(1,0){70}}
\put(120,50){\vector(0,1){15}}
\put(120,50){\qbezier(-35,0)(0,0)(35,0)}

\put(100,47){\qbezier(4,-8)(-2,-18)(-8,3)}
\put(100,47){\qbezier(4,-8)(10,-2)(16,-10)}
\put(100,47){\qbezier(26,-8)(22,-16)(16,-10)}
\put(100,47){\qbezier(26,-8)(30,0)(33,-11)}
\put(100,47){\qbezier(45,3)(40,-18)(33,-11)}

\put(100,3){\qbezier(10,12)(2,10)(-3,-3)}
\put(100,3){\qbezier(10,12)(15,14)(20,6)}
\put(100,3){\qbezier(28,12)(24,2)(20,6)}
\put(100,3){\qbezier(28,12)(34,20)(40,-3)}

\put(109,23){subsonic}
\put(104,45){supersonic}
\put(104,1.5){supersonic}

\put(113,63){$\psi$}
\put(150,-5){$\varphi$}
\put(88,-14){Figure of Theorem \ref{soniclocationtaylor}}

\put(100.25,7){\vector(1,2){1}}

\put(142.5,40){\vector(-1,2){1}}

\end{picture}
\vskip15mm

\begin{theorem}
\label{soniclocation}
Let $S$ be the sonic curve of a $C^2(\overline G)$ transonic flow of Meyer type.
Then, $S$ is a disjoint union of three connected parts (may be empty): $S_{e}$, $S_{+}$, $S_{-}$,
where $S_{e}$ is the set of exceptional points,
while $S_{+}$ and $S_{-}$ denote the other two connected parts approaching the upper and lower walls, 
respectively. Furthermore,

{\rm(i)} $S_{e}$ is a closed segment parallel to the $\psi$-axis;

{\rm(ii)} $S_{+}$ and $S_{-}$ are two graphes of function with respect to $\varphi$, respectively.
Along the positive direction of $S$,
$\varphi$ is strictly decreasing on $S_{+}$ while strictly increasing on $S_{-}$.

{\rm(iii)} If $S_{e}$ is empty, then $S=S_{+}$ or $S=S_{-}$.
\end{theorem}

\Proof
Denote
$$
S: \varphi=S_1(s),\,\psi=S_2(s),\,(S'_1(s))^2+(S'_2(s))^2>0,\quad 0\le s\le1,\quad S_1,S_2\in C^2([0,1]).
$$
Without loss of generality, it is assumed that $S_2(0)=\psi_1$ and $S_2(1)=\psi_2$,
which is equivalent to that the subsonic region is located on the left of $S$.
Note
$$
\pd q\varphi(S_1(s),S_2(s))S'_1(s)+\pd q\psi(S_1(s),S_2(s))S'_2(s)=0,\quad0\le s\le1.
$$
It is clear that
\begin{align}
\label{sss-1}
\mbox{if $S'_1(s)=0$ for some $s\in[0,1]$, then $\pd q\psi(S_1(s),S_2(s))=0$.}
\end{align}
Lemma \ref{soniccurveform} and Proposition \ref{lemma-aaa3} show that
\begin{align}
\label{sss-2}
\mbox{if $\pd q\psi(S_1(s),S_2(s))=0$ for some $s\in(0,1)$, then $S'_1(s)=0$.}
\end{align}
Owing to \eqref{sss-1} and Proposition \ref{soniccurve3},
$S'_1(0)\not=0$ if $(S_1(0),S_2(0))$ is not exceptional.
There are four cases to be considered.

(i) The case that $(S_1(0),S_2(0))$ is exceptional and $S'_1(0)\ge0$.
Then, it follows from Propositions \ref{soniccurve3}, \ref{lemma-aaa3}
and Lemma \ref{soniccurveform2} that $S'_1(0)=0$. Set
$$
s^*=\sup\{0\le s\le 1:S'_1(s)=0\}.
$$
Due to \eqref{sss-1} and Lemmas \ref{lemma-aaa1}, \ref{lemma-aaa2}, one gets that
\begin{align}
\label{sss-3}
S'_1(s)=0,\quad0\le s\le s^*.
\end{align}
Indeed, if \eqref{sss-3} is false, there exist $s_1,s_2\in[0,s^*]$ such that $s_1<s_2$,
$S'_1(s_1)=S'_1(s_2)=0$ and
\begin{align}
\label{sss-3zz1}
S'_1(s)\not=0,\quad s_1<s<s_2.
\end{align}
It follows from \eqref{sss-1} and Lemma \ref{lemma-aaa2} that
\begin{gather}
\label{sss-4}
S_1(s)\le S_1(s_1),\quad s_1\le s\le s_1+\varepsilon,
\\
\label{sss-5}
S_1(s)\le S_1(s_2),\quad s_2-\varepsilon\le s\le s_2
\end{gather}
for some positive number $\varepsilon\,(\le s_2-s_1)$. Additionally, Lemma \ref{lemma-aaa1} yields
\begin{align}
\label{sss-6}
S_1(s)\ge\min\{S_1(s_1),S_1(s_2)\},\quad s_1\le s\le s_2.
\end{align}
Due to \eqref{sss-4}--\eqref{sss-6}, $S'_1=0$ either on $[s_1,s_1+\varepsilon]$
or on $[s_2-\varepsilon,s_2]$,
which contradicts \eqref{sss-3zz1}.
Owing to \eqref{sss-3}, \eqref{sss-1} and Proposition \ref{soniccurve3},
$(S_1(s),S_2(s))$ is exceptional for each $s\in[0,s^*]$.
If $s^*=1$, $S_e=S$
is a closed segment parallel to the $\psi$-axis.
Otherwise, $0\le s^*<1$. Lemma \ref{lemma-aaa2} and the definition of $s^*$ imply
$$
S'_1(s)<0,\quad s^*<s\le 1.
$$
It follows from Proposition \ref{soniccurve3} and \eqref{sss-2}
that $(S_1(s),S_2(s))$ is not exceptional for each $s\in(s^*,1)$.
Furthermore, Propositions \ref{soniccurve3}, \ref{lemma-aaa3} and Lemma \ref{soniccurveform2}
show that $(S_1(1),S_2(1))$ is not exceptional either.
That is to say,
$S_e=\{(S_1(s),S_2(s)):0\le s\le s^*\}$ and $S_+=\{(S_1(s),S_2(s)):s^*<s\le1\}$.

(ii) The case that $(S_1(0),S_2(0))$ is exceptional and $S'_1(0)<0$.
Let us prove that $S'_1(s)<0$ for each $s\in[0,1]$
by contradiction. Otherwise, there exists $s_0\in(0,1]$ such that $S'_1(s_0)=0$ and
$S'_1(s)<0$ for each $s\in[0,s_0)$. Then, \eqref{sss-1}
and Lemma \ref{lemma-aaa2} show that $S_1(s_0-\varepsilon)\le S_1(s_0)$
for sufficiently small $\varepsilon\in(0,s_0)$,
which contradicts that $S'_1(s)<0$ for each $s\in[0,s_0)$. Therefore,
$$
S'_1(s)<0,\quad 0\le s\le 1.
$$
It follows from Proposition \ref{soniccurve3} and \eqref{sss-2}
that $(S_1(s),S_2(s))$ is not exceptional for each $s\in(0,1)$.
Furthermore, $(S_1(1),S_2(1))$ is not exceptional either
due to Propositions \ref{soniccurve3}, \ref{lemma-aaa3} and Lemma \ref{soniccurveform2}.
That is to say,
$S_e=\{(S_1(0),S_2(0))\}$ and $S_+=\{(S_1(s),S_2(s)):0<s\le1\}$.

(iii) The case that $(S_1(0),S_2(0))$ is not exceptional and $S'_1(0)>0$.
Set
$$
s_*=\sup\{0\le s\le 1:S'_1(s)>0\}.
$$
Then $0<s_*\le1$.
Let us prove
\begin{align}
\label{sss-7}
S'_1(s)>0,\quad0\le s<s_*
\end{align}
by contradiction. Otherwise, there exist $s_0\in(0,s_*)$ and 
$\varepsilon_0\in(0,s_*-s_0]$
such that $S'_1(s_0)=0$ and $S'_1(s)>0$ for each $s\in(s_0,s_0+\varepsilon_0)$.
Then, \eqref{sss-1} and Lemma \ref{lemma-aaa2} show that $S_1(s_0+\varepsilon)\le S_1(s_0)$
for sufficiently small $\varepsilon\in(0,\varepsilon_0)$,
which contradicts that $S'_1(s)>0$ for each $s\in(s_0,s_0+\varepsilon_0)$.
Thus $(S_1(s),S_2(s))$ is not exceptional for each $s\in(0,s_*)$
due to \eqref{sss-7}, \eqref{sss-2} and Proposition \ref{soniccurve3}.
If $s_*=1$, then either $S_-=S$, or $S_-=\{(S_1(s),S_2(s)):0\le s<1\}$ and $S_e=\{(S_1(1),S_2(1))\}$.
Otherwise, $0<s_*<1$ and $S'_1(s_*)=0$.
Set
$$
s^*=\sup\{s_*\le s\le 1:S'_1(s)=0\}.
$$
Then, $0<s_*\le s^*\le1$. Similar to the discussion in (i),
one can get that
$$
S'_1(s)=0,\quad s_*\le s\le s^*;
$$
furthermore, $S_-=\{(S_1(s),S_2(s)):0\le s<s_*\}$ and $S_e=\{(S_1(1),S_2(1)):s_*\le s\le 1\}$
if $s^*=1$, while
$S_-=\{(S_1(s),S_2(s)):0\le s<s_*\}$, $S_e=\{(S_1(1),S_2(1)):s_*\le s\le s^*\}$
and $S_+=\{(S_1(s),S_2(s)):s^*<s\le1\}$ if $s^*<1$.

(iv) The case that $(S_1(0),S_2(0))$ is not exceptional and $S'_1(0)<0$.
As the proof of (ii), one can prove that $S'_1(s)<0$
and $(S_1(s),S_2(s))$ is not exceptional for each $s\in[0,1]$.
Then $S_+=S$.
$\hfill\Box$\vskip 4mm

\begin{theorem}
\label{soniclocationtaylor}
For any $C^2(\overline G)$ transonic flow of Taylor type,
each sonic point is not exceptional and
$\varphi$ is strictly monotone along the sonic curve.
\end{theorem}

\Proof
Without loss of generality, we assume that $S$ is the sonic curve intersecting the lower wall.
Denote
$$
S: \varphi=S_1(s),\,\psi=S_2(s),\,(S'_1(s))^2+(S'_2(s))^2>0,\quad 0\le s\le1,\quad S_1,S_2\in C^2([0,1])
$$
with $S_2(0)=S_2(1)=\psi_2$ and $S_1(0)<S_1(1)$.
Lemma \ref{lemma-aaa1} implies that $S'_1(0)\ge0$ and $S'_1(1)\le0$.

First we show that there is no exceptional point in the nozzle.
Otherwise, it is assumed that $(S_1(s_0),S_2(s_0))$ is an exceptional point
for some $s_0\in(0,1)$.
Owing to Remark \ref{structure1},
$S'_2(s_0)>0$ and $S_1(s_0)$
achieves a local maximum if the left of $(S_1(s_0),S_2(s_0))$ is subsonic,
while $S'_2(s_0)<0$ and $S_1(s_0)$
achieves a local minimum if the right of $(S_1(s_0),S_2(s_0))$ is subsonic.
Note that $S_2(0)=S_2(1)$ and $S_1(0)<S_1(1)$.
If the left of $(S_1(s_0),S_2(s_0))$ is subsonic,
then Lemma \ref{lemma-aaa2} yields that 
there exists $s_1\in(s_0,1)$ such that $S'_2(s_1)>0$, $S_1(s_1)$
achieves a local strict minimum and the left of $(S_1(s_1),S_2(s_1))$ is subsonic.
Similarly, if the right of $(S_1(s_0),S_2(s_0))$ is subsonic,
then there exists $s_2\in(0,s_0)$ such that $S'_2(s_2)<0$, $S_1(s_2)$
achieves a local strict maximum and the right of $(S_1(s_2),S_2(s_2))$ is subsonic.
Each case contradicts Lemma \ref{lemma-aaa1}.
Therefore, $(S_1(s),S_2(s))$ is not exceptional for each $s\in(0,1)$.
Furthermore, Proposition \ref{soniccurve3} implies
that $S_1$ is a strictly increasing function on $[0,1]$.

\hskip40mm
\setlength{\unitlength}{0.6mm}
\begin{picture}(250,45)

\put(0,0){\qbezier(-45,0)(0,0)(45,0)}

\put(-20,0){\cbezier(-5,0)(0,2)(8,8)(5,12)}
\put(-20,0){\cbezier(10,27)(5,25)(-2,18)(5,12)}
\put(-20,0){\cbezier(10,27)(15,30)(25,32)(40,0)}

\put(-14.5,10.5){\circle*{1.5}}
\put(-14.5,10.5){\qbezier[12](0,-6)(0,0)(0,6)}

\put(-18,18){\circle*{1.5}}
\put(-18,18){\qbezier[12](0,-6)(0,0)(0,6)}

\put(-18,4.6){\vector(1,1){1}}
\put(-14,24.7){\vector(3,2){1}}

\put(-13,2){supersonic}
\put(-43,8){subsonic}
\put(-13,30){subsonic}
\put(15,12){subsonic}

\put(-13.5,9){$s=s_0$}
\put(-37,16){$s=s_1$}

\put(-52,-14){the left of $(S_1(s_0),S_2(s_0))$ is subsonic}

\put(120,0){\qbezier(-45,0)(0,0)(45,0)}

\put(135,0){\cbezier(5,0)(0,2)(-8,8)(-5,12)}
\put(135,0){\cbezier(-10,27)(-5,25)(2,18)(-5,12)}
\put(135,0){\cbezier(-10,27)(-15,30)(-25,32)(-40,0)}

\put(129.5,10.5){\circle*{1.5}}
\put(129.5,10.5){\qbezier[12](0,-6)(0,0)(0,6)}

\put(133,18){\circle*{1.5}}
\put(133,18){\qbezier[12](0,-6)(0,0)(0,6)}

\put(133,4.6){\vector(3,-2){1}}
\put(129,24.7){\vector(3,-2){1}}

\put(99,2){supersonic}
\put(133,8){subsonic}
\put(110,30){subsonic}
\put(75,12){subsonic}

\put(110.5,9){$s=s_0$}
\put(134,16){$s=s_2$}

\put(70,-14){the right of $(S_1(s_0),S_2(s_0))$ is subsonic}

\end{picture}
\vskip15mm

Let us show that $(S_1(0),S_2(0))$ is not exceptional by contradiction.
Otherwise, $S'_1(0)\ge0$ and Lemma \ref{soniccurveform2} yield that
$S$ near $(S_1(0),S_2(0))$ is a graph of function with respect to $\psi$.
Then, Propositions \ref{soniccurve3}, \ref{lemma-aaa3} imply
that $S_1(s)\le S_1(0)$ for sufficiently small $s>0$,
which contradicts that $S_1$ is strictly increasing on $[0,1]$.
Similarly, one can prove that $(S_1(1),S_2(1))$ is not exceptional either.
$\hfill\Box$\vskip 4mm

Theorems \ref{soniclocation} and \ref{soniclocationtaylor} can be described in the physical plane.
Assume that the nozzle in the physical plane is of the following general form
\begin{align*}
%\label{soniccurvef}
{\mathcal G}=\big\{(x,y):f_1(x)<y<f_2(x),l_1< x< l_2\big\},
\end{align*}
where $l_1<l_2$ and
$$
f_1(x)<f_2(x),\quad l_1\le x\le l_2.
$$

The counter part of Theorem \ref{soniclocation} is the following

\begin{theorem}
\label{soniclocationphysical}
For any $C^2$ transonic flow of Meyer type in the nozzle ${\mathcal G}$,
its sonic curve ${\mathcal S}$ is a disjoint union of three connected parts (may be empty):
${\mathcal S}_{e}$, ${\mathcal S}_{+}$, ${\mathcal S}_{-}$,
where ${\mathcal S}_{e}$ is the set of exceptional points,
which is a closed line segment,
while ${\mathcal S}_{+}$ and ${\mathcal S}_{-}$ denote the other two connected parts
approaching the upper and lower walls, respectively.
Moreover, if $S_{e}$ is empty, then $S=S_{+}$ or $S=S_{-}$.
\end{theorem}

Similarly, in the physical plan, Theorem \ref{soniclocationtaylor} becomes

\begin{theorem}
\label{soniclocationtaylorphysical}
For any $C^2$ transonic flow of Taylor type in the nozzle ${\mathcal G}$,
each point on its sonic curve is not exceptional.
\end{theorem}

\hskip35mm
\setlength{\unitlength}{0.6mm}
\begin{picture}(250,65)
\put(-55,0){\vector(1,0){110}}
\put(0,40){\vector(0,1){17}}
\put(51,-4){$x$} \put(-4,53){$y$}
\put(-40,5){${\mathcal G}$}

\put(-60,48){\cbezier(14,2)(30,-5)(40,-7)(60,-8)}
\put(60,48){\cbezier(-14,2)(-30,-5)(-40,-7)(-60,-8)}

\put(-32,50){\cbezier(11,-8)(14,-16)(16,-19)(18.5,-20.8)}
\put(-32,50){\qbezier(18.5,-20.8)(20,-22)(21,-20)}
\put(-32,50){\qbezier(25,-15)(22,-16)(21,-20)}
\put(-32,50){\qbezier(25,-15)(28,-16)(30,-19)}
\put(-32,50){\qbezier(36,-32)(35,-25)(30,-19)}
\put(-32,50){\qbezier(38,-45)(37,-38.5)(36,-32)}
\put(-30,50){\qbezier(36,-45)(34,-57)(21,-63.8)}

\put(6,5){\circle*{1.5}}

\put(4,18){\circle*{1.5}}
\put(-30,16){subsonic}
\put(10,16){supersonic}
\put(-9.5,-9.5){${\mathcal S}_-$}
\put(-16,34){${\mathcal S}_+$}
\put(-2,9){${\mathcal S}_e$}

\put(10,0){\qbezier(-56,-25)(-40,-18)(-20,-14)}
\put(10,0){\qbezier(-20,-14)(10,-10)(20,-13)}
\put(10,0){\qbezier(36,-20)(30,-15)(20,-13)}

\put(-1.3,30.15){\vector(-1,1){1}}
\put(5.9,8){\vector(-1,4){1}}
\put(1.42,-5.2){\vector(1,1){1}}

\put(-36,-36){Figure of Theorem \ref{soniclocationphysical}}

\put(75,0){\vector(1,0){110}}
\put(130,40){\vector(0,1){17}}
\put(181,-4){$x$} \put(126,53){$y$}
\put(140,0){\qbezier(-56,-25)(-40,-18)(-20,-14)}
\put(140,0){\qbezier(-20,-14)(10,-10)(20,-13)}
\put(140,0){\qbezier(36,-20)(30,-15)(20,-13)}

\put(70,48){\cbezier(14,2)(30,-5)(40,-7)(60,-8)}
\put(190,48){\cbezier(-14,2)(-30,-5)(-40,-7)(-60,-8)}

\put(112,40){\cbezier(-8,3)(-2,-16)(4,-12)(6,-9)}
\put(112,40){\cbezier(14,-11)(12,-8)(8,-5)(6,-9)}
\put(112,40){\cbezier(14,-11)(15,-13)(20,-16)(21,-12)}
\put(112,40){\cbezier(31,-12)(26,-6)(23,-8)(21,-12)}
\put(112,40){\cbezier(31,-12)(32,-14)(35,-16)(38,-11)}
\put(112,40){\cbezier(38,-11)(40,-8)(42,-3)(44,3)}

\put(112,-9){\qbezier(6,6)(0,5)(-2,-7.5)}
\put(112,-9){\cbezier(22,6)(18,-3)(12,10)(6,6)}
\put(112,-9){\qbezier(22,6)(26,12)(29,10)}
\put(112,-9){\qbezier(38,-3)(34,8)(29,10)}

\put(104,43.3){\circle*{1.5}}
\put(156,43.3){\circle*{1.5}}
\put(110,-16.4){\circle*{1.5}}
\put(150,-11.7){\circle*{1.5}}

\put(152,32.5){\vector(-1,-2){1}}
\put(114,-5.1){\vector(1,1){1}}

\put(118,15){subsonic}
\put(114,36){supersonic}
\put(114,-11){supersonic}
\put(90,5){${\mathcal G}$}

\put(94,-36){Figure of Theorem \ref{soniclocationtaylorphysical}}

\end{picture}
\vskip30mm

\subsection{Characteristics from sonic points}

Characteristics play an important role in understanding
the behavior of supersonic flows near sonic curves.
So we first study characteristics from sonic points.

As will been shown in $\S$ 3.3, in the supersonic region,
it holds that
\begin{gather}
\label{nov-3}
\pd{W}{\psi}(\varphi,\psi)+\beta(\varphi,\psi)\pd{W}{\varphi}(\varphi,\psi)
+\omega(\varphi,\psi) W(\varphi,\psi)=0,
\\
\label{nov-4}
\pd{Z}{\psi}(\varphi,\psi)-\beta(\varphi,\psi)\pd{Z}{\varphi}(\varphi,\psi)
+\omega(\varphi,\psi) Z(\varphi,\psi)=0,
\end{gather}
where
\begin{align*}
\beta(\varphi,\psi)=&
\Big(\frac{-A'(q(\varphi,\psi))}{B'(q(\varphi,\psi))}\Big)^{1/2}
=\Big(\frac{\gamma+1}{2}q^2(\varphi,\psi)-1\Big)^{1/2}
\Big(1-\frac{\gamma-1}{2}q^2(\varphi,\psi)\Big)^{-1/(\gamma-1)-1/2},
\\
\omega(\varphi,\psi)=&\frac{\gamma+1}{2}q^4(\varphi,\psi)
\Big(\frac{\gamma+1}{2}q^2(\varphi,\psi)-1\Big)^{-2}
\Big(1-\frac{\gamma-1}{2}q^2(\varphi,\psi)\Big)^{1/(\gamma-1)}
A'(q(\varphi,\psi))\pd{q}{\psi}(\varphi,\psi),
\\
W(\varphi,\psi)=&A'(q(\varphi,\psi))\pd{q}{\varphi}(\varphi,\psi)
-\frac{A'(q(\varphi,\psi))}{\beta(\varphi,\psi)}\pd{q}{\psi}(\varphi,\psi),
\\
Z(\varphi,\psi)=&-A'(q(\varphi,\psi))\pd{q}{\varphi}(\varphi,\psi)
-\frac{A'(q(\varphi,\psi))}{\beta(\varphi,\psi)}\pd{q}{\psi}(\varphi,\psi).
\end{align*}
Note that $\beta$ is positive in the supersonic region and vanishes at the sonic curve.
The positive and negative characteristics of the system \eqref{nov-3}, \eqref{nov-4} are governed by
$$
\Phi'_{+}(\psi)=\beta(\Phi'_{+}(\psi),\psi)
\quad\mbox{and}\quad
\Phi'_{-}(\psi)=-\beta(\Phi'_{-}(\psi),\psi),
$$
respectively.

The following a priori estimates hold along characteristics.

\begin{lemma}
\label{lemma-sc4}
Assume that $\Sigma_{+}\,(\Sigma_{-})$ is a positive (negative) characteristic.
Then, $\theta(\Phi_{\pm}(\psi),\psi)\mp H(q(\Phi_{\pm}(\psi),\psi))$
is invariant on $\Sigma_{\pm}$, where
$$
H(q)=\int_{c_*}^q\big(-A'(s)B'(s)\big)^{1/2}ds,\quad c_*\le q<q_{max}.
$$
\end{lemma}

\Proof
On $\Sigma_{\pm}$, it holds that
\begin{align*}
\frac{d}{d\psi}\theta(\Phi_{\pm}(\psi),\psi)
=&\pm\beta(\Phi_{\pm}(\psi),\psi)\pd{\theta}{\varphi}(\Phi_{\pm}(\psi),\psi)
+\pd{\theta}{\psi}(\Phi_{\pm}(\psi),\psi)
\\
=&\pm\big(-A'(q(\Phi_{\pm}(\psi),\psi))B'(q(\Phi_{\pm}(\psi),\psi))\big)^{1/2}
\pd{q}{\psi}(\Phi_{\pm}(\psi),\psi)
\\
&\qquad-A'(q(\Phi_{\pm}(\psi),\psi))
\pd{q}{\varphi}(\Phi_{\pm}(\psi),\psi)
\\
=&\pm\Big(\pd{H(q)}{\psi}(\Phi_{\pm}(\psi),\psi)
\pm\beta(\Phi_{\pm}(\psi),\psi)\pd{H(q)}{\varphi}(\Phi_{\pm}(\psi),\psi)\Big)
\\
=&\pm\frac{d}{d\psi}H(q(\Phi_{\pm}(\psi),\psi)).
\end{align*}
So, $\theta(\Phi_{\pm}(\psi),\psi)\mp H(q(\Phi_{\pm}(\psi),\psi))$
is invariant on $\Sigma_{\pm}$.
$\hfill\Box$\vskip 4mm

Consider a $C^2$ transonic flow of Meyer type in ${\mathcal G}$
whose sonic curve ${\mathcal S}$ intersects the upper wall
at $(x_1,f_1(x_1))$.
Assume that the subsonic region is located on the left of ${\mathcal S}$,
${\mathcal S}={\mathcal S}_{+}\cup {\mathcal S}_{e}$,
${\mathcal S}_{+}$ is not empty and ${\mathcal S}_{e}$ is a line segment with nonempty interior.
In the potential plane, the lower wall corresponds to $\psi=0$
and the sonic curve is given by
$$
S: \varphi=S_1(s),\,\psi=S_2(s),\,(S'_1(s))^2+(S'_2(s))^2>0,\quad 0\le s\le1,\quad S_1,S_2\in C^2([0,1])
$$
with $S_2(0)=0$ and $S_2(1)=m$.
Set
$$
s^*=\sup\big\{s:S'_1(s)=0,0\le s\le1\big\}.
$$
Then, as proved in Theorem \ref{soniclocation}, $0<s^*<1$,
\begin{align}
\label{nov-1}
S'_1(s)=0,\quad\pd{q}{\psi}(S_1(s),S_2(s))=0,
\quad0\le s\le s^*,
\\
\label{nov-2}
S'_1(s)<0,\quad\pd{q}{\psi}(S_1(s),S_2(s))>0,
\quad s^*<s\le 1.
\end{align}

Let us first study characteristics from $(S_1(s),S_2(s))$ for $0\le s<s^*$.

\begin{lemma}
\label{nov-lemma000}
Assume that $0\le s<s^*$. Then the positive and negative characteristics from $(S_1(s),S_2(s))$
coincide and are given by
\begin{align}
\label{nov-lemma000-1}
\Phi_{+}(\psi)=\Phi_{-}(\psi)=S_1(s),\quad 0\le\psi\le S_2(s^*).
\end{align}
\end{lemma}

\Proof
It is not hard to check that the functions given by \eqref{nov-lemma000-1}
are positive and negative characteristics from $(S_1(s),S_2(s))$, respectively.
Let us verify the uniqueness of the positive characteristic.
Assume that $\tilde\Phi$ is a positive characteristic from $(S_1(s),S_2(s))$.
Since $\tilde\Phi$ is a nondecreasing function, it suffices to
verify that $\tilde\Phi$ coincide with the functions given by \eqref{nov-lemma000-1}
in a right neighborhood of $S_2(s)$.
Owing to \eqref{eq2or} and \eqref{nov-1}, one gets that
\begin{align*}
A''(c_*)\Big(\pd{q}{\varphi}(S_1(s),\psi)\Big)^2
+B'(c_*)\pd{^2q}{\psi^2}(S_1(s),\psi)
+B''(c_*)\Big(\pd{q}{\psi}(S_1(s),\psi)\Big)^2=0,
\quad 0\le\psi\le S_2(s^*)
\end{align*}
and
\begin{align*}
\pd{q}{\psi}(S_1(s),\psi)=0,\quad\pd{^2q}{\psi^2}(S_1(s),\psi)=0,
\quad 0\le\psi\le S_2(s^*),
\end{align*}
which imply
\begin{align}
\label{nov-lemma0-0}
\pd{q}{\varphi}(S_1(s),\psi)=0,\quad 0\le\psi\le S_2(s^*).
\end{align}
Fix a small positive number $\delta$ such that
$(S_1(s),S_1(s)+\delta]\times[S_2(s),S_2(s)+\delta]$ belongs to the supersonic region.
Due to \eqref{nov-lemma0-0}, there exists $M>0$ such that
\begin{align}
\label{nov-lemma0-1}
0\le\beta(\varphi,\psi)\le M(\varphi-S_1(s)),\quad
S_1(s)\le\varphi\le S_1(s)+\delta,\, S_2(s)\le\psi\le S_2(s)+\delta.
\end{align}
Choose a number $\varepsilon\in(0,\delta)$ such that $\tilde\Phi(S_2(s)+\varepsilon)\le S_1(s)+\delta$.
Owing to \eqref{nov-lemma0-1}, $\tilde\Phi$ satisfies
\begin{align}
\label{nov-lemma0-2}
0\le\tilde\Phi'(\psi)\le M(\tilde\Phi(\psi)-S_1(s)),
\quad S_2(s)\le\psi\le S_2(s)+\varepsilon;\qquad \tilde\Phi(S_2(s))=S_1(s).
\end{align}
The unique solution to \eqref{nov-lemma0-2} is
\begin{align*}
\tilde\Phi(\psi)=S_1(s),\quad S_2(s)\le\psi\le S_2(s)+\varepsilon.
\end{align*}
The uniqueness of the negative characteristic can be proved similarly.
$\hfill\Box$\vskip 4mm

We now turn to characteristics from $(S_1(s),S_2(s))$ for $s^*<s<1$.

\begin{lemma}
\label{nov-lemma1}
Assume that $s^*<s<1$. Then, there exist uniquely
a positive and a negative characteristics from $(S_1(s),S_2(s))$,
which contain no other sonic point except $(S_1(s),S_2(s))$.
\end{lemma}

\Proof
According to \eqref{nov-2}, the sonic curve from $(S_1(s^*),S_2(s^*))$ to $(S_1(1),S_2(1))$ is
a graph of function with respect to $\varphi$, which is denoted by
$$
\psi=S(\varphi),\quad S_1(1)\le\varphi\le S_1(s^*)
$$
with
\begin{align}
\label{nov-lemma1-1000}
S'(\varphi)=-\pd{q}{\varphi}(\varphi,S(\varphi))\Big(\pd{q}{\psi}(\varphi,S(\varphi))\Big)^{-1}
\in(-\infty,+\infty),\quad S_1(1)\le\varphi<S_1(s^*).
\end{align}

First we show that $(S_1(s),S_2(s))$ is the unique intersecting point between
any characteristic from $(S_1(s),S_2(s))$ and the sonic curve.
Assume that $\Phi$ is a characteristic from $(S_1(s),S_2(s))$,
which is defined in a right neighborhood of $S_2(s)$.
Owing to \eqref{nov-lemma1-1000} and $\Phi'(S_2(s))=0$,
the path of $\Phi$ leaves for the supersonic region after $\psi=S_2(s)$ locally.
Therefore, $\Phi$ is strictly monotone near $\psi=S_2(s)$
and its inverse function is denoted by $\Psi$, which solves
\begin{align*}
\Psi'(\varphi)=\frac1{\beta(\varphi,\Psi(\varphi))},\quad\varphi>S_1(s),
\qquad\Psi(S_1(s))=S_2(s)
\end{align*}
if $\Phi$ is a positive characteristic, while
\begin{align*}
\Psi'(\varphi)=-\frac1{\beta(\varphi,\Psi(\varphi))},\quad\varphi<S_1(s),
\qquad\Psi(S_1(s))=S_2(s)
\end{align*}
if $\Phi$ is a negative characteristic.
Since $\Psi'=\infty$ at sonic points, it follows from \eqref{nov-lemma1-1000} that
the path of $\Psi$ never approaches $\{(S_1(\tilde s),S_2(\tilde s)):s^*<\tilde s\le1\}$
after it leaves $(S_1(s),S_2(s))$.
Let $s_0\in[0,(s+s^*)/2]$ satisfy
$$
S_2(s_0)=\max\Big\{S_2(\tilde s):0\le\tilde s\le\frac{s+s^*}2\Big\}.
$$
Then, $s_0>s^*$ since ${\mathcal S}\in C^2$ and ${\mathcal S}_{e}$ is a segment
parallel to the $\psi$-axis.
Owing to the monotonicity of $\Psi$,
the path of $\Psi$ never approaches
$\{(S_1(\tilde s),S_2(\tilde s)):0\le\tilde s\le s_0\}$
if $\Phi$ is a positive characteristic, while
never approaches $\{(S_1(\tilde s),S_2(\tilde s)):0\le\tilde s<s\}$
if $\Phi$ is a negative characteristic.
Summing up, $(S_1(s),S_2(s))$ is the unique intersecting point between
the path of $\Psi$ and the sonic curve.

Below we prove the local existence and the uniqueness
of negative characteristics under the additional assumption $S'(S_1(s))<0$,
and the other cases can be discussed similarly.
Since $q(S_1(s),S_2(s))=c_*$, there exists a positive number $\tau_0$
such that
$$
G_0=\{(\varphi,\psi):S_1(s)-\tau_0<\varphi<S_1(s),
S(\varphi)<\psi<S(\varphi)+\tau_0\}
$$
belongs to the supersonic region and
\begin{align}
\label{nov-lemma1-1}
\frac1{\beta(\varphi,\psi)}\ge1-\min_{[S_1(s)-\tau_0,S_1(s)]}S',\quad (\varphi,\psi)\in G_0.
\end{align}
For $0<\varepsilon\le\tau_0/2$, $q(S_1(s),S_2(s)+\varepsilon)>c_*$.
Therefore, there exists uniquely a solution $\Psi_{\varepsilon}$
to the following problem
\begin{align}
\label{nov-lemma1-2}
\Psi'_{\varepsilon}(\varphi)=-\frac1{\beta(\varphi,\Psi_{\varepsilon}(\varphi))},
\quad\varphi<S_1(s);
\qquad\Psi_{\varepsilon}(S_1(s))=S_2(s)+\varepsilon.
\end{align}
Let $[\varphi_\varepsilon,S_1(s)]$ be the maximum interval of existence
for the solution to the problem \eqref{nov-lemma1-2} on $\overline G_0$.
It follows from \eqref{nov-lemma1-1} that
\begin{align}
\label{nov-lemma1-3}
S(\varphi)+(S_1(s)-\varphi)\le
S_2(s)+\varepsilon+\Big(1-\min_{[S_1(s)-\tau_0,S_1(s)]}S'\Big)(S_1(s)-\varphi)\le\Psi_{\varepsilon}(\varphi),
\quad \varphi_\varepsilon\le\varphi\le S_1(s),
\end{align}
which shows that the path of $\Psi_{\varepsilon}$ does not intersect $\partial G_0\cap S$.
Therefore, for $0<\varepsilon_1<\varepsilon_2\le\tau_0/2$,
\begin{align}
\label{nov-lemma1-4}
\varphi_{\varepsilon_1}\le\varphi_{\varepsilon_2}\le\varphi_{\tau_0/2}
\quad\mbox{and}\quad
\Psi_{\varepsilon_1}(\varphi)<\Psi_{\varepsilon_2}(\varphi)\mbox{ for each }
\varphi_{\tau_0/2}\le\varphi\le S_1(s).
\end{align}
Let
$$
\Psi(\varphi)=\lim_{\varepsilon\to0^+}\Psi_{\varepsilon}(\varphi),
\quad \varphi_{\tau_0/2}\le\varphi\le S_1(s).
$$
Then, it follows from \eqref{nov-lemma1-2}--\eqref{nov-lemma1-4} that $\Psi$
solves the following problem
\begin{align}
\label{nov-lemma1-2-001}
\Psi'(\varphi)=-\frac1{\beta(\varphi,\Psi(\varphi))},
\quad\varphi_{\tau_0/2}<\varphi<S_1(s);
\qquad\Psi(S_1(s))=S_2(s)
\end{align}
and satisfies
\begin{align}
\label{nov-lemma1-2-002}
S(\varphi)+(S_1(s)-\varphi)\le\Psi(\varphi)<\Psi_{\tau_0/2}(\varphi),
\quad\varphi_{\tau_0/2}\le\varphi\le S_1(s).
\end{align}
Note that \eqref{nov-lemma1-2-001} and \eqref{nov-lemma1-2-002} imply
$$
\Psi'(\varphi)<0,\quad\varphi_{\tau_0/2}\le\varphi<S_1(s).
$$
Thus $\Psi$ is strictly decreasing and its inverse function is a
negative characteristic from $(S_1(s),S_2(s))$ due to \eqref{nov-lemma1-2-001}.

Let us turn to the uniqueness.
Assume that $\Phi$ and $\tilde\Phi$ are two negative characteristics from $(S_1(s),S_2(s))$,
which are defined in a right neighborhood of $S_2(s)$.
Since $\pd{q}{\psi}(S_1(s),S_2(s))>0$ from \eqref{nov-2}, there exists a positive number 
$\tau_1\,(\le\tau_0)$
such that
\begin{align}
\label{nov-lemma1-5}
\frac12\pd{q}{\psi}(S_1(s),S_2(s))\le\pd{q}{\psi}(\varphi,\psi)\le\frac32\pd{q}{\psi}(S_1(s),S_2(s)),
\quad(\varphi,\psi)\in G_1,
\end{align}
where
$$
G_1=\{(\varphi,\psi):
S_1(s)-\tau_1<\varphi<S_1(s),
S(\varphi)<\psi<S(\varphi)+\tau_1\}.
$$
Since $\Phi'(S_2(s))=\tilde\Phi'(S_2(s))=0$,
there exists a positive number $\delta$ so small that
\begin{align}
\label{nov-lemma1-6}
(\Phi(\psi),\psi),(\tilde\Phi(\psi),\psi)\in G_1,\quad S_2(s)<\psi<S_2(s)+\delta
\end{align}
and
\begin{gather}
\label{nov-lemma1-7-0}
0\le\sup_{G_1}\Big|\pd q\varphi\Big|(S_1(s)-\Phi(\psi))
\le\frac{1}{4}\pd{q}{\psi}(S_1(s),S_2(s))(\psi-S_2(s)),
\quad S_2(s)\le\psi\le S_2(s)+\delta,
\\
\label{nov-lemma1-7}
0\le\sup_{G_1}\Big|\pd q\varphi\Big|(S_1(s)-\tilde\Phi(\psi))
\le\frac{1}{4}\pd{q}{\psi}(S_1(s),S_2(s))(\psi-S_2(s)),
\quad S_2(s)\le\psi\le S_2(s)+\delta.
\end{gather}
Set
$$
h(\psi)=(\Phi(\psi)-\tilde\Phi(\psi))^2,\quad S_2(s)\le\psi\le S_2(s)+\delta.
$$
Then, $h$ solves
\begin{align}
\label{nov-lemma1-8}
h'(\psi)=\xi(\psi)(q(\zeta(\psi),\psi)-c_*)^{-1/2}h(\psi),
\quad S_2(s)<\psi<S_2(s)+\delta;\qquad h(S_2(s))=0,
\end{align}
where $\xi\in C([S_2(s),S_2(s)+\delta])$, while
\begin{align}
\label{nov-lemma1-9}
\min\{\Phi(\psi),\tilde\Phi(\psi)\}\le\zeta(\psi)\le\max\{\Phi(\psi),\tilde\Phi(\psi)\},
\quad S_2(s)<\psi<S_2(s)+\delta.
\end{align}
It follows from \eqref{nov-lemma1-5}--\eqref{nov-lemma1-7} and \eqref{nov-lemma1-9} that
\begin{align}
\label{nov-lemma1-10}
q(\zeta(\psi),\psi)-c_*\le&\sup_{G_1}\Big|\pd q\varphi\Big|(S_1(s)-\zeta(\psi))
+\frac{3}{2}\pd{q}{\psi}(S_1(s),S_2(s))(\psi-S_2(s))
\nonumber
\\
\le&\frac{7}{4}\pd{q}{\psi}(S_1(s),S_2(s))(\psi-S_2(s)),
\quad S_2(s)<\psi<S_2(s)+\delta,
\end{align}
\begin{align}
\label{nov-lemma1-11}
q(\zeta(\psi),\psi)-c_*\ge&-\sup_{G_1}\Big|\pd q\varphi\Big|(S_1(s)-\zeta(\psi))
+\frac{1}{2}\pd{q}{\psi}(S_1(s),S_2(s))(\psi-S_2(s))
\nonumber
\\
\ge&\frac{1}{4}\pd{q}{\psi}(S_1(s),S_2(s))(\psi-S_2(s)),
\quad S_2(s)<\psi<S_2(s)+\delta.
\end{align}
Substituting \eqref{nov-lemma1-10} and \eqref{nov-lemma1-11} into \eqref{nov-lemma1-8}
yields
\begin{align}
\label{nov-lemma1-12}
|h'(\psi)|\le M(\psi-S_2(s))^{-1/2}h(\psi),
\quad S_2(s)<\psi<S_2(s)+\delta;\qquad h(S_2(s))=0
\end{align}
for some positive constant $M$. Since \eqref{nov-lemma1-12} admits only a trivial solution,
the uniqueness is proved.
$\hfill\Box$\vskip 4mm

For fixed $s\in(s^*,1)$, according to Lemma \ref{nov-lemma1},
there exist uniquely a positive and a negative characteristics from $(S_1(s),S_2(s))$;
furthermore, for each characteristic,
the lower endpoint is $(S_1(s),S_2(s))$, the upper endpoint is located at the upper wall
or the outlet of the nozzle.
Without loss of generality, the nozzle is assumed to be so long that
all upper endpoints of the characteristics from $(S_1(s),S_2(s))$ for $s\in(s^*,1)$
are located at the upper wall.
Therefore, the positive and the negative characteristics from $(S_1(s),S_2(s))$ are governed by
$$
\Sigma_{s,+}:\varphi=\Phi_{s,+}(\psi):
\left\{
\begin{aligned}
&\frac{d}{d\psi}\Phi_{s,+}(\psi)=\beta(\Psi_{s,+}(\psi),\psi),
\quad S_2(s)<\psi<m,
\\
&\Phi_{s,+}(S_2(s))=S_1(s),\,\Phi_{s,+}(m)=\varphi_{s,+}
\end{aligned}
\right.
$$
and
$$
\Sigma_{s,-}:\varphi=\Phi_{s,-}(\psi):
\left\{
\begin{aligned}
&\frac{d}{d\psi}\Phi_{s,-}(\psi)=-\beta(\Psi_{s,-}(\psi),\psi),
\quad S_2(s)<\psi<m,
\\
&\Phi_{s,-}(S_2(s))=S_1(s),\,\Phi_{s,-}(m)=\varphi_{s,-},
\end{aligned}
\right.
$$
respectively, where $s\in(s^*,1)$.

For $s^*<s_1<s_2<1$,
the uniqueness of the characteristic shows
$$
\varphi_{s_1,+}>\varphi_{s_2,+},\quad \varphi_{s_1,-}>\varphi_{s_2,-}
$$
and
\begin{gather*}
\Phi_{s_1,+}(\psi)>\Phi_{s_2,+}(\psi),\quad
\Phi_{s_1,-}(\psi)>\Phi_{s_2,-}(\psi),\quad S_2(s_2)\le\psi\le m.
\end{gather*}
Set
\begin{align}
\label{varphipm}
\varphi_-=\lim_{s\to (s^*)^+}\varphi_{s,-},\quad \varphi_+=\lim_{s\to (s^*)^+}\varphi_{s,+}.
\end{align}

\begin{remark}
\label{nov-lemma1-re}
It is not clear at this moment
whether the positive and negative characteristics from $(S_1(s^*),S_2(s^*))$
are unique. However,
$$
\left\{
\begin{aligned}
&S_1(s^*),\quad&&0\le\psi\le S_2(s^*),
\\
&\lim_{s\to (s^*)^+}\Phi_{s,+}(\psi),\quad &&S_2(s^*)<\psi\le m
\end{aligned}
\right.
\quad\mbox{ and }\quad
\left\{
\begin{aligned}
&S_1(s^*),\quad&&0\le\psi\le S_2(s^*),
\\
&\lim_{s\to (s^*)^+}\Phi_{s,-}(\psi),\quad &&S_2(s^*)<\psi\le m
\end{aligned}
\right.
$$
are the maximal positive and the minimal negative characteristics from $(S_1(s^*),S_2(s^*))$,
respectively.
\end{remark}

\hskip70mm
\setlength{\unitlength}{0.6mm}
\begin{picture}(250,65)
\put(-45,0){\vector(1,0){90}}
\put(0,50){\vector(0,1){15}}
\put(0,50){\qbezier(-45,0)(0,0)(45,0)}

\put(-32,50){\qbezier(6,-14)(0.6,-6)(-2,0)}
\put(-32,50){\cbezier(6,-14)(8.7,-16.5)(11,-16.8)(12,-17)}
\put(-32,50){\qbezier(18,-14)(15,-16.5)(12,-17)}
\put(-32,50){\qbezier(18,-14)(19,-13.5)(19.5,-12)}
\put(-32,50){\qbezier(19.5,-12)(21,-10.4)(22,-8)}
\put(-32,50){\cbezier(28,-12)(26,-8)(23.5,-6)(22,-8)}
\put(-32,50){\cbezier(28,-12)(30,-14)(33,-22)(36,-32)}
\put(-32,50){\qbezier(36,-50)(36,-40)(36,-32)}

\put(-32,50){\qbezier(6,-14)(5.8,-12)(4,-5)}
\put(-32,50){\qbezier(1,0)(3,-2)(4,-5)}
\put(-32,50){\qbezier(6,-14)(6.2,-12)(8,-5)}
\put(-32,50){\qbezier(11,0)(9,-2)(8,-5)}
\put(-32,50){\qbezier[9](6,-14)(6,-9.5)(6,-5)}

\put(-32,50){\qbezier(21,-10)(20.8,-5)(19,-0)}
\put(-32,50){\qbezier(21,-10)(21.2,-5)(23,-0)}
\put(-32,50){\qbezier[7](21,-10)(21,-6.5)(21,-3)}

\put(4,18){\qbezier[16](0,0)(0,8)(0,16)}
\put(4,18){\qbezier(0,0)(-1,14)(-5,23)}
\put(4,18){\qbezier(-5,23)(-6,25)(-9,32)}
\put(-5,50){\circle*{1.5}}
\put(1,52){$m$}
\put(-10,-5){$\varphi_{-}$}
\put(-5,0){\qbezier[50](0,0)(0,25)(0,50)}

\put(4,18){\qbezier(0,0)(0.3,4)(2.3,12)}
\put(4,18){\qbezier(6,23)(4,18)(2.3,12)}
\put(4,18){\qbezier(6,23)(7,25)(12,32)}
\put(16,50){\circle*{1.5}}
\put(16,-5){$\varphi_{+}$}
\put(16,0){\qbezier[50](0,0)(0,25)(0,50)}

\put(4,18){\circle*{1.5}}
\put(6,16){$(S_1(s^*),S_2(s^*))$}
\put(-30,6){subsonic}
\put(10,6){supersonic}
\put(-7,36){$S$}

\put(-30,29){\vector(1,1){17}}
\put(-30,29){\vector(1,4){3}}
\put(-30,29){\vector(4,1){30}}
\put(-35,25){$\Sigma_-$}

\put(-9,23){\vector(0,1){25}}
\put(-9,23){\vector(-2,3){15}}
\put(-9,23){\vector(1,1){19}}
\put(-13,18){$\Sigma_+$}

\put(-7,63){$\psi$}
\put(40,-5){$\varphi$}
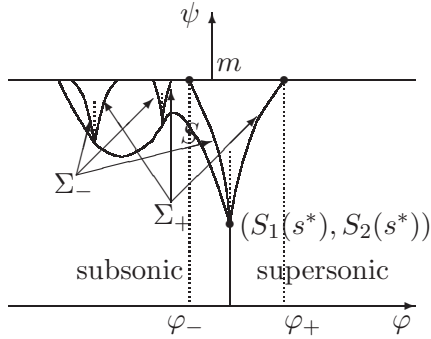
\put(-95,-14){Figure: characteristics for the transonic flow in the potential plane}
\end{picture}
\vskip15mm

Let $x^*$ be the point such that $(x^*,f_1(x^*))$ is transformed to
$(\varphi_{-},m)$ in the coordinates transformation
from the physical plane to the potential plane.
That is to say,
\begin{align}
\label{xupper}
x^*=\sup\big\{x\in[l_1,l_2]:(x,f_1(x)) \mbox{ is the endpoint of the negative characteric
from a point on }{\mathcal S}_+\big\}.
\end{align}

\begin{proposition}
\label{lemma-sc3}
For each $x\in[x_1,x^*]$, $f''_1(x)>0$.
\end{proposition}

\Proof
Fix $s\in(s^*,1)$. It follows from the definition of $W$ that
\begin{align*}
\frac{W(\varphi,\Psi_{s,+}(\varphi))}{A'(q(\varphi,\Psi_{s,+}(\varphi)))}=
\pd{q}{\varphi}(\varphi,\Psi_{s,+}(\varphi))
-\frac1{\beta(\varphi,\Psi_{s,+}(\varphi))}\pd{q}{\psi}(\varphi,\Psi_{s,+}(\varphi),
\quad S_1(s)<\varphi\le\varphi_{s,+}.
\end{align*}
Since $\beta(S_1(s),S_2(s))=0$ and $\pd{q}{\psi}(S_1(s),S_2(s))>0$,
there exists a positive number $\delta\,(\le\varphi_{s,+}-S_1(s))$ such that
$W(S_1(s)+\delta,\Psi_{s,+}(S_1(s)+\delta))>0$. Then, \eqref{nov-3} yields
\begin{align}
\label{nov-5}
W(\varphi_{s,+},m)>0.
\end{align}
Similarly, one can get
\begin{align}
\label{nov-6}
Z(\varphi_{s,-},m)>0.
\end{align}
Due to to the arbitrariness of $s\in(s^*,1)$, it follows from \eqref{nov-5} and \eqref{nov-6} that
$$
W(\varphi,m)>0,\quad S_1(1)<\varphi<\varphi_{+},\quad
W(\varphi_{+},m)\ge0
$$
and
$$
Z(\varphi,m)>0,\quad S_1(1)<\varphi<\varphi_{-},\quad
Z(\varphi_{-},m)\ge0.
$$
Thus,
\begin{align}
\label{lemma-sc3-2-1}
\pd{q}{\psi}(\varphi,m)
=-\frac{\beta(\varphi,m)}{2A'(q(\varphi,m))}
(W(\varphi,m)+Z(\varphi,m))>0,\quad S_1(1)<\varphi\le\varphi_{-}.
\end{align}
Moreover, \eqref{nov-2} gives
\begin{align}
\label{lemma-sc3-2-2}
\pd{q}{\psi}(S_1(1),m)>0.
\end{align}
Assume that $(x,f_1(x))$ is transformed into $(\Phi(x),m)$ for $l_1<x<l_2$ in the coordinates transformation
from the physical plane to the potential plane.
Then,
\begin{align}
\label{lemma-sc3-2-3-000}
\Phi'(x)=\frac{q(x,f_1(x))}{\cos\theta(x,f_1(x))}
=q(x,f_1(x))\sqrt{1+(f'_1(x))^2}>0,
\quad l_1<x<l_2
\end{align}
and
\begin{align}
\label{lemma-sc3-2-3}
\Phi'(x)B'(q(\Phi(x),m))\pd{q}{\psi}(\Phi(x),m)=\frac{d\theta(x,f_1(x))}{dx}
=\frac{f''_1(x)}{1+(f'_1(x))^2},\quad l_1<x<l_2.
\end{align}
The conclusion of the proposition follows from \eqref{lemma-sc3-2-1}--\eqref{lemma-sc3-2-3}.
$\hfill\Box$\vskip 4mm

\hskip35mm
\setlength{\unitlength}{0.6mm}
\begin{picture}(250,65)
\put(-55,0){\vector(1,0){110}}
\put(0,40){\vector(0,1){17}}
\put(51,-4){$x$} \put(-4,53){$y$}
\put(-40,5){${\mathcal G}$}

\put(-60,48){\cbezier(14,2)(30,-5)(40,-7)(60,-8)}
\put(60,48){\cbezier(-14,2)(-30,-5)(-40,-7)(-60,-8)}

\put(4,46){$y=f_1(x)$}

\put(-32,50){\cbezier(11,-8)(14,-16)(16,-19)(18.5,-20.8)}
\put(-32,50){\qbezier(18.5,-20.8)(20,-22)(21,-20)}
\put(-32,50){\qbezier(25,-15)(22,-16)(21,-20)}
\put(-32,50){\cbezier(25,-15)(30,-16)(35,-24)(36,-32)}
\put(-32,50){\qbezier(39,-62)(37.5,-47)(36,-32)}

\put(4,18){\qbezier[20](0,0)(-2,14)(-8,23)}
\put(-4,40.5){\circle*{1.5}}
\put(-6,-5){$x^*$}
\put(-4,0){\qbezier[40](0,0)(0,20)(0,40)}
\put(-21,42.3){\circle*{1.5}}
\put(-24,-5){$x_1$}
\put(-21,0){\qbezier[43](0,0)(0,21.5)(0,43)}

\put(4,18){\circle*{1.5}}
\put(-2,34){$\Sigma_-$}
\put(-30,12){subsonic}
\put(10,12){supersonic}

\put(10,0){\qbezier(-56,-25)(-40,-18)(-20,-14)}
\put(10,0){\qbezier(-20,-14)(10,-10)(20,-13)}
\put(10,0){\qbezier(36,-20)(30,-15)(20,-13)}

\put(-11,-18){$y=f_2(x)$}
\put(-16,34){${\mathcal S}_+$}
\put(-2,4){${\mathcal S}_e$}

\put(-36,-36){Figure of Proposition \ref{lemma-sc3}}

\put(75,0){\vector(1,0){110}}
\put(130,40){\vector(0,1){17}}
\put(181,-4){$x$} \put(126,53){$y$}
\put(90,5){${\mathcal G}$}

\put(70,48){\cbezier(14,2)(30,-5)(40,-7)(60,-8)}
\put(190,48){\cbezier(-14,2)(-30,-5)(-40,-7)(-60,-8)}

\put(134,46){$y=f_1(x)$}

\put(140,-12){\qbezier(0,0)(-3,24)(-13,36)}
\put(140,-12){\qbezier(-23,36)(-18,44)(-13,36)}
\put(140,-12){\qbezier(-23,36)(-27,30)(-31,36)}
\put(140,-12){\qbezier(-39,56)(-35,40)(-31,36)}

\put(140,-12){\qbezier[22](0,0)(-1,24)(-4,32)}
\put(140,-12){\qbezier[18](-15,52.5)(-9,48)(-4,32)}

\put(101,44){\circle*{1.5}}
\put(122,-5){$x^*$}
\put(125,0){\qbezier[40](0,0)(0,20)(0,40)}
\put(125,40.6){\circle*{1.5}}
\put(99,-5){$x_1$}
\put(101,0){\qbezier[42](0,0)(0,21)(0,42)}

\put(130,34){$\Sigma_-$}
\put(100,12){subsonic}
\put(128,24){supersonic}

\put(140,0){\qbezier(-56,-25)(-40,-18)(-20,-14)}
\put(140,0){\qbezier(-20,-14)(10,-10)(20,-13)}
\put(140,0){\qbezier(36,-20)(30,-15)(20,-13)}

\put(119,-18){$y=f_2(x)$}
\put(109,25){${\mathcal S}_+$}

\put(94,-36){Figure of Remark \ref{lemma-sc3-re}}

\end{picture}
\vskip30mm

\begin{remark}
\label{lemma-sc3-re}
Proposition \ref{lemma-sc3} still holds if ${\mathcal S}_e$ is empty or a single point set.
\end{remark}

Similar to Lemmas \ref{nov-lemma000}, \ref{nov-lemma1}, 
Proposition \ref{lemma-sc3} and Remarks \ref{nov-lemma1-re}, \ref{lemma-sc3-re}, one can prove
the following two general conclusions. One is on characteristics from sonic points,
and the other is on the geometry of walls.

\hskip35mm
\setlength{\unitlength}{0.6mm}
\begin{picture}(250,65)
\put(-55,0){\vector(1,0){110}}
\put(0,40){\vector(0,1){17}}
\put(51,-4){$x$} \put(-4,53){$y$}
\put(-40,5){${\mathcal G}$}

\put(-60,48){\cbezier(14,2)(30,-5)(40,-7)(60,-8)}
\put(60,48){\cbezier(-14,2)(-30,-5)(-40,-7)(-60,-8)}

\put(4,46){$y=f_1(x)$}

\put(-32,50){\cbezier(19,-25)(18,-24)(15,-20)(10,-8)}
\put(-32,50){\cbezier(28,-22)(26,-22.3)(23,-30)(19,-25)}
\put(-32,50){\cbezier(36,-32)(35,-28)(30,-20)(28,-22)}

\put(-32,50){\qbezier(38,-45)(37,-38.5)(36,-32)}

\put(-30,50){\cbezier(36,-45)(35,-55)(30,-62)(26,-63.3)}

\put(-30,50){\qbezier(36,-45)(35.6,-62.5)(33,-62.7)}

\put(6,4){\qbezier[10](0,0)(1,-6.5)(1.5,-9.75)}
\put(6,4){\cbezier(0,0)(1.5,-7)(3.5,-12)(6,-16)}

\put(4,18){\cbezier(0,0)(-3,12)(-8,18)(-15,22.8)}
\put(4,18){\qbezier[10](0,0)(-1,6.5)(-1.5,9.75)}
\put(4,18){\cbezier(0,0)(-1,10)(2,20)(4,22.6)}

\put(-8,24.5){\qbezier(0,0)(-3,12)(-10,17)}
\put(-8,24.5){\qbezier[10](0,0)(-1,6.5)(-1.5,9.75)}
\put(-8,24.5){\cbezier(0,0)(-0.1,4)(0,12)(2,15.7)}

\put(3,-6){\cbezier(0,0)(-0.5,-5)(-2,-6)(-3,-6.6)}
\put(3,-6){\qbezier[10](0,0)(0.5,-3.25)(1,-6.5)}
\put(3,-6){\cbezier(0,0)(1.5,-4.5)(2,-5)(4,-6.2)}

\put(6,4){\circle*{1.5}}

\put(4,18){\circle*{1.5}}
\put(-30,16){subsonic}
\put(10,16){supersonic}
\put(-8,-9.5){${\mathcal S}_-$}
\put(-22,26){${\mathcal S}_+$}
\put(-2,9){${\mathcal S}_e$}

\put(10,0){\qbezier(-56,-25)(-40,-18)(-20,-14)}
\put(10,0){\qbezier(-20,-14)(10,-10)(20,-13)}
\put(10,0){\qbezier(36,-20)(30,-15)(20,-13)}

\put(-11,-18){$y=f_2(x)$}

\put(-28,3){$\Sigma_-$}
\put(-22,5){\vector(1,3){10}}
\put(-22,5){\vector(2,3){19}}
\put(-22,5){\vector(3,-2){23.5}}
\put(-22,5){\vector(3,-1){27}}

\put(23,3){$\Sigma_+$}
\put(22,5){\vector(-1,2){16}}
\put(22,5){\vector(-1,1){29.5}}
\put(22,5){\vector(-1,-1){16}}
\put(22,5){\vector(-2,-3){10.5}}

\put(-36,-36){Figure of Theorem \ref{thm-aaa}}

\put(75,0){\vector(1,0){110}}
\put(130,40){\vector(0,1){17}}
\put(181,-4){$x$} \put(-4,53){$y$}
\put(90,5){${\mathcal G}$}

\put(70,48){\cbezier(14,2)(30,-5)(40,-7)(60,-8)}
\put(190,48){\cbezier(-14,2)(-30,-5)(-40,-7)(-60,-8)}

\put(134,46){$y=f_1(x)$}

\put(98,50){\cbezier(19,-25)(15,-22)(11,-12)(9.8,-8)}
\put(98,50){\cbezier(28,-22)(26,-22.3)(23,-30)(19,-25)}
\put(98,50){\cbezier(36,-32)(35,-28)(30,-20)(28,-22)}

\put(98,50){\qbezier(38,-45)(37,-38.5)(36,-32)}
\put(100,50){\cbezier(36,-45)(35,-55)(30,-62)(26,-63.3)}

\put(100,50){\qbezier(36,-45)(35.6,-62.5)(32,-62.7)}
\put(136,4){\circle*{1.5}}
\put(133,-12.3){\qbezier[12](0,0)(0,6.15)(0,12.3)}
\put(129,1.5){$x_*$}
\put(127,-12.6){\qbezier[12](0,0)(0,6.3)(0,12.6)}
\put(122,1.5){$x_2$}
\put(127,-12.9){\circle*{1.5}}
\put(133,-12.5){\circle*{1.5}}

\put(134,18){\cbezier(0,0)(-3,12)(-8,18)(-15,23)}
\put(119,41.1){\circle*{1.5}}
\put(114,-5){$x^*$}
\put(119,0){\qbezier[40](0,0)(0,20)(0,40)}
\put(107.5,42.4){\circle*{1.5}}
\put(103,-5){$x_1$}
\put(107.5,0){\qbezier[42](0,0)(0,21)(0,42)}

\put(134,18){\circle*{1.5}}
\put(128,34){$\Sigma_-$}
\put(134.5,-9){$\Sigma_+$}
\put(100,16){subsonic}
\put(140,16){supersonic}
\put(122,-9.5){${\mathcal S}_-$}
\put(108,26){${\mathcal S}_+$}
\put(128,9){${\mathcal S}_e$}

\put(140,0){\qbezier(-56,-25)(-40,-18)(-20,-14)}
\put(140,0){\qbezier(-20,-14)(10,-10)(20,-13)}
\put(140,0){\qbezier(36,-20)(30,-15)(20,-13)}

\put(119,-18){$y=f_2(x)$}

\put(94,-36){Figure of Theorem \ref{lemma-sc3-negative}}

\end{picture}
\vskip30mm

\begin{theorem}
\label{thm-aaa}
Let ${\mathcal S}={\mathcal S}_+\cup {\mathcal S}_e\cup {\mathcal S}_-$
be the sonic curve of a $C^2$ transonic flow of Meyer type in ${\mathcal G}$.

{\rm (i)} There exist uniquely a positive and a negative characteristics
from each interior point of ${\mathcal S}_+\cup {\mathcal S}_-$.

{\rm (ii)} All positive and negative characteristics
from points of ${\mathcal S}_e\setminus\partial({\mathcal S}_+\cup {\mathcal S}_-)$
are the same line ${\mathcal S}_e$.

{\rm (iii)} There are positive and negative characteristics
from points of ${\mathcal S}_e\cap\partial({\mathcal S}_+\cup {\mathcal S}_-)$.
Particularly, the maximal positive and the minimal negative characteristics
are unique.
\end{theorem}

Characteristics from boundary sonic points are as follows.

\begin{remark}
Let ${\mathcal S}={\mathcal S}_+\cup {\mathcal S}_e\cup {\mathcal S}_-$
be the sonic curve of a $C^2$ transonic flow of Meyer type in ${\mathcal G}$.
Assume that the subsonic region is located on the left.
For the sonic point on the upper wall,

{\rm (i)} If it belongs to ${\mathcal S}_+$, then there is not any characteristic
from this point.

{\rm (ii)} If it belongs to ${\mathcal S}_e\setminus\partial{\mathcal S}_-$,
then the positive and negative characteristics
from this point are the same line ${\mathcal S}_e$.

{\rm (iii)} If it belongs to ${\mathcal S}_e\cap\partial{\mathcal S}_-$,
then there are positive and negative characteristics
from this point. Particularly, the maximal positive and the minimal negative characteristics
are unique.

{\rm (iv)} If it belongs to ${\mathcal S}_-$,
then there exist uniquely a positive and a negative characteristics
from this point.

For the sonic point on the lower wall,

{\rm (i)} If it belongs to ${\mathcal S}_-$, then there is not any characteristic
from this point.

{\rm (ii)} If it belongs to ${\mathcal S}_e\setminus\partial{\mathcal S}_+$,
then the positive and negative characteristics
from this point are the same line ${\mathcal S}_e$.

{\rm (iii)} If it belongs to ${\mathcal S}_e\cap\partial{\mathcal S}_+$,
then there are positive and negative characteristics
from this point. Particularly, the maximal positive and the minimal negative characteristics
are unique.

{\rm (iv)} If it belongs to ${\mathcal S}_+$,
then there exist uniquely a positive and a negative characteristics
from this point.
\end{remark}

\begin{theorem}
\label{lemma-sc3-negative}
Let ${\mathcal S}={\mathcal S}_+\cup {\mathcal S}_e\cup {\mathcal S}_-$
be the sonic curve of a $C^2$ transonic flow of Meyer type in ${\mathcal G}$.
Assume that
the subsonic region is located on the left and
${\mathcal S}$ intersects the upper and lower walls
at $(x_1,f_1(x_1))$ and $(x_2,f_2(x_2))$,
respectively.

{\rm (i)} If ${\mathcal S}_+\not=\emptyset$, then
$f''_1(x)>0$ for each $x\in[x_1,x^*]$, where
$x^*$ is given by \eqref{xupper}.

{\rm (ii)} If ${\mathcal S}_-\not=\emptyset$, then
$f''_2(x)<0$ for each $x\in[x_2,x_*]$, where
\begin{align}
\label{xlower}
x_*=\sup\big\{x\in[l_1,l_2]:(x,f_2(x)) \mbox{ is the endpoint of the positive characteric
from a point on }{\mathcal S}_-\big\}.
\end{align}
\end{theorem}

As a consequence of Theorem \ref{soniclocationtaylorphysical}
and the proof of Proposition \ref{lemma-sc3},
we can obtain

\begin{theorem}
\label{transonictaylortype}
Consider a $C^2$ transonic flow of Taylor type in ${\mathcal G}$.
Assume that
the sonic curve intersects the upper wall at $(x_{1,-},f_1(x_{1,-}))$ and $(x_{1,+},f_1(x_{1,+}))$,
while the lower wall at $(x_{2,-},f_2(x_{2,-}))$ and $(x_{2,+},f_2(x_{2,+}))$
with $x_{k,-}<x_{k,+}\,(k=1,2)$. Then,
$$
f''_1(x)>0,\quad x_{1,-}\le x\le x_{1,+}
$$
and
$$
f''_2(x)<0,\quad x_{2,-}\le x\le x_{2,+}.
$$
\end{theorem}

\setlength{\unitlength}{0.6mm}
\begin{picture}(250,65)

\put(75,0){\vector(1,0){110}}
\put(130,40){\vector(0,1){17}}
\put(181,-4){$x$} \put(126,53){$y$}
\put(140,0){\qbezier(-56,-25)(-40,-18)(-20,-14)}
\put(140,0){\qbezier(-20,-14)(10,-10)(20,-13)}
\put(140,0){\qbezier(36,-20)(30,-15)(20,-13)}

\put(70,48){\cbezier(14,2)(30,-5)(40,-7)(60,-8)}
\put(190,48){\cbezier(-14,2)(-30,-5)(-40,-7)(-60,-8)}

\put(112,40){\cbezier(-8,3)(-2,-16)(4,-12)(6,-9)}
\put(112,40){\cbezier(14,-11)(12,-8)(8,-5)(6,-9)}
\put(112,40){\cbezier(14,-11)(15,-13)(20,-16)(21,-12)}
\put(112,40){\cbezier(31,-12)(26,-6)(23,-8)(21,-12)}
\put(112,40){\cbezier(31,-12)(32,-14)(35,-16)(38,-11)}
\put(112,40){\cbezier(38,-11)(40,-8)(42,-3)(44,3)}

\put(112,-9){\qbezier(6,6)(0,5)(-2,-7.5)}
\put(112,-9){\cbezier(22,6)(18,-3)(12,10)(6,6)}
\put(112,-9){\qbezier(22,6)(26,12)(29,10)}
\put(112,-9){\qbezier(38,-3)(34,8)(29,10)}

\put(104,43.3){\circle*{1.5}}
\put(156,43.3){\circle*{1.5}}
\put(110,-16.1){\circle*{1.5}}
\put(150,-11.6){\circle*{1.5}}

\put(96,-5){$x_{1,-}$}
\put(104,0){\qbezier[44](0,0)(0,22)(0,44)}
\put(154,-5){$x_{1,+}$}
\put(156,0){\qbezier[44](0,0)(0,22)(0,44)}
\put(108,1.5){$x_{2,-}$}
\put(110,0){\qbezier[16](0,0)(0,-8)(0,-16)}
\put(142,1.5){$x_{2,+}$}
\put(150,0){\qbezier[12](0,0)(0,-6)(0,-12)}

\put(118,15){subsonic}
\put(114,36){supersonic}
\put(114,-11){supersonic}
\put(90,5){${\mathcal G}$}
\put(134,46){$y=f_1(x)$}
\put(119,-18){$y=f_2(x)$}

\put(94,-36){Figure of Theorem \ref{transonictaylortype}}

\end{picture}
\vskip30mm

\begin{remark}
For a smooth transonic flow of Taylor type, it has been mentioned in \cite{Bers} that
the boundary of a supersonic enclosure on a wall cannot
contain a straight segment.
Theorem \ref{transonictaylortype} strengthen it to be that
the curvature of this boundary must be nonzero.
\end{remark}

\subsection{Instability of transonic flows with nonexceptional points}

According to Theorems \ref{soniclocation} and \ref{lemma-sc3-negative},
we are ready to show that $C^2$ transonic flows with nonexceptional points
are unstable with respect to
small changes in the shape of the nozzle.

\hskip75mm
\setlength{\unitlength}{0.6mm}
\begin{picture}(250,65)
\put(-55,0){\vector(1,0){110}}
\put(0,40){\vector(0,1){17}}
\put(51,-4){$x$} \put(-4,53){$y$}
\put(-40,5){${\mathcal G}$}

\put(-60,48){\cbezier(14,2)(30,-5)(40,-7)(60,-8)}
\put(60,48){\cbezier(-14,2)(-30,-5)(-40,-7)(-60,-8)}

\put(-60,43){\qbezier[30](14,5)(20,1)(40,-3)}
\put(-60,43){\qbezier[26](40,-3)(50,-5)(60,-5.5)}

\put(60,43){\qbezier[30](-14,5)(-20,1)(-40,-3)}
\put(60,43){\qbezier[26](-40,-3)(-50,-5)(-60,-5.5)}

\put(14,52){$y=f_1(x)$}

\put(-38,52){$\tilde f''_1(\hat x_1)=0$}
\put(-15,39.3){\circle*{1.5}}
\put(-26,50){\vector(1,-1){10}}

\put(12,34.5){$y=\tilde f_1(x)$}

\put(-32,50){\cbezier(19,-25)(18,-24)(15,-20)(10,-8)}
\put(-32,50){\cbezier(28,-22)(26,-22.3)(23,-30)(19,-25)}
\put(-32,50){\cbezier(36,-32)(35,-28)(30,-20)(28,-22)}

\put(-32,50){\qbezier(38,-45)(37,-38.5)(36,-32)}
\put(-30,50){\cbezier(36,-45)(35,-55)(30,-62)(26,-63)}

\put(-30,50){\qbezier(36,-45)(35.6,-62.5)(32,-62.7)}
\put(6,4){\circle*{1.5}}
\put(3,-12.3){\qbezier[12](0,0)(0,6.15)(0,12.3)}
\put(-1,1.5){$x_*$}
\put(-3,-12.6){\qbezier[12](0,0)(0,6.3)(0,12.6)}
\put(-8,1.5){$x_2$}
\put(-3,-12.9){\circle*{1.5}}
\put(3,-12.3){\circle*{1.5}}

\put(4,18){\cbezier(0,0)(-3,12)(-8,18)(-15,22.7)}
\put(-11,41.1){\circle*{1.5}}
\put(-16,-5){$x^*$}
\put(-11,0){\qbezier[40](0,0)(0,20)(0,40)}
\put(-22,42.6){\circle*{1.5}}
\put(-27,-5){$x_1$}
\put(-22,0){\qbezier[42](0,0)(0,21)(0,42)}

\put(4,18){\circle*{1.5}}
\put(0,30){$\Sigma_-$}
\put(5.5,-9){$\Sigma_+$}
\put(-30,16){subsonic}
\put(10,16){supersonic}
\put(-10,-9.5){${\mathcal S}_-$}
\put(-22,26){${\mathcal S}_+$}
\put(-2,9){${\mathcal S}_e$}

\put(10,0){\qbezier(-56,-25)(-40,-18)(-20,-14)}
\put(10,0){\qbezier(-20,-14)(10,-10)(20,-13)}
\put(10,0){\qbezier(36,-20)(30,-15)(20,-13)}

\put(10,-2){\qbezier[44](-56,-25.5)(-40,-18)(-20,-14)}
\put(10,-2){\qbezier[52](-20,-14)(10,-10)(20,-13)}
\put(10,-2){\qbezier[22](36,-20.5)(30,-15)(20,-13)}

\put(16,-25){$y=\tilde f_2(x)$}

\put(28,-10){$y=f_2(x)$}

\put(-20,-28){$\tilde f''_2(\hat x_2)=0$}
\put(-1,-15){\circle*{1.5}}
\put(-8,-22){\vector(1,1){6}}

\put(-35,-40){Figure of Theorem \ref{soniclocation2}}

\end{picture}
\vskip30mm

\begin{theorem}
\label{soniclocation2}
Consider a $C^2(\overline{\mathcal G})$ transonic flow of Meyer type,
whose subsonic region is located on the left.
Assume that the sonic curve ${\mathcal S}={\mathcal S}_+\cup{\mathcal S}_e\cup{\mathcal S}_-$
intersects the upper and lower walls
at $(x_1,f_1(x_1))$ and $(x_2,f_2(x_2))$, respectively.
If there is a nonexceptional point on ${\mathcal S}$, then $\varphi$
is unstable with respect to small changes in the shape of the nozzle
in the following sense.

{\rm (i)} The case ${\mathcal S}_+\not=\emptyset$.
Let $x^*$ be given by \eqref{xupper}.
If $\tilde f_1\in C^\infty([l_1,l_2])$ is a $C^1$ small perturbation of $f_1$
and satisfies $\tilde f''_1(\hat x_1)=0$ with some $x_1<\hat x_1<x^*$,
then there is not a $C^2(\tilde{\mathcal G})$ transonic flow of Meyer type
with $\tilde{\mathcal G}$ being the corresponding domain when the upper wall is replaced by $\tilde f_1$,
which is a $C^0$ small perturbation of the background flow.

{\rm (ii)} The case ${\mathcal S}_-\not=\emptyset$.
Let $x_*$ be given by \eqref{xlower}.
If $\tilde f_2\in C^\infty([l_1,l_2])$ is a $C^1$ small perturbation of $f_2$
and satisfies $\tilde f''_2(\hat x_2)=0$ with some $x_2<\hat x_2<x_*$,
then there is not a $C^2(\tilde{\mathcal G})$ transonic flow of Meyer type
with $\tilde{\mathcal G}$ being the corresponding domain when the lower wall is replaced by $\tilde f_2$,
which is a $C^0$ small perturbation of the background flow.
\end{theorem}

\Proof
We prove (i) only and the proof of (ii) is similar.
Consider a $C^2(\tilde{\mathcal G})$ transonic flow of Meyer type, whose
sonic curve is denoted by $\tilde{\mathcal S}=\tilde{\mathcal S}_+
\cup\tilde{\mathcal S}_e\cup\tilde{\mathcal S}_-$
according to Theorem \ref{soniclocation}.
If the transonic flow is close enough to the background flow in $C^0$-norm,
then $\tilde{\mathcal S}_+\not=\emptyset$ and its two endpoints are small perturbations of
the two endpoints of ${\mathcal S}_+$, respectively. Thus,
$\tilde x_1$ and $\tilde x^*$ are small perturbations of $x_1$ and $x^*$, respectively,
and it follows from Theorem \ref{lemma-sc3-negative} that
$$
\tilde f''_1(x)>0,\quad\tilde x_1<x<\tilde x^*,
$$
where $\tilde x_1$ and $\tilde x^*$ are defined in a similar way as for $x_1$ and $x^*$, respectively.
This contradicts that $\tilde f''_1(\tilde x_1)=0$.
$\hfill\Box$\vskip 4mm

\begin{remark}
The instability in Theorem \ref{soniclocation2} is weak. It is unknown
whether the flow is unstable if the nozzle wall is perturbed in $C^2$
or other smooth spaces.
\end{remark}

\begin{remark}
\label{soniclocation2-re}
Theorem \ref{soniclocation2} motives us to seek a smooth transonic flow of Meyer type
whose sonic points are exceptional.
Indeed, we do prove the existence of such a smooth transonic flow
for the symmetric de Laval nozzle $\Omega$ in the paper.
\end{remark}

Moreover, for $C^2$ transonic flows of Meyer type in $\Omega$, one can prove

\begin{proposition}
\label{streamlineconvex}
Assume that $f''(0)=0$.
Consider a $C^2$ transonic flow of Meyer type in $\Omega$
whose each streamline in the supersonic region is a graph of function with respect to $x$.
If there is a nonexceptional point at the sonic curve,
then not all streamlines are convex.
\end{proposition}

\Proof
It is assumed that the subsonic region is located on the left without loss of generality.
We transform the flow from the physical plane
to the potential plane.
It is assumed that $(0,f(0))$ is transformed to $(0,m)$
and the sonic curve is transformed to
$$
S: \varphi=S_1(s),\,\psi=S_2(s),\, (S'_1(s))^2+(S'_2(s))^2>0,\quad 0\le s\le1
$$
with $S_2(0)=0$ and $S_2(1)=m$.
Set
$$
s^*=\sup\big\{s:S'_1(s)=0,0\le s\le1\big\}.
$$
The symmetry of the nozzle and Proposition \ref{soniccurve3} show that
$(S_1(0),S_2(0))$ is an exceptional point
and $\theta(S_1(0),S_2(0))=0$.
Similar to the discussion in $\S 2.5$, one can prove that $0\le s^*<1$
and
\begin{align*}
S'_1(s)=0,\quad\pd{q}{\psi}(S_1(s),S_2(s))=0,\quad\theta(S_1(s),S_2(s))=0,
\quad0\le s\le s^*,
\\
S'_1(s)<0,\quad\pd{q}{\psi}(S_1(s),S_2(s))>0,\quad\theta(S_1(s),S_2(s))<0,
\quad s^*<s\le 1.
\end{align*}
Lemma \ref{lemma-sc4} implies that
$$
\theta(\varphi_-,m)+H(q(\varphi_-,m))=0,\quad
\theta(\varphi_+,m)-H(q(\varphi_+,m))=0,
$$
where $\varphi_\pm$ are defined by \eqref{varphipm} and $H$ is given in Lemma \ref{lemma-sc4}.
Thanks to $q(\varphi_\pm,m)>c_*$, one gets that
$\theta(\varphi_-,m)<0<\theta(\varphi_+,m)$ and thus $\varphi_-<0<\varphi_+$.
As shown in Proposition \ref{lemma-sc3}, one can get $W(0,m)>0$ and thus $Z(0,m)=-W(0,m)<0$.
Denote $\Sigma_-$ to be the negative characteristic from $(0,m)$.
Due to Lemma \ref{lemma-sc4}, $\theta+H(q)$ is invariant on $\Sigma_-$,
which equals identically to $\theta(0,m)+H(q(0,m))>0$.
Note that
$$
\theta(S_1(s),S_2(s))+H(q(S_1(s),S_2(s)))=\theta(S_1(s),S_2(s))\le0,
\quad 0\le s\le 1.
$$
Thus $\Sigma_-$ never approach the sonic curve.
Assume that $\Sigma_-$ intersects the lower wall at $(\varphi_1,0)$.
Then $S_1(0)<\varphi_1$.
Since the flow is sonic at $(S_1(0),0)$ and supersonic at $(\varphi_1,0)$,
there exists $\varphi_2\in(S_1(0),\varphi_1)$ such that $W(\varphi_2,0)<0$.
Let $\Sigma_+$ be the positive characteristic from $(\varphi_2,0)$
and denote $(\varphi_*,\psi_*)$ to be the intersecting point of $\Sigma_+$ and $\Sigma_-$.
It follows from \eqref{nov-3} and \eqref{nov-4} that
$W(\varphi_*,\psi_*)<0$ and $Z(\varphi_*,\psi_*)<0$, which imply
$$
\pd{q}{\psi}(\varphi_*,\psi_*)
=-\frac{\beta(\varphi_*,\psi_*)}{2A'(q(\varphi_*,\psi_*))}
(W(\varphi_*,\psi_*)+Z(\varphi_*,\psi_*))<0.
$$
Assume that $y=f_*(x)$ is the streamline across $(\varphi_*,\psi_*)$
and $(x_*,f_*(x_*))$ corresponds to $(\varphi_*,\psi_*)$ in the coordinates transformation.
Then, as shown in \eqref{lemma-sc3-2-3-000} and \eqref{lemma-sc3-2-3},
$f''_*(x_*)$ and $\pd{q}{\psi}(\varphi_*,\psi_*)$ are of the same sign.
Thus $f''_*(x_*)<0$.
$\hfill\Box$\vskip 4mm

\hskip-10mm
\setlength{\unitlength}{0.6mm}
\begin{picture}(250,65)

\put(80,0){\vector(1,0){120}}
\put(140,40){\vector(0,1){17}}
\put(197,-4){$x$} \put(136,53){$y$}
\put(97,5){$\Omega$}

\put(80,48){\cbezier(14,2)(30,-5)(40,-7)(60,-8)}
\put(200,48){\cbezier(-14,2)(-30,-5)(-40,-7)(-60,-8)}

\put(108,50){\cbezier(13,-8.5)(16,-20)(19,-22)(20,-20)}
\put(108,50){\cbezier(30,-18)(25,-10)(22,-19)(20,-20)}
\put(108,50){\qbezier(36,-32)(35,-25)(30,-18)}
\put(108,50){\qbezier(36,-50)(36,-40)(36,-32)}

\put(140,40){\qbezier[8](0,0)(-3.6,-2)(-4,-5.6)}
\put(140,40){\qbezier[56](0,0)(15,-10)(20,-40)}
\put(165,43.4){\qbezier[56](0,-0.6)(-11,-10)(-16,-43.4)}

\put(110,12){subsonic}
\put(160,12){supersonic}
\put(156,32){$\Sigma_+$}
\put(145,32){$\Sigma_-$}
\put(154,22){\circle*{1.5}}
\put(156,20){$(x_*,f_*(x_*))$}
\put(118,30){${\mathcal S}_+$}
\put(137,10){${\mathcal S}_e$}

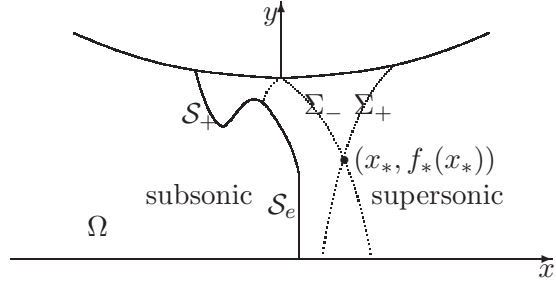
\put(90,-14){Figure of the proof of Proposition \ref{streamlineconvex}}

\end{picture}
\vskip30mm

\section{Formulation of the smooth transonic flow problem and main results}

As mentioned in Remark \ref{soniclocation2-re},
we seek a smooth transonic flow of Meyer type,
whose sonic points are exceptional, in the symmetric de Laval nozzle $\Omega$.
In this section, let us formulate this smooth transonic flow problem
both in the physical plane and in the potential plane,
and state the existence and uniqueness theorems.

\subsection{Formulation of the transonic flow problem in the de Laval nozzle}

To seek a smooth transonic flow in $\Omega$ whose sonic points are exceptional,
the first step is to determine the location of the sonic curve.
It is not hard to verify from Lemma \ref{lemmasoniccuver} and
Proposition \ref{soniccurve1} that

\begin{proposition}
\label{soniccurve2}
For a $C^2(\overline\Omega)$ transonic flow of Meyer type,
the following statements are equivalent

{\rm(i)} The endpoint of the sonic curve on the wall is located at the throat point;

{\rm(ii)} The flow angle always equals to zero at the sonic curve;

{\rm(iii)} Every point at the sonic curve is exceptional;

{\rm(iv)} The velocity vector is along the normal direction at the sonic curve;

{\rm(v)} The potential is a constant at the sonic curve;

{\rm(vi)} The sonic curve is located at the throat of the nozzle.
\end{proposition}

Therefore, the sonic curve of the sought smooth transonic flow should be
located at the throat. That is to say,
the flow should be subsonic in the convergent part, then sonic at the throat
and then supersonic in the divergent part.
For such a smooth transonic flow,
the velocity vector is along the normal direction at the sonic curve
and thus the mass flux is
\begin{align}
\label{mass-m}
m=f(0) c_*^{1+2/(\gamma-1)}.
\end{align}
Furthermore, the potential of the flow is a constant at the sonic curve and we normalize it to be zero
in the paper.
Then, the potential of the flow at the inlet, which is a constant, is free.
So is the potential at the outlet.

Thus, the problem of a transonic flow in $\Omega$,
whose velocity vector is along the normal direction at the inlet and the outlet,
which satisfies the slip condition on the wall
and whose sonic curve is located at the throat, can be formulated as
\begin{align}
\label{phytran1}
&\mbox{div}(\rho(|\nabla\varphi|^2)\nabla\varphi)=0,\quad&&(x,y)\in\Omega,
\\
\label{phytran2}
&\varphi(g(y),y)=C_{\text{\rm in}},\quad&&0<y<f(l_-),
\\
\label{phytran3}
&\pd{\varphi}y(x,0)=0,\quad&&g(0)<x<t(0),
\\
\label{phytran4}
&\pd{\varphi}y(x,f(x))-f'(x)\pd{\varphi}x(x,f(x))=0,\quad&&l_-<x<l_+,
\\
\label{phytran5}
&\varphi(t(y),y)=C_{\text{\rm out}},\quad&&0<y<f(l_+),
\\
\label{phytran6}
&|\nabla\varphi(0,y)|=c_*,\,\varphi(0,y)=0,\quad&&0<y<f(0),
\\
\label{phytran7}
&|\nabla\varphi(x,y)|<c_*,\quad&&(x,y)\in\Omega_-,
\\
\label{phytran8}
&|\nabla\varphi(x,y)|>c_*,\quad&&(x,y)\in\Omega_+,
\end{align}
where $C_{\text{\rm in}},C_{\text{\rm out}}\in\mathbb R$
are free, the outlet $\Gamma_{\text{\rm out}}:x=t(y)\,(0\le y\le f(l_+))$ is free,
$\Omega_-$ and $\Omega_+$ are the convergent part and the divergent part of the nozzle, respectively,
i.e.
$$
\Omega_-=\big\{(x,y)\in\Omega:x<0\big\},\quad\Omega_+=\big\{(x,y)\in\Omega:x>0\big\}.
$$
Since the sonic curve of the smooth transonic flow is located at the throat,
we can decompose the transonic flow problem \eqref{phytran1}--\eqref{phytran8} into
a subsonic-sonic flow problem and a sonic-supersonic flow problem as follows
\begin{align}
\label{phytran9}
&\mbox{div}(\rho(|\nabla\varphi|^2)\nabla\varphi)=0,\quad&&(x,y)\in\Omega_-,
\\
\label{phytran10}
&\varphi(g(y),y)=C_{\text{\rm in}},\quad&&0<y<f_-(l_-),
\\
\label{phytran11}
&\pd{\varphi}y(x,0)=0,\quad&&g(0)<x<0,
\\
\label{phytran12}
&\pd{\varphi}y(x,f_-(x))-f'_-(x)\pd{\varphi}x(x,f_-(x))=0,\quad&&l_-<x<0,
\\
\label{phytran13}
&|\nabla\varphi(0,y)|=c_*,\,\varphi(0,y)=0,\quad&&0<y<f_-(0),
\\
\label{phytran14}
&|\nabla\varphi(x,y)|<c_*,\quad&&(x,y)\in\Omega_-
\end{align}
and
\begin{align}
\label{phytran15}
&\mbox{div}(\rho(|\nabla\varphi|^2)\nabla\varphi)=0,\quad&&(x,y)\in\Omega_+,
\\
\label{phytran16}
&|\nabla\varphi(0,y)|=c_*,\,\varphi(0,y)=0,\quad&&0<y<f_+(0),
\\
\label{phytran17}
&\pd{\varphi}y(x,0)=0,\quad&&0<x<t(0),
\\
\label{phytran18}
&\pd{\varphi}y(x,f_+(x))-f'_+(x)\pd{\varphi}x(x,f_+(x))=0,\quad&&0<x<l_+,
\\
\label{phytran19}
&\varphi(t(y),y)=C_{\text{\rm out}},\quad&&0<y<f_+(l_+),
\\
\label{phytran20} 
&|\nabla\varphi(x,y)|>c_*,\quad&&(x,y)\in\Omega_+,
\end{align}
where $f_-$ and $f_+$ are the walls of the convergent part and divergent part of the nozzle, respectively, i.e.
\begin{gather*}
f_-(x)=f(x),\quad l_-\le x\le0
\quad\mbox{ and }\quad
f_+(x)=f(x),\quad 0\le x\le l_+.
\end{gather*}

\begin{remark}
For the subsonic-sonic flow problem \eqref{phytran9}--\eqref{phytran14},
there are two boundary conditions at the sonic curve (see \eqref{phytran13}).
Indeed, the potential at the sonic curve is a free constant,
which is normalized to be zero and thus the potential at the inlet is free
(see \eqref{phytran10}, where $C_{\text{\rm in}}$ is a free parameter).
\end{remark}

It should be noted that to study a subsonic-sonic flow problem and a sonic-supersonic flow problem,
one must determine and control precisely
the speed of the flow near the sonic state,
which plays a very essential role in the mathematical analysis.
In \eqref{phytran1}, the speed of the flow is the absolute value of the gradient
of a solution, which is very hard to estimate precisely.
It turns out to be more convenient to solve
the subsonic-sonic flow problem \eqref{phytran9}--\eqref{phytran14}
and the sonic-supersonic flow problem \eqref{phytran15}--\eqref{phytran20} in the potential plane,
where the speed of the flow is a solution to \eqref{eq2or}.
We will thus formulate these two problems in the potential plane in next two subsections.

\subsection{Formulation of the subsonic-sonic flow problem in the potential plane}

We will show that there exists a subsonic-sonic flow in the divergent part of the nozzle $\Omega_-$.
In this subsection, we formulate this subsonic-sonic flow problem in the potential plane.

Assume that $f_-\in C^{3}([l_-,0])\,(-1\le l_-<0)$ satisfies
\begin{align}
\label{s-2}
\delta_{1,-}(-x)^{\lambda_{-}}\le f''_-(x)\le\delta_{2,-}(-x)^{\lambda_{-}},
\quad l_-\le x\le0
\end{align}
with positive constants $\lambda_{-}$, $\delta_{1,-}$ and $\delta_{2,-}$
such that $\lambda_{-}>2$ and $\delta_{1,-}\le\delta_{2,-}$.
Choose the inlet $\Gamma_{\text{\rm in}}$ as a small perturbation
of the arc
$$
\Gamma_0:x=g_0(y),\quad 0\le y\le f_-(l_-),
$$
where
$$
g_0(y)=l_--\frac{f_-(l_-)}{f'_-(l_-)}-\sqrt{R_0^2-y^2},
\quad 0\le y\le f_-(l_-),
\qquad R_0=\frac{f_-(l_-)\sqrt{(f'_-(l_-))^2+1}}{-f'_-(l_-)}.
$$
More precisely,
$$
\Gamma_{\text{\rm in}}:x=g(y),\quad 0\le y\le f_-(l_-)
$$
with $g\in C^{3}([0,f_-(l_-)])$ satisfying
\begin{align}
\label{s-3}
g'(0)=0,\quad
g(f_-(l_-))=l_-,\quad
g'(f_-(l_-))=-f'_-(l_-)
\end{align}
and
\begin{align}
\label{s-4}
\frac{1}{2R_0}\le\Big(\frac{g'(y)}{\sqrt{(g'(y))^2+1}}\Big)'
\le\frac{3}{2R_0},
\quad\Big|\Big(\frac{g'(y)}{\sqrt{(g'(y))^2+1}}\Big)''\Big|
\le\epsilon(l_-)(-l_-)^{3\lambda_{-}/2},
\quad 0\le y\le f_-(l_-),
\end{align}
where $0<\epsilon(s)\le1$ for $-1\le s<0$ and $\lim_{s\to0^-}\epsilon(s)=0$.

Since the velocity vector of the flow is along the normal direction at the inlet
and the wall of the nozzle is solid,
the flow angle at the inlet and at the upper wall
can be expressed as
$$
\Theta_{\text{\rm in}}(y)=-\arctan g'(y),\quad 0\le y\le f_{-}(l_{-})
$$
and
$$
\Theta_{-}(x)=\arctan f'_-(x),\quad l_-\le x\le 0,
$$
respectively.
Let the speed of the flow at the inlet
and at the upper wall be
denoted by
$$
q(g(y),y)={\mathscr Q}_{\text{\rm in}}(y),\quad0\le y\le f_{-}(l_{-})
$$
and
$$
q(x,f(x))={\mathscr Q}_{-}(x),\quad l_{-}\le x\le 0,
$$
respectively.
The incoming mass flux is given by
\begin{align}
\label{ass-1}
m_{\text{\rm in}}=\int_{0}^{f_{-}(l_{-})}
\frac{{\mathscr Q}_{\text{\rm in}}(y)\rho({\mathscr Q}_{\text{\rm in}}^2(y))}
{\cos\Theta_{\text{\rm in}}(y)}dy.
\end{align}
Since the potential at the sonic curve is normalized to be zero, the
potential at the inlet is given by
\begin{align}
\label{ass-2}
\zeta_{-}=\int^{l_{-}}_0
\frac{{\mathscr Q}_{-}(x)}{\cos\Theta_{-}(x)}dx.
\end{align}
At the inlet, the stream function is
$$
\psi(g(y),y)=\Psi_{\text{\rm in}}(y),\quad0\le y\le f_{-}(l_{-})
$$
satisfying $\Psi_{\text{\rm in}}(0)=0$, $\Psi_{\text{\rm in}}(f_{-}(l_{-}))=m_{\text{\rm in}}$ and
$$
\Psi'_{\text{\rm in}}(y)=
\frac{{\mathscr Q}_{\text{\rm in}}(y)\rho({\mathscr Q}_{\text{\rm in}}^2(y))}
{\cos\Theta_{\text{\rm in}}(y)}>0,
\quad0\le y\le f_{-}(l_{-}).
$$
Thus, $\Psi_{\text{\rm in}}$ is expressed as
\begin{align}
\label{Psiin}
\Psi_{\text{\rm in}}(y)=\int_{0}^{y}
\frac{{\mathscr Q}_{\text{\rm in}}(s)\rho({\mathscr Q}_{\text{\rm in}}^2(s))}
{\cos\Theta_{\text{\rm in}}(s)}ds,\quad
0\le y\le f_{-}(l_{-}).
\end{align}
Denote by $Y_{\text{\rm in}}$ the inverse function of
$\Psi_{\text{\rm in}}$, i.e.
\begin{align}
\label{Yin}
Y_{\text{\rm in}}(\psi)=\Psi_{\text{\rm in}}^{-1}(\psi),\quad
0\le\psi\le m_{\text{\rm in}}.
\end{align}
At the upper wall, the potential function is
$$
\varphi(x,f_{-}(x))=\Phi_{-}(x),\quad l_{-}\le x\le 0
$$
satisfying $\Phi_{-}(l_{-})=\zeta_{-}$, $\Phi_{-}(0)=0$
and
$$
\Phi'_{-}(x)=
\frac{{\mathscr Q}_{-}(x)}{\cos\Theta_{-}(x)}>0, \quad l_{-}\le x\le0.
$$
Thus, $\Phi_{-}$ is expressed as
\begin{align}
\label{Phiub}
\Phi_{-}(x)=\int_{0}^{x}
\frac{{\mathscr Q}_{-}(s)}{\cos\Theta_{-}(s)}ds,\quad  l_{-}\le x\le0.
\end{align}
Denote by $X_{-}$ the inverse function of $\Phi_{-}$, i.e.
\begin{align}
\label{Xub}
X_{-}(\varphi)=\Phi_{-}^{-1}(\varphi),\quad
\zeta_-\le\varphi\le0.
\end{align}

Therefore, the subsonic-sonic flow problem in the potential plane
is formulated as follows
\begin{align}
\label{p-eq} &\pd{^2A(q)}{\varphi^2}+\pd{^2B(q)}{\psi^2}=0,
\quad&&(\varphi,\psi)\in(\zeta_{-},0)\times(0,m_{\text{\rm in}}),
\\
\label{p-inbc} &\pd{A(q)}{\varphi}(\zeta_-,\psi)=
-\frac{\Theta'_{\text{\rm in}}(y)\cos\Theta_{\text{\rm in}}(y)}
{{\mathscr Q}_{\text{\rm in}}(y)\rho({\mathscr Q}_{\text{\rm in}}^2(y))}
\Big|_{y=Y_{\text{\rm in}}(\psi)},\quad&&\psi\in(0,m_{\text{\rm in}}),
\\
\label{p-lbbc}
&\pd{q}{\psi}(\varphi,0)=0,\quad&&\varphi\in(\zeta_{-},0),
\\
\label{p-ubbc} &\pd{B(q)}{\psi}(\varphi,m_{\text{\rm in}})=
\frac{\Theta'_{-}(x)\cos\Theta_{-}(x)}{{\mathscr Q}_{-}(x)}\Big|_{x=X_{-}(\varphi)},
\quad&&\varphi\in(\zeta_{-},0),
\\
\label{p-outbc} &q(0,\psi)=c_*,\quad&&\psi\in(0,m_{\text{\rm in}}),
\\
\label{Qinq}
&{\mathscr Q}_{\text{\rm in}}(y)=q(0,\psi)\Big|_{\psi=\Psi_{\text{\rm in}}(y)},
\quad&&y\in[0,f_{-}(l_{-})],
\\
\label{Qubq}
&{\mathscr Q}_{-}(x)=q(\varphi,m_{\text{\rm in}})\Big|_{\varphi=\Phi_{-}(x)},
\quad&&x\in[l_{-},0].
\end{align}

\begin{remark}
A natural boundary condition at the inlet seems to be
\begin{align}
\label{p-reinbc}
\pd{A(q)}{\varphi}(\zeta_-,\psi)=
-\frac{\Theta'_{\text{\rm in}}(Y_{\text{\rm in}}(\psi))
\cos\Theta_{\text{\rm in}}(Y_{\text{\rm in}}(\psi))}
{q(\zeta_-,\psi)\rho(q^2(\zeta_-,\psi))},\quad\psi\in(0,m_{\text{\rm in}})
\end{align}
instead of \eqref{p-inbc}. However, since
$\frac{d}{dq}\Big(\frac1{q\rho(q^2)}\Big)<0$ in $q\in(0,c_*)$ and
$-\Theta'_{\text{\rm in}}(Y_{\text{\rm in}}(\psi))\cos\Theta_{\text{\rm in}}(Y_{\text{\rm in}}(\psi))
\ge1/(2R_0)$ from \eqref{s-4}, it seems difficult to obtain the well-posedness of
the problem \eqref{p-eq}, \eqref{p-reinbc}, \eqref{p-lbbc}--\eqref{p-outbc}.
Therefore, we prescribe the boundary condition at the inlet by \eqref{p-inbc}.
\end{remark}

The problem \eqref{p-eq}--\eqref{Qubq} is a boundary value problem
for a second-order quasilinear degenerate equation in a rectangle
with two free parameters.
Furthermore, the degeneracy is characteristic and
the boundary conditions \eqref{p-inbc} and \eqref{p-ubbc} are nonlinear, nonlocal and implicit.
As far as we know, there is no known theory to show the existence and uniqueness of solutions
even if the problem is uniformly elliptic.
Thus we  will use a fixed point argument
to prove the existence of solutions to the problem \eqref{p-eq}--\eqref{Qubq} as follows:
For given ${\mathscr Q}_{\text{\rm in}}$ and ${\mathscr Q}_{-}$,
we solve the problem \eqref{p-eq}--\eqref{p-outbc}, then define a mapping by
\eqref{Qinq} and \eqref{Qubq}.
To this end, one must determine and control precisely
the rate of the solution tending to $c_*$ in order to solve the problem \eqref{p-eq}--\eqref{Qubq}.
However, there is no a background solution to this problem to suggest what the rate is.
And we will carry out some  precise elliptic estimates in this paper
to determine and control this rate.

Since \eqref{p-eq} is degenerate, it is more convenient to 
introduce weak solutions to the problem \eqref{p-eq}--\eqref{p-outbc}.
After determining the rate of the solution tending to $c_*$,
one can establish the regularity of weak solutions.

\begin{definition}
\label{defsolution}
A function $q\in L^\infty((\zeta_{-},0)\times(0,m_{\text{\rm in}}))$ is said to be a
weak solution to the problem
\eqref{p-eq}--\eqref{p-outbc}, if
$$
0<\inf_{(\zeta_{-},0)\times(0,m_{\text{\rm in}})}q \le\sup_{(\zeta_{-},0)\times(0,m_{\text{\rm in}})}q\le
c_*
$$
and
\begin{align*}
&\int^{0}_{\zeta_-}\int_0^{m_{\text{\rm in}}}\Big(
A(q(\varphi,\psi))\pd{^2\xi}{\varphi^2}(\varphi,\psi)
+B(q(\varphi,\psi))\pd{^2\xi}{\psi^2}(\varphi,\psi)\Big)d\varphi
d\psi
\\
&\qquad
+\int_0^{m_{\text{\rm in}}}\frac{\Theta'_{\text{\rm in}}(y)\cos\Theta_{\text{\rm in}}(y)}
{{\mathscr Q}_{\text{\rm in}}(y)\rho({\mathscr Q}_{\text{\rm in}}^2(y))}
\Big|_{y=Y_{\text{\rm in}}(\psi)}\xi(\zeta_-,\psi)d\psi
+\int^0_{\zeta_-}\frac{\Theta'_{-}(x)\cos\Theta_{-}(x)}{{\mathscr Q}_{-}(x)}
\Big|_{x=X_{-}(\varphi)}\xi(\varphi,m_{\text{\rm in}})d\varphi=0
\end{align*}
for any $\xi\in C^2([\zeta_{-},0]\times[0,m_{\text{\rm in}}])$
with
$$
\pd{\xi}{\psi}(\cdot,0)\Big|_{(\zeta_{-},0)}
=\pd{\xi}{\psi}(\cdot,m_{\text{\rm in}})\Big|_{(\zeta_{-},0)}=0\quad\mbox{ and
}\quad\pd{\xi}{\varphi}(\zeta_-,\cdot)\Big|_{(0,m_{\text{\rm in}})}
=\xi(0,\cdot)\Big|_{(0,m_{\text{\rm in}})}=0.
$$
\end{definition}

Though the condition \eqref{s-2} seems quite strict, yet, 
it is almost necessary for the existence of $C^2$ solutions
as shown in the following remark.

\begin{remark}
\label{rem-aaa}
If the problem \eqref{p-eq}--\eqref{p-outbc} admits
a $C^2([\zeta_{-},0]\times[0,m_{\text{\rm in}}])$ solution,
then
\begin{align}
\label{rem-aaa-1}
f''_-(x)=O(x^2),\quad l_-<x<0.
\end{align}
\end{remark}

\Proof
It follows from \eqref{p-eq} and \eqref{p-outbc} that
\begin{align*}
A''(c_*)\Big(\pd{q}{\varphi}(0,\psi)\Big)^2
+B'(c_*)\pd{^2q}{\psi^2}(0,\psi)
+B''(c_*)\Big(\pd{q}{\psi}(0,\psi)\Big)^2=0,
\quad 0\le\psi\le m_{\text{\rm in}}
\end{align*}
and
\begin{align*}
\pd{q}{\psi}(0,\psi)=0,\quad\pd{^2q}{\psi^2}(0,\psi)=0,
\quad 0\le\psi\le m_{\text{\rm in}},
\end{align*}
which imply
\begin{align*}
\pd{q}{\varphi}(0,\psi)=0,\quad 0\le\psi\le m_{\text{\rm in}}
\end{align*}
and thus
\begin{align*}
\pd{q}{\varphi}(\varphi,\psi)=O(\varphi),\quad
{q}(\varphi,\psi)=c_*-O(\varphi^2),\quad(\varphi,\psi)\in (\zeta_{-},0)\times(0,m_{\text{\rm in}}).
\end{align*}
Hence
\begin{align*}
\pd{^2A(q)}{\varphi^2}(\varphi,\psi)=A'(q(\varphi,\psi))\pd{^2q}{\varphi^2}(\varphi,\psi)
+A''(q(\varphi,\psi))\Big(\pd{q}{\varphi}(\varphi,\psi)\Big)^2
=O(\varphi^2),\quad
\\
(\varphi,\psi)\in (\zeta_{-},0)\times(0,m_{\text{\rm in}}),
\end{align*}
which, together with \eqref{p-eq}, \eqref{p-lbbc} and \eqref{p-ubbc}, yields
\begin{align}
\label{rem-aaa-2}
\frac{\Theta'_{-}(x)\cos\Theta_{-}(x)}{{\mathscr Q}_{-}(x)}\Big|_{x=X_{-}(\varphi)}
=-\int_0^{m_{\text{\rm in}}}\pd{^2A(q)}{\varphi^2}(\varphi,\psi)d\psi
=O(\varphi^2),
\quad\varphi\in(\zeta_-,0).
\end{align}
Then, \eqref{rem-aaa-1} follows from \eqref{rem-aaa-2},
$\Theta_{-}(0)=0$, ${\mathscr Q}_{-}(0)=c_*>0$, $X_{-}(0)=0$ and $X'_{-}(0)=1/{c_*}>0$.
$\hfill\Box$\vskip 4mm

\subsection{Formulation of the sonic-supersonic flow problem in the potential plane}

A sonic-supersonic flow is expected in the divergent part of the nozzle $\Omega_+$.
In this subsection, we formulate this sonic-supersonic flow problem in the potential plane.

Assume that $f_+\in C^{3}([0,l_+])\,(0<l_+\le1)$
satisfies
\begin{align}
\label{a-a-0}
\delta_{1,+}x^{\lambda_{+}}\le f''_+(x)\le\delta_{2,+}x^{\lambda_{+}},
\quad f'''(x)\ge0,
\quad0\le x\le{l_+}
\end{align}
with positive constants $\lambda_{+}$, $\delta_{1,+}$ and $\delta_{2,+}$
such that $\lambda_{+}>2$ and $\delta_{1,+}\le\delta_{2,+}$.

\begin{remark}
\label{qqq5}
The second condition in \eqref{a-a-0} can be relaxed by
\begin{align*}
f'''(x)\ge-\omega x^{\lambda_{+}-1},
\quad0\le x\le{l_+},
\end{align*}
where $\omega$ is a positive constant depending only on $\gamma$, $m$,
$\lambda_+$, $\delta_{1,+}$ and $\delta_{2,+}$ (see Remark \ref{remqqq} for the detail).
\end{remark}

Set
$$
Q(\varphi,\psi)=A(q(\varphi,\psi)).
$$
Then, \eqref{eq2or} takes the form
\begin{align}
\label{a-a-2}
Q_{\varphi\varphi}+(K'(Q)Q_{\psi})_\psi=0,
\end{align}
where
\begin{align*}
K(s)=B(A^{-1}_+(s)),\quad s<0
\end{align*}
and
\begin{align*}
K'(s)=&\Big(1-\frac{\gamma+1}{2}q^2\Big)^{-1}
\Big(1-\frac{\gamma-1}{2}q^2\Big)^{2/(\gamma-1)+1}\Big|_{q=A^{-1}_+(s)}<0,
\qquad &&s<0,
\\[2mm]
K''(s)=&(\gamma+1)q^4\Big(1-\frac{\gamma+1}{2}q^2\Big)^{-3}
\Big(1-\frac{\gamma-1}{2}q^2\Big)^{3/(\gamma-1)+1}
\Big|_{q=A^{-1}_+(s)}<0,\qquad &&s<0.
\end{align*}
Since $\lim_{s\to 0^-}K'(s)=-\infty$, \eqref{a-a-2} is a singular nonlinear wave equation
in the sonic-supersonic region.
We now transform \eqref{a-a-2} into a first order system to use the method of characteristics.
Set
$$
U(\varphi,\psi)=Q_\varphi(\varphi,\psi),\quad
V(\varphi,\psi)=Q_\psi(\varphi,\psi).
$$
Then, \eqref{a-a-2} is transformed into
the following system
\begin{align*}
&U_\varphi-(b(Q)V)_\psi=0,
\\
&V_\varphi-U_\psi=0,
\end{align*}
i.e.
\begin{align}
\label{a-a-3}
\left(
\begin{aligned}
U
\\
V
\end{aligned}
\right)_\varphi
+\left(
\begin{array}{cc}
0&-b(Q)
\\[2mm]
-1&0
\end{array}
\right)
\left(
\begin{aligned}
U
\\
V
\end{aligned}
\right)_\psi
+\left(
\begin{array}{cc}
0&-p(Q)Q_\psi
\\[2mm]
0&0
\end{array}
\right) \left(
\begin{aligned}
U
\\
V
\end{aligned}
\right)=\left(
\begin{aligned}
0
\\
0
\end{aligned}
\right),
\end{align}
where
$$
b(Q)=-K'(Q)>0,\quad p(Q)=-K''(Q)>0,\qquad Q<0.
$$
Set
$$
R(Q)=-\frac{1}{2\sqrt{b(Q)}}\left(
\begin{array}{cc}
-\sqrt{b(Q)}&\sqrt{b(Q)}
\\[2mm]
1&1
\end{array}
\right),\quad R^{-1}(Q)=\left(
\begin{array}{cc}
1&-\sqrt{b(Q)}
\\[2mm]
-1&-\sqrt{b(Q)}
\end{array}
\right).
$$
Multiplying \eqref{a-a-3} by $R^{-1}(Q)$ on the left side, one gets
$$
R^{-1}(Q)\left(
\begin{aligned}
U
\\
V
\end{aligned}
\right)_\varphi +\left(
\begin{array}{cc}
\sqrt{b(Q)}&0
\\[2mm]
0&-\sqrt{b(Q)}
\end{array}
\right)R^{-1}(Q)\left(
\begin{aligned}
U
\\
V
\end{aligned}
\right)_\psi+R^{-1}(Q)\left(
\begin{array}{cc}
0&-p(Q)Q_\psi
\\[2mm]
0&0
\end{array}
\right) \left(
\begin{aligned}
U
\\
V
\end{aligned}
\right)=\left(
\begin{aligned}
0
\\
0
\end{aligned}
\right),
$$
which is equivalent to
\begin{align*}
\left(
\begin{aligned}
W
\\
Z
\end{aligned}
\right)_\varphi +\left(
\begin{array}{cc}
\sqrt{b(Q)}&0
\\[2mm]
0&-\sqrt{b(Q)}
\end{array}
\right)\left(
\begin{aligned}
W
\\
Z
\end{aligned}
\right)_\psi =&-R^{-1}(Q)\left(
\begin{array}{cc}
0&-p(Q)Q_\psi
\\[2mm]
0&0
\end{array}
\right) \left(
\begin{aligned}
U
\\
V
\end{aligned}
\right) +\pd{R^{-1}(Q)}{\varphi}\left(
\begin{aligned}
U
\\
V
\end{aligned}
\right)
\\
&\qquad+\left(
\begin{array}{cc}
\sqrt{b(Q)}&0
\\[2mm]
0&-\sqrt{b(Q)}
\end{array}
\right)\pd{R^{-1}(Q)}{\psi} \left(
\begin{aligned}
U
\\
V
\end{aligned}
\right),
\end{align*}
where
$$
\left(
\begin{aligned}
W
\\
Z
\end{aligned}
\right)=R^{-1}(Q)\left(
\begin{aligned}
U
\\
V
\end{aligned}
\right)=\left(
\begin{aligned}
U-\sqrt{b(Q)}V
\\
-U-\sqrt{b(Q)}V
\end{aligned}
\right)=\left(
\begin{aligned}
Q_\varphi-\sqrt{b(Q)}Q_\psi
\\
-Q_\varphi-\sqrt{b(Q)}Q_\psi
\end{aligned}
\right).
$$
Therefore, \eqref{a-a-3} becomes
\begin{align*}
%\label{s3}
&W_\varphi+b^{1/2}(Q)W_\psi
=\frac14b^{-1}(Q)p(Q)(Q_\varphi-{b^{1/2}(Q)}Q_\psi)(W+Z),
\\
%\label{s4}
&Z_\varphi-b^{1/2}(Q)Z_\psi
=\frac14b^{-1}(Q)p(Q)(Q_\varphi+{b^{1/2}(Q)}Q_\psi)(W+Z).
\end{align*}

The flow angle on the upper wall is
$$
\Theta_{+}(x)=\arctan f'_+(x),\quad0\le x\le l_+.
$$
Let the speed of the flow on the upper wall be given by
$$
q(x,f_+(x))={\mathscr Q}_{+}(x),\quad0\le x\le l_+.
$$
Then, the potential function on the upper wall is
$$
\varphi(x,f_+(x))=\Phi_{+}(x),\quad0\le x\le l_+
$$
satisfying $\Phi_{+}(0)=0$ and
$$
\Phi'_{+}(x)=\frac{{\mathscr Q}_{+}(x)}{\cos\Theta_{+}(x)}>0, \quad0\le x\le l_+.
$$
Thus, $\Phi_{+}$ is expressed as
\begin{align}
\label{Phiub+} \Phi_{+}(x)=\int_{0}^{x}
\frac{{\mathscr Q}_{+}(s)}{\cos\Theta_{+}(s)}ds,\quad 0\le x\le l_+.
\end{align}
Denote by $X_{+}$ the inverse function of $\Phi_{+}$, i.e.
\begin{align}
\label{Xub+}
X_{+}(\varphi)=\Phi_{+}^{-1}(\varphi),\quad
0\le\varphi\le\zeta_+,
\end{align}
where
\begin{align}
\label{zeta+}
\zeta_+=\Phi_{+}(l_+).
\end{align}

Therefore, the sonic-supersonic flow problem in the potential plane
can be formulated as
\begin{align}
\label{supp-eq} &Q_{\varphi\varphi}-(b(Q)Q_{\psi})_\psi=0,
\quad&&(\varphi,\psi)\in(0,\zeta_{+})\times(0,m),
\\
\label{supp-inbc1}
&Q(0,\psi)=0,\quad&&\psi\in(0,m),
\\
\label{supp-inbc2}
&Q_\varphi(0,\psi)=0,\quad&&\psi\in(0,m),
\\
\label{supp-lbbc}
&Q_{\psi}(\varphi,0)=0,\quad&&\varphi\in(0,\zeta_{+}),
\\
\label{supp-ubbc}
&Q_{\psi}(\varphi,m)=-\frac{\Theta'_{+}(X_{+}(\varphi))\cos\Theta_{+}(X_{+}(\varphi))}
{b(Q(\varphi,m)))A^{-1}_+(Q(\varphi,m))},
\quad&&\varphi\in(0,\zeta_{+}),
\\
\label{supp-q}
&{\mathscr Q}_{+}(x)=A^{-1}_+(Q(\Phi_{+}(x),m))
\quad&&x\in[0,l_+]
\end{align}
or equivalently as
\begin{align}
\label{supp-neq1}
&W_\varphi+b^{1/2}(Q)W_\psi
=\frac14b^{-1}(Q)p(Q)(Q_\varphi-{b^{1/2}(Q)}Q_\psi)(W+Z),
\quad&&(\varphi,\psi)\in(0,\zeta_{+})\times(0,m),
\\
\label{supp-neq2}
&Z_\varphi-b^{1/2}(Q)Z_\psi
=\frac14b^{-1}(Q)p(Q)(Q_\varphi+{b^{1/2}(Q)}Q_\psi)(W+Z),
\quad&&(\varphi,\psi)\in(0,\zeta_{+})\times(0,m),
\\
\label{supp-ninbc1} &W(0,\psi)=0,\quad&&\psi\in(0,m),
\\
\label{supp-ninbc2} &Z(0,\psi)=0,\quad&&\psi\in(0,m),
\\
\label{supp-nlbbc}
&W(\varphi,0)+Z(\varphi,0)=0,\quad&&\varphi\in(0,\zeta_{+}),
\\
\label{supp-nubbc}
&W(\varphi,m)+Z(\varphi,m)=\frac{2\Theta'_{+}(X_{+}(\varphi))\cos\Theta_{+}(X_{+}(\varphi))}
{b^{1/2}(Q(\varphi,m))A^{-1}_+(Q(\varphi,m))},
\quad&&\varphi\in(0,\zeta_{+}),
\\
\label{supp-q1}
&Q_\varphi(\varphi,\psi)=\frac{1}{2}\big(W(\varphi,\psi)-Z(\varphi,\psi)\big),
\quad&&(\varphi,\psi)\in(0,\zeta_{+})\times(0,m),
\\
\label{supp-q3}
&Q(0,\psi)=0,\quad&&\psi\in(0,m),
\\
\label{supp-q4}
&{\mathscr Q}_{+}(x)=A^{-1}_+(Q(\Phi_{+}(x),m)),
\quad&&x\in[0,l_+].
\end{align}

The problems \eqref{supp-eq}--\eqref{supp-ubbc} and \eqref{supp-neq1}--\eqref{supp-q3}
should be regarded as initial boundary value problems, where $\varphi$-direction
plays the role of the time.
Since $\lim_{Q\to0^-}b(Q)=+\infty$, \eqref{supp-eq}, \eqref{supp-neq1} and \eqref{supp-neq2}
are singular at sonic points.
Furthermore, as shown in Theorem \ref{thm-aaa},
all positive and negative characteristics
from any sonic point are the same line $\{0\}\times[0,m]$.
That is to say, the singularity at the sonic curve is so strong
in the sense that there is not any characteristic moving from a sonic point
to the supersonic region.
The source terms in \eqref{supp-neq1} and \eqref{supp-neq2} are also
singular at sonic points because $\lim_{Q\to0^-}b^{-1}(Q)p(Q)=+\infty$.
Moreover, the boundary conditions \eqref{supp-ubbc} and \eqref{supp-nubbc}
are nonlinear, nonlocal and implicit.
We will use a fixed point argument
to prove the existence of solutions to the problem \eqref{supp-neq1}--\eqref{supp-q4} as
follows:
For given $Q$ and ${\mathscr Q}_{-}$,
we first solve the problem \eqref{supp-neq1}--\eqref{supp-nubbc},
then define a new $Q$ by solving the problem \eqref{supp-q1}, \eqref{supp-q3}
and a new ${\mathscr Q}_{-}$ by \eqref{supp-q4}.
To get a fixed point, one must choose a space with a precise behavior near the sonic curve
for $Q$ and then show that the new $Q$ also belongs to this space
by some elaborate calculations and optimal estimates.

Since \eqref{supp-eq}, \eqref{supp-neq1} and \eqref{supp-neq2}
are singular, we first seek weak solutions, which are defined as follows.
By more complicated and precise estimates, one can get smooth solutions.

\begin{definition}
A function $Q$ is said to be a weak solution
to the problem \eqref{supp-eq}--\eqref{supp-ubbc}
if $Q\in C^{0,1}([0,\zeta_{+}]\times[0,m])$ satisfies \eqref{supp-inbc1}
and
\begin{align*}
&\int_0^{\zeta_+}\int_0^m
\big(Q_\varphi(\varphi,\psi)\xi_\varphi(\varphi,\psi)
-b(Q(\varphi,\psi))Q_\psi(\varphi,\psi)\xi_\psi(\varphi,\psi)\big)d\varphi d\psi
\\
&\qquad-\int_0^{\zeta_+}\frac{\Theta'_{+}(X_{+}(\varphi))\cos\Theta_{+}(X_{+}(\varphi))}
{A^{-1}_+(Q(\varphi,m))}\xi(\varphi,m) d\varphi=0
\end{align*}
for any $\xi\in C^1([0,\zeta_{+}]\times[0,m])$ with
$\xi(\zeta_{+},\cdot)\Big|_{(0,m)}=0$.
\end{definition}

\begin{definition}
A triad of functions $(W,Z,Q)$ is said to be a weak solution
to the problem \eqref{supp-neq1}--\eqref{supp-q3}
if $W,Z\in L^\infty((0,\zeta_{+})\times(0,m))$
and $Q\in C^{0,1}([0,\zeta_{+}]\times[0,m])$ satisfy

{\rm (i)} For any $\xi,\eta\in C^1([0,\zeta_{+}]\times[0,m])$ with
$$
\xi(\zeta_{+},\cdot)\Big|_{(0,m)}=\eta(\zeta_{+},\cdot)\Big|_{(0,m)}=0,
\quad
\xi(\cdot,0)\Big|_{(0,\zeta_{+})}=-\eta(\cdot,0)\Big|_{(0,\zeta_{+})},
\quad
\xi(\cdot,m)\Big|_{(0,\zeta_{+})}=-\eta(\cdot,m)\Big|_{(0,\zeta_{+})},
$$
it holds that
\begin{align*}
&\int_0^{\zeta_+}\int_0^m
\Big(W\xi_\varphi
+b^{1/2}(Q)W\xi_\psi
+\frac12b^{-1/2}(Q)p(Q)Q_\psi W\xi
\\
&\qquad\qquad+\frac14b^{-1}(Q)p(Q)
(Q_\varphi-{b^{1/2}(Q)}Q_\psi)
(W+Z)\xi\Big)d\varphi d\psi
\\
&\qquad+\int_0^{\zeta_+}\int_0^m
\Big(Z\eta_\varphi
-b^{1/2}(Q)Z\eta_\psi
-\frac12b^{-1/2}(Q)p(Q)Q_\psi Z\eta
\\
&\qquad\qquad+\frac14b^{-1}(Q)p(Q)
(Q_\varphi+{b^{1/2}(Q)}Q_\psi)
(W+Z)\eta\Big)d\varphi d\psi
\\
&\qquad-2\int_0^{\zeta_+}\frac{\Theta'_{+}(X_{+}(\varphi))\cos\Theta_{+}(X_{+}(\varphi))}
{A^{-1}_+(Q(\varphi,m))}\xi(\varphi,m) d\varphi=0;
\end{align*}

{\rm (ii)} \eqref{supp-q1} holds in the sense of distribution;

{\rm (iii)} $Q$ satisfies \eqref{supp-q3}.
\end{definition}

As the discussion for the subsonic-sonic flow problem in Remark \ref{rem-aaa},
the condition \eqref{a-a-0} on the wall of the nozzle is almost necessary for the existence of $C^2$ solutions.

\begin{remark}
\label{rem-aaaa}
If the problem \eqref{supp-eq}--\eqref{supp-ubbc} admits
a $C^2([0,\zeta_{+}]\times[0,m])$ solution,
then
\begin{align*}
f''_+(x)=O(x^2),\quad 0<x<l_+.
\end{align*}
\end{remark}

\subsection{Main results of the transonic flow problem in the de Laval nozzle}

We now state the main results on the transonic flow problem.
Since the sonic curve of the transonic flow studied in the paper is
located at the throat of the nozzle,
the results of the transonic flow can be stated for the subsonic-sonic part
and the sonic-supersonic part separately.

For the subsonic-sonic flow in the convergent part of the nozzle,
the first result concerns the existence and regularity estimates of such a solution.

\begin{theorem}
\label{theoremm1}
Assume that $f_-\in C^{3}([l_-,0])$ satisfies \eqref{s-2}
and $g\in C^{3}([0,f_-(l_-)])$
satisfies \eqref{s-3} and \eqref{s-4}.
There exists a positive constant $\delta_{0,-}$ depending only on $\gamma$, $m$,
$\lambda_-$, $\delta_{1,-}$, $\delta_{2,-}$ and $\epsilon(\cdot)$, such that
for any $-\delta_{0,-}\le l_-<0$,
the problem \eqref{p-eq}--\eqref{Qubq} admits a weak solution
$q\in C^\infty((\zeta_{-},0)\times(0,m_{\text{\rm in}})) \cap
C^1([\zeta_{-},0)\times[0,m_{\text{\rm in}}]) \cap
C([\zeta_{-},0]\times[0,m_{\text{\rm in}}])$ with $m_{\text{\rm in}}=m$ satisfying
\begin{gather}
\label{rr-aaa4}
\Big|\pd{q}{\psi}(\varphi,\psi)\Big|\le
C_{1,-}(-\zeta_-)^{\lambda-1/2}(-\varphi)^{1/2},\quad
(\varphi,\psi)\in(\zeta_{-},0)\times(0,m),
\\
\label{rr-aaa5}
\Big|A(q(\varphi',\psi'))-A(q(\varphi'',\psi''))\Big|\le
C_{1,-}(|\varphi'-\varphi''|^{1/2}+|\psi'-\psi''|),\quad
(\varphi',\psi'),(\varphi'',\psi'')\in(\zeta_{-},0)\times(0,m),
\\
\label{rr-aaa6}
c_*-C_{3,-}(-\varphi)^{{\lambda_-}/2+1}\le q(\varphi,\psi)
\le c_*-C_{2,-}(-\varphi)^{{\lambda_-}/2+1},\quad
(\varphi,\psi)\in(\zeta_{-},0)\times(0,m),
\end{gather}
where $C_{i,-}\,(i=1,2,3)\,(C_{2,-}\le C_{3,-})$
are positive constants depending only on $\gamma$, $m$,
$\lambda_-$, $\delta_{1,-}$ and $\delta_{2,-}$.
Moreover,

{\rm (i)} If $f_-$ satisfies
\begin{align}
\label{ssrr-1}
|f'''_-(x)|\le\delta_{3,-}(-x)^{\lambda_-/4+1/2},
\quad l_-\le x\le0
\end{align}
with a positive constant $\delta_{3,-}$
additionally, then any weak solution satisfies $q\in C^1([\zeta_{-},0]\times[0,m_{\text{\rm in}}])$
and
\begin{gather}
\label{rr-aaa7}
\Big|\pd{q}\varphi(\varphi,\psi)\Big|\le C_{4,-}(-\varphi)^{{\lambda_-}/4+1/2},\quad
\Big|\pd{q}\psi(\varphi,\psi)\Big|\le C_{4,-}(-\varphi)^{{\lambda_-}/2+1},
\quad
(\varphi,\psi)\in(\zeta_{-},0)\times(0,m),
\end{gather}
where the positive constant $C_{4,-}$ depends only on $\gamma$, $m$,
$\lambda_-$, $\delta_{1,-}$, $\delta_{2,-}$ and $\delta_{3,-}$.

{\rm (ii)} If $f_-\in C^{4}([l_-,0])$ satisfies \eqref{ssrr-1} and
\begin{align}
\label{ssrr-2}
|f^{(4)}_-(x)|\le\delta_{4,-},\quad l_-\le x\le0
\end{align}
with a positive constant $\delta_{4,-}$ additionally,
then any weak solution satisfies $q\in C^{1,1}((\zeta_{-},0]\times[0,m])$ and
\begin{align}
\label{rr-aaa8}
\Big|\pd{^2q}{\varphi^2}(\varphi,\psi)\Big|\le C_{5,-},\quad
\Big|\pd{^2q}{\varphi\partial\psi}(\varphi,\psi)\Big|\le C_{5,-}(-\varphi)^{{\lambda_-}/4+1/2},
&\quad\Big|\pd{^2q}{\psi^2}(\varphi,\psi)\Big|\le C_{5,-}(-\varphi)^{{\lambda_-}/2+1},
\nonumber
\\
&\hskip-20mm
(\varphi,\psi)\in(\zeta_{-}/2,0)\times(0,m)
\end{align}
with a positive constant $C_{5,-}$ depending only on $\gamma$, $m$, $\lambda_-$,
$\delta_{1,-}$, $\delta_{2,-}$,
$\delta_{3,-}$ and $\delta_{4,-}$.

{\rm (iii)} If $g\in C^{3,\alpha}([0,f_-(l_-)])$
and $f_-\in C^{3,\alpha}([l_-,0))$ for a number $\alpha\in(0,1)$ additionally,
then any weak solution satisfies $q\in C^{2,\alpha}([\zeta_{-},0)\times[0,m])$.
\end{theorem}

The next result yields the uniqueness.

\begin{theorem}
\label{theoremm2}
There exists at most one solution $\varphi\in C^{1,1}(\overline\Omega_-)$
to the problem \eqref{phytran9}--\eqref{phytran14}.
\end{theorem}

For the sonic-supersonic flow in the divergent part of the nozzle,
we start with the existence and regularity estimates.

\begin{theorem}
\label{theoremm3}
Assume that $f_+\in C^{3}([0,l_+])$ satisfies \eqref{a-a-0}.
There exists a positive constant $\delta_{0,+}>0$ depending only on $\gamma$, $m$,
$\lambda_+$, $\delta_{1,+}$ and $\delta_{2,+}$, such that
for any $0<l_+\le\delta_{0,+}$,
the problem \eqref{supp-eq}--\eqref{supp-q} admits at least one weak solution
$Q\in C^{0,1}([0,\zeta_{+}]\times[0,m])$ satisfying
\begin{gather}
\label{rr-aaa1}
-C_{2,+}\varphi^{\lambda_++2}\le Q(\varphi,\psi)\le-C_{1,+}\varphi^{\lambda_++2},
\quad(\varphi,\psi)\in(0,\zeta_{+})\times(0,m),
\\
\label{rr-aaa2}
-C_{2,+}\varphi^{\lambda_++1}\le Q_\varphi(\varphi,\psi)
\le-C_{1,+}\varphi^{\lambda_++1},
\quad
|Q_\psi(\varphi,\psi)|\le C_{2,+}\varphi^{3{\lambda_+}/2+1},
\quad(\varphi,\psi)\in(0,\zeta_{+})\times(0,m)
\end{gather}
with positive constants $C_{1,+}$ and $C_{2,+}$
$(C_{1,+}\le C_{2,+})$ depending only on $\gamma$, $m$,
$\lambda_+$, $\delta_{1,+}$ and $\delta_{2,+}$.
Furthermore, assume that $f_+$ satisfies
\begin{align}
\label{rr-1}
|f'''_+(x)|\le\delta_{3,+}x^{\lambda_+-1},
\quad0\le x\le{l_+}
\end{align}
with a positive constant $\delta_{3,+}$
additionally. Then there exist positive constants $\tilde\delta_{0,+}$ and $C_{3,+}$,
both depending only on $\gamma$, $m$,
$\lambda_+$, $\delta_{1,+}$, $\delta_{2,+}$ and $\delta_{3,+}$,
such that
if $0<l_+\le\tilde\delta_{0,+}$, then the problem \eqref{supp-eq}--\eqref{supp-q} admits at least
a solution $Q\in C^{1,1}([0,\zeta_+]\times[0,m])$, which satisfies
\eqref{rr-aaa1}, \eqref{rr-aaa2} and
\begin{align}
\label{rr-aaa3}
|Q_{\varphi\varphi}(\varphi,\psi)|\le C_{3,+}\varphi^{\lambda_+},
\quad
|Q_{\varphi\psi}(\varphi,\psi)|\le C_{3,+}\varphi^{5{\lambda_+}/4+1/2},
&\quad
|Q_{\psi\psi}(\varphi,\psi)|\le C_{3,+}\varphi^{3{\lambda_+}/2+1},
\nonumber
\\
&\hskip-20mm
\quad(\varphi,\psi)\in(0,\zeta_{+})\times(0,m).
\end{align}
\end{theorem}

The following theorem shows the uniqueness of weak solutions.

\begin{theorem}
\label{theoremm4}
Assume that $f_+\in C^{4}([0,l_+])$ satisfies \eqref{a-a-0} and \eqref{rr-1}.
Then the problem \eqref{supp-eq}--\eqref{supp-q} admits at most one weak solution
$Q\in C^{0,1}([0,\zeta_+]\times[0,m])$
satisfying \eqref{rr-aaa1} and \eqref{rr-aaa2}
with any given positive constants $C_{i,+}\,(i=1,2)$.
\end{theorem}

Let us connect the subsonic-sonic flow and the sonic-supersonic flow to
get a global smooth transonic flow.

\begin{theorem}
\label{theoremm5}
Assume that $f_-\in C^{4}([l_-,0])$ satisfies \eqref{s-2}, \eqref{ssrr-1} and \eqref{ssrr-2},
$g\in C^{3,\alpha}([0,f_-(l_-)])$ satisfies \eqref{s-3} and \eqref{s-4},
while $f_+\in C^{4}([0,l_+])$ satisfies \eqref{a-a-0} and \eqref{rr-1}.
If $-l_-$ and $l_+$ are sufficiently small,
the problem \eqref{p-eq}--\eqref{Qubq}
admits a unique solution $q_-\in C^{1,1}([\zeta_{-},0]\times[0,m])$
satisfying \eqref{rr-aaa4}--\eqref{rr-aaa6}, \eqref{rr-aaa7} and \eqref{rr-aaa8},
while the problem \eqref{supp-eq}--\eqref{supp-q}
admits a unique solution $Q_+\in C^{1,1}([0,\zeta_+]\times[0,m])$
satisfying \eqref{rr-aaa1}, \eqref{rr-aaa2} and \eqref{rr-aaa3}.
Connect $q_-$ and $Q_+$ in the following way
$$
q(\varphi,\psi)=\left\{
\begin{aligned}
&q_-(\varphi,\psi),\quad&&(\varphi,\psi)\in [\zeta_{-},0]\times[0,m],
\\
&A^{-1}_+(Q_+(\varphi,\psi)),\quad&&(\varphi,\psi)\in [0,\zeta_{+}]\times[0,m].
\end{aligned}
\right.
$$
Then, $q\in C^{1,1}([\zeta_{-},\zeta_{+}]\times[0,m])$ is a solution to \eqref{eq2or}.
Furthermore, $\pd{q}{\varphi}=\pd{q}{\psi}=0$ at the sonic curve
and $\pd{^2q}{\psi^2}$ is continuous at sonic points.
Therefore, $q$ satisfies \eqref{eq2or} at the sonic curve, i.e.
$$
\lim_{\varphi\to0}\Big(\pd{^2A(q(\varphi,\psi))}{\varphi^2}+\pd{^2B(q(\varphi,\psi))}{\psi^2}\Big)=0,
\quad0\le\psi\le m.
$$
\end{theorem}

In the physical plane, Theorem \ref{theoremm5} can be stated as follows.

\begin{theorem}
\label{theoremm6}
Assume that $f_-\in C^{4}([l_-,0])$ satisfies \eqref{s-2}, \eqref{ssrr-1} and \eqref{ssrr-2},
$g\in C^{3,\alpha}([0,f_-(l_-)])$ satisfies \eqref{s-3} and \eqref{s-4},
while $f_+\in C^{4}([0,l_+])$ satisfies \eqref{a-a-0} and \eqref{rr-1}.
If $-l_-$ and $l_+$ are sufficiently small,
the transonic flow problem \eqref{phytran1}--\eqref{phytran8}
admits a unique classical solution $\varphi\in C^{2,1}(\overline\Omega)$,
which satisfies \eqref{rr-aaa4}--\eqref{rr-aaa6}, \eqref{rr-aaa7} and \eqref{rr-aaa8}
for $q=|\nabla\varphi|$ in $\Omega_-$,
while satisfies \eqref{rr-aaa1}, \eqref{rr-aaa2} and \eqref{rr-aaa3}
for $Q=A(|\nabla\varphi|)$ in $\Omega_+$.
\end{theorem}

\hskip75mm
\setlength{\unitlength}{0.6mm}
\begin{picture}(250,65)
\put(-70,0){\vector(1,0){140}}
\put(0,40){\vector(0,1){17}}
\put(67,-4){$x$} \put(-4,53){$y$}
\put(-2,25){$\Omega$}
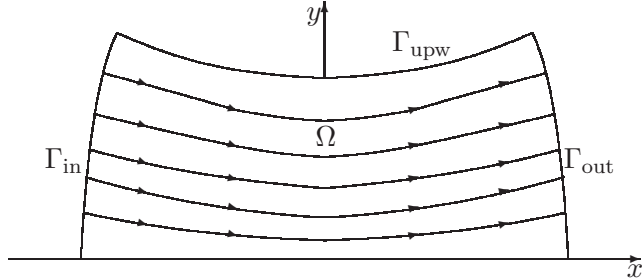
\put(-80,-15){Figure: the smooth transonic flow in the de Laval nozzle}

\put(-60,48){\cbezier(14,2)(30,-5)(40,-7)(60,-8)}
\put(60,48){\cbezier(-14,2)(-30,-5)(-40,-7)(-60,-8)}

\put(-40,48){\qbezier(-6,2)(-12,-10)(-14,-48)}
\put(40,48){\qbezier(6,2)(12,-10)(14,-48)}

\put(15,46){$\Gamma_{\text{upw}}$}
\put(-62,20){$\Gamma_{\text{in}}$}
\put(53,20){$\Gamma_{\text{out}}$}

\put(-60,30){\qbezier(9,2)(30,-3)(40,-5)}
\put(-60,30){\qbezier(40,-5)(50,-7)(60,-7.5)}
\put(-60,30){\qbezier(8.5,-6)(30,-12)(60,-14.4)}
\put(-60,30){\qbezier(7.1,-12)(33,-19.5)(60,-20.9)}
\put(-60,25){\qbezier(7.2,-15)(30,-20)(60,-21)}
\put(-40,21.3){\vector(3,-1){1}}
\put(-40,14.9){\vector(4,-1){1}}
\put(-40,7.7){\vector(4,-1){1}}
\put(-20,17.6){\vector(4,-1){1}}
\put(-20,11.2){\vector(4,-1){1}}
\put(-20,5.3){\vector(4,-1){1}}
\put(-40,29.5){\vector(3,-1){1}}
\put(-20,25.1){\vector(3,-1){1}}

\put(60,30){\qbezier(-9,2)(-30,-3)(-40,-5)}
\put(60,30){\qbezier(-40,-5)(-50,-7)(-60,-7.5)}
\put(60,30){\qbezier(-8.5,-6)(-30,-12)(-60,-14.4)}
\put(60,30){\qbezier(-7.1,-12)(-33,-19.5)(-60,-20.9)}
\put(60,25){\qbezier(-7.2,-15)(-30,-20)(-60,-21)}
\put(40,21.4){\vector(4,1){1}}
\put(40,15){\vector(4,1){1}}
\put(40,7.8){\vector(4,1){1}}
\put(20,17.9){\vector(4,1){1}}
\put(20,11.3){\vector(4,1){1}}
\put(20,5.4){\vector(4,1){1}}
\put(40,29.8){\vector(4,1){1}}
\put(20,25.3){\vector(4,1){1}}

\put(-60,38){\qbezier(11.3,3)(30,-2)(40,-5)}
\put(-60,38){\qbezier(40,-5)(50,-7)(60,-7.5)}
\put(60,38){\qbezier(-11.3,3)(-30,-2)(-40,-5)}
\put(60,38){\qbezier(-40,-5)(-50,-7)(-60,-7.5)}
\put(-40,38.9){\vector(3,-1){1}}
\put(-20,33.1){\vector(3,-1){1}}
\put(40,38.9){\vector(3,1){1}}
\put(20,33.3){\vector(3,1){1}}

\end{picture}
\vskip15mm

\begin{remark}
In Theorem \ref{theoremm6} we get a $C^{1,1}$ transonic flow.
However, this transonic flow pattern is strongly singular in the sense that
the sonic curve is a characteristic degenerate boundary in the subsonic-sonic region,
while in the sonic-supersonic region all characteristics from sonic points coincide,
which are the sonic curve and never approach the supersonic region.
\end{remark}

\begin{remark}
In the de Laval nozzle $\Omega$, there is also a subsonic-sonic-subsonic flow,
whose sonic curve is located at the throat of the nozzle.
\end{remark}

\begin{remark}
For the nozzle $\Omega$, if there is a flat part between the convergent part and the divergent part,
there is also a transonic flow, which is subsonic in the convergent part, sonic in the flat part,
while supersonic in the divergent part.
\end{remark}

\begin{remark}
The geometry of nozzles is important for smooth transonic flows.
As mentioned in Remarks \ref{rem-aaa} and \ref{rem-aaaa},
$$
f''(x)=O(x^2),\quad l_-<x<l_+
$$
is necessary for a $C^2$ transonic flow
whose sonic curve is located at the throat of the nozzle.
Moreover, in \cite{WX2}, we get a continuous subsonic-sonic flow
in a convergent nozzle with straight solid walls.
However, this subsonic-sonic flow is singular in the sense that while the speed is
continuous yet the acceleration blows up at the sonic state,
and there is no way to extend it to be a transonic flow or a subsonic-sonic-subsonic flow.
\end{remark}

\section{Subsonic-sonic flows in the convergent part of the nozzle}

In this section, we will establish the well-posedness of the subsonic-sonic flow problem in the
convergent part of the nozzle.
For convenience, we abbreviate
$\lambda_{-}$ and $\delta_{k,-}$ by $\lambda$ and $\delta_k$ for $k=1,2,3,4$, respectively,
and use $C_i$, $\underline c$, $\mu_i$, $\kappa_i$, $M_i$ to denote generic positive constants.
Furthermore, a parenthesis after a generic constant means that
this constant depends only on the variables in the parentheses.

Owing to \eqref{s-2}, $f_-(l_-)$ and $R_0$ satisfy
$$
m c_*^{-1-2/(\gamma-1)}<f_-(l_-)<m c_*^{-1-2/(\gamma-1)}+\delta_{2},
\quad
-\delta_{2}<f'_-(l_-)<0
$$
and
\begin{align*}
%\label{s-nn}
\frac{\delta_{1}c_*^{1+2/(\gamma-1)}}{(\lambda+1)(m+\delta_{2}c_*^{1+2/(\gamma-1)})
\sqrt{\delta_{2}^2+1}}(-l_-)^{\lambda+1}
<\frac{1}{R_0}<
\frac{\delta_{2}c_*^{1+2/(\gamma-1)}}{(\lambda+1)m}(-l_-)^{\lambda+1}
\end{align*}
with $m$ given by \eqref{mass-m}.

In order to solve the problem \eqref{p-eq}--\eqref{p-outbc} by a fixed point argument,
one needs to specify ${\mathscr Q}_{\text{\rm in}}$ and ${\mathscr Q}_{-}$
in advance as follows:
${\mathscr Q}_{\text{\rm in}}\in C^{0,1}([0,f_{-}(l_{-})])$ satisfies
\begin{align}
\label{qin}
\max\Big\{\frac{c_*}{2},c_*-C_2(-l_-)^{\lambda/2+1}\Big\}\le {\mathscr Q}_{\text{\rm in}}(y)
\le \max\Big\{\frac{c_*}{2},c_*-C_1(-l_-)^{\lambda/2+1}\Big\},
\quad|{\mathscr Q}'_{\text{\rm in}}(y)|\le 1,
\quad
\nonumber
\\
0\le y\le f_{-}(l_{-})
\end{align}
and ${\mathscr Q}_{-}\in C([l_{-},0])\cap
C^{1/2}([l_{-},l_{-}/2])$ satisfies
\begin{align}
\label{qub}
\underline c\le {\mathscr Q}_{-}(x)\le c_*,
\quad l_{-}\le x\le 0,\qquad
[{\mathscr Q}_{-}]_{1/2;(l_{-},l_{-}/2)}\le C_0,
\end{align}
where $C_i\,(i=0,1,2)$ and $\underline c\,(C_1\le C_2,\,\underline c<c_*)$ will be determined.
Then, $m_{\text{\rm in}}$, $\zeta_-$, $\Psi'_{\text{\rm in}}$ and $\Phi'_{-}$
satisfy
\begin{align*}
m_1=\frac{1}2m c_*^{-2/(\gamma-1)}\rho(c_*^2/4)\le m_{\text{\rm in}}\le m_2
=c_*\rho(c_*^2)(m c_*^{-1-2/(\gamma-1)}+\delta_{2})\sqrt{\delta_{2}^2+1},
\end{align*}
\begin{align*}
%\label{zeta-}
c_*\sqrt{\delta_{2}^2+1}l_-
\le\zeta_-\le\underline c l_-,
\end{align*}
\begin{align*}
%\label{Psiine}
\frac{1}2c_*\rho(c_*^2/4)
\le\Psi'_{\text{\rm in}}(y)
\le c_*\rho(c_*^2)\sqrt{\delta_{2}^2+1},\quad0\le y\le f_{-}(l_{-}),
\end{align*}
\begin{align*}
%\label{Phiube}
\underline c\le\Phi'_{-}(x)\le c_*\sqrt{\delta_{2}^2+1},\quad
l_{-}\le x\le0.
\end{align*}
Moreover, one can verify that
\begin{gather}
\label{inestimate1}
\mu_{1}(-\zeta_-)^{\lambda+1}
\le-\frac{\Theta'_{\text{\rm in}}(y)\cos\Theta_{\text{\rm in}}(y)}
{{\mathscr Q}_{\text{\rm in}}(y)\rho({\mathscr Q}_{\text{\rm in}}^2(y))}
\Big|_{y=Y_{\text{\rm in}}(\psi)}
\le\mu_{2}(-\zeta_-)^{\lambda+1},
\quad\psi\in(0,m_{\text{\rm in}}),
\\
\label{inestimate2}
\Big|\Big(\frac{\Theta'_{\text{\rm in}}(y)\cos\Theta_{\text{\rm in}}(y)}
{{\mathscr Q}_{\text{\rm in}}(y)\rho({\mathscr Q}_{\text{\rm in}}^2(y))}
\Big|_{y=Y_{\text{\rm in}}(\psi)}\Big)'\Big|
\le\mu_{2}(C_2+1)\big(\epsilon(\zeta_-)+\zeta_-^2\big)(-\zeta_-)^{3\lambda/2},
\quad\psi\in(0,m_{\text{\rm in}}),
\\
\label{inestimate}
\mu_{1}(-\varphi)^{\lambda}\le
\frac{\Theta'_{-}(x)\cos\Theta_{-}(x)}{{\mathscr Q}_{-}(x)}\Big|_{x=X_{-}(\varphi)}
\le\mu_{2}(-\varphi)^{\lambda},
\quad\varphi\in(\zeta_{-},0),
\end{gather}
where $\mu_{1}$ and $\mu_{2}\,(\mu_{1}\le\mu_{2})$ depend only on
$\gamma$, $m$, $\lambda$, $\delta_1$, $\delta_2$ and $\underline c$.

\subsection{A comparison principle}

We will state a comparison principle
for weak solutions to the problem \eqref{p-eq}--\eqref{p-outbc}.
More generally, instead of \eqref{p-outbc}, one imposes the following boundary condition
\begin{align}
\label{p-outbc1}
q(0,\psi)=c,\quad\psi\in(0,m_{\text{\rm in}}),
\end{align}
where $0<c\le c_*$.
Weak solutions, supersolutions and subsolutions to the problem
\eqref{p-eq}--\eqref{p-ubbc}, \eqref{p-outbc1}
are defined in the following sense.

\begin{definition}
A function $q\in
L^\infty((\zeta_{-},0)\times(0,m_{\text{\rm in}}))$ is said to be a
weak supersolution (subsolution) to the problem
\eqref{p-eq}--\eqref{p-ubbc}, \eqref{p-outbc1}, if
$$
0<\inf_{(\zeta_{-},0)\times(0,m_{\text{\rm in}})}q \le\sup_{(\zeta_{-},0)\times(0,m_{\text{\rm in}})}q\le
c_*
$$
and
\begin{align*}
&\int^{0}_{\zeta_-}\int_0^{m_{\text{\rm in}}}\Big(
A(q(\varphi,\psi))\pd{^2\xi}{\varphi^2}(\varphi,\psi)
+B(q(\varphi,\psi))\pd{^2\xi}{\psi^2}(\varphi,\psi)\Big)d\varphi
d\psi
\\
&\qquad\qquad-A(c)\int_0^{m_{\text{\rm in}}}\pd{\xi}{\varphi}(0,\psi)d\psi
+\int_0^{m_{\text{\rm in}}}\frac{\Theta'_{\text{\rm in}}(y)\cos\Theta_{\text{\rm in}}(y)}
{{\mathscr Q}_{\text{\rm in}}(y)\rho({\mathscr Q}_{\text{\rm in}}^2(y))}
\Big|_{y=Y_{\text{\rm in}}(\psi)}\xi(\zeta_-,\psi)d\psi
\\
&\qquad\qquad
+\int^0_{\zeta_-}\frac{\Theta'_{-}(x)\cos\Theta_{-}(x)}{{\mathscr Q}_{-}(x)}
\Big|_{x=X_{-}(\varphi)}\xi(\varphi,m_{\text{\rm in}})d\varphi\le(\ge)0
\end{align*}
for any nonnegative function $\xi\in C^2([\zeta_{-},0]\times[0,m_{\text{\rm in}}])$
with
$$
\pd{\xi}{\psi}(\cdot,0)\Big|_{(\zeta_{-},0)}
=\pd{\xi}{\psi}(\cdot,m_{\text{\rm in}})\Big|_{(\zeta_{-},0)}=0\quad\mbox{ and
}\quad\pd{\xi}{\varphi}(\zeta_-,\cdot)\Big|_{(0,m_{\text{\rm in}})}
=\xi(0,\cdot)\Big|_{(0,m_{\text{\rm in}})}=0.
$$
Furthermore, $q\in L^\infty((\zeta_{-},0)\times(0,m_{\text{\rm in}}))$ is said to be a
weak solution to the problem \eqref{p-eq}--\eqref{p-ubbc}, \eqref{p-outbc1}
if $q$ is both a weak supersolution and a weak
subsolution.
\end{definition}

Let $q\in L^\infty((\zeta_{-},0)\times(0,m_{\text{\rm in}}))$ be a
weak solution to the problem \eqref{p-eq}--\eqref{p-ubbc}, \eqref{p-outbc1}.
If $q$ satisfies
\begin{align}
\label{class}
\sup_{(\zeta_{-},0)\times(0,m_{\text{\rm in}})}q<c_*,
\end{align}
then $q\in C^\infty((\zeta_{-},0)\times(0,m_{\text{\rm in}}))
\cap C^1([\zeta_{-},0]\times[0,m_{\text{\rm in}}))
\cap C([\zeta_{-},0]\times[0,m_{\text{\rm in}}])$
due to the classical elliptic theory.

It follows from the proof of Proposition 3.1 in \cite{WX2} that the following comparison principle holds.

\begin{lemma}
\label{lemma-comparison}
Assume that $q_{+}\in
L^\infty((\zeta_{-},0)\times(0,m_{\text{\rm in}}))$ and $q_{-}\in
L^\infty((\zeta_{-},0)\times(0,m_{\text{\rm in}}))$ are a weak supersolution and a weak
subsolution to the problem \eqref{p-eq}--\eqref{p-ubbc}, \eqref{p-outbc1},
respectively. Then
$$
q_{+}(\varphi,\psi)\ge q_{-}(\varphi,\psi),\quad
(\varphi,\psi)\in(\zeta_{-},0)\times(0,m_{\text{\rm in}}).
$$
\end{lemma}

\subsection{Solutions to elliptic boundary problems}

We now study solutions to the problem \eqref{p-eq}--\eqref{p-ubbc}, \eqref{p-outbc1}
with $0<c<c_*$.
Due to Lemma \ref{lemma-comparison}, such a solution is unique.
First we show the existence for sufficiently small $c$.

\begin{lemma}
\label{s-lemma1}
There exists $\kappa_1\in(0,1]$ depending only on $\gamma$,
$m$, $\lambda$, $\delta_1$ and $\delta_2$ such that for any $-\kappa_1\le\zeta_-<0$,
the problem \eqref{p-eq}--\eqref{p-ubbc}, \eqref{p-outbc1}
admits a unique solution with \eqref{class} and
\begin{align}
\label{class1}
\sup_{(0,m_{\text{\rm in}})}q(\zeta_{-},\cdot)\le c_*-
\big(\epsilon(\zeta_-)+\zeta_-^2\big)(-\zeta_{-})^{\lambda/2+1}
\end{align}
if $0<c\le c_*/2$.
\end{lemma}

\Proof
Set
$$
\varepsilon_0=\min\Big\{\frac{c_*}{12m_{2}^2},
\inf_{c_*/2<s<5c_*/6}\frac{A'(s)}{B'(s)}\Big\}.
$$
Choose $\kappa_1\in(0,1]$ such that
$$
\kappa_1^2\le\frac{c_*}6,\quad
\kappa_1\le\frac{c_*}{24},\quad
\kappa_1^\lambda\le\frac{2\varepsilon_0m_1}
{\mu_{2}}\inf_{c_*/2<s<5c_*/6}B'(s).
$$
Then,
$$
q_+(\varphi,\psi)=\frac23c_*-2\varphi-\varphi^2+\varepsilon_0\psi^2,
\quad(\varphi,\psi)\in [\zeta_{-},0]\times[0,m_{\text{\rm in}}]
$$
is a supersolution to the problem \eqref{p-eq}--\eqref{p-ubbc}, \eqref{p-outbc1} if $0<c\le c_*/2$.
In fact,
\begin{align}
\label{qqq1}
\frac12c_*\le q_+(\varphi,\psi)\le\frac56c_*,
\quad(\varphi,\psi)\in [\zeta_{-},0]\times[0,m_{\text{\rm in}}],
\end{align}
\begin{align*}
\pd{^2A(q_+)}{\varphi^2}(\varphi,\psi)+\pd{^2B(q_+)}{\psi^2}(\varphi,\psi)
\le& A'(q_+(\varphi,\psi))\pd{^2q_+}{\varphi^2}(\varphi,\psi)
+B'(q_+(\varphi,\psi))\pd{^2q_+}{\psi^2}(\varphi,\psi)
\\
=&-2A'(q_+(\varphi,\psi))+2\varepsilon_0B'(q_+(\varphi,\psi))
\\
\le&0,\quad(\varphi,\psi)\in(\zeta_{-},0)\times(0,m_{\text{\rm in}}),
\end{align*}
\begin{gather*}
\pd{A(q_+)}{\varphi}(\zeta_-,\psi)=
-2(1+\zeta_-)A'(q_+(\zeta_-,\psi))
\le0,\quad\psi\in(0,m_{\text{\rm in}}),
\\
\pd{B(q_+)}{\psi}(\varphi,0)=0,\quad\varphi\in(\zeta_{-},0),
\\
\pd{B(q_+)}{\psi}(\varphi,m_{\text{\rm in}})
=2\varepsilon_0m_{\text{\rm in}}B'(q_+(\varphi,m_{\text{\rm in}}))
\ge\mu_{2}(-\varphi)^{\lambda},\quad\varphi\in(\zeta_{-},0).
\end{gather*}
Additionally, it is not hard to verify that
$$
q_-(\varphi,\psi)=A_{-}^{-1}\big(A(c)+\mu_{2}\varphi\big),
\quad(\varphi,\psi)\in [\zeta_{-},0]\times[0,m_{\text{\rm in}}]
$$
is a subsolution to the problem \eqref{p-eq}--\eqref{p-ubbc},
\eqref{p-outbc1}. Then, one can get the existence and the uniqueness of the solution
to the problem \eqref{p-eq}--\eqref{p-ubbc} by a standard argument in the classical elliptic theory.
Finally, \eqref{class1} follows from \eqref{qqq1} and
$\big(\epsilon(\zeta_-)+\zeta_-^2\big)(-\zeta_{-})^{\lambda/2+1}\le2\kappa_1\le{c_*}/{12}$.
$\hfill\Box$\vskip 4mm

Below, we investigate properties of the solution
to the problem \eqref{p-eq}--\eqref{p-ubbc}, \eqref{p-outbc1}
with $0<c<c_*$.
As ithe proof of Proposition 3.2 in \cite{WX2}, we introduce two functions.
The first one is
\begin{align*}
E(s)=A(B^{-1}(s)),\quad s<0.
\end{align*}
Clearly, $E\in C^\infty((-\infty,0))$ satisfies
\begin{align*}
E'(s)=& \Big(1-\frac{\gamma+1}{2}q^2\Big)
\Big(1-\frac{\gamma-1}{2}q^2\Big)^{-2/(\gamma-1)-1}\Big|_{q=B^{-1}(s)}
>0,&&\quad s<0,
\\
E''(s)=&-(\gamma+1)q^4\Big(1-\frac{\gamma-1}{2}q^2\Big)^{-3/(\gamma-1)-2}
\Big|_{q=B^{-1}(s)} <0,&&\quad s<0,
\\
E'''(s)=&-(\gamma+1)q^4(4+3q^2)\Big(1-\frac{\gamma-1}{2}q^2\Big)^{-4/(\gamma-1)-3}
\Big|_{q=B^{-1}(s)} <0,&&\quad s<0.
\end{align*}
The other function is
\begin{align*}
G(s)=E'(E^{-1}(s)),\quad s<0,
\end{align*}
which satisfies $G\in C^\infty((-\infty,0))$ and
\begin{align*}
G'(s)=&\frac{E''(E^{-1}(s))}{E'(E^{-1}(s))}
=\frac{E''(t)}{E'(t)}\Big|_{t=E^{-1}(s)}<0,&&\quad s<0,
\\[4mm]
G''(s)=&\frac{E'''(t)E'(t)-(E''(t))^2}{(E'(t))^3}\Big|_{t=E^{-1}(s)}
<-\frac{(E''(t))^2}{(E'(t))^3}\Big|_{t=E^{-1}(s)}<0,&&\quad s<0.
\end{align*}

\begin{lemma}
\label{s-lemma2}
Assume that $c_*/3\le c<c_*$ and
$q$, satisfying \eqref{class} and \eqref{class1},
is the solution to the problem \eqref{p-eq}--\eqref{p-ubbc}, \eqref{p-outbc1}.
Then
\begin{gather}
\label{s-lemma2-0}
q(\varphi,\psi)\ge A_{-}^{-1}\big(A(c_*/3)+\mu_{2}\varphi\big),\quad
(\varphi,\psi)\in[\zeta_{-},0]\times[0,m_{\text{\rm in}}],
\\
\label{s-lemma2-1}
\Big|\pd{q}{\psi}(\varphi,\psi)\Big|\le
\mu_{3}(C_2+1)(-\zeta_-)^{\lambda-1/2}(-\varphi)^{1/2},\quad
(\varphi,\psi)\in[\zeta_{-},0]\times[0,m_{\text{\rm in}}]
\end{gather}
and
\begin{align}
\label{s-lemma2-2}
\Big|A(q(\varphi',\psi'))-A(q(\varphi'',\psi''))\Big|\le
\mu_{4}(C_2+1)(|\varphi'-\varphi''|^{1/2}+|\psi'-\psi''|),\quad
\nonumber
\\
(\varphi',\psi'),(\varphi'',\psi'')\in[\zeta_{-},0]\times[0,m_{\text{\rm in}}],
\end{align}
where $\mu_{i}=\mu_{i}(\gamma, m, \lambda, \delta_1, \delta_2, \underline c)$ for $i=3,4$.
\end{lemma}

\Proof
First, \eqref{s-lemma2-0} follows from Lemma \ref{lemma-comparison} and the fact that
$$
q_-(\varphi,\psi)=A_{-}^{-1}\big(A(c_*/3)+\mu_{2}\varphi\big),
\quad(\varphi,\psi)\in [\zeta_{-},0]\times[0,m_{\text{\rm in}}]
$$
is a subsolution to the problem \eqref{p-eq}--\eqref{p-ubbc}, \eqref{p-outbc1}.
To prove \eqref{s-lemma2-1}, one can adapt the proof of Proposition 3.2 in \cite{WX2}
as follows.
Set
$$
z(\varphi,\psi)=\pd{B(q)}{\psi}(\varphi,\psi),
\quad(\varphi,\psi)\in[\zeta_-,0]\times[0,m_{\text{\rm in}}].
$$
Then $z$ solves the following problem
\begin{align}
\label{prop2-1}
&j_{1}(\varphi,\psi)\pd{^2z}{\varphi^2}
+\pd{^2z}{\psi^2}
+j_{2}(\varphi,\psi)\pd{z}{\varphi}
+j_{3}(\varphi,\psi)\pd{z}{\psi}
+j_{4}(\varphi,\psi){z}=0, \hskip-40mm&& \nonumber
\\
&\quad&&\hskip-25mm(\varphi,\psi)\in(\zeta_{-},0)\times(0,m_{\text{\rm in}}),
\\
\label{prop2-2}
&\pd{z}{\varphi}(\zeta_-,\psi)+e(\psi)z(\zeta_-,\psi)=
-\frac{B'(q(\zeta_-,\psi))}{A'(q(\zeta_-,\psi))}
\Big(\frac{\Theta'_{\text{\rm in}}(y)\cos\Theta_{\text{\rm in}}(y)}
{{\mathscr Q}_{\text{\rm in}}(y)\rho({\mathscr Q}_{\text{\rm in}}^2(y))}
\Big|_{y=Y_{\text{\rm in}}(\psi)}\Big)',
&&\psi\in(0,m_{\text{\rm in}}),
\\
\label{prop2-3}
&z(\varphi,0)=0, &&\varphi\in(\zeta_{-},0),
\\
\label{prop2-4}
&z(\varphi,m_{\text{\rm in}})=\frac{\Theta'_{-}(x)\cos\Theta_{-}(x)}{{\mathscr Q}_{-}(x)}
\Big|_{x=X_{-}(\varphi)},
\quad&&\varphi\in(\zeta_{-},0),
\\
\label{prop2-5}
&z(0,\psi)=0, &&\psi\in(0,m_{\text{\rm in}}),
\end{align}
where for $(\varphi,\psi)\in(\zeta_{-},0)\times(0,m_{\text{\rm in}})$,
\begin{gather*}
j_{1}(\varphi,\psi)=G(A(q(\varphi,\psi))),
\quad
j_{2}(\varphi,\psi)=2G'(A(q(\varphi,\psi)))
\pd{A(q(\varphi,\psi))}{\varphi},
\\
j_{3}(\varphi,\psi)=-G'(A(q(\varphi,\psi)))
{\pd{B(q(\varphi,\psi))}{\psi}},
\quad
j_{4}(\varphi,\psi)=G''(A(q(\varphi,\psi)))
\Big(\pd{A(q(\varphi,\psi))}{\varphi}\Big)^2,
\end{gather*}
and
\begin{align*}
e(\psi)=-\frac{E''(q(\zeta_-,\psi))}{(E'(q(\zeta_-,\psi)))^2}
\frac{\Theta'_{\text{\rm in}}(y)\cos\Theta_{\text{\rm in}}(y)}
{{\mathscr Q}_{\text{\rm in}}(y)\rho({\mathscr Q}_{\text{\rm in}}^2(y))}
\Big|_{y=Y_{\text{\rm in}}(\psi)}<0,
\quad\psi\in(0,m_{\text{\rm in}}).
\end{align*}
It follows from \eqref{s-lemma2-0}, \eqref{class1} and \eqref{inestimate2} that
there exists a positive number $\tilde\mu_{3}$ depending only on
$\gamma$, $m$, $\lambda$, $\delta_1$, $\delta_2$ and $\underline c$
such that
$$
\Big|\frac{B'(q(\zeta_-,\psi))}{A'(q(\zeta_-,\psi))}
\Big(\frac{\Theta'_{\text{\rm in}}(y)\cos\Theta_{\text{\rm in}}(y)}
{{\mathscr Q}_{\text{\rm in}}(y)\rho({\mathscr Q}_{\text{\rm in}}^2(y))}
\Big|_{y=Y_{\text{\rm in}}(\psi)}\Big)'\Big|
\le\frac12\tilde\mu_{3}(C_2+1)(-\zeta_-)^{\lambda-1},
\quad\psi\in(0,m_{\text{\rm in}}).
$$
Note
\begin{align*}
\frac14j_{1}(\varphi,\psi)(-\varphi)^{-3/2}
-j_{4}(\varphi,\psi)(-\varphi)^{1/2}
\ge&\sqrt{-j_{1}(\varphi,\psi)j_{4}(\varphi,\psi)}
(-\varphi)^{-1/2}
\\
\ge&-\frac12j_{2}(\varphi,\psi)(-\varphi)^{-1/2},
\quad(\varphi,\psi)\in(\zeta_{-},0)\times(0,m_{\text{\rm in}}).
\end{align*}
One can verify that
$$
z_{\pm}(\varphi,\psi)=\pm
\tilde\mu_{3}(C_2+1)(-\zeta_-)^{\lambda-1/2}(-\varphi)^{1/2},
\quad(\varphi,\psi)\in[\zeta_{-},0]\times[0,m_{\text{\rm in}}]
$$
are super and sub solutions to the problem
\eqref{prop2-1}--\eqref{prop2-5}.
Then, \eqref{s-lemma2-1} follows from the classical comparison principle.
Finally, one can get \eqref{s-lemma2-2} from
\eqref{s-lemma2-0} and \eqref{s-lemma2-1} by a standard process
as in the proof of Proposition 3.2 in \cite{WX2}.
$\hfill\Box$\vskip 4mm

\subsection{The weak solution to the degenerate elliptic boundary problem}

We are now ready to solve the problem \eqref{p-eq}--\eqref{p-outbc}.
First we determine $\underline c$ by
\begin{align}
\label{underlinec}
\underline c=A_{-}^{-1}\big(A(c_*/3)-\mu_{2}\big)
\in(0,c_*/3).
\end{align}

\begin{proposition}
\label{s-prop1}
There exists $\kappa_2\in(0,\kappa_1]$ depending only on $\gamma$,
$m$, $\lambda$, $\delta_1$, $\delta_2$, $\epsilon(\cdot)$ and $C_2$ such that for any
$-\kappa_2\le\zeta_-<0$,
the problem \eqref{p-eq}--\eqref{p-outbc}
admits a weak solution $q\in C^\infty((\zeta_{-},0)\times(0,m_{\text{\rm in}}))
\cap C^1([\zeta_{-},0)\times[0,m_{\text{\rm in}}))
\cap C([\zeta_{-},0)\times[0,m_{\text{\rm in}}])$
satisfying
\begin{align}
\label{s-prop1-1}
\Big|\pd{q}{\psi}(\varphi,\psi)\Big|\le
\mu_{3}(C_2+1)(-\zeta_-)^{\lambda-1/2}(-\varphi)^{1/2},\quad
(\varphi,\psi)\in[\zeta_{-},0]\times[0,m_{\text{\rm in}}],
\end{align}
\begin{align}
\label{s-prop1-2}
\Big|A(q(\varphi',\psi'))-A(q(\varphi'',\psi''))\Big|\le
\mu_{4}(C_2+1)(|\varphi'-\varphi''|^{1/2}+|\psi'-\psi''|),\quad
\nonumber
\\
(\varphi',\psi'),(\varphi'',\psi'')\in[\zeta_{-},0]\times[0,m_{\text{\rm in}}],
\end{align}
\begin{align}
\label{s-prop1-3}
c_*-\mu_{6}(-\varphi)^{\lambda/2+1}\le q(\varphi,\psi)
\le c_*-\mu_{5}(-\varphi)^{\lambda/2+1},\quad
(\varphi,\psi)\in[\zeta_{-},0]\times[0,m_{\text{\rm in}}],
\end{align}
where $\mu_{5}$ and $\mu_{6}\,(\mu_{5}\le\mu_{6})$ depend only on
$\gamma$, $m$, $\lambda$, $\delta_1$ and $\delta_2$.
\end{proposition}

\Proof
Set
\begin{align*}
{\mathscr E}=\Big\{c_*/3<c<c_*:
\mbox{the problem \eqref{p-eq}--\eqref{p-ubbc}, \eqref{p-outbc1}
admits a solution } q_c
\mbox{ with \eqref{class} and \eqref{class1}}
\Big\}.
\end{align*}
According to Lemmas \ref{lemma-comparison} and \ref{s-lemma1},
${\mathscr E}\,(\not=\emptyset)$ is an interval.
We claim that ${\mathscr E}=(c_*/3,c_*)$ if $-\kappa_2\le\zeta_-<0$,
where $\kappa_2\in(0,\kappa_1]$
will be determined.

Assume that $c\in{\mathscr E}$.
It follows from Lemma \ref{s-lemma2} that
$q_c$ satisfies \eqref{s-lemma2-1} and \eqref{s-lemma2-2}.
Integrating \eqref{p-eq} and \eqref{p-inbc} over $(0,m_{\text{\rm in}})$ with respect to $\psi$ and using
\eqref{p-lbbc} and \eqref{p-ubbc}, one gets that
\begin{align}
\label{s-prop1-4}
&\frac{d^2}{d\varphi^2}\int_0^{m_{\text{\rm in}}}A(q_c(\varphi,\psi))d\psi
=-\frac{\Theta'_{-}(x)\cos\Theta_{-}(x)}{{\mathscr Q}_{-}(x)}\Big|_{x=X_{-}(\varphi)},
\quad\varphi\in(\zeta_{-},0),
\\
\nonumber
&\Big(\frac{d}{d\varphi}\int_0^{m_{\text{\rm in}}}A(q_c(\varphi,\psi))d\psi\Big)\Big|_{\varphi=\zeta_-}=
-\int_0^{m_{\text{\rm in}}}\frac{\Theta'_{\text{\rm in}}(y)\cos\Theta_{\text{\rm in}}(y)}
{{\mathscr Q}_{\text{\rm in}}(y)\rho({\mathscr Q}_{\text{\rm in}}^2(y))}
\Big|_{y=Y_{\text{\rm in}}(\psi)}d\psi.
\end{align}
Direction calculations yield
\begin{align*}
%\label{s-prop1-5-1}
\int_{\zeta_-}^0\frac{\Theta'_{-}(x)\cos\Theta_{-}(x)}{{\mathscr Q}_{-}(x)}
\Big|_{x=X_{-}(\varphi)}d\varphi
=\int_{l_-}^0\frac{\Theta'_{-}(x)\cos\Theta_{-}(x)}{{\mathscr Q}_{-}(x)}
\Phi'_{-}(x)dx=\int_{l_-}^0\Theta'_{-}(x)dx=-\Theta_{-}(-l_-)
\end{align*}
and
\begin{align*}
%\label{s-prop1-5-2}
\int_0^{m_{\text{\rm in}}}\frac{\Theta'_{\text{\rm in}}(y)\cos\Theta_{\text{\rm in}}(y)}
{{\mathscr Q}_{\text{\rm in}}(y)\rho({\mathscr Q}_{\text{\rm in}}^2(y))}
\Big|_{y=Y_{\text{\rm in}}(\psi)}d\psi
=&\int_0^{f_-(l_-)}\frac{\Theta'_{\text{\rm in}}(y)\cos\Theta_{\text{\rm in}}(y)}
{{\mathscr Q}_{\text{\rm in}}(y)\rho({\mathscr Q}_{\text{\rm in}}^2(y))}
\Psi'_{\text{\rm in}}(y)dy
\\
=&\int_0^{f_-(l_-)}\Theta'_{\text{\rm in}}(y)dy=\Theta_{\text{\rm in}}(f_-(l_-)).
\end{align*}
Therefore,
\begin{align}
\label{s-prop1-5}
\Big(\frac{d}{d\varphi}\int_0^{m_{\text{\rm in}}}A(q_c(\varphi,\psi))d\psi\Big)\Big|_{\varphi=0}=0.
\end{align}
It follows from \eqref{s-prop1-4}, \eqref{s-prop1-5} and \eqref{p-outbc1} that
$$
\int_0^{m_{\text{\rm in}}}A(q_c(\varphi,\psi))d\psi
={m_{\text{\rm in}}}A(c)-\int_{\varphi}^0\int_{\tilde\varphi}^0
\frac{\Theta'_{-}(x)\cos\Theta_{-}(x)}{{\mathscr Q}_{-}(x)}
\Big|_{x=X_{-}(s)}ds d\tilde\varphi,\quad\varphi\in[\zeta_-,0],
$$
which, together with \eqref{inestimate}, leads to
\begin{align*}
%\label{s-prop1-6}
\int_0^{m_{\text{\rm in}}}A(q_c(\varphi,\psi))d\psi
\le{m_{\text{\rm in}}}A(c)-\frac{\mu_{1}}{(\lambda+1)(\lambda+2)}(-\varphi)^{\lambda+2},
\quad\varphi\in[\zeta_-,0].
\end{align*}
Therefore, there exists a Lipschitz continuous curve $L:\psi=\omega(\varphi)\,(\zeta_-\le\varphi\le0)$
in $[\zeta_{-},0]\times(0,m_{\text{\rm in}})$ satisfying
\begin{align*}
A(q_c(\varphi,\omega(\varphi)))\le
-\frac{\mu_{1}}{{m_{\text{\rm in}}}(\lambda+1)(\lambda+2)}(-\varphi)^{\lambda+2},
\quad\varphi\in[\zeta_-,0].
\end{align*}
Thus, there exists $M_0=M_0(\gamma,m,\lambda,\delta_1,\delta_2)$ such that
\begin{align}
\label{s-prop1-7-0}
q_c(\varphi,\omega(\varphi))\le
c_*-M_0(-\varphi)^{\lambda/2+1},
\quad\varphi\in[\zeta_-,0],
\end{align}
which, together with \eqref{s-lemma2-1}, leads to
\begin{align}
\label{s-prop1-7}
q_c(\zeta_-,\psi)\le
c_*-M_0(-\zeta_-)^{\lambda/2+1}+\mu_{3}(C_2+1)m_2(-\zeta_-)^{\lambda},
\quad
\psi\in[0,m_{\text{\rm in}}].
\end{align}
Decompose $(\zeta_{-},0)\times(0,m_{\text{\rm in}})$ into $\Omega_1$ and $\Omega_2$ by the curve $L$,
where
$$
\Omega_1=\big\{(\varphi,\psi)\in(\zeta_{-},0)\times(0,m_{\text{\rm in}}):
\psi<\omega(\varphi)\big\},\quad
\Omega_2=\big\{(\varphi,\psi)\in(\zeta_{-},0)\times(0,m_{\text{\rm in}}):
\psi>\omega(\varphi)\big\}.
$$
Then, $q_c$ solves the following two problems
\begin{align*}
&\pd{^2A(q_c)}{\varphi^2}+\pd{^2B(q_c)}{\psi^2}=0,
\quad&&(\varphi,\psi)\in\Omega_1,
\\
&q_c(\zeta_-,\psi)\le
c_*-M_0(-\zeta_-)^{\lambda/2+1}+\mu_{3}(C_2+1)m_2(-\zeta_-)^{\lambda},
\quad&&\psi\in(0,\omega(\zeta_{-})),
\\
&\pd{q_c}{\psi}(\varphi,0)=0,\quad&&\varphi\in(\zeta_{-},0),
\\
&q_c(\varphi,\omega(\varphi))\le c_*-M_0(-\varphi)^{\lambda/2+1},
\quad&&\varphi\in(\zeta_-,0),
\\
&q_c(0,\psi)=c,\quad&&\psi\in(0,\omega(0))
\end{align*}
and
\begin{align*}
&\pd{^2A(q_c)}{\varphi^2}+\pd{^2B(q_c)}{\psi^2}=0,
\quad&&(\varphi,\psi)\in\Omega_2,
\\
&q_c(\zeta_-,\psi)\le
c_*-M_0(-\zeta_-)^{\lambda/2+1}+\mu_{3}(C_2+1)m_2(-\zeta_-)^{\lambda},
\quad&&\psi\in(\omega(\zeta_{-}),m_{\text{\rm in}}),
\\
&q_c(\varphi,\omega(\varphi))\le c_*-M_0(-\varphi)^{\lambda/2+1},
\quad&&\varphi\in(\zeta_-,0),
\\
&\pd{B(q_c)}{\psi}(\varphi,m_{\text{\rm in}})=
\frac{\Theta'_{-}(x)\cos\Theta_{-}(x)}{{\mathscr Q}_{-}(x)}\Big|_{x=X_{-}(\varphi)},
\quad&&\varphi\in(\zeta_{-},0),
\\
&q_c(0,\psi)=c,\quad&&\psi\in(\omega(0),m_{\text{\rm in}}).
\end{align*}
Set
$$
\bar q_1(\varphi,\psi)=c_*-\mu_{5}(-\varphi)^{\lambda/2+1},\quad
(\varphi,\psi)\in\overline\Omega_1
$$
and
$$
\bar q_2(\varphi,\psi)=c_*-\mu_{5}(-\varphi)^{\lambda/2+1}
(1+m_2-\psi),\quad
(\varphi,\psi)\in\overline\Omega_2
$$
with $\mu_{5}=\frac{M_0}{2(1+m_2)}$.
Note that $\bar q_1$ and $\bar q_2$ satisfy
\begin{align*}
%\label{subsonicq}
\pd{^2A(\bar q_k)}{\varphi^2}(\varphi,\psi)+\pd{^2B(\bar q_k)}{\psi^2}(\varphi,\psi)
\le0,\quad(\varphi,\psi)\in\Omega_k,\qquad k=1,2.
\end{align*}
Choose a number $\tilde\kappa_2\in(0,\kappa_1]$ depending only on
$\gamma$, $m$, $\lambda$, $\delta_1$, $\delta_2$, $\epsilon(\cdot)$ and $C_2$,
such that
\begin{align}
\label{s-prop1-8-1}
\mu_{3}(C_2+1)m_2\tilde\kappa_2^{\lambda/2-1}\le\frac{M_0}2,
\quad
\sup_{0<s<\tilde\kappa_2}\epsilon(-s)+\tilde\kappa_2^2\le\frac{M_0}{4},
\quad
\mu_{2}\tilde\kappa_2^{\lambda/2-1}\le
\mu_{5}\inf_{\underline c<s<c_*}B'(s).
\end{align}
Owing to \eqref{s-prop1-8-1} and \eqref{inestimate},
one can get from the classical comparison principle that
\begin{align*}
q_c(\varphi,\psi)\le\bar q_k(\varphi,\psi),\quad
(\varphi,\psi)\in\overline\Omega_k,\qquad k=1,2.
\end{align*}
Therefore,
\begin{align}
\label{s-prop1-9}
q_c(\varphi,\psi)\le c_*-\mu_{5}(-\varphi)^{\lambda/2+1},\quad
(\varphi,\psi)\in[\zeta_{-},0]\times[0,m_{\text{\rm in}}].
\end{align}
Owing to \eqref{s-lemma2-2} and \eqref{p-outbc1}, one also gets that
\begin{align}
\label{s-prop1-9-1}
A(q_c(\varphi,\psi))\le A(c)+\mu_{4}(C_2+1)(-\varphi)^{1/2},\quad
(\varphi,\psi)\in[\zeta_{-},0]\times[0,m_{\text{\rm in}}].
\end{align}
Moreover, it follows from \eqref{s-prop1-7} and \eqref{s-prop1-8-1} that
\begin{align}
\label{s-prop1-10}
q_c(\zeta_-,\psi)\le c_*-2\big(\epsilon(\zeta_-)+\zeta_-^2\big)(-\zeta_{-})^{\lambda/2+1},
\quad 0\le\psi\le m_{\text{\rm in}}.
\end{align}

Set $c_0=\sup{\mathscr E}$. We show that $c_0=c_*$ by contradiction.
Otherwise, $c_0<c_*$. Then, for any $c_*/3<c<c_0$,
$c\in{\mathscr E}$ and $q_c$ satisfies \eqref{s-prop1-9}--\eqref{s-prop1-10},
where $q_c$ is the solution to the problem
\eqref{p-eq}--\eqref{p-ubbc}, \eqref{p-outbc1}. Set
$$
q_{c_0}(\varphi,\psi)=\lim_{c\to c_0^-}q_{c}(\varphi,\psi),
\quad(\varphi,\psi)\in[\zeta_{-},0]\times[0,m_{\text{in}}].
$$
Then, it follows from \eqref{s-prop1-9}--\eqref{s-prop1-10} that
\begin{align}
\label{s-prop1-10ww}
\sup_{(\zeta_{-},0)\times(0,m_{\text{in}})}q_{c_0}<c_*,
\quad
\sup_{(0,m_{\text{in}})}q_{c_0}(\zeta_{-},\cdot)
\le c_*-2\big(\epsilon(\zeta_-)+\zeta_-^2\big)(-\zeta_{-})^{\lambda/2+1},
\end{align}
and it is not hard to show that $q_{c_0}$ is the solution
to the problem \eqref{p-eq}--\eqref{p-ubbc}, \eqref{p-outbc1}
with $c=c_0$. Hence $c_0\in{\mathscr E}$.
It follows from \eqref{s-prop1-10ww}
and the stability theory of uniformly elliptic problems
that the problem \eqref{p-eq}--\eqref{p-ubbc}, \eqref{p-outbc1}
admits a unique solution with \eqref{class} and \eqref{class1}
if $c_0<c<c_0+\varepsilon$ with a sufficiently small positive number $\varepsilon$.
This contradicts that $c_0=\sup{\mathscr E}<c_*$.
Therefore, ${\mathscr E}=(c_*/3,c_*)$.

For any $c_*/3<c<c_*$, let $q_c$ be the solution to the
problem \eqref{p-eq}--\eqref{p-ubbc}, \eqref{p-outbc1}  with \eqref{class} and \eqref{class1}. 
Then, $q_c$
satisfies \eqref{s-lemma2-1}, \eqref{s-lemma2-2}, \eqref{s-prop1-9}
and \eqref{s-prop1-10}. Set
$$
q(\varphi,\psi)=\lim_{c\to c_*^-}q_{c}(\varphi,\psi),
\quad(\varphi,\psi)\in[\zeta_{-},0]\times[0,m_{\text{in}}].
$$
One can show that $q$ is the weak solution to
the problem \eqref{p-eq}--\eqref{p-outbc}
satisfying \eqref{s-prop1-1}, \eqref{s-prop1-2} and the second inequality in \eqref{s-prop1-3}.
Moreover, $q\in C^\infty((\zeta_{-},0)\times(0,m_{\text{\rm in}}))
\cap C^1([\zeta_{-},0)\times[0,m_{\text{\rm in}}))
\cap C([\zeta_{-},0]\times[0,m_{\text{\rm in}}])$ thanks to 
the Schauder theory on uniformly elliptic equations.

It remains to verify the first inequality in \eqref{s-prop1-3}.
Similar to the proof of \eqref{s-prop1-7-0},
there exists a Lipschitz continuous curve $\tilde L:\psi=\tilde\omega(\varphi)\,(\zeta_-\le\varphi\le0)$
in $[\zeta_{-},0]\times(0,m_{\text{\rm in}})$ satisfying
\begin{align*}
q(\varphi,\tilde\omega(\varphi))\ge c_*-M_1(-\varphi)^{\lambda/2+1},
\quad\varphi\in[\zeta_-,0]
\end{align*}
with $M_1=M_1(\gamma,m,\lambda,\delta_1,\delta_2)$.
Decompose $(\zeta_{-},0)\times(0,m_{\text{\rm in}})$ into $\tilde\Omega_1$
and $\tilde\Omega_2$ by the curve $\tilde L$,
where
$$
\tilde\Omega_1=\big\{(\varphi,\psi)\in(\zeta_{-},0)\times(0,m_{\text{\rm in}}):
\psi<\tilde\omega(\varphi)\big\},\quad
\tilde\Omega_2=\big\{(\varphi,\psi)\in(\zeta_{-},0)\times(0,m_{\text{\rm in}}):
\psi>\tilde\omega(\varphi)\big\}.
$$
Prescribe a Dirichlet condition for $q$ on $\tilde L$. 
Then the problem \eqref{p-eq}--\eqref{p-outbc} can be regarded as two problems on $\tilde\Omega_1$
and $\tilde\Omega_2$, respectively.
Set
\begin{gather*}
\underline q_1(\varphi,\psi)=c_*-M_2(-\varphi)^{\lambda/2+1}(1+m_{\text{\rm in}}^2-\psi^2),\quad
(\varphi,\psi)\in\overline{\tilde\Omega}_1,
\\
\underline q_2(\varphi,\psi)=c_*-M_2(-\varphi)^{\lambda/2+1}
(1+m_{\text{\rm in}}^2-\psi^2+2m_{\text{\rm in}}\psi),\quad
(\varphi,\psi)\in\overline{\tilde\Omega}_2.
\end{gather*}
Owing to $\lambda>2$, \eqref{inestimate1} and \eqref{inestimate},
there exist sufficiently small $\kappa_2\in(0,\tilde\kappa_2]$
and sufficiently large $M_2$
such that $\underline q_k$ is a subsolution to the problem of $q$ on $\tilde\Omega_k$ for $k=1,2$,
where $M_2$ depends on $\gamma$, $m$, $\lambda$, $\delta_1$ and $\delta_2$,
while $\kappa_2$ also on $\epsilon(\cdot)$ and $C_2$.
Then, the first inequality in \eqref{s-prop1-3} follows from
a comparison principle. Here, the comparison principle
is not the classical one since \eqref{p-eq} is degenerate at the sonic state
and it can be proved in a similar way as for  Proposition 3.1 in \cite{WX2}.
$\hfill\Box$\vskip 4mm

\subsection{Existence of subsonic-sonic flows}

We now choose
\begin{align}
\label{c12}
C_1=\mu_{5}\underline c^{\lambda/2+1},\quad
C_2=\mu_{6}\Big(c_*\sqrt{\delta_2^2+1}\Big)^{\lambda/2+1}.
\end{align}
Assume that $f_-\in C^{3}([l_-,0])$ satisfies \eqref{s-2}
with $-\delta_0\le l_-<0$,
where
\begin{align}
\label{delta0}
\delta_0=\frac{1}{c_*\sqrt{\delta_2^2+1}}\min\Big\{{\kappa_2},
\Big(\frac{1}{\mu_{3}(C_2+1)c_*\rho(c_*^2)\sqrt{\delta_2^2+1}}\Big)^{1/\lambda},
\Big(\frac{c_*}{2\mu_{6}}\Big)^{2/(\lambda+2)}\Big\},
\end{align}
which depends only on $\gamma$,
$m$, $\lambda$, $\delta_1$, $\delta_2$ and $\epsilon(\cdot)$.
Assume that $g\in C^{3}([0,f_-(l_-)])$ satisfies \eqref{s-3} and
\eqref{s-4}.
Choose
\begin{align}
\label{c0}
C_0=\mu_{4}(C_2+1)\Big(c_*\sqrt{\delta_2^2+1}\Big)^{1/2}
\sup\Big\{\frac{1}{A'(s)}:\underline c<s<c_*-\mu_{5}
\Big(\frac{-\underline c l_-}{2}\Big)^{\lambda/2+1}\Big\}.
\end{align}
Set
\begin{align*}
{\mathscr S}=\Big\{&({\mathscr Q}_{\text{\rm in}},{\mathscr Q}_{-})\in
C^{0,1}([0,f_{-}(l_{-})])\times C([l_{-},0])\cap C^{1/2}([l_{-},l_-/2]):
\\
&\qquad
{\mathscr Q}_{\text{\rm in}}\mbox{ satisfies \eqref{qin} with \eqref{c12}},
\mbox{ while } {\mathscr Q}_{-}\mbox{ satisfies \eqref{qub}
with \eqref{underlinec} and \eqref{c0}}\Big\}
\end{align*}
with the norm
$$
\|({\mathscr Q}_{\text{\rm in}},{\mathscr Q}_{-})\|_{\mathscr S}
=\max\Big\{\|{\mathscr Q}_{\text{\rm in}}\|_{L^\infty((0,f_{-}(l_{-})))},
\|{\mathscr Q}_{-}\|_{L^\infty((l_{-},0))}\Big\},\quad
({\mathscr Q}_{\text{\rm in}},{\mathscr Q}_{-})\in {\mathscr S}.
$$
For given $({\mathscr Q}_{\text{\rm in}},{\mathscr Q}_{-})\in {\mathscr S}$,
it follows from Lemma \ref{lemma-comparison} and
Proposition \ref{s-prop1} that
the problem \eqref{p-eq}--\eqref{p-outbc}
admits a unique weak solution $q\in C^\infty((\zeta_{-},0)\times(0,m_{\text{\rm in}}))
\cap C^1([\zeta_{-},0)\times[0,m_{\text{\rm in}}))
\cap C([\zeta_{-},0]\times[0,m_{\text{\rm in}}])$
satisfying the estimates \eqref{s-prop1-1}--\eqref{s-prop1-3}.
Set
\begin{align*}
\hat{\mathscr Q}_{\text{\rm in}}(y)=q(\zeta_-,\Psi_{\text{\rm in}}(y)),\quad
0\le y\le f_{-}(l_{-})
\end{align*}
and
\begin{align*}
\hat{\mathscr Q}_{-}(x)=q(\Phi_{-}(x),m_{\text{\rm in}}),
\quad l_{-}\le x\le 0.
\end{align*}
It follows from \eqref{underlinec}--\eqref{s-prop1-3} and \eqref{c12}--\eqref{c0}
that $\hat{\mathscr Q}_{\text{\rm in}}\in C^{0,1}([0,f_{-}(l_{-})])$ satisfies \eqref{qin},
and $\hat {\mathscr Q}_{-}\in C([l_{-},0])\cap C^{1/2}([l_{-},l_-/2])$ satisfies \eqref{qub}.
Therefore, we can define a mapping $J$ from ${\mathscr S}$ to itself
as follows
\begin{align}
\label{s-mapj}
J(({\mathscr Q}_{\text{\rm in}},{\mathscr Q}_{-}))
=(\hat{\mathscr Q}_{\text{\rm in}},\hat{\mathscr Q}_{-}),
\quad({\mathscr Q}_{\text{\rm in}},{\mathscr Q}_{-})\in {\mathscr S}.
\end{align}

\begin{theorem}
\label{theoremsub}
Assume that $f_-\in C^{3}([l_-,0])$ satisfies
\eqref{s-2}, and $g\in C^{3}([0,f_-(l_-)])$
satisfies \eqref{s-3} and \eqref{s-4}. If $-\delta_0\le l_-<0$,
then, the problem \eqref{p-eq}--\eqref{Qubq} admits a weak solution
$q\in C^\infty((\zeta_{-},0)\times(0,m_{\text{\rm in}})) \cap
C^1([\zeta_{-},0)\times[0,m_{\text{\rm in}}]) \cap
C([\zeta_{-},0]\times[0,m_{\text{\rm in}}])$ with the following estimates
\begin{align*}
%\label{theoremsub-1}
\Big|\pd{q}{\psi}(\varphi,\psi)\Big|\le
\mu_{3}(C_2+1)(-\zeta_-)^{\lambda-1/2}(-\varphi)^{1/2},\quad
(\varphi,\psi)\in(\zeta_{-},0)\times(0,m_{\text{\rm in}}),
\end{align*}
\begin{align*}
%\label{theoremsub-2}
\Big|A(q(\varphi',\psi'))-A(q(\varphi'',\psi''))\Big|\le
\mu_{4}(C_2+1)(|\varphi'-\varphi''|^{1/2}+|\psi'-\psi''|),\quad
\\
(\varphi',\psi'),(\varphi'',\psi'')\in(\zeta_{-},0)\times(0,m_{\text{\rm in}})
\end{align*}
and
\begin{align*}
%\label{theoremsub-3}
c_*-\mu_{6}(-\varphi)^{\lambda/2+1}
\le q(\varphi,\psi)\le c_*-\mu_{5}(-\varphi)^{\lambda/2+1},
\quad(\varphi,\psi)\in(\zeta_{-},0)\times(0,m_{\text{\rm in}}),
\end{align*}
where $\delta_0$ and $C_2$ are defined in \eqref{delta0} and \eqref{c12}, respectively,
while $\mu_{i}\,(i=3,4,5,6)$ are given in Lemma \ref{s-lemma2} and Proposition \ref{s-prop1}.
Furthermore,
\begin{align}
\label{theoremsub-5}
m_{\text{\rm in}}=m.
\end{align}
\end{theorem}

\Proof
As mentioned above, the mapping $J$ defined by \eqref{s-mapj}
is from ${\mathscr S}$ to itself.
It follows from \eqref{s-prop1-1}--\eqref{s-prop1-3} that $J$ is compact.
Therefore, the first part of the theorem follows from
Proposition \ref{s-prop1} and
the Schauder fixed point theorem provided that $J$ is also continuous.

Now let us show that $J$ is continuous.
Assume that $\{({\mathscr Q}_{\text{\rm in}}^{(n)},{\mathscr Q}_{-}^{(n)})\}_{n=0}^\infty
\subset{\mathscr S}$ satisfies
\begin{align}
\label{subt2-1}
\lim_{n\to\infty}\|{\mathscr Q}_{\text{\rm in}}^{(n)}
-{\mathscr Q}_{\text{\rm in}}^{(0)}\|_{L^\infty((0,f_-(l_-)))}=0,
\quad
\lim_{n\to\infty}\|{\mathscr Q}_{-}^{(n)}
-{\mathscr Q}_{-}^{(0)}\|_{L^\infty((l_-,0))}=0.
\end{align}
Let $m_{\text{\rm in}}^{(n)}$, $\zeta_{-}^{(n)}$,
$\Psi_{\text{\rm in}}^{(n)}$, $Y_{\text{\rm in}}^{(n)}$,
$\Phi_-^{(n)}$ and $X_-^{(n)}$ be
defined by \eqref{ass-1}--\eqref{Xub} with
${\mathscr Q}_{\text{\rm in}}={\mathscr Q}_{\text{\rm in}}^{(n)}$ and
${\mathscr Q}_{-}={\mathscr Q}_{-}^{(n)}$ for $n=0,1,2,\cdots$.
Then, \eqref{subt2-1} yields
\begin{gather}
\label{subt2-1-1}
\lim_{n\to\infty}m_{\text{\rm in}}^{(n)}=m_{\text{\rm in}}^{(0)},
\quad
\lim_{n\to\infty}\zeta_{-}^{(n)}=\zeta_{-}^{(0)},
\\
\label{subt2-1-2}
\lim_{n\to\infty}\|\Psi_{\text{\rm in}}^{(n)}-\Psi_{\text{\rm in}}^{(n)}\|_{0,1;(0,f_-(l_-))}=0,
\quad
Y_{\text{\rm in}}^{(n)}\longrightarrow Y_{\text{\rm in}}^{(0)}\mbox{ uniformly in }
(0,m_{\text{\rm in}}^{(0)}),
\\
\label{subt2-1-3}
\lim_{n\to\infty}\|\Phi_-^{(n)}-\Phi_-^{(0)}\|_{0,1;(l_-,0)}=0,
\quad
X_-^{(n)}\longrightarrow X_-^{(0)}\mbox{ uniformly in }(\zeta_{-}^{(0)},0).
\end{gather}
It follows from the definition of $J$ that
$$
J(({\mathscr Q}_{\text{\rm in}}^{(n)},{\mathscr Q}_{-}^{(n)}))
=(q^{(n)}(\zeta_-^{(n)},\Psi_{\text{\rm in}}^{(n)}(\cdot)),
q^{(n)}(\Phi_{-}^{(n)}(\cdot),m_{\text{\rm in}}^{(n)})),
\quad n=0,1,2,\cdots,
$$
where $q^{(n)}\in C^\infty((\zeta^{(n)}_{-},0)\times(0,m^{(n)}_{\text{\rm in}}))
\cap C^1([\zeta^{(n)}_{-},0)\times[0,m^{(n)}_{\text{\rm in}}))
\cap C([\zeta^{(n)}_{-},0]\times[0,m^{(n)}_{\text{\rm in}}])$
is the unique weak solution to the following problem
\begin{align}
\label{subt2-p-eq}
&\pd{^2A(q^{(n)})}{\varphi^2}+\pd{^2B(q^{(n)})}{\psi^2}=0,
\quad&&(\varphi,\psi)\in(\zeta_{-}^{(n)},0)\times(0,m_{\text{\rm in}}^{(n)}),
\\
\label{subt2-p-inbc}
&\pd{A(q^{(n)})}{\varphi}(\zeta_-^{(n)},\psi)=
-\frac{\Theta'_{\text{\rm in}}(y)\cos\Theta_{\text{\rm in}}(y)}
{{\mathscr Q}_{\text{\rm in}}^{(n)}(y)\rho(({\mathscr Q}_{\text{\rm in}}^{(n)}(y))^2)}
\Big|_{y=Y_{\text{\rm in}}^{(n)}(\psi)},\quad&&\psi\in(0,m_{\text{\rm in}}^{(n)}),
\\
\label{subt2-p-lbbc}
&\pd{q^{(n)}}{\psi}(\varphi,0)=0,\quad&&\varphi\in(\zeta_{-}^{(n)},0),
\\
\label{subt2-p-ubbc}
&\pd{B(q^{(n)})}{\psi}(\varphi,m^{(n)}_{\text{\rm in}})=
\frac{\Theta'_{-}(x)\cos\Theta_{-}(x)}{{\mathscr Q}_{-}^{(n)}(x)}\Big|_{x=X_{-}^{(n)}(\varphi)},
\quad&&\varphi\in(\zeta_{-}^{(n)},0),
\\
\label{subt2-p-outbc}
&q^{(n)}(0,\psi)=c_*,\quad&&\psi\in(0,m_{\text{\rm in}}^{(n)})
\end{align}
for $n=0,1,2,\cdots$.
Set
$$
\hat q^{(n)}(\hat\varphi,\hat\psi)
=q^{(n)}(\zeta_{-}^{(n)}\hat\varphi,m_{\text{\rm in}}^{(n)}\hat\psi),
\quad(\hat\varphi,\hat\psi)\in(0,1)\times(0,1),
\qquad n=0,1,2,\cdots.
$$
It remains to verify that
\begin{align}
\label{subt2-2}
\lim_{n\to\infty}\big(\|\hat q^{(n)}(1,\cdot)-\hat q^{(0)}(1,\cdot)\|_{L^\infty((0,1))}
+\|\hat q^{(n)}(\cdot,1)-\hat q^{(0)}(\cdot,1)\|_{L^\infty((0,1))}\big)=0,
\end{align}
which, together with \eqref{subt2-1-1}--\eqref{subt2-1-3}, implies
$$
\lim_{n\to\infty}\|J(({\mathscr Q}_{\text{\rm in}}^{(n)},{\mathscr Q}_{-}^{(n)}))
-J(({\mathscr Q}_{\text{\rm in}}^{(0)},{\mathscr Q}_{-}^{(0)}))\|_{\mathscr S}=0.
$$
We prove \eqref{subt2-2} by contradiction. Otherwise,
there exist  a positive number $\varepsilon_0$ and a subsequence of
$\{q^{(n)}\}_{n=1}^\infty$, denoted by itself for convenience, such that
for each $n=1,2,\cdots$,
\begin{align}
\label{subt2-3}
\|\hat q^{(n)}(1,\cdot)-\hat q^{(0)}(1,\cdot)\|_{L^\infty((0,1))}
+\|\hat q^{(n)}(\cdot,1)-\hat q^{(0)}(\cdot,1)\|_{L^\infty((0,1))}\ge\varepsilon_0.
\end{align}
Since $q^{(n)}$ satisfies \eqref{s-prop1-2} and \eqref{s-prop1-3}
with $m_{\text{\rm in}}=m_{\text{\rm in}}^{(n)}$ and $\zeta_{-}=\zeta_{-}^{(n)}$
for each $n=1,2,\cdots$,
there exists a subsequence of $\{q^{(n)}\}_{n=1}^\infty$, denoted by itself again for convenience,
such that
\begin{align}
\label{subt2-4}
&\lim_{n\to\infty}\|\hat q^{(n)}-\hat q^*\|_{L^\infty((0,1)\times(0,1))}=0,
\end{align}
where
$$
\hat q^*(\hat\varphi,\hat\psi)
=q^*(\zeta_{-}^{(0)}\hat\varphi,m_{\text{\rm in}}^{(0)}\hat\psi),
\quad(\hat\varphi,\hat\psi)\in(0,1)\times(0,1)
$$
and $q^*\in L^\infty((\zeta^{(0)}_{-},0)\times(0,m^{(0)}_{\text{\rm in}}))$
satisfies \eqref{s-prop1-2} and \eqref{s-prop1-3}
with $m_{\text{\rm in}}=m_{\text{\rm in}}^{(0)}$ and $\zeta_{-}=\zeta_{-}^{(0)}$.
Letting $n\to\infty$ in \eqref{subt2-p-eq}--\eqref{subt2-p-outbc}
and using \eqref{subt2-1-1}--\eqref{subt2-1-3} and \eqref{subt2-4}, one can get
that $q^*$ solves the following problem
\begin{align*}
&\pd{^2A(q^{*})}{\varphi^2}+\pd{^2B(q^{*})}{\psi^2}=0,
\quad&&(\varphi,\psi)\in(\zeta_{-}^{(0)},0)\times(0,m_{\text{\rm in}}^{(0)}),
\\
&\pd{A(q^{(*)})}{\varphi}(\zeta_-^{(0)},\psi)=
-\frac{\Theta'_{\text{\rm in}}(y)\cos\Theta_{\text{\rm in}}(y)}
{{\mathscr Q}_{\text{\rm in}}^{(0)}(y)\rho(({\mathscr Q}_{\text{\rm in}}^{(0)}(y))^2)}
\Big|_{y=Y_{\text{\rm in}}^{(0)}(\psi)},\quad&&\psi\in(0,m_{\text{\rm in}}^{(0)}),
\\
&\pd{q^{*}}{\psi}(\varphi,0)=0,\quad&&\varphi\in(\zeta_{-}^{(0)},0),
\\
&\pd{B(q^{*})}{\psi}(\varphi,m^{(0)}_{\text{\rm in}})=
\frac{\Theta'_{-}(x)\cos\Theta_{-}(x)}{{\mathscr Q}_{-}^{(0)}(x)}\Big|_{x=X_{-}^{(0)}(\varphi)},
\quad&&\varphi\in(\zeta_{-}^{(0)},0),
\\
&q^{*}(0,\psi)=c_*,\quad&&\psi\in(0,m_{\text{\rm in}}^{(0)}).
\end{align*}
It follows from Lemma \ref{lemma-comparison} that
$$
q^*(\varphi,\psi)=q^{(0)}(\varphi,\psi),\quad
(\varphi,\psi)\in(\zeta_{-}^{(0)},0)\times(0,m_{\text{\rm in}}^{(0)}),
$$
which contradicts \eqref{subt2-3} and \eqref{subt2-4}.
Hence \eqref{subt2-2} holds.
Finally, it follows from \eqref{phytran13} that
the sonic curve of the flow lies at the throat of the nozzle 
and the velocity vector is along the normal direction at the sonic curve
in the physical plane. Thus \eqref{theoremsub-5} holds.
$\hfill\Box$\vskip 4mm

\subsection{Uniqueness of subsonic-sonic flows}

We show the uniqueness of the subsonic-sonic flow in the physical plane.
Note that the boundary condition \eqref{phytran13} means that
the velocity vector of the flow is along the normal direction at the sonic curve.
Thus, for the problem \eqref{phytran9}--\eqref{phytran14},
the solution satisfies both a Dirichlet and a Neumann boundary conditions at the sonic curve
although there is a free parameter at the inlet.
Hence the uniqueness theorem follows easily.

\begin{theorem}
There exists at most one solution $\varphi\in C^{1,1}(\overline\Omega_-)$
to the problem \eqref{phytran9}--\eqref{phytran14}.
\end{theorem}

\Proof Let $\varphi_1,\varphi_2\in C^{1,1}(\overline\Omega_-)$ be two solutions
to the problem \eqref{phytran9}--\eqref{phytran14}. Assume that
$$
\varphi_k(g(y),y)=C_{\text{\rm in},k},\quad0<y<f_-(l_-),\,k=1,2.
$$
For $k=1,2$, multiplying the equation with respect to $\varphi_k$ by
$(\varphi_1-\varphi_2)-(C_{\text{\rm in},1}-C_{\text{\rm in},2})$
and then integrating over $\Omega_-$ by parts, one can get that
\begin{align*}
\int_{\Omega_-}\rho(|\nabla\varphi_k(x,y)|^2)\nabla\varphi_k(x,y)\cdot
\nabla(\varphi_1(x,y)-\varphi_2(x,y))dxdy
+(C_{\text{\rm in},1}-C_{\text{\rm in},2})\rho(c_*^2)c_*f_-(0)=0.
\end{align*}
Hence
\begin{align*}
\int_{\Omega_-}\big(\rho(|\nabla\varphi_1(x,y)|^2)\nabla\varphi_1(x,y)
-\rho(|\nabla\varphi_2(x,y)|^2)\nabla\varphi_2(x,y)\big)
\cdot(\nabla\varphi_1(x,y)-\nabla\varphi_2(x,y))dxdy=0,
\end{align*}
which, together with \eqref{phytran13} and \eqref{phytran14}, leads to
$$
\varphi_1(x,y)=\varphi_2(x,y),\quad (x,y)\in\Omega_-.
$$
$\hfill\Box$\vskip 4mm

\subsection{Regularity of the subsonic-sonic flow}

The subsonic-sonic flow in Theorem \ref{theoremsub}
is only continuous near the sonic state. We investigate its regularity
in this subsection.

\begin{theorem}
Let $q$ be the weak solution to the problem \eqref{p-eq}--\eqref{Qubq}
under the assumptions of Theorem \ref{theoremsub}.

{\rm(i)} If $f_-$ satisfies \eqref{ssrr-1} additionally,
then $q\in C^1([\zeta_{-},0]\times[0,m])$ satisfies
\begin{align}
\label{subreg6}
\Big|\pd{q}\varphi(\varphi,\psi)\Big|\le\mu_{7}(-\varphi)^{\lambda/4+1/2},\quad
\Big|\pd{q}\psi(\varphi,\psi)\Big|\le\mu_{7}(-\varphi)^{\lambda/2+1},\quad
(\varphi,\psi)\in(\zeta_{-},0)\times(0,m)
\end{align}
with $\mu_{7}=\mu_{7}(\gamma,m,\lambda,\delta_1,\delta_2,\delta_3)$.

{\rm(ii)} If $f_-\in C^4([l_-,0])$ satisfies \eqref{ssrr-1} and \eqref{ssrr-2} additionally,
then $q\in C^{1,1}((\zeta_{-},0]\times[0,m])$ satisfies
\begin{align}
\label{subreg7}
\Big|\pd{^2q}{\varphi^2}(\varphi,\psi)\Big|\le\mu_{8},\quad
\Big|\pd{^2q}{\varphi\partial\psi}(\varphi,\psi)\Big|\le\mu_{8}(-\varphi)^{\lambda/4+1/2},
&\quad
\Big|\pd{^2q}{\psi^2}(\varphi,\psi)\Big|\le\mu_{8}(-\varphi)^{\lambda/2+1},
\nonumber
\\
&(\varphi,\psi)\in(\zeta_{-}/2,0)\times(0,m)
\end{align}
with $\mu_{8}=\mu_{8}(\gamma,m,\lambda,\delta_1,\delta_2,\delta_3,\delta_4)$.

{\rm (iii)} If $g\in C^{3,\alpha}([0,f_-(l_-)])$
and $f_-\in C^{3,\alpha}([l_-,0))$ for a number $\alpha\in(0,1)$ additionally,
then $q\in C^{2,\alpha}([\zeta_{-},0)\times[0,m])$.
\end{theorem}

\Proof
Thanks to the Schauder theory on uniformly elliptic equations,
one can get that \eqref{subreg6} holds
in $(\zeta_{-},\zeta_{-}/2]\times(0,m)$ in the case (i),
$q\in C^{1,1}((\zeta_{-},0)\times[0,m])$ in the case (ii)
and $q\in C^{2,\alpha}([\zeta_{-},0)\times[0,m])$ in the case (iii).
So it suffices to verify \eqref{subreg6} and \eqref{subreg7} in $(\zeta_{-}/2,0)\times(0,m)$.

For any $0<\varepsilon<-\zeta_-/2\le1/2$, $q\in C^2([-2\varepsilon,-\varepsilon/2]\times(0,m))
\cap C^1([-2\varepsilon,-\varepsilon/2]\times[0,m])$ satisfies
\begin{align*}
&\pd{}{\varphi}\Big(a(\varphi,\psi)\pd{q}\varphi\Big)+\pd{}{\psi}\Big(h(\varphi,\psi)\pd{q}\psi\Big)=0,
\quad&&(\varphi,\psi)\in(-2\varepsilon,-\varepsilon/2)\times(0,m),
\\
&\pd{q}{\psi}(\varphi,0)=0,\quad&&\varphi\in(-2\varepsilon,-\varepsilon/2),
\\
&\pd{q}{\psi}(\varphi,m)=e(\varphi),\quad&&\varphi\in(-2\varepsilon,-\varepsilon/2),
\end{align*}
where
\begin{align*}
a(\varphi,\psi)=A'(q(\varphi,\psi)),\quad
h(\varphi,\psi)=B'(q(\varphi,\psi)),\quad
\quad(\varphi,\psi)\in[-2\varepsilon,-\varepsilon/2]\times[0,m]
\end{align*}
and
$$
e(\varphi)=\frac{\Theta'_{-}(X_{-}(\varphi))\cos\Theta_{-}(X_{-}(\varphi))}
{B'(q(\varphi,m))q(\varphi,m)},
\quad\varphi\in[-2\varepsilon,-\varepsilon/2].
$$
Shifting and rescaling the potential as
$$
\tilde\varphi=\varepsilon^{-\lambda/4-1/2}(\varphi+\varepsilon),\quad\varphi\in[-2\varepsilon,-\varepsilon/2]
$$
and setting
$$
\tilde q(\tilde\varphi,\psi)=
q(-\varepsilon+\varepsilon^{\lambda/4+1/2}\tilde\varphi,\psi)-c_*,\quad
(\tilde\varphi,\psi)\in[-\varepsilon^{-\lambda/4+1/2},\varepsilon^{-\lambda/4+1/2}/2]\times[0,m],
$$
one can check that
$\tilde q\in C^2([-\varepsilon^{-\lambda/4+1/2},\varepsilon^{-\lambda/4+1/2}/2]\times(0,m))
\cap C^1([-\varepsilon^{-\lambda/4+1/2},\varepsilon^{-\lambda/4+1/2}/2]\times[0,m])$ solves
\begin{align*}
&\pd{}{\tilde\varphi}\Big(\varepsilon^{-\lambda/2-1}\tilde a(\tilde\varphi,\psi)
\pd{\tilde q}{\tilde\varphi}\Big)+\pd{}{\psi}\Big(\tilde h(\tilde\varphi,\psi)\pd{\tilde q}\psi\Big)=0,
\quad&&(\tilde\varphi,\psi)\in(-\varepsilon^{-\lambda/4+1/2},\varepsilon^{-\lambda/4+1/2}/2)\times(0,m),
\\
&\pd{\tilde q}{\psi}(\tilde\varphi,0)=0,\quad&&\tilde\varphi\in(-\varepsilon^{-\lambda/4+1/2},\varepsilon^{-\lambda/4+1/2}/2),
\\
&\pd{\tilde q}{\psi}(\tilde\varphi,m)=\tilde e(\tilde\varphi),
\quad&&\tilde\varphi\in(-\varepsilon^{-\lambda/4+1/2},\varepsilon^{-\lambda/4+1/2}/2),
\end{align*}
where $\tilde a$ and $\tilde h$ are defined on
$[-\varepsilon^{-\lambda/4+1/2},\varepsilon^{-\lambda/4+1/2}/2]\times[0,m]$
by
$$
\tilde a(\tilde\varphi,\psi)=a(-\varepsilon+\varepsilon^{\lambda/4+1/2}\tilde\varphi,\psi),\quad
\tilde h(\tilde\varphi,\psi)=h(-\varepsilon+\varepsilon^{\lambda/4+1/2}\tilde\varphi,\psi),
$$
and
$$
\tilde e(\tilde\varphi)=e(-\varepsilon+\varepsilon^{\lambda/4+1/2}\tilde\varphi),\quad
\tilde\varphi\in[-\varepsilon^{-\lambda/4+1/2},\varepsilon^{-\lambda/4+1/2}/2].
$$
It follows from the H\"older continuity estimates for uniformly elliptic equations
that there exists a number $\beta\in(0,1)$ such that
\begin{align}
\label{sssa1}
[\tilde q]_{\beta;(-1/4,1/4)\times(0,m)}\le M_1\|\tilde q\|_{L^\infty((-1/2,1/2)\times(0,m))}
\le M_2\varepsilon^{\lambda/2+1},
\end{align}
where $M_i=M_i(\gamma,m,\lambda,\delta_1,\delta_2)$ for $i=1,2$.
Then, \eqref{sssa1} and \eqref{ssrr-1} yield
$$
\varepsilon^{-\lambda/2-1}\|\tilde a\|_{\beta;(-1/4,1/4)\times(0,m)}
+\|\tilde h\|_{\beta;(-1/4,1/4)\times(0,m)}\le M_3,
\quad
\|\tilde e\|_{\beta;(-1/4,1/4)}\le M_3\varepsilon^{\lambda/2+1}
$$
with $M_3=M_3(\gamma,m,\lambda,\delta_1,\delta_2,\delta_3)$.
So one gets from the Schauder estimates on uniformly elliptic equations that
\begin{align}
\label{subreg8}
\|\tilde q\|_{1,\beta;(-1/8,1/8)\times(0,m)}\le M_4\big(\|\tilde q\|_{L^\infty((-1/4,1/4)\times(0,m))}
+\|\tilde e\|_{\beta;(-1/4,1/4)}\big)\le M_5\varepsilon^{\lambda/2+1},
\end{align}
where $M_i=M_i(\gamma,m,\lambda,\delta_1,\delta_2,\delta_3)$ for $i=4,5$.
It follows from \eqref{subreg8} that
$$
\Big|\pd{q}\varphi(-\varepsilon,\psi)\Big|
\le M_6\varepsilon^{\lambda/4+1/2},\quad
\Big|\pd{q}\psi(-\varepsilon,\psi)\Big|\le
M_6\varepsilon^{\lambda/2+1},\quad\psi\in[0,m]
$$
with $M_6=M_6(\gamma,m,\lambda,\delta_1,\delta_2,\delta_3)$,
which lead to \eqref{subreg6} owing to the arbitrariness of $\varepsilon\in(0,-\zeta_-/2)$.
Furthermore, we get from \eqref{subreg8}, \eqref{ssrr-1} and \eqref{ssrr-2} that
$$
\varepsilon^{-\lambda/2-1}\|\tilde a\|_{1,\beta;(-1/8,1/8)\times(0,m)}
+\|\tilde h\|_{1,\beta;(-1/8,1/8)\times(0,m)}\le M_7,
\quad
\|\tilde e\|_{1,\beta;(-1/8,1/8)}\le M_7\varepsilon^{\lambda/2+1}
$$
with $M_7=M_7(\gamma,m,\lambda,\delta_1,\delta_2,\delta_3,\delta_4)$.
Using the Schauder estimates again yields
$$
\|\tilde q\|_{2,\beta;(-1/{16},1/{16})\times(0,m)}\le
M_8\big(\|\tilde q\|_{L^\infty((-1/8,1/8)\times(0,m))}
+\|\tilde e\|_{1,\beta;(-1/8,1/8)}\big)\le M_9\varepsilon^{\lambda/2+1},
$$
where $M_i=M_i(\gamma,m,\lambda,\delta_1,\delta_2,\delta_3,\delta_4)$ for $i=8,9$.
Therefore,
$$
\Big|\pd{^2q}{\varphi^2}(-\varepsilon,\psi)\Big|\le M_{10},\quad
\Big|\pd{^2q}{\varphi\partial\psi}(-\varepsilon,\psi)\Big|
\le M_{10}\varepsilon^{\lambda/4+1/2},\quad
\Big|\pd{^2q}{\psi^2}(-\varepsilon,\psi)\Big|\le M_{10}\varepsilon^{\lambda/2+1},\quad
\psi\in[0,m]
$$
with $M_{10}=M_{10}(\gamma,m,\lambda,\delta_1,\delta_2,\delta_3,\delta_4)$,
which yield \eqref{subreg7} owing to the arbitrariness of $\varepsilon\in(0,-\zeta_-/2)$.
$~\hfill\Box$\vskip 4mm

\section{Sonic-supersonic flows in the divergent part of the nozzle}

In this section, we will establish the well-posedness of sonic-supersonic flows
in the divergent part of the nozzle.
For convenience,
$\lambda_{+}$ and $\delta_{k,+}$ will be abbreviated by
$\lambda$ and $\delta_k$ for $k=1,2,3,4$, respectively.

We will solve the problem \eqref{supp-eq}--\eqref{supp-q} by a fixed point argument.
Let ${\mathscr Q}_{+}$ be given in advance such that
${\mathscr Q}_{+}\in C([0,l_+])$ satisfies
\begin{align}
\label{qub+}
c_*\le {\mathscr Q}_{+}(x)\le\overline c=\frac{c_*+\sqrt{2/(\gamma-1)}}{2},\quad0\le x\le l_+.
\end{align}
Due to \eqref{qub+} and the first formula in \eqref{a-a-0}, it holds that
\begin{gather*}
%\label{a-a-4}
{c_*}\le\Phi'_{+}(x)\le\overline c\sqrt{\delta_{2}^2+1},
\quad0\le x\le l_+,
\\
%\label{a-0-0}
\frac{\delta_{1}}{{\overline c}^{\lambda}(\delta_{2}^2+1)^{\lambda/2+1}}
\varphi^{\lambda}\le\Theta'_{+}(X_{+}(\varphi))
\le\frac{\delta_{2}}{c_*^{\lambda}}\varphi^{\lambda},\quad0\le\varphi\le\zeta_+.
\end{gather*}

\subsection{An iteration scheme}

The existence of solutions
to the problem \eqref{supp-eq}--\eqref{supp-q} will be proved 
by using the Schauder fixed point theorem.
For given ${\mathscr Q}_{+}$ and $\tilde Q$ in some suitable set,
we solve the following problem
\begin{align}
\label{itsupp-eq} &Q_{\varphi\varphi}-(b(\tilde Q)Q_{\psi})_\psi=0,
\quad&&(\varphi,\psi)\in(0,\zeta_{+})\times(0,m),
\\
\label{itsupp-inbc1}
&Q(0,\psi)=0,\quad&&\psi\in(0,m),
\\
\label{itsupp-inbc2}
&Q_\varphi(0,\psi)=0,\quad&&\psi\in(0,m),
\\
\label{itsupp-lbbc}
&Q_{\psi}(\varphi,0)=0,\quad&&\varphi\in(0,\zeta_{+}),
\\
\label{itsupp-ubbc}
&Q_{\psi}(\varphi,m)=-\frac{\Theta'_{+}(X_{+}(\varphi))\cos\Theta_{+}(X_{+}(\varphi))}
{b(\tilde Q(\varphi,m)))A^{-1}_+(\tilde Q(\varphi,m))},
\quad&&\varphi\in(0,\zeta_{+}).
\end{align}
If the problem \eqref{itsupp-eq}--\eqref{itsupp-ubbc} admits a unique
solution $Q$, then we can define a mapping as
\begin{align*}
J(({\mathscr Q}_{+},\tilde Q))=\big(A^{-1}_+(Q(\Phi_{+}(\cdot),m)),
Q(\cdot,\cdot)\big).
\end{align*}
A fixed point of the mapping $J$ is a desired solution
to the problem \eqref{supp-eq}--\eqref{supp-q}.
Note that the problem \eqref{itsupp-eq}--\eqref{itsupp-ubbc}
is equivalent to
\begin{align}
\label{b-1}
&W_\varphi+b^{1/2}(\tilde Q)W_\psi
=\frac14b^{-1}(\tilde Q)p(\tilde Q)\tilde W(W+Z),
\quad&&(\varphi,\psi)\in(0,\zeta_{+})\times(0,m),
\\
\label{b-2}
&Z_\varphi-b^{1/2}(\tilde Q)Z_\psi
=-\frac14b^{-1}(\tilde Q)p(\tilde Q)\tilde Z(W+Z),
\quad&&(\varphi,\psi)\in(0,\zeta_{+})\times(0,m),
\\
\label{b-3}
&W(0,\psi)=0,\quad&&\psi\in(0,m),
\\
\label{b-4}
&Z(0,\psi)=0,\quad&&\psi\in(0,m),
\\
\label{b-5}
&W(\varphi,0)+Z(\varphi,0)=0,\quad&&\varphi\in(0,\zeta_{+}),
\\
\label{b-6}
&W(\varphi,m)+Z(\varphi,m)=\tilde h(\varphi),
\quad&&\varphi\in(0,\zeta_{+}),
\\
\label{b-7}
&Q_\varphi(\varphi,\psi)=\frac{1}{2}\big(W(\varphi,\psi)-Z(\varphi,\psi)\big),
\quad&&(\varphi,\psi)\in(0,\zeta_{+})\times(0,m),
\\
\label{b-9}
&Q(0,\psi)=0,\quad&&\psi\in(0,m),
\end{align}
where
\begin{align}
\label{tildew}
&\tilde W(\varphi,\psi)=
\tilde Q_\varphi(\varphi,\psi)-b^{1/2}(\tilde Q(\varphi,\psi))
\tilde Q_\psi(\varphi,\psi),
&&\quad(\varphi,\psi)\in(0,\zeta_{+})\times(0,m),
\\
\label{tildez}
&\tilde Z(\varphi,\psi)=
-\tilde Q_\varphi(\varphi,\psi)-b^{1/2}(\tilde Q(\varphi,\psi))
\tilde Q_\psi(\varphi,\psi),
&&\quad(\varphi,\psi)\in(0,\zeta_{+})\times(0,m),
\\
\label{tildeh}
&\tilde h(\varphi)=\frac{2\Theta'_{+}(X_{+}(\varphi))\cos\Theta_{+}(X_{+}(\varphi))}
{b^{1/2}(\tilde Q(\varphi,m))A^{-1}_+(\tilde Q(\varphi,m))},
&&\quad\varphi\in(0,\zeta_{+}).
\end{align}

\begin{remark}
For the system \eqref{b-1}, \eqref{b-2} and \eqref{b-7},
the compatibility condition is of the form
\begin{align}
\label{b-8}
Q_\psi(\varphi,\psi)=-\frac{1}{2}b^{-1/2}(\tilde Q(\varphi,\psi))
\big(W(\varphi,\psi)+Z(\varphi,\psi)\big),
\quad(\varphi,\psi)\in(0,\zeta_{+})\times(0,m).
\end{align}
\end{remark}

It turns out to be more convenient to solve the problem \eqref{b-1}--\eqref{b-9}
instead of the problem \eqref{itsupp-eq}--\eqref{itsupp-ubbc}.
This will be done by solving  first the
problem \eqref{b-1}--\eqref{b-6} and then the problem
\eqref{b-7}, \eqref{b-9}. The problem \eqref{b-1}--\eqref{b-6} can be solved
by using the
contraction mapping theorem. For given $(w,z)$ in some suitable set,
we solve the following problem
\begin{align}
\label{d-1}
&W_\varphi+b^{1/2}(\tilde Q)W_\psi
=\frac14b^{-1}(\tilde Q)p(\tilde Q)\tilde W(W+z),
\quad&&(\varphi,\psi)\in(0,\zeta_{+})\times(0,m),
\\
\label{d-2}
&Z_\varphi-b^{1/2}(\tilde Q)Z_\psi
=-\frac14b^{-1}(\tilde Q)p(\tilde Q)\tilde Z(w+Z),
\quad&&(\varphi,\psi)\in(0,\zeta_{+})\times(0,m),
\\
\label{d-3}
&W(0,\psi)=0,\quad&&\psi\in(0,m),
\\
\label{d-4}
&Z(0,\psi)=0,\quad&&\psi\in(0,m),
\\
\label{d-5}
&W(\varphi,0)+Z(\varphi,0)=0,\quad&&\varphi\in(0,\zeta_{+}),
\\
\label{d-6}
&W(\varphi,m)+Z(\varphi,m)=\tilde h(\varphi),
\quad&&\varphi\in(0,\zeta_{+})
\end{align}
and obtain a solution $(W,Z)={\mathscr T}((w,z))$.
Then, the unique fixed point of the contraction mapping ${\mathscr T}$ is the unique solution
to the problem \eqref{b-1}--\eqref{b-6}.

In the following two subsections, we will solve the problem \eqref{b-1}--\eqref{b-9},
and consequently, give a solution to the problem \eqref{supp-eq}--\eqref{supp-q}.

\subsection{A hyperbolic system with singularity}

Fix ${\mathscr Q}_{+}\in C([0,l_+])$ satisfying \eqref{qub+}.
Assume that $\tilde Q\in C^{0,1}([0,\zeta_{+}]\times[0,m])$ is given and satisfies
\begin{gather}
\label{a-0-1}
-\mu_2\varphi^{\lambda+2}\le\tilde Q(\varphi,\psi)\le-\mu_1\varphi^{\lambda+2}
\quad(\varphi,\psi)\in(0,\zeta_{+})\times(0,m),
\\
\label{a-0-2}
-\beta_2\varphi^{\lambda+1}\le\tilde Q_\varphi(\varphi,\psi)
\le-\beta_1\varphi^{\lambda+1},
\quad
|\tilde Q_\psi(\varphi,\psi)|\le\beta_3\varphi^{3\lambda/2+1},
\quad(\varphi,\psi)\in(0,\zeta_{+})\times(0,m),
\end{gather}
where $0<\mu_1\le\mu_2$, $0<\beta_1\le\beta_2$, $\beta_3>0$.
Here and what follows, $\mu_i$, $\beta_i$, $\tau$, $\kappa$ and $\tau_i$ will 
denote generic positive constants.
Owing to $\lambda>2$, \eqref{a-a-0}, \eqref{a-0-1} and \eqref{a-0-2},
there exists $\tau\in(0,1]$ such that if $0<\zeta_+\le\tau$, then
\begin{gather}
\label{h-1}
\mu_3\varphi^{5\lambda/4+1/2}\le\tilde h(\varphi)
\le\mu_4\varphi^{5\lambda/4+1/2},
\quad\tilde h'(\varphi)\ge0,
\quad\varphi\in(0,\zeta_{+}),
\\
-{\mu_5\beta_2}\varphi^{\lambda+1}\le\tilde W(\varphi,\psi)\le0,
\quad
0\le\tilde Z(\varphi,\psi)\le{\mu_5\beta_2}\varphi^{\lambda+1},
\quad(\varphi,\psi)\in(0,\zeta_{+})\times(0,m),
\\
\label{qqq3}
-{\mu_5\beta_2}\varphi^{-1}
\le\frac14b^{-1}(\tilde Q(\varphi,\psi))p(\tilde Q(\varphi,\psi))\tilde W(\varphi,\psi)
\le0,
\quad(\varphi,\psi)\in(0,\zeta_{+})\times(0,m),
\\
\label{qqq4}
0\le\frac14b^{-1}(\tilde Q(\varphi,\psi))p(\tilde Q(\varphi,\psi))\tilde Z(\varphi,\psi)
\le{\mu_5\beta_2}\varphi^{-1},
\quad(\varphi,\psi)\in(0,\zeta_{+})\times(0,m),
\end{gather}
where the constants $\mu_i\,(i=3,4,5)\,(\mu_3\le\mu_4)$ depend on $\gamma$, $m$,
$\lambda$, $\delta_1$ and $\delta_2$, while $\tau$ also on $\beta_1$ and $\beta_3$.
The problem \eqref{b-1}--\eqref{b-9} will be solved 
under the assumptions \eqref{a-0-1}--\eqref{qqq4} in this subsection.

\begin{remark}
\label{remqqq}
The constants $\mu_1$, $\mu_2$, $\beta_1$, $\beta_2$ and $\beta_3$,
which will be determined in the proof of the existence theorem
(Theorem \ref{theoremsuper}), depend only on $\gamma$, $m$,
$\lambda$, $\delta_{1}$ and $\delta_{2}$.
Therefore, the second condition in \eqref{a-a-0},
which is only used to guarantee the second estimate in \eqref{h-1},
can be relaxed as noted in Remark \ref{qqq5}.
\end{remark}

As mentioned in the end of the last subsection,
we first consider the problem \eqref{d-1}--\eqref{d-6}.
Note that the problem \eqref{d-1}--\eqref{d-6} can be decomposed into two problems:
one is of homogeneous source terms and the other is of homogeneous boundary condtions,
which are solved by the following three lemmas.
For $(w,z)\in\mathbb R^2$, we always use $|(w,z)|$ to denote $\max\{|w|,|z|\}$ in this section.

\begin{lemma}
\label{lemma-0}
The problem
\begin{align}
\label{hd-1}
&W_\varphi+b^{1/2}(\tilde Q)W_\psi
=0,
\quad&&(\varphi,\psi)\in(0,\zeta_{+})\times(0,m),
\\
\label{hd-2}
&Z_\varphi-b^{1/2}(\tilde Q)Z_\psi
=0,
\quad&&(\varphi,\psi)\in(0,\zeta_{+})\times(0,m),
\\
\label{hd-3}
&W(0,\psi)=0,\quad&&\psi\in(0,m),
\\
\label{hd-4}
&Z(0,\psi)=0,\quad&&\psi\in(0,m),
\\
\label{hd-5}
&W(\varphi,0)+Z(\varphi,0)=0,\quad&&\varphi\in(0,\zeta_{+}),
\\
\label{hd-6}
&W(\varphi,m)+Z(\varphi,m)=\tilde h(\varphi),
\quad&&\varphi\in(0,\zeta_{+})
\end{align}
admits a unique weak solution
$(W,Z)\in L^\infty((0,\zeta_{+})\times(0,m))
\times L^\infty((0,\zeta_{+})\times(0,m))$ satisfying
\begin{align}
\label{hd-un}
\lim_{\varphi\to0^+}\big\||(W,Z)(\varphi,\cdot)|\big\|_{L^\infty((0,m))}
=0.
\end{align}
Furthermore, $(W,Z)$ satisfies
\begin{align*}
-\mu_7\varphi^{\lambda+1}
\le W(\varphi,\psi)\le-\mu_6\varphi^{\lambda+1},
\quad
\mu_6\varphi^{\lambda+1}
\le Z(\varphi,\psi)\le\mu_7\varphi^{\lambda+1},
\quad(\varphi,\psi)\in(0,\zeta_{+})\times(0,m)
\end{align*}
and
\begin{align*}
0\le W(\varphi,\psi)+Z(\varphi,\psi)\le\mu_4\varphi^{5\lambda/4+1/2},
\quad(\varphi,\psi)\in(0,\zeta_{+})\times(0,m),
\end{align*}
where $\mu_6$ and $\mu_7\,(\mu_6\le\mu_7)$ depend only on $\gamma$, $m$,
$\lambda$, $\mu_1$, $\mu_2$, $\mu_3$ and $\mu_4$.
\end{lemma}

\Proof
First consider the uniqueness.
Let $(W_1,Z_1)$ and $(W_2,Z_2)$ be two weak solutions to
the problem \eqref{hd-1}--\eqref{hd-6} satisfying \eqref{hd-un}.
Set
$$
(W,Z)(\varphi,\psi)=(W_1-W_2,Z_1-Z_2)(\varphi,\psi),\quad
(\varphi,\psi)\in(0,\zeta_{+})\times(0,m).
$$
Then, $(W,Z)$ satisfies
$$
\lim_{\varepsilon\to0^+}\big\||(W,Z)(\varepsilon,\cdot)|\big\|_{L^\infty((0,m))}
=0
$$
and is a weak solution to the following problem
\begin{align*}
&U_\varphi+b^{1/2}(\tilde Q)U_\psi
=0,
\quad&&(\varphi,\psi)\in(\varepsilon,\zeta_{+})\times(0,m),
\\
&V_\varphi-b^{1/2}(\tilde Q)V_\psi
=0,
\quad&&(\varphi,\psi)\in(\varepsilon,\zeta_{+})\times(0,m),
\\
&U(\varepsilon,\psi)=W(\varepsilon,\psi),
\quad&&\psi\in(0,m),
\\
&V(\varepsilon,\psi)=Z(\varepsilon,\psi),
\quad&&\psi\in(0,m),
\\
&U(\varphi,0)+V(\varphi,0)=0,
\quad&&\varphi\in(\varepsilon,\zeta_{+}),
\\
&U(\varphi,m)+V(\varphi,m)=0,
\quad&&\varphi\in(\varepsilon,\zeta_{+})
\end{align*}
for any $0<\varepsilon<\zeta_+$.
Hence
\begin{gather*}
|(W,Z)(\varphi,\psi)|
\le\big\||(W,Z)(\varepsilon,\cdot)|\big\|_{L^\infty((0,m))},
\quad 0<\varepsilon<\varphi<\zeta_{+},\,\psi\in(0,m).
\end{gather*}
Letting $\varepsilon\to0^+$ yields
$$
(W,Z)(\varphi,\psi)=0,\quad
(\varphi,\psi)\in(0,\zeta_{+})\times(0,m),
$$
i.e.
$$
(W_1,Z_1)(\varphi,\psi)=(W_2,Z_2)(\varphi,\psi),
\quad(\varphi,\psi)\in(0,\zeta_{+})\times(0,m).
$$

We now turn to the existence.
For any $0<\varepsilon<\zeta_+$, it follows from
the classical theory for strictly hyperbolic systems that
the problem
\begin{align}
\label{r-1}
&W^\varepsilon_\varphi+b^{1/2}(\tilde Q)W^\varepsilon_\psi
=0,
\quad&&(\varphi,\psi)\in(\varepsilon,\zeta_{+})\times(0,m),
\\
\label{r-2}
&Z^\varepsilon_\varphi-b^{1/2}(\tilde Q)Z^\varepsilon_\psi
=0,
\quad&&(\varphi,\psi)\in(\varepsilon,\zeta_{+})\times(0,m),
\\
\label{r-3}
&W^\varepsilon(\varepsilon,\psi)=0,\quad&&\psi\in(0,m),
\\
&Z^\varepsilon(\varepsilon,\psi)=0,\quad&&\psi\in(0,m),
\label{r-4}
\\
&W^\varepsilon(\varphi,0)+Z^\varepsilon(\varphi,0)=0,
\quad&&\varphi\in(\varepsilon,\zeta_{+}),
\label{r-5}
\\
\label{r-6}
&W^\varepsilon(\varphi,m)+Z^\varepsilon(\varphi,m)=\tilde h(\varphi),
\quad&&\varphi\in(\varepsilon,\zeta_{+})
\end{align}
admits a unique weak solution
$(W^\varepsilon,Z^\varepsilon)\in C^{0,1}([\varepsilon,\zeta_{+}]\times[0,m])
\times C^{0,1}([\varepsilon,\zeta_{+}]\times[0,m])$.
Moreover, it is clear from \eqref{h-1} that
\begin{gather}
\label{hd-9}
W^\varepsilon(\varphi,\psi)\le0,\quad
W^\varepsilon_\varphi(\varphi,\psi)\le0,\quad
W^\varepsilon_\psi(\varphi,\psi)\ge0,\quad
(\varphi,\psi)\in[\varepsilon,\zeta_{+}]\times[0,m],
\\
\label{hd-10}
Z^\varepsilon(\varphi,\psi)\ge0,\quad
Z^\varepsilon_\varphi(\varphi,\psi)\ge0,\quad
Z^\varepsilon_\psi(\varphi,\psi)\ge0,\quad
(\varphi,\psi)\in[\varepsilon,\zeta_{+}]\times[0,m].
\end{gather}
For any $0<\varepsilon_1<\varepsilon_2<\zeta_{+}$,
it is not hard to show that
\begin{align*}
W^{\varepsilon_1}(\varphi,\psi)
\le W^{\varepsilon_2}(\varphi,\psi),\quad
Z^{\varepsilon_1}(\varphi,\psi)
\ge Z^{\varepsilon_2}(\varphi,\psi),\quad
(\varphi,\psi)\in[\varepsilon_2,\zeta_{+}]\times[0,m].
\end{align*}
Set
\begin{align*}
(W,Z)(\varphi,\psi)=\lim_{\varepsilon\to0^+}(W^{\varepsilon},Z^{\varepsilon})(\varphi,\psi),\quad
(\varphi,\psi)\in(0,\zeta_{+})\times(0,m).
\end{align*}

\vskip10mm

\hskip45mm
\setlength{\unitlength}{0.6mm}
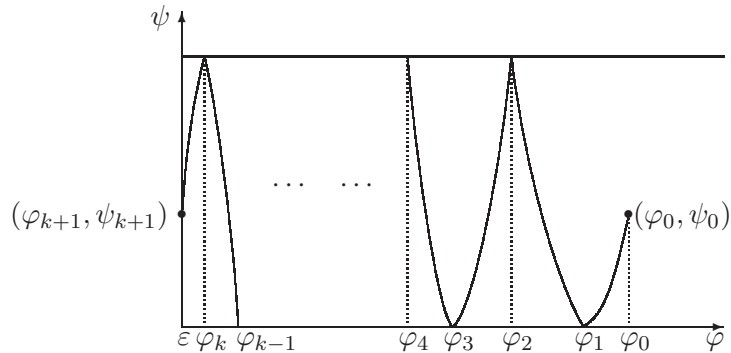
\begin{picture}(250,65)
\put(0,0){\vector(1,0){120}}
\put(0,0){\vector(0,1){70}}
\put(116,-4){$\varphi$} \put(-7,67){$\psi$}
\put(0,60){\qbezier(0,0)(55,0)(120,0)}

\put(0,25){\cbezier(0,0)(1,15)(2.5,25)(5,35)}

\put(5,60){\cbezier(0,0)(2.5,-10)(6.5,-40)(7.5,-60)}

\put(20,30){$\cdots\quad\cdots$}

\put(60,0){\cbezier(0,0)(-3,2)(-8,35)(-10,60)}
\put(60,0){\cbezier(0,0)(4,3)(10,30)(13,60)}

\put(73,60){\cbezier(0,0)(3,-25)(10,-50)(16,-60)}

\put(89,0){\cbezier(0,0)(4,3)(6,6)(10,25)}

\put(99,0){\qbezier[25](0,0)(0,12.5)(0,25)}
\put(73,0){\qbezier[60](0,0)(0,30)(0,60)}
\put(50,0){\qbezier[60](0,0)(0,30)(0,60)}
\put(5,0){\qbezier[60](0,0)(0,30)(0,60)}

\put(99,25){\circle*{1.5}}
\put(100,23){$(\varphi_0,\psi_0)$}

\put(97,-4){$\varphi_0$}
\put(87,-4){$\varphi_1$}
\put(71,-4){$\varphi_2$}
\put(58,-4){$\varphi_3$}
\put(48,-4){$\varphi_4$}

\put(12,-4){$\varphi_{k-1}$}
\put(3,-4){$\varphi_k$}
\put(-1,-4){$\varepsilon$}

\put(0,25){\circle*{1.5}}
\put(-38,23){$(\varphi_{k+1},\psi_{k+1})$}

\put(30,-17){the case that $k$ is even}

\end{picture}
\vskip23mm

\hskip45mm
\setlength{\unitlength}{0.6mm}
\begin{picture}(250,65)
\put(0,0){\vector(1,0){120}}
\put(0,0){\vector(0,1){70}}
\put(116,-4){$\varphi$} \put(-7,67){$\psi$}
\put(0,60){\qbezier(0,0)(55,0)(120,0)}

\put(0,25){\cbezier(0,0)(1,-6)(2,-24)(4,-25)}

\put(4,0){\cbezier(0,0)(3,2)(6,25)(9,60)}

\put(20,30){$\cdots\quad\cdots$}

\put(60,0){\cbezier(0,0)(-3,2)(-8,35)(-10,60)}
\put(60,0){\cbezier(0,0)(4,3)(10,30)(13,60)}

\put(73,60){\cbezier(0,0)(3,-25)(10,-50)(16,-60)}

\put(89,0){\cbezier(0,0)(4,3)(6,6)(10,25)}

\put(99,0){\qbezier[25](0,0)(0,12.5)(0,25)}
\put(73,0){\qbezier[60](0,0)(0,30)(0,60)}
\put(50,0){\qbezier[60](0,0)(0,30)(0,60)}
\put(13,0){\qbezier[60](0,0)(0,30)(0,60)}

\put(99,25){\circle*{1.5}}
\put(100,23){$(\varphi_0,\psi_0)$}

\put(97,-4){$\varphi_0$}
\put(87,-4){$\varphi_1$}
\put(71,-4){$\varphi_2$}
\put(58,-4){$\varphi_3$}
\put(48,-4){$\varphi_4$}

\put(12,-4){$\varphi_{k-1}$}
\put(3,-4){$\varphi_k$}
\put(-1,-4){$\varepsilon$}

\put(0,25){\circle*{1.5}}
\put(-38,23){$(\varphi_{k+1},\psi_{k+1})$}

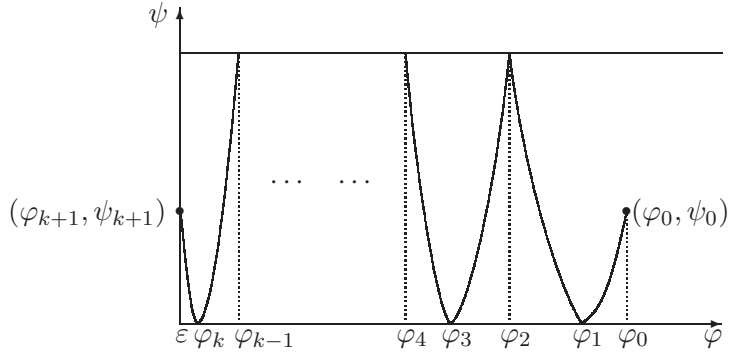
\put(30,-17){the case that $k$ is odd}

\end{picture}
\vskip17mm

Fix $(\varphi_0,\psi_0)\in(\varepsilon,\zeta_{+})\times(0,m)$.
We estimate $W(\varphi_0,\psi_0)$ by the method of characteristics.
Let $\psi=\Psi_1(\varphi)$ be the positive characteristic across $(\varphi_0,\psi_0)$,
which approaches either the initial boundary $\{\varepsilon\}\times[0,m]$
or the lower boundary $(\varepsilon,\zeta_{+})\times\{0\}$ at a point
$(\varphi_1,\psi_1)$.
If $\varphi_1>\varepsilon$, then there exists a negative characteristic
$\psi=\Psi_2(\varphi)$ from $(\varphi_1,\psi_1)$,
which approaches either the initial boundary $\{\varepsilon\}\times[0,m]$
or the upper boundary $(\varepsilon,\zeta_{+})\times\{m\}$ at a point
$(\varphi_2,\psi_2)$.
If $\varphi_2>\varepsilon$, then there exists a positive characteristic
$\psi=\Psi_3(\varphi)$ from $(\varphi_2,\psi_2)$,
which approaches either the initial boundary $\{\varepsilon\}\times[0,m]$
or the lower boundary $(\varepsilon,\zeta_{+})\times\{0\}$ at a point
$(\varphi_3,\psi_3)$. Since the system \eqref{r-1}, \eqref{r-2} is strictly hyperbolic,
there exists a nonnegative integer $k$ such that
$$
\varphi_0>\varphi_1>\cdots>\varphi_k>\varphi_{k+1}=\varepsilon,
$$
$$
\psi_j=\left\{
\begin{aligned}
&0,\quad1\le j\le k\mbox{ and $j$ is odd},
\\
&m,\quad1\le j\le k\mbox{ and $j$ is even},
\end{aligned}
\right.
\quad0\le\psi_{k+1}\le m,
$$
$$
\left\{
\begin{aligned}
&\Psi'_{j}(\varphi)=b^{1/2}(\tilde Q(\varphi,\Psi_{j}(\varphi))),
\quad\varphi_{j}<\varphi<\varphi_{j-1},
\\
&\Psi_{j}(\varphi_j)=\psi_{j},\quad \Psi_{j}(\varphi_{j-1})=\psi_{j-1},
\end{aligned}
\right.
\quad1\le j\le k+1\mbox{ and $j$ is odd},
$$
$$
\left\{
\begin{aligned}
&\Psi'_{j}(\varphi)=-b^{1/2}(\tilde Q(\varphi,\Psi_{j}(\varphi))),
\quad\varphi_{j}<\varphi<\varphi_{j-1},
\\
&\Psi_{j}(\varphi_j)=\psi_{j},\quad \Psi_{j}(\varphi_{j-1})=\psi_{j-1},
\end{aligned}
\right.
\quad1\le j\le k+1\mbox{ and $j$ is even}.
$$
Then,
\begin{align*}
W^\varepsilon(\varphi_0,\psi_0)=
-\sum_{2\le2j\le k}\tilde h(\varphi_{2j}).
\end{align*}
For $1\le j\le k+1$, it holds that
$$
\int_{\varphi_{j}}^{\varphi_{j-1}}
b^{1/2}(\tilde Q(s,\Psi_{j}(s)))ds=|\psi_{j-1}-\psi_{j}|
\left\{
\begin{aligned}
&\in(0,m),&&\mbox{ if }j=1,
\\
&=m,&&\mbox{ if }2\le j\le k,
\\
&\in(0,m],&&\mbox{ if }j=k+1.
\end{aligned}
\right.
$$
Thus, it follows from this and \eqref{a-0-1} that
$$
\varphi_{2j}^{-\lambda/4-1/2}(\varphi_{2j}-\varphi_{2(j+1)})
\le\int_{\varphi_{2(j+1)}}^{\varphi_{2j}}s^{-\lambda/4-1/2}ds
\le\hat\mu_6m,
\quad 0\le 2j\le k-1,
$$
$$
\varphi_{2j}^{-\lambda/4-1/2}(\varphi_{2(j-1)}-\varphi_{2j})
\ge\int_{\varphi_{2j}}^{\varphi_{2(j-1)}}s^{-\lambda/4-1/2}ds
\ge\hat\mu_7m,\quad 2\le 2j\le k+1,
$$
where $\hat\mu_6$ and $\hat\mu_7$ are positive constants depending only on $\gamma$,
$\lambda$, $\mu_1$ and $\mu_2$.
Hence
\begin{align*}
W^\varepsilon(\varphi_0,\psi_0)\le-\mu_3
\sum_{2\le2j\le k}\varphi_{2j}^{5\lambda/4+1/2}
\le-\frac{\mu_3}{\hat\mu_6m}
\sum_{2\le2j\le k-1}\varphi_{2j}^{\lambda}(\varphi_{2j}-\varphi_{2(j+1)})
\le-\frac{\mu_3}{\hat\mu_6m}\int_{\varphi_{k}}^{\varphi_2}t^{\lambda}dt
\\
\mbox{if }k\ge2,
\end{align*}
\begin{align*}
W^\varepsilon(\varphi_0,\psi_0)\ge-\mu_4
\sum_{2\le2j\le k}\varphi_{2j}^{5\lambda/4+1/2}
\ge-\frac{\mu_4}{\hat\mu_7m}
\sum_{2\le2j\le k}\varphi_{2j}^{\lambda}(\varphi_{2(j-1)}-\varphi_{2j})
\ge-\frac{\mu_4}{\hat\mu_7m}\int_{\varepsilon}^{\varphi_0}t^{\lambda}dt
\\
\mbox{if }k\ge0.
\end{align*}
These estimates lead to
\begin{align*}
-\mu_7\varphi_0^{\lambda+1}
\le W(\varphi_0,\psi_0)\le-\mu_6\varphi_0^{\lambda+1}
\end{align*}
with $\mu_6$ and $\mu_7\,(\mu_6\le\mu_7)$ depending only on $\gamma$, $m$,
$\lambda$, $\mu_1$, $\mu_2$, $\mu_3$ and $\mu_4$.
Therefore,
\begin{align}
\label{hd-11}
-\mu_7\varphi^{\lambda+1}
\le W(\varphi,\psi)\le-\mu_6\varphi^{\lambda+1},
\quad(\varphi,\psi)\in(0,\zeta_{+})\times(0,m).
\end{align}
Similar, one can get that
\begin{align}
\label{hd-12}
\mu_6\varphi^{\lambda+1}
\le Z(\varphi,\psi)\le\mu_7\varphi^{\lambda+1},
\quad(\varphi,\psi)\in(0,\zeta_{+})\times(0,m).
\end{align}
Moreover, it follows from \eqref{hd-9} and \eqref{hd-10} that
$$
0\le W^\varepsilon(\varphi,\psi)+Z^\varepsilon(\varphi,\psi)\le
\tilde h(\varphi),
\quad(\varphi,\psi)\in[\varepsilon,\zeta_{+}]\times[0,m],
$$
which implies
\begin{align}
\label{hd-13}
0\le W(\varphi,\psi)+Z(\varphi,\psi)\le\mu_4\varphi^{5\lambda/4+1/2},
\quad(\varphi,\psi)\in(0,\zeta_{+})\times(0,m).
\end{align}
It is not hard to show that $(W,Z)\in L^\infty((0,\zeta_{+})\times(0,m))
\times L^\infty((0,\zeta_{+})\times(0,m))$ with \eqref{hd-11}--\eqref{hd-13}
is a weak solution to the problem \eqref{hd-1}--\eqref{hd-6}
satisfying \eqref{hd-un}.
$\hfill\Box$\vskip 4mm

We now turn to the problem with nonhomogeneous source terms and homogeneous boundary conditions,
and consider first the strictly hyperbolic case and then the singular case.

\begin{lemma}
\label{lemma-0a}
Assume that $0<\varepsilon<\zeta_{+}$, $W_0,Z_0\in L^\infty((0,m))$,
and $w,z\in L^\infty((\varepsilon,\zeta_{+})\times(0,m))$ satisfy
\begin{align*}
|(w,z)(\varphi,\psi)|\le\kappa\varphi^{5\lambda/4+1/2},\quad
(\varphi,\psi)\in(\varepsilon,\zeta_+)\times(0,m).
\end{align*}
Then, the problem
\begin{align}
\label{a-1}
&W_\varphi+b^{1/2}(\tilde Q)W_\psi
=\frac14b^{-1}(\tilde Q)p(\tilde Q)\tilde W(W+z),
\quad&&(\varphi,\psi)\in(\varepsilon,\zeta_{+})\times(0,m),
\\
\label{a-2}
&Z_\varphi-b^{1/2}(\tilde Q)Z_\psi
=-\frac14b^{-1}(\tilde Q)p(\tilde Q)\tilde Z(w+Z),
\quad&&(\varphi,\psi)\in(\varepsilon,\zeta_{+})\times(0,m),
\\
\label{a-3}
&W(\varepsilon,\psi)=W_0(\psi),\quad&&\psi\in(0,m),
\\
\label{a-4}
&Z(\varepsilon,\psi)=Z_0(\psi),\quad&&\psi\in(0,m),
\\
\label{a-5}
&W(\varphi,0)+Z(\varphi,0)=0,\quad&&\varphi\in(\varepsilon,\zeta_{+}),
\\
\label{a-6}
&W(\varphi,m)+Z(\varphi,m)=0,
\quad&&\varphi\in(\varepsilon,\zeta_{+})
\end{align}
admits a unique solution $(W,Z)\in L^\infty((\varepsilon,\zeta_{+})\times(0,m))
\times L^\infty((\varepsilon,\zeta_{+})\times(0,m))$.
Moreover,
$(W,Z)$ satisfies
\begin{align*}
|(W,Z)(\varphi,\psi)|\le\frac{4\mu_5\beta_2}{5\lambda+2+4\mu_5\beta_2}\kappa\varphi^{5\lambda/4+1/2}
+\big\||(W_0,Z_0)|\big\|_{L^\infty((0,m))},\quad
(\varphi,\psi)\in(\varepsilon,\zeta_{+})\times(0,m).
\end{align*}
\end{lemma}

\Proof
According to the classical theory for strictly hyperbolic systems, the problem \eqref{a-1}--\eqref{a-6}
admits a unique weak solution $(W,Z)\in L^\infty((\varepsilon,\zeta_{+})\times(0,m))
\times L^\infty((\varepsilon,\zeta_{+})\times(0,m))$.

We first estimate $W$ along a positive characteristic.
Assume that
$$
\Sigma_+:\psi=\Psi_{+}(\varphi),\quad
\hat\varphi\le\varphi\le\check\varphi
\qquad(\varepsilon\le\hat\varphi<\check\varphi\le\zeta_{+})
$$
is a positive characteristic of \eqref{a-1}, i.e.
$$
\left\{
\begin{aligned}
&\Psi'_{+}(\varphi)=b^{1/2}(\tilde Q(\varphi,\Psi_{+}(\varphi))),
\quad\hat\varphi<\varphi<\check\varphi,
\\
&0<\Psi_{+}(\varphi)<m,\quad\hat\varphi<\varphi<\check\varphi.
\end{aligned}
\right.
$$
On $\Sigma_+$, $W$ satisfies
\begin{align*}
&\frac{d}{d\varphi}W(\varphi,\Psi_{+}(\varphi))
-\frac14b^{-1}(\tilde Q(\varphi,\Psi_{+}(\varphi)))
p(\tilde Q(\varphi,\Psi_{+}(\varphi)))\tilde W(\varphi,\Psi_{+}(\varphi))
W(\varphi,\Psi_{+}(\varphi))
\\
=&\frac14b^{-1}(\tilde Q(\varphi,\Psi_{+}(\varphi)))
p(\tilde Q(\varphi,\Psi_{+}(\varphi)))\tilde W(\varphi,\Psi_{+}(\varphi))
z(\varphi,\Psi_{+}(\varphi)),
\quad\hat\varphi<\varphi<\check\varphi,
\end{align*}
i.e.
\begin{align*}
&\Big(W(\varphi,\Psi_{+}(\varphi))
\mbox{exp}\Big\{-\frac14\int_{\hat\varphi}^{\varphi}
b^{-1}(\tilde Q(s,\Psi_{+}(s)))p(\tilde Q(s,\Psi_{+}(s)))
\tilde W(s,\Psi_{+}(s))ds\Big\}\Big)'
\\
=&-z(\varphi,\Psi_{+}(\varphi))\Big(\mbox{exp}\Big\{-\frac14\int_{\hat\varphi}^{\varphi}
b^{-1}(\tilde Q(s,\Psi_{+}(s)))p(\tilde Q(s,\Psi_{+}(s)))
\tilde W(s,\Psi_{+}(s))ds\Big\}\Big)',
\quad\hat\varphi<\varphi<\check\varphi.
\end{align*}
This leads to
\begin{align*}
&\Big|\Big(W(\varphi,\Psi_{+}(\varphi))
\mbox{exp}\Big\{-\frac14\int_{\hat\varphi}^{\varphi}
b^{-1}(\tilde Q(s,\Psi_{+}(s)))p(\tilde Q(s,\Psi_{+}(s)))
\tilde W(s,\Psi_{+}(s))ds\Big\}\Big)'\Big|
\\
\le&\kappa\varphi^{5\lambda/4+1/2}
\Big(\mbox{exp}\Big\{-\frac14\int_{\hat\varphi}^{\varphi}
b^{-1}(\tilde Q(s,\Psi_{+}(s)))p(\tilde Q(s,\Psi_{+}(s)))
\tilde W(s,\Psi_{+}(s))ds\Big\}\Big)',
\quad\hat\varphi<\varphi<\check\varphi.
\end{align*}
Then, for any $\hat\varphi\le\varphi\le\check\varphi$, it holds that
\begin{align*}
&\Big|W(\varphi,\Psi_{+}(\varphi))
\mbox{exp}\Big\{-\frac14\int_{\hat\varphi}^{\varphi}
b^{-1}(\tilde Q(s,\Psi_{+}(s)))p(\tilde Q(s,\Psi_{+}(s)))
\tilde W(s,\Psi_{+}(s))ds\Big\}-W(\hat\varphi,\Psi_{+}(\hat\varphi))\Big|
\\
\le&\kappa\int_{\hat\varphi}^\varphi
t^{5\lambda/4+1/2}\Big(\mbox{exp}\Big\{-\frac14\int_{\hat\varphi}^{t}
b^{-1}(\tilde Q(s,\Psi_{+}(s)))p(\tilde Q(s,\Psi_{+}(s)))
\tilde W(s,\Psi_{+}(s))ds\Big\}\Big)'dt
\\
=&\kappa\varphi^{5\lambda/4+1/2}\mbox{exp}\Big\{-\frac14\int_{\hat\varphi}^{\varphi}
b^{-1}(\tilde Q(s,\Psi_{+}(s)))p(\tilde Q(s,\Psi_{+}(s)))
\tilde W(s,\Psi_{+}(s))ds\Big\}
-\kappa{\hat\varphi}^{5\lambda/4+1/2}
\\
&\qquad-\Big(\frac{5}{4}\lambda+\frac12\Big)
\kappa\int_{\hat\varphi}^\varphi
t^{5\lambda/4-1/2}\mbox{exp}\Big\{-\frac14\int_{\hat\varphi}^{t}
b^{-1}(\tilde Q(s,\Psi_{+}(s)))p(\tilde Q(s,\Psi_{+}(s)))
\tilde W(s,\Psi_{+}(s))ds\Big\}dt,
\end{align*}
which yields
\begin{align}
\label{a-9}
&(\kappa\varphi^{5\lambda/4+1/2}-|W(\varphi,\Psi_{+}(\varphi))|)
\mbox{exp}\Big\{-\frac14\int_{\hat\varphi}^{\varphi}
b^{-1}(\tilde Q(s,\Psi_{+}(s)))p(\tilde Q(s,\Psi_{+}(s)))
\tilde W(s,\Psi_{+}(s))ds\Big\}
\nonumber
\\
\ge&\kappa{\hat\varphi}^{5\lambda/4+1/2}-|W(\hat\varphi,\Psi_{+}(\hat\varphi))|
\nonumber
\\
&\qquad+\Big(\frac{5}{4}\lambda+\frac12\Big)\kappa\int_{\hat\varphi}^\varphi
t^{5\lambda/4-1/2}\mbox{exp}\Big\{-\frac14\int_{\hat\varphi}^{t}
b^{-1}(\tilde Q(s,\Psi_{+}(s)))p(\tilde Q(s,\Psi_{+}(s)))
\tilde W(s,\Psi_{+}(s))ds\Big\}dt.
\end{align}
Similarly, we can estimate $Z$ along a negative characteristic.
Assume that
$$
\Sigma_-:\psi=\Psi_{-}(\varphi),\quad
\hat\varphi\le\varphi\le\check\varphi
\qquad(\varepsilon\le\hat\varphi<\check\varphi\le\zeta_{+})
$$
is a negative characteristic of \eqref{a-2}, i.e.
$$
\left\{
\begin{aligned}
&\Psi'_{-}(\varphi)=-b^{1/2}(\tilde Q(\varphi,\Psi_{-}(\varphi))),
\quad\hat\varphi<\varphi<\check\varphi,
\\
&0<\Psi_{-}(\varphi)<m,\quad\hat\varphi<\varphi<\check\varphi.
\end{aligned}
\right.
$$
On $\Sigma_-$, $Z$ satisfies
\begin{align*}
&\Big(Z(\varphi,\Psi_{-}(\varphi))
\mbox{exp}\Big\{\frac14\int_{\hat\varphi}^{\varphi}
b^{-1}(\tilde Q(s,\Psi_{-}(s)))p(\tilde Q(s,\Psi_{-}(s)))
\tilde Z(s,\Psi_{-}(s))ds\Big\}\Big)'
\\
=&-w(\varphi,\Psi_{-}(\varphi))
\Big(\mbox{exp}\Big\{\frac14\int_{\hat\varphi}^{\varphi}
b^{-1}(\tilde Q(s,\Psi_{-}(s)))p(\tilde Q(s,\Psi_{-}(s)))
\tilde Z(s,\Psi_{-}(s))ds\Big\}\Big)',
\quad\hat\varphi<\varphi<\check\varphi,
\end{align*}
which yields
\begin{align*}
&\Big|\Big(Z(\varphi,\Psi_{-}(\varphi))
\mbox{exp}\Big\{\frac14\int_{\hat\varphi}^{\varphi}
b^{-1}(\tilde Q(s,\Psi_{-}(s)))p(\tilde Q(s,\Psi_{-}(s)))
\tilde Z(s,\Psi_{-}(s))ds\Big\}\Big)'\Big|
\\
\le&\kappa\varphi^{5\lambda/4+1/2}
\Big(\mbox{exp}\Big\{\frac14\int_{\hat\varphi}^{\varphi}
b^{-1}(\tilde Q(s,\Psi_{-}(s)))p(\tilde Q(s,\Psi_{-}(s)))
\tilde Z(s,\Psi_{-}(s))ds\Big\}\Big)',
\quad\hat\varphi<\varphi<\check\varphi.
\end{align*}
Therefore, for any $\hat\varphi\le\varphi\le\check\varphi$, it holds that
\begin{align}
\label{a-10}
&(\kappa\varphi^{5\lambda/4+1/2}-|Z(\varphi,\Psi_{-}(\varphi))|)
\mbox{exp}\Big\{\frac14\int_{\hat\varphi}^{\varphi}
b^{-1}(\tilde Q(s,\Psi_{-}(s)))p(\tilde Q(s,\Psi_{-}(s)))
\tilde Z(s,\Psi_{-}(s))ds\Big\}
\nonumber
\\
\ge&\kappa{\hat\varphi}^{5\lambda/4+1/2}-|Z(\hat\varphi,\Psi_{-}(\hat\varphi))|
\nonumber
\\
&\qquad+\Big(\frac{5}{4}\lambda+\frac12\Big)\kappa\int_{\hat\varphi}^\varphi
t^{5\lambda/4-1/2}\mbox{exp}\Big\{\frac14\int_{\hat\varphi}^{t}
b^{-1}(\tilde Q(s,\Psi_{-}(s)))p(\tilde Q(s,\Psi_{-}(s)))
\tilde Z(s,\Psi_{-}(s))ds\Big\}dt.
\end{align}

Next we estimate $W$ by the method of characteristics.
As shown in the proof of Lemma \ref{lemma-0},
for any $(\varphi_0,\psi_0)\in(\varepsilon,\zeta_{+})\times(0,m)$,
there exists a nonnegative integer $k$ such that
$$
\varphi_0>\varphi_1>\cdots>\varphi_k>\varphi_{k+1}=\varepsilon,
$$
$$
\psi_j=\left\{
\begin{aligned}
&0,\quad1\le j\le k\mbox{ and $j$ is odd},
\\
&m,\quad1\le j\le k\mbox{ and $j$ is even},
\end{aligned}
\right.
\quad0\le\psi_{k+1}\le m,
$$
$$
\left\{
\begin{aligned}
&\Psi'_{j}(\varphi)=b^{1/2}(\tilde Q(\varphi,\Psi_{j}(\varphi))),
\quad\varphi_{j}<\varphi<\varphi_{j-1},
\\
&\Psi_{j}(\varphi_j)=\psi_{j},\quad \Psi_{j}(\varphi_{j-1})=\psi_{j-1},
\end{aligned}
\right.
\quad1\le j\le k+1\mbox{ and $j$ is odd},
$$
$$
\left\{
\begin{aligned}
&\Psi'_{j}(\varphi)=-b^{1/2}(\tilde Q(\varphi,\Psi_{j}(\varphi))),
\quad\varphi_{j}<\varphi<\varphi_{j-1},
\\
&\Psi_{j}(\varphi_j)=\psi_{j},\quad \Psi_{j}(\varphi_{j-1})=\psi_{j-1},
\end{aligned}
\right.
\quad1\le j\le k+1\mbox{ and $j$ is even}.
$$
One can see the figures in the proof of Lemma \ref{lemma-0}.
Define a function $r(s)\,(\varepsilon\le s<\varphi)$ as follows
$$
r(s)=\left\{
\begin{aligned}
&-\frac14b^{-1}(\tilde Q(s,\Psi_{j}(s)))p(\tilde Q(s,\Psi_{j}(s)))
\tilde W(s,\Psi_{j}(s)),\quad
\\
&\qquad\qquad\qquad\qquad\qquad\qquad\varphi_{j}\le s<\varphi_{j-1},\, 1\le j\le k+1\mbox{ and $j$ is odd},
\\
&\frac14b^{-1}(\tilde Q(s,\Psi_{j}(s)))p(\tilde Q(s,\Psi_{j}(s)))
\tilde Z(s,\Psi_{j}(s)),\quad
\\
&\qquad\qquad\qquad\qquad\qquad\qquad\varphi_{j}\le s<\varphi_{j-1},\,1\le j\le k+1\mbox{ and $j$ is even}.
\end{aligned}
\right.
$$
Then, for $1\le j\le k+1$, it follows from \eqref{a-9} and \eqref{a-10} that
\begin{align*}
&(\kappa\varphi_{j-1}^{5\lambda/4+1/2}-|W(\varphi_{j-1},\psi_{j-1})|)
\mbox{exp}\Big\{\int_{\varphi_j}^{\varphi_{j-1}}r(s)ds\Big\}
\\
\ge&\kappa{\varphi}_{j}^{5\lambda/4+1/2}-|W(\varphi_j,\psi_{j})|
+\Big(\frac{5}{4}\lambda+\frac12\Big)\kappa\int_{\varphi_j}^{\varphi_{j-1}}
t^{5\lambda/4-1/2}\mbox{exp}\Big\{\int_{\varphi_j}^{t}r(s)ds\Big\}dt,
\quad\mbox{ if $j$ is odd},
\end{align*}
\begin{align*}
&(\kappa\varphi_{j-1}^{5\lambda/4+1/2}-|Z(\varphi_{j-1},\psi_{j-1})|)
\mbox{exp}\Big\{\int_{\varphi_j}^{\varphi_{j-1}}r(s)ds\Big\}
\\
\ge&\kappa{\varphi}_{j}^{5\lambda/4+1/2}-|Z(\varphi_j,\psi_{j})|
+\Big(\frac{5}{4}\lambda+\frac12\Big)\kappa\int_{\varphi_j}^{\varphi_{j-1}}
t^{5\lambda/4-1/2}\mbox{exp}\Big\{\int_{\varphi_j}^{t}r(s)ds\Big\}dt,
\quad\mbox{ if $j$ is even},
\end{align*}
which are equivalent to
\begin{align}
\label{refin1}
&\kappa\varphi_{j-1}^{5\lambda/4+1/2}-|W(\varphi_{j-1},\psi_{j-1})|
-\Big(\frac{5}{4}\lambda+\frac12\Big)\kappa\int_{\varepsilon}^{\varphi_{j-1}}
t^{5\lambda/4-1/2}\mbox{exp}\Big\{-\int^{\varphi_{j-1}}_{t}r(s)ds\Big\}dt
\nonumber
\\
\ge&\Big(\kappa{\varphi}_{j}^{5\lambda/4+1/2}-|W(\varphi_j,\psi_{j})|
-\Big(\frac{5}{4}\lambda+\frac12\Big)\kappa\int_{\varepsilon}^{\varphi_{j}}
t^{5\lambda/4-1/2}\mbox{exp}\Big\{-\int^{\varphi_j}_{t}r(s)ds\Big\}dt\Big)
\nonumber
\\
&\qquad\cdot\mbox{exp}\Big\{-\int_{\varphi_j}^{\varphi_{j-1}}r(s)ds\Big\},
\qquad\qquad\quad\mbox{ if $j$ is odd},
\end{align}
\begin{align}
\label{refin2}
&\kappa\varphi_{j-1}^{5\lambda/4+1/2}-|Z(\varphi_{j-1},\psi_{j-1})|
-\Big(\frac{5}{4}\lambda+\frac12\Big)\kappa\int_{\varepsilon}^{\varphi_{j-1}}
t^{5\lambda/4-1/2}\mbox{exp}\Big\{-\int^{\varphi_{j-1}}_{t}r(s)ds\Big\}dt
\nonumber
\\
\ge&\Big(\kappa{\varphi}_{j}^{5\lambda/4+1/2}-|Z(\varphi_j,\psi_{j})|
-\Big(\frac{5}{4}\lambda+\frac12\Big)\kappa\int_{\varepsilon}^{\varphi_{j}}
t^{5\lambda/4-1/2}\mbox{exp}\Big\{-\int^{\varphi_j}_{t}r(s)ds\Big\}dt\Big)
\nonumber
\\
&\qquad\cdot\mbox{exp}\Big\{-\int_{\varphi_j}^{\varphi_{j-1}}r(s)ds\Big\},
\qquad\qquad\quad\mbox{ if $j$ is even}.
\end{align}
If $k$ is odd, one can get from \eqref{refin1}, \eqref{refin2}, \eqref{a-3}--\eqref{a-6} that
\begin{align}
\label{a-11}
&\kappa\varphi_{0}^{5\lambda/4+1/2}-|W(\varphi_{0},\psi_{0})|
-\Big(\frac{5}{4}\lambda+\frac12\Big)\kappa\int_{\varepsilon}^{\varphi_{0}}
t^{5\lambda/4-1/2}\mbox{exp}\Big\{-\int^{\varphi_{0}}_{t}r(s)ds\Big\}dt
\nonumber
\\
\ge&\Big(\kappa\varphi_{k+1}^{5\lambda/4+1/2}-|W(\varphi_{k+1},\psi_{k+1})|
-\Big(\frac{5}{4}\lambda+\frac12\Big)\kappa\int_{\varepsilon}^{\varphi_{k+1}}
t^{5\lambda/4-1/2}\mbox{exp}\Big\{-\int^{\varphi_{k+1}}_{t}r(s)ds\Big\}dt\Big)
\nonumber
\\
&\qquad\cdot\mbox{exp}\Big\{-\int_{\varphi_{k+1}}^{\varphi_{0}}r(s)ds\Big\}
\nonumber
\\
=&(\kappa\varepsilon^{5\lambda/4+1/2}-|W_0(\psi_{k+1})|)
\mbox{exp}\Big\{-\int_{\varepsilon}^{\varphi_0}r(s)ds\Big\}
\nonumber
\\
\ge&\kappa\varepsilon^{5\lambda/4+1/2}\mbox{exp}\Big\{-\int_{\varepsilon}^{\varphi_0}r(s)ds\Big\}
-|W_0(\psi_{k+1})|.
\end{align}
Similarly, if $k$ is even, one can get that
\begin{align}
\label{a-12}
&\kappa\varphi^{5\lambda/4+1/2}_0-|W(\varphi_0,\psi_0)|
-\Big(\frac{5}{4}\lambda+\frac12\Big)\kappa\int_{\varepsilon}^{\varphi_0}
t^{5\lambda/4-1/2}\mbox{exp}\Big\{-\int^{\varphi_0}_{t}r(s)ds\Big\}dt
\nonumber
\\
\ge&\kappa\varepsilon^{5\lambda/4+1/2}\mbox{exp}\Big\{-\int_{\varepsilon}^{\varphi_0}r(s)ds\Big\}
-|Z_0(\psi_{k+1})|.
\end{align}
Then, it follows from \eqref{a-11} and \eqref{a-12} that
\begin{align*}
|W(\varphi_0,\psi_0)|
\le&\kappa\varphi^{5\lambda/4+1/2}_0+\big\||(W_0,Z_0)|\big\|_{L^\infty((0,m))}
-\kappa\varepsilon^{5\lambda/4+1/2}\mbox{exp}\Big\{-\int_{\varepsilon}^{\varphi_0}r(s)ds\Big\}
\\
&\qquad-\Big(\frac{5}{4}\lambda+\frac12\Big)\kappa\int_{\varepsilon}^{\varphi_0}
t^{5\lambda/4-1/2}\mbox{exp}\Big\{-\int^{\varphi_0}_{t}r(s)ds\Big\}dt
\\
\le&\kappa\varphi^{5\lambda/4+1/2}_0+\big\||(W_0,Z_0)|\big\|_{L^\infty((0,m))}
-\kappa\varphi^{-\mu_5\beta_2}_0\varepsilon^{5\lambda/4+1/2+\mu_5\beta_2}
\\
&\qquad-\Big(\frac{5}{4}\lambda+\frac12\Big)\kappa\int_{\varepsilon}^{\varphi_0}
t^{5\lambda/4-1/2}\Big(\frac{\varphi_0}{t}\Big)^{-\mu_5\beta_2}dt
\\
\le&\kappa\varphi^{5\lambda/4+1/2}_0+\big\||(W_0,Z_0)|\big\|_{L^\infty((0,m))}
-\kappa\varphi^{-\mu_5\beta_2}_0\varepsilon^{5\lambda/4+1/2+\mu_5\beta_2}
\\
&\qquad-\frac{5\lambda+2}{5\lambda+2+4\mu_5\beta_2}
\kappa(\varphi^{5\lambda/4+1/2}_0
-\varphi^{-\mu_5\beta_2}_0\varepsilon^{5\lambda/4+1/2+\mu_5\beta_2})
\\
\le&\frac{4\mu_5\beta_2}{5\lambda+2+4\mu_5\beta_2}\kappa\varphi^{5\lambda/4+1/2}_0
+\big\||(W_0,Z_0)|\big\|_{L^\infty((0,m))}.
\end{align*}
The estimate of $Z$ can be proved similarly.
$\hfill\Box$\vskip 4mm

\begin{lemma}
\label{lemma-0ann}
Assume that $(w,z)\in L^\infty((0,\zeta_{+})\times(0,m))$ satisfies
\begin{align}
\label{qqq8}
|(w,z)(\varphi,\psi)|\le\kappa\varphi^{5\lambda/4+1/2},\quad
(\varphi,\psi)\in(0,\zeta_+)\times(0,m).
\end{align}
Then, the problem
\begin{align}
\label{ann-1}
&W_\varphi+b^{1/2}(\tilde Q)W_\psi
=\frac14b^{-1}(\tilde Q)p(\tilde Q)\tilde W(W+z),
\quad&&(\varphi,\psi)\in(0,\zeta_{+})\times(0,m),
\\
\label{ann-2}
&Z_\varphi-b^{1/2}(\tilde Q)Z_\psi
=-\frac14b^{-1}(\tilde Q)p(\tilde Q)\tilde Z(w+Z),
\quad&&(\varphi,\psi)\in(0,\zeta_{+})\times(0,m),
\\
\label{ann-3}
&W(0,\psi)=0,\quad&&\psi\in(0,m),
\\
\label{ann-4}
&Z(0,\psi)=0,\quad&&\psi\in(0,m),
\\
\label{ann-5}
&W(\varphi,0)+Z(\varphi,0)=0,\quad&&\varphi\in(0,\zeta_{+}),
\\
\label{ann-6}
&W(\varphi,m)+Z(\varphi,m)=0,
\quad&&\varphi\in(0,\zeta_{+})
\end{align}
admits a unique solution $(W,Z)\in L^\infty((0,\zeta_{+})\times(0,m))
\times L^\infty((0,\zeta_{+})\times(0,m))$ satisfying
\begin{align}
\label{hd-unnn}
\lim_{\varphi\to0^+}\big\||(W,Z)(\varphi,\cdot)|\big\|_{L^\infty((0,m))}
=0.
\end{align}
Moreover,
$(W,Z)$ satisfies
\begin{align*}
|(W,Z)(\varphi,\psi)|\le\frac{4\mu_5\beta_2}{5\lambda+2+4\mu_5\beta_2}
\kappa\varphi^{5\lambda/4+1/2},\quad
(\varphi,\psi)\in(0,\zeta_{+})\times(0,m).
\end{align*}
\end{lemma}

\Proof
We start with the uniqueness.
Let $(W_1,Z_1)$ and $(W_2,Z_2)$ be two weak solutions to
the problem \eqref{ann-1}--\eqref{ann-6} satisfying \eqref{hd-unnn}.
Set
$$
(W,Z)(\varphi,\psi)=(W_1-W_2,Z_1-Z_2)(\varphi,\psi),\quad
(\varphi,\psi)\in(0,\zeta_{+})\times(0,m).
$$
Then, $(W,Z)$ satisfies
$$
\lim_{\varphi\to0^+}\big\||(W,Z)(\varphi,\cdot)|\big\|_{L^\infty((0,m))}
=0
$$
and is a weak solution to the following problem
\begin{align*}
&U_\varphi+b^{1/2}(\tilde Q)U_\psi
=\frac14b^{-1}(\tilde Q)p(\tilde Q)\tilde W U,
\quad&&(\varphi,\psi)\in(\varepsilon,\zeta_{+})\times(0,m),
\\
&V_\varphi-b^{1/2}(\tilde Q)V_\psi
=-\frac14b^{-1}(\tilde Q)p(\tilde Q)\tilde Z V,
\quad&&(\varphi,\psi)\in(\varepsilon,\zeta_{+})\times(0,m),
\\
&U(0,\psi)=W(\varepsilon,\psi),\quad&&\psi\in(0,m),
\\
&V(0,\psi)=Z(\varepsilon,\psi),\quad&&\psi\in(0,m),
\\
&U(\varphi,0)+V(\varphi,0)=0,\quad&&\varphi\in(\varepsilon,\zeta_{+}),
\\
&U(\varphi,m)+V(\varphi,m)=0,\quad&&\varphi\in(\varepsilon,\zeta_{+})
\end{align*}
for any $0<\varepsilon<\zeta_+$.
It follows from Lemma \ref{lemma-0a} that
\begin{gather*}
|(W,Z)(\varphi,\psi)|
\le\big\||(W,Z)(\varepsilon,\cdot)|\big\|_{L^\infty((0,m))},
\quad 0<\varepsilon<\varphi<\zeta_{+},\,\psi\in(0,m).
\end{gather*}
Letting $\varepsilon\to0^+$ yields
$$
(W,Z)(\varphi,\psi)=(0,0),\quad
(\varphi,\psi)\in(0,\zeta_{+})\times(0,m),
$$
i.e.
$$
(W_1,Z_1)(\varphi,\psi)=(W_2,Z_2)(\varphi,\psi),\quad
(\varphi,\psi)\in(0,\zeta_{+})\times(0,m).
$$

We now turn to the existence.
For any $0<\varepsilon<\zeta_+$, it follows from Lemma \ref{lemma-0a} that
the following problem
\begin{align*}
&W_\varphi^\varepsilon+b^{1/2}(\tilde Q)W_\psi^\varepsilon
=\frac14b^{-1}(\tilde Q)p(\tilde Q)\tilde W(W^\varepsilon+z),
\quad&&(\varphi,\psi)\in(\varepsilon,\zeta_{+})\times(0,m),
\\
&Z_\varphi^\varepsilon-b^{1/2}(\tilde Q)Z_\psi^\varepsilon
=-\frac14b^{-1}(\tilde Q)p(\tilde Q)\tilde Z(w+Z^\varepsilon),
\quad&&(\varphi,\psi)\in(\varepsilon,\zeta_{+})\times(0,m),
\\
&W^\varepsilon(\varepsilon,\psi)=0,\quad&&\psi\in(0,m),
\\
&Z^\varepsilon(\varepsilon,\psi)=0,\quad&&\psi\in(0,m),
\\
&W^\varepsilon(\varphi,0)+Z^\varepsilon(\varphi,0)=0,
\quad&&\varphi\in(\varepsilon,\zeta_{+}),
\\
&W^\varepsilon(\varphi,m)+Z^\varepsilon(\varphi,m)=0,
\quad&&\varphi\in(\varepsilon,\zeta_{+})
\end{align*}
admits a unique solution $(W^\varepsilon,Z^\varepsilon)
\in L^\infty((\varepsilon,\zeta_{+})\times(0,m))
\times L^\infty((\varepsilon,\zeta_{+})\times(0,m))$;
furthermore, $(W^\varepsilon,Z^\varepsilon)$ satisfies
\begin{align}
\label{nn1}
|(W^\varepsilon,Z^\varepsilon)(\varphi,\psi)|\le\frac{4\mu_5\beta_2}{5\lambda+2+4\mu_5\beta_2}
\kappa\varphi^{5\lambda/4+1/2},\quad
(\varphi,\psi)\in(\varepsilon,\zeta_{+})\times(0,m).
\end{align}
Define
$$
(W^\varepsilon,Z^\varepsilon)(\varphi,\psi)=(0,0),
\quad(\varphi,\psi)\in(0,\varepsilon)\times(0,m).
$$
Then, there exists a decreasing sequence
$\{\varepsilon_n\}_{n=1}^\infty\subset(0,\zeta_{+})$,
which tends to $0$ as $n\to\infty$,
and $W,Z\in L^\infty((0,\zeta_{+})\times(0,m))$ such that
$$
(W^{\varepsilon_n},Z^{\varepsilon_n})\rightharpoonup (W,Z)
\mbox{ weakly * in }L^\infty((0,\zeta_{+})\times(0,m))\times L^\infty((0,\zeta_{+})\times(0,m))
\mbox{ as }n\to\infty.
$$
Moreover, it follows from \eqref{nn1} that
\begin{align}
\label{nn3}
|(W,Z)(\varphi,\psi)|\le\frac{4\mu_5\beta_2}{5\lambda+2+4\mu_5\beta_2}
\kappa\varphi^{5\lambda/4+1/2},\quad
(\varphi,\psi)\in(0,\zeta_{+})\times(0,m).
\end{align}
It is not hard to show that $(W,Z)$ with \eqref{nn3}
is a weak solution to the problem \eqref{ann-1}--\eqref{ann-6}
satisfying \eqref{hd-unnn}.
$\hfill\Box$\vskip 4mm

\begin{remark}
\label{qqq9}
In Lemma \ref{lemma-0ann}, if \eqref{qqq8} is replaced by
\begin{align*}
|w(\varphi,\psi)|\le\kappa\varphi^{\alpha},\quad
|z(\varphi,\psi)|\le\kappa\varphi^{\alpha},\quad
(\varphi,\psi)\in(0,\zeta_+)\times(0,m)
\end{align*}
with a positive constant $\alpha$.
then the solution $(W,Z)$ to the problem \eqref{ann-1}--\eqref{ann-6} satisfies
\begin{gather*}
|(W,Z)(\varphi,\psi)|\le\frac{\mu_5\beta_2}{\alpha+\mu_5\beta_2}\kappa\varphi^{\alpha},\quad
(\varphi,\psi)\in(0,\zeta_{+})\times(0,m).
\end{gather*}
\end{remark}

To solve the problem \eqref{d-1}--\eqref{d-6},
we introduce a norm
\begin{align*}
\|(w,z)\|_{{\mathscr B}}=\sup\Big\{
\varphi^{-5\lambda/4-1/2}|(w,z)(\varphi,\psi)|:(\varphi,\psi)\in(0,\zeta_{+})\times(0,m)\Big\}
\end{align*}
and define
$$
{\mathscr B}=\Big\{(w,z)\in L^\infty((0,\zeta_{+})\times(0,m))
\times L^\infty((0,\zeta_{+})\times(0,m)):
\|(w-W^*,z-Z^*)\|_{\mathscr B}<+\infty\Big\},
$$
where $(W^*,Z^*)$ is the unique weak solution to the problem \eqref{hd-1}--\eqref{hd-6}
satisfying \eqref{hd-un}.
It follows from Lemma \ref{lemma-0} that
\begin{align}
\label{h-2}
-\mu_7\varphi^{\lambda+1}
\le W^*(\varphi,\psi)\le-\mu_6\varphi^{\lambda+1},
\quad
\mu_6\varphi^{\lambda+1}
\le Z^*(\varphi,\psi)\le\mu_7\varphi^{\lambda+1},
\quad(\varphi,\psi)\in(0,\zeta_{+})\times(0,m)
\end{align}
and
\begin{align}
\label{h-3}
0\le W^*(\varphi,\psi)+Z^*(\varphi,\psi)\le\mu_4\varphi^{5\lambda/4+1/2},
\quad(\varphi,\psi)\in(0,\zeta_{+})\times(0,m).
\end{align}
Then we have

\begin{proposition}
\label{proposition-a}
For given $(w,z)\in{\mathscr B}$,
the problem \eqref{d-1}--\eqref{d-6} admits a unique weak solution
$(W,Z)\in{\mathscr B}$.
Furthermore, if $\|(w-W^*,z-Z^*)\|_{\mathscr B}\le\mu_8\beta_2$ additionally,
then $\|(W-W^*,Z-Z^*)\|_{\mathscr B}\le\mu_8\beta_2$, where
$\mu_8={4\mu_4\mu_5}/(5\lambda+2)$.
\end{proposition}

\Proof
Since the uniqueness follows from a similar proof in Lemma \ref{lemma-0ann},
it suffices to  consider the existence.
Consider the following problem
\begin{align}
\label{scrd-1}
&{\mathscr W}_\varphi+b^{1/2}(\tilde Q){\mathscr W}_\psi
=\frac14b^{-1}(\tilde Q)p(\tilde Q)\tilde W
\big({\mathscr W}+(W^*+Z^*)+(z-Z^*)\big),
&&(\varphi,\psi)\in(0,\zeta_{+})\times(0,m),
\\
\label{scrd-2}
&{\mathscr Z}_\varphi-b^{1/2}(\tilde Q){\mathscr Z}_\psi
=-\frac14b^{-1}(\tilde Q)p(\tilde Q)\tilde Z
\big((w-W^*)+(W^*+Z^*)+{\mathscr Z}),
&&(\varphi,\psi)\in(0,\zeta_{+})\times(0,m),
\\
\label{scrd-3}
&{\mathscr W}(0,\psi)=0,\quad&&\psi\in(0,m),
\\
\label{scrd-4}
&{\mathscr Z}(0,\psi)=0,\quad&&\psi\in(0,m),
\\
\label{scrd-5}
&{\mathscr W}(\varphi,0)+{\mathscr Z}(\varphi,0)=0,
\quad&&\varphi\in(0,\zeta_{+}),
\\
\label{scrd-6}
&{\mathscr W}(\varphi,m)+{\mathscr Z}(\varphi,m)=0,
\quad&&\varphi\in(0,\zeta_{+}).
\end{align}
Thanks to Lemma \ref{lemma-0ann},
the problem \eqref{scrd-1}--\eqref{scrd-6} admits a unique weak solution
$({\mathscr W},{\mathscr Z})\in L^\infty((0,\zeta_{+})\times(0,m))
\times L^\infty((0,\zeta_{+})\times(0,m))$ satisfying
\begin{align}
\label{scrd-9}
|({\mathscr W},{\mathscr Z})(\varphi,\psi)|\le
\frac{4\mu_5\beta_2}{5\lambda+2+4\mu_5\beta_2}
\big(\mu_4+\|(w-W^*,z-Z^*)\|_{\mathscr B}\big)\varphi^{5\lambda/4+1/2},
\nonumber
\\
(\varphi,\psi)\in(0,\zeta_{+})\times(0,m).
\end{align}
Therefore, $({\mathscr W}+W^*,{\mathscr Z}+Z^*)\in{\mathscr B}$
and it is not hard to verify that
$({\mathscr W}+W^*,{\mathscr Z}+Z^*)$
is a weak solution to the problem \eqref{d-1}--\eqref{d-6}.
Finally, the rest of the proposition follows from
\eqref{scrd-9} directly.
$~\hfill\Box$\vskip 4mm

We are now ready to solve the problem \eqref{b-1}--\eqref{b-9}.

\begin{proposition}
\label{prop-b}
For given ${\mathscr Q}_{+}\in C([0,l_+])$ with \eqref{qub+}
and given $\tilde Q\in C^{0,1}([0,\zeta_{+}]\times[0,m])$
with \eqref{a-0-1}--\eqref{qqq4},
the problem \eqref{b-1}--\eqref{b-9}
admits a unique weak solution $(W,Z,Q)$ with $(W,Z)\in{\mathscr B}$.
Furthermore, the following estimates hold:
\begin{gather}
\label{b-11}
-\mu_7\varphi^{\lambda+1}-\mu_8\beta_2\varphi^{5\lambda/4+1/2}
\le Q_\varphi(\varphi,\psi)\le-\mu_6\varphi^{\lambda+1}+\mu_8\beta_2\varphi^{5\lambda/4+1/2},
\quad(\varphi,\psi)\in(0,\zeta_{+})\times(0,m),
\\
\label{b-12}
|Q_\psi(\varphi,\psi)|\le
\mu_9(\beta_2+1)\varphi^{3\lambda/2+1},
\quad(\varphi,\psi)\in(0,\zeta_{+})\times(0,m)
\end{gather}
and
\begin{align}
\label{b-13}
-\sigma_2\varphi^{\lambda+2}
-m\mu_9(\beta_2+1)\varphi^{3\lambda/2+1}
\le Q(\varphi,\psi)\le
-\sigma_1\varphi^{\lambda+2}
+m\mu_9(\beta_2+1)\varphi^{3\lambda/2+1},
\nonumber
\\
\quad
(\varphi,\psi)\in(0,\zeta_+)\times(0,m),
\end{align}
where $\mu_9$ depends only on
$\gamma$, $\lambda$, $\mu_1$, $\mu_2$, $\mu_4$ and $\mu_8$,
while $\sigma_1$ and $\sigma_2\,(\sigma_1\le\sigma_2)$ depend only on $\gamma$, $m$,
$\lambda$, $\delta_1$ and $\delta_2$.
\end{proposition}

\Proof
Let $(w,z)\in{\mathscr B}$ be given. It follows from Proposition \ref{proposition-a} that
the problem \eqref{d-1}--\eqref{d-6} admits a unique weak solution
$(W,Z)\in{\mathscr B}$. Define the mapping
$$
{\mathscr T}:{\mathscr B}\longrightarrow{\mathscr B},
\quad (w,z)\longmapsto(W,Z).
$$
For given $(w_1,z_1),(w_2,z_2)\in{\mathscr B}$,
it follows from the definition of ${\mathscr T}$ that
$$
{\mathscr T}((w_1,z_1))(\varphi,\psi)-{\mathscr T}((w_2,z_2))(\varphi,\psi)
=(U,V)(\varphi,\psi),
\quad(\varphi,\psi)\in(0,\zeta_{+})\times(0,m),
$$
where $(U,V)\in L^\infty((0,\zeta_{+})\times(0,m))
\times L^\infty((0,\zeta_{+})\times(0,m))$ solves the following problem
\begin{align*}
&U_\varphi+b^{1/2}(\tilde Q)U_\psi
=\frac14b^{-1}(\tilde Q)p(\tilde Q)\tilde W(U+(z_1-z_2)),
\quad&&(\varphi,\psi)\in(0,\zeta_{+})\times(0,m),
\\
&V_\varphi-b^{1/2}(\tilde Q)V_\psi
=-\frac14b^{-1}(\tilde Q)p(\tilde Q)\tilde Z((w_1-w_2)+V),
\quad&&(\varphi,\psi)\in(0,\zeta_{+})\times(0,m),
\\
&U(0,\psi)=0,\quad&&\psi\in(0,m),
\\
&V(0,\psi)=0,\quad&&\psi\in(0,m),
\\
&U(\varphi,0)+V(\varphi,0)=0,
\quad&&\varphi\in(0,\zeta_{+}),
\\
&U(\varphi,m)+V(\varphi,m)=0,
\quad&&\varphi\in(0,\zeta_{+})
\end{align*}
satisfying
$$
\lim_{\varphi\to0^+}\big\||(U,V)(\varphi,\cdot)|\big\|_{L^\infty((0,m))}
=0.
$$
Moreover, Lemma \ref{lemma-0ann} leads to that
$(U,V)$ satisfies
\begin{align*}
|(U,V)(\varphi,\psi)|\le
\frac{4\mu_5\beta_2}{5\lambda+2+4\mu_5\beta_2}
\|(w_1-w_2,z_1-z_2)\|_{\mathscr B}\varphi^{5\lambda/4+1/2},\quad
(\varphi,\psi)\in(0,\zeta_{+})\times(0,m).
\end{align*}
Therefore,
$$
\|{\mathscr T}((w_1,z_1))-{\mathscr T}((w_2,z_2))\|_{\mathscr B}\le
\frac{4\mu_5\beta_2}{5\lambda+2+4\mu_5\beta_2}
\|(w_1-w_2,z_1-z_2)\|_{\mathscr B},
$$
which shows that ${\mathscr T}$ is a contraction mapping on ${\mathscr B}$.
Thus, the problem \eqref{b-1}--\eqref{b-6}
admits a unique weak solution $(W,Z)\in{\mathscr B}$;
furthermore, one gets from Proposition \ref{proposition-a} that
\begin{align}
\label{nn5}
\|(W-W^*,Z-Z^*)\|_{\mathscr B}\le\mu_8\beta_2.
\end{align}
Then, we can get the weak solution $Q\in C^{0,1}([0,\zeta_{+}]\times[0,m])$
to the problem \eqref{b-7}, \eqref{b-9}.
And this $(W,Z,Q)$ is the unique weak solution
to the problem \eqref{b-1}--\eqref{b-9}
with $(W,Z)\in{\mathscr B}$.
It follows from \eqref{nn5}, \eqref{h-2} and \eqref{h-3} that
\begin{align*}
-\mu_7\varphi^{\lambda+1}-\mu_8\beta_2\varphi^{5\lambda/4+1/2}
\le Q_\varphi(\varphi,\psi)\le-\mu_6\varphi^{\lambda+1}+\mu_8\beta_2\varphi^{5\lambda/4+1/2},
\quad(\varphi,\psi)\in(0,\zeta_{+})\times(0,m),
\end{align*}
\begin{align*}
-\frac12\mu_4\varphi^{5\lambda/4+1/2}-\mu_8\beta_2\varphi^{5\lambda/4+1/2}
\le b^{1/2}(\tilde Q)Q_\psi(\varphi,\psi)\le
\mu_8\beta_2\varphi^{5\lambda/4+1/2},
\quad(\varphi,\psi)\in(0,\zeta_{+})\times(0,m),
\end{align*}
which yield \eqref{b-11} and \eqref{b-12} immediately.

Finally, we prove \eqref{b-13}.
Owing to \eqref{b-1}, \eqref{b-2}, \eqref{b-7} and \eqref{b-8},
it holds in the sense of distribution that
\begin{align*}
Q_{\varphi\varphi}(\varphi,\psi)
=&\frac{1}{2}\big(W_\varphi(\varphi,\psi)-Z_\varphi(\varphi,\psi)\big)
\\
=&-\frac12{b^{1/2}(\tilde
Q(\varphi,\psi))}\big(W_\psi(\varphi,\psi)+Z_\psi(\varphi,\psi)\big)
\\
&\qquad-\frac14b^{-1/2}(\tilde Q(\varphi,\psi))p(\tilde
Q(\varphi,\psi))\tilde
Q_\psi(\varphi,\psi)\big(W(\varphi,\psi)+Z(\varphi,\psi)\big)
\\
=&b^{1/2}(\tilde Q(\varphi,\psi))\Big({b^{1/2}(\tilde
Q(\varphi,\psi))}Q_\psi(\varphi,\psi)\Big)_\psi
+\Big(b^{1/2}(\tilde Q(\varphi,\psi))\Big)_\psi{b^{1/2}(\tilde
Q(\varphi,\psi))}Q_\psi(\varphi,\psi)
\\
=&\big({b(\tilde Q(\varphi,\psi))}Q_\psi(\varphi,\psi)\big)_\psi,
\quad(\varphi,\psi)\in(0,\zeta_{+})\times(0,m).
\end{align*}
It follows from this formula, \eqref{b-5} and \eqref{b-6} that
\begin{align*}
\frac{d^2}{d\varphi^2}\int_0^m Q(\varphi,\psi)d\psi
=-\frac{\Theta'_{+}(X_{+}(\varphi))\cos\Theta_{+}(X_{+}(\varphi))}
{A^{-1}_+(\tilde Q(\varphi,m))},
\quad0<\varphi<\zeta_+,
\end{align*}
which, together with \eqref{b-11} and \eqref{b-9}, yields
\begin{align*}
\int_0^m Q(\varphi,\psi)d\psi
=-\int_0^\varphi\int_0^t\frac{\Theta'_{+}(X_{+}(s))\cos\Theta_{+}(X_{+}(s))}
{A^{-1}_+(\tilde Q(s,m))}dsdt,
\quad0<\varphi<\zeta_+.
\end{align*}
Then, there exist $\sigma_1$ and $\sigma_2\,(\sigma_1\le\sigma_2)$ depending only on
$\gamma$, $m$, $\lambda$, $\delta_1$ and $\delta_2$ such that
\begin{align*}
-\sigma_2\varphi^{\lambda+2}\le\frac1m\int_0^m Q(\varphi,\psi)d\psi\le
-\sigma_1\varphi^{\lambda+2},\quad0<\varphi<\zeta_+,
\end{align*}
which, together with \eqref{b-12}, leads to \eqref{b-13}.
$\hfill\Box$\vskip 4mm

\subsection{Existence of sonic-supersonic flows}

Based on the well-posedness of the linearized problem \eqref{b-1}--\eqref{b-9},
which is equivalent to \eqref{itsupp-eq}--\eqref{itsupp-ubbc},
we are going to prove the existence of weak solutions to the nonlinear problem \eqref{supp-eq}--\eqref{supp-q}.

Assume that $f_+\in C^{3}([0,l_+])$ satisfies \eqref{a-a-0}
with $0<l_+\le\delta_0$,
where $\delta_0$ is a positive constant to be determined below.
Let ${\mathscr Q}_{+}\in C([0,l_+])$ and ${\mathscr Q}\in C^{0,1}([0,1]\times[0,m])$
be given such that \eqref{qub+} holds and
\begin{gather}
\label{c-1}
-2\sigma_2\Big(\overline c\sqrt{\delta_2^2+1}l_+\Big)^{\lambda+2}\phi^{\lambda+2}
\le{\mathscr Q}(\phi,\psi)
\le-\frac12\sigma_1(c_*l_+)^{\lambda+2}\phi^{\lambda+2},
\quad(\phi,\psi)\in(0,1)\times(0,m),
\\
\label{c-2-0}
-\beta_2\Big(\overline c\sqrt{\delta_2^2+1}l_+\Big)^{\lambda+2}\phi^{\lambda+1}
\le{\mathscr Q}_\phi(\phi,\psi)\le-\beta_1(c_*l_+)^{\lambda+2}\phi^{\lambda+1},
\quad(\phi,\psi)\in(0,1)\times(0,m),
\\
\label{c-2}
|{\mathscr Q}_\psi(\phi,\psi)|\le\beta_3
\Big(\overline c\sqrt{\delta_2^2+1}l_+\Big)^{3\lambda/2+1}
\phi^{3\lambda/2+1},
\quad(\phi,\psi)\in(0,1)\times(0,m),
\end{gather}
where $\sigma_1$ and $\sigma_2$ are defined in Proposition
\ref{prop-b}, which depend only on $\gamma$, $m$,
$\lambda$, $\delta_1$ and $\delta_2$,
while $\beta_i\,(i=1,2,3)\,(\beta_1\le\beta_2)$ will be determined below.
Define
$$
\tilde Q(\varphi,\psi)={\mathscr Q}\Big(\frac{\varphi}{\zeta_+},\psi\Big),
\quad(\varphi,\psi)\in(0,\zeta_{+})\times(0,m),
$$
where $\zeta_+$ is given in \eqref{zeta+}.
Then, $\tilde Q\in C^{0,1}([0,\zeta_{+}]\times[0,m])$ satisfies
\begin{gather*}
-2\sigma_2\Big(\frac{\overline c\sqrt{\delta_2^2+1}}{c_*}\Big)^{\lambda+2}\varphi^{\lambda+2}
\le\tilde Q(\varphi,\psi)\le
-\frac12\sigma_1\Big(\frac{c_*}{\overline c\sqrt{\delta_2^2+1}}\Big)^{\lambda+2}\varphi^{\lambda+2},
\quad(\varphi,\psi)\in(0,\zeta_{+})\times(0,m),
\\
-\beta_2\Big(\frac{\overline c\sqrt{\delta_2^2+1}}{c_*}\Big)^{\lambda+2}
\varphi^{\lambda+1}\le\tilde Q_\varphi(\varphi,\psi)
\le-\beta_1\Big(\frac{c_*}{\overline c\sqrt{\delta_2^2+1}}\Big)^{\lambda+2}\varphi^{\lambda+1},
\quad(\varphi,\psi)\in(0,\zeta_{+})\times(0,m),
\\
|\tilde Q_\psi(\varphi,\psi)|\le
\beta_3\Big(\frac{\overline c\sqrt{\delta_2^2+1}}{c_*}\Big)^{3\lambda/2+1}\varphi^{3\lambda/2+1},
\quad(\varphi,\psi)\in(0,\zeta_{+})\times(0,m).
\end{gather*}
It follows from these estimates and \eqref{a-a-0} that
there exists $\tau_1\in(0,1]$ such that if $0<\zeta_+\le\tau_1$, then
\begin{gather*}
\mu_3\varphi^{5\lambda/4+1/2}\le\tilde h(\varphi)
\le\mu_4\varphi^{5\lambda/4+1/2},
\quad\tilde h'(\varphi)\ge0,
\quad\varphi\in(0,\zeta_{+}),
\\
-{\mu_5\beta_2}\varphi^{\lambda+1}\le\tilde W(\varphi,\psi)\le0,
\quad
0\le\tilde Z(\varphi,\psi)\le{\mu_5\beta_2}\varphi^{\lambda+1},
\quad(\varphi,\psi)\in(0,\zeta_{+})\times(0,m),
\\
-{\mu_5\beta_2}\varphi^{-1}
\le\frac14b^{-1}(\tilde Q(\varphi,\psi))p(\tilde Q(\varphi,\psi))\tilde W(\varphi,\psi)
\le0,
\quad(\varphi,\psi)\in(0,\zeta_{+})\times(0,m),
\\
0\le\frac14b^{-1}(\tilde Q(\varphi,\psi))p(\tilde Q(\varphi,\psi))\tilde Z(\varphi,\psi)
\le{\mu_5\beta_2}\varphi^{-1},
\quad(\varphi,\psi)\in(0,\zeta_{+})\times(0,m),
\end{gather*}
where $\tilde h$, $\tilde W$ and $\tilde Z$ are defined by 
\eqref{tildeh}, \eqref{tildew} and \eqref{tildez}, respectively,
$\mu_i\,(i=3,4,5)\,(\mu_3\le\mu_4)$ depend on $\gamma$, $m$,
$\lambda$, $\delta_1$ and $\delta_2$, while $\tau_1$ also on
$\beta_1$ and $\beta_3$.
According to Proposition \ref{prop-b}, the problem \eqref{b-1}--\eqref{b-9}
admits a unique weak solution $(W,Z,Q)$ with $(W,Z)\in{\mathscr B}$
if $0<\zeta_+\le\tau_1$;
furthermore, $Q$ satisfies \eqref{b-11}--\eqref{b-13}.
That is to say, the problem \eqref{itsupp-eq}--\eqref{itsupp-ubbc}
admits a unique weak solution $Q$ with $(W,Z)\in{\mathscr B}$
if $0<\zeta_+\le\tau_1$, and $Q$ satisfies
\begin{gather}
\label{c-4}
-\mu_7\varphi^{\lambda+1}-\mu_8\beta_2\varphi^{5\lambda/4+1/2}
\le Q_\varphi(\varphi,\psi)\le-\mu_6\varphi^{\lambda+1}+\mu_8\beta_2\varphi^{5\lambda/4+1/2},
\quad(\varphi,\psi)\in(0,\zeta_{+})\times(0,m),
\\
\label{c-5}
|Q_\psi(\varphi,\psi)|\le
\mu_9(\beta_2+1)\varphi^{3\lambda/2+1},
\quad(\varphi,\psi)\in(0,\zeta_{+})\times(0,m)
\end{gather}
and
\begin{align}
\label{c-6}
-\sigma_2\varphi^{\lambda+2}
-m\mu_9(\beta_2+1)\varphi^{3\lambda/2+1}
\le Q(\varphi,\psi)\le
-\sigma_1\varphi^{\lambda+2}
+m\mu_9(\beta_2+1)\varphi^{3\lambda/2+1},
\nonumber
\\
\quad
(\varphi,\psi)\in(0,\zeta_+)\times(0,m),
\end{align}
where $\mu_i\,(i=6,7,8,9)\,(\mu_6\le\mu_7)$ depend only on
$\gamma$, $m$,
$\lambda$, $\delta_1$ and $\delta_2$.
Choose
\begin{align}
\label{beta}
\beta_1=\frac12\mu_6,\quad
\beta_2=2\mu_7,\quad
\beta_3=\mu_9(\beta_2+1)
\end{align}
and
\begin{align}
\label{tau2}
\tau_2=\min\Big\{\Big(\frac{\mu_6}{2\mu_8\beta_2}\Big)^{4/(\lambda-2)},
\Big(\frac{\sigma_1}{2m\mu_9(\beta_2+1)}\Big)^{2/(\lambda-2)},
\Big(\frac{-A(\overline c)}{2\sigma_2}\Big)^{1/(\lambda+2)}\Big\},
\end{align}
which all depend only on $\gamma$, $m$,
$\lambda$, $\delta_1$ and $\delta_2$.
Then, if $0<\zeta_+\le\tau_2$, it follows from \eqref{c-4}--\eqref{tau2} that
\begin{gather}
\label{c-4-1}
-\beta_2\varphi^{\lambda+1}
\le Q_\varphi(\varphi,\psi)\le-\beta_1\varphi^{\lambda+1},
\quad(\varphi,\psi)\in(0,\zeta_{+})\times(0,m),
\\
\label{c-5-1}
|Q_\psi(\varphi,\psi)|\le\beta_3\varphi^{3\lambda/2+1},
\quad(\varphi,\psi)\in(0,\zeta_{+})\times(0,m),
\\
\label{c-6-1}
-2\sigma_2\varphi^{\lambda+2}
\le Q(\varphi,\psi)\le
-\frac12\sigma_1\varphi^{\lambda+2},
\quad
(\varphi,\psi)\in(0,\zeta_+)\times(0,m)
\end{gather}
and
\begin{align}
\label{c-7}
c_*\le A^{-1}_+(Q(\varphi,\psi))\le\overline c,
\quad(\varphi,\psi)\in(0,\zeta_+)\times(0,m).
\end{align}
Set
$$
\hat {\mathscr Q}_{+}(x)=A^{-1}_+(Q(\Phi_{+}(x),m)),\quad0\le x\le l_+
$$
and
$$
\hat {\mathscr Q}(\phi,\psi)=Q(\zeta_+\phi,\psi),\quad
(\phi,\psi)\in(0,1)\times(0,m).
$$
Due to \eqref{c-4-1}--\eqref{c-7}, one can get that
$\hat {\mathscr Q}_{+}\in C([0,l_+])$ satisfies
\begin{align}
\label{c-8}
c_*\le \hat {\mathscr Q}_{+}(x)\le\overline c,\quad0\le x\le l_+,
\end{align}
and
$\hat {\mathscr Q}\in C^{0,1}([0,1]\times[0,m])$ satisfies
\begin{gather}
\label{c-6-2}
-2\sigma_2\Big(\overline c\sqrt{\delta_2^2+1}l_+\Big)^{\lambda+2}\phi^{\lambda+2}
\le\hat {\mathscr Q}(\phi,\psi)\le
-\frac12\sigma_1(c_*l_+)^{\lambda+2}\phi^{\lambda+2},
\quad
(\phi,\psi)\in(0,1)\times(0,m),
\\
\label{c-4-2}
-\beta_2\Big(\overline c\sqrt{\delta_2^2+1}l_+\Big)^{\lambda+2}\phi^{\lambda+1}
\le\hat {\mathscr Q}_\phi(\phi,\psi)\le
-\beta_1(c_*l_+)^{\lambda+2}\phi^{\lambda+1},
\quad(\phi,\psi)\in(0,1)\times(0,m),
\\
\label{c-5-2}
\Big|\hat {\mathscr Q}_\psi(\phi,\psi)\Big|\le
\beta_3\Big(\overline c\sqrt{\delta_2^2+1}l_+\Big)^{3\lambda/2+1}\phi^{3\lambda/2+1},
\quad(\phi,\psi)\in(0,1)\times(0,m).
\end{gather}

Now we choose
\begin{align}
\label{delta0+}
\delta_0=\frac1{\overline c\sqrt{\delta_2^2+1}}\min\{\tau_1,\tau_2\},
\end{align}
which depends only on $\gamma$, $m$,
$\lambda$, $\delta_1$ and $\delta_2$.
Due to \eqref{c-8}--\eqref{c-5-2}, we get that
if $0<l_+\le\delta_0$, then $0<\zeta_+\le\min\{\tau_1,\tau_2\}$,
$\hat {\mathscr Q}_{+}\in C([0,l_+])$ satisfies \eqref{qub+}
and $\hat {\mathscr Q}\in C^{0,1}([0,1]\times[0,m])$ satisfies
\eqref{c-1}--\eqref{c-2} with \eqref{beta}.

For $0<l_+<\delta_0$, define
\begin{align*}
{\mathscr S}=\Big\{&({\mathscr Q}_{+},{\mathscr Q})\in
C([0,l_+])\times C^{0,1}([0,1]\times[0,m]):
\\
&\qquad
{\mathscr Q}_{+} \mbox{ satisfies }\eqref{qub+},
\mbox{ while } {\mathscr Q}\mbox{ satisfies \eqref{c-1}--\eqref{c-2}
with \eqref{beta}}\Big\}
\end{align*}
with the norm
$$
\|({\mathscr Q}_{+},{\mathscr Q})\|_{\mathscr S}
=\max\Big\{\|{\mathscr Q}_{+}\|_{L^\infty((0,l_+))},
\|{\mathscr Q}\|_{L^\infty((0,1)\times(0,m))}\Big\},\quad
({\mathscr Q}_{+},{\mathscr Q})\in {\mathscr S}.
$$
Owing to the discussion above, we can
define a mapping $J$ from ${\mathscr S}$ to itself
as follows
\begin{align}
\label{mapj}
J(({\mathscr Q}_{+},{\mathscr Q}))=(\hat {\mathscr Q}_{+},
\hat {\mathscr Q}),\quad ({\mathscr Q}_{+},{\mathscr Q})\in {\mathscr S}.
\end{align}

Now the existence of sonic-supersonic flows can be stated as

\begin{theorem}
\label{theoremsuper}
Assume that $f_+\in C^{3}([0,l_+])$ satisfies \eqref{a-a-0}.
If $0<l_+\le\delta_0$,
then, the problem \eqref{supp-eq}--\eqref{supp-q} admits at least one weak solution
$Q$ satisfying
\begin{align}
\label{theoremsuper-1}
-2\sigma_2\varphi^{\lambda+2}\le Q(\varphi,\psi)\le-\frac12\sigma_1\varphi^{\lambda+2},
\quad
-\beta_2\varphi^{\lambda+1}\le Q_\varphi(\varphi,\psi)
\le-\beta_1\varphi^{\lambda+1},
&\quad
|Q_\psi(\varphi,\psi)|\le\beta_3\varphi^{3\lambda/2+1},
\nonumber
\\
&\hskip-15mm
(\varphi,\psi)\in(0,\zeta_{+})\times(0,m),
\end{align}
where $\delta_0$, $\sigma_i\,(i=1,2)$ and $\beta_i\,(i=1,2,3)$
are defined in \eqref{delta0+}, Proposition \ref{prop-b} and \eqref{beta}, respectively,
which all depend only on $\gamma$, $m$,
$\lambda$, $\delta_1$ and $\delta_2$.
\end{theorem}

\Proof
As mentioned above, the mapping $J$ defined by \eqref{mapj}
is from ${\mathscr S}$ to itself. It follows from \eqref{c-2-0}, \eqref{c-2}
and the embedding theorem that $J$ is compact.
Therefore, Theorem \ref{theoremsuper} follows from
\eqref{c-4-1}--\eqref{c-6-1} and
the Schauder fixed point theorem provided that $J$ is also continuous.

It remains to show that $J$ is continuous.
Assume that $\{({\mathscr Q}_{+}^{(n)},{\mathscr Q}^{(n)})\}_{n=0}^\infty
\subset{\mathscr S}$ satisfies
\begin{align}
\label{t2-1}
\lim_{n\to\infty}\|{\mathscr Q}_{+}^{(n)}
-{\mathscr Q}_{+}^{(0)}\|_{L^\infty((0,l_+))}=0,
\quad
\lim_{n\to\infty}\|{\mathscr Q}^{(n)}
-{\mathscr Q}^{(0)}\|_{L^\infty((0,1)\times(0,m))}=0.
\end{align}
Let $\Phi_+^{(n)}$, $X_+^{(n)}$ and $\zeta_{+}^{(n)}$ be
defined by \eqref{Phiub+}--\eqref{zeta+} with
${\mathscr Q}_{+}={\mathscr Q}_{+}^{(n)}$ for $n=0,1,2,\cdots$ and set
$$
\tilde Q^{(n)}(\varphi,\psi)={\mathscr Q}^{(n)}
\Big(\frac{\varphi}{\zeta_{+}^{(n)}},\psi\Big),
\quad(\varphi,\psi)\in(0,\zeta_{+}^{(n)})\times(0,m),
\quad n=0,1,2,\cdots.
$$
Then, \eqref{t2-1} leads to
\begin{gather}
\label{t2-1-1}
\lim_{n\to\infty}\|\Phi_+^{(n)}-\Phi_+^{(0)}\|_{0,1;(0,l_+)}=0,
\quad
\lim_{n\to\infty}\zeta_{+}^{(n)}=\zeta_{+}^{(0)},
\\
\label{t2-1-3}
X_+^{(n)}\longrightarrow X_+^{(0)}\mbox{ uniformly in }(0,\zeta_{+}^{(0)}),
\quad
\tilde Q^{(n)}\longrightarrow \tilde Q^{(0)}
\mbox{ uniformly in } (0,\zeta_{+}^{(0)})\times(0,m).
\end{gather}
It follows from the definition of $J$ that
$$
J(({\mathscr Q}_{+}^{(n)},{\mathscr Q}^{(n)}))
=(\hat {\mathscr Q}_{+}^{(n)},\hat {\mathscr Q}^{(n)}),
\quad n=0,1,2,\cdots,
$$
where
$$
\hat {\mathscr Q}_{+}^{(n)}(x)=A^{-1}_+(Q^{(n)}(\Phi_{+}^{(n)}(x),m)),\quad0\le x\le l_+,
\quad n=0,1,2,\cdots
$$
and
$$
\hat {\mathscr Q}^{(n)}(\phi,\psi)=Q^{(n)}(\zeta_+^{(n)}\phi,\psi),\quad
(\phi,\psi)\in(0,1)\times(0,m),
\quad n=0,1,2,\cdots
$$
with $Q^{(n)}$ being the unique weak solution to the following problem
\begin{align}
\label{t2itsupp-eq} &Q_{\varphi\varphi}^{(n)}
-(b(\tilde Q^{(n)})Q_{\psi}^{(n)})_\psi=0,
\quad&&(\varphi,\psi)\in(0,\zeta_{+}^{(n)})\times(0,m),
\\
\label{t2itsupp-inbc1}
&Q^{(n)}(0,\psi)=0,\quad&&\psi\in(0,m),
\\
\label{t2itsupp-inbc2}
&Q^{(n)}_\varphi(0,\psi)=0,\quad&&\psi\in(0,m),
\\
\label{t2itsupp-lbbc}
&Q^{(n)}_{\psi}(\varphi,0)=0,\quad&&\varphi\in(0,\zeta_{+}^{(n)}),
\\
\label{t2itsupp-ubbc}
&Q^{(n)}_{\psi}(\varphi,m)=-\frac{\Theta'_{+}(X^{(n)}_{+}(\varphi))
\cos\Theta_{+}(X^{(n)}_{+}(\varphi))}
{b(\tilde Q^{(n)}(\varphi,m)))A^{-1}_+(\tilde Q^{(n)}(\varphi,m))},
\quad&&\varphi\in(0,\zeta_{+}^{(n)})
\end{align}
for $n=0,1,2,\cdots$.
So, to show show the continuity of $J$, it suffices to prove
\begin{align}
\label{t2-2}
\lim_{n\to\infty}\|\hat {\mathscr Q}^{(n)}
-\hat {\mathscr Q}^{(0)}\|_{L^\infty((0,1)\times(0,m))}=0,
\end{align}
since this limit, together with \eqref{t2-1-1}, also implies
$$
\lim_{n\to\infty}\|\hat {\mathscr Q}^{(n)}_+
-\hat {\mathscr Q}^{(0)}_+\|_{L^\infty((0,l_+))}=0.
$$
We will prove \eqref{t2-2} by contradiction. Otherwise,
there exists a positive number $\varepsilon_0$ and a subsequence of
$\{\hat {\mathscr Q}^{(n)}\}_{n=1}^\infty$, denoted by itself for convenience, such that
\begin{align}
\label{t2-3}
\|\hat {\mathscr Q}^{(n)}
-\hat {\mathscr Q}^{(0)}\|_{L^\infty((0,1)\times(0,m))}\ge\varepsilon_0,
\quad n=1,2,\cdots.
\end{align}
Since $\hat {\mathscr Q}^{(n)}$ satisfies \eqref{c-1}--\eqref{c-2}
for each $n=1,2,\cdots$, there
exist a subsequence of
$\{\hat {\mathscr Q}^{(n)}\}_{n=1}^\infty$, denoted by itself again for convenience,
and a function ${\mathscr Q}^*$ with \eqref{c-1}--\eqref{c-2}
such that
\begin{gather}
\label{t2-4}
\hat {\mathscr Q}^{(n)}\longrightarrow {\mathscr Q}^{*}
\mbox{ in }L^\infty((0,1)\times(0,m))
\mbox{ as } n\to\infty,
\\
\label{t2-5}
\hat {\mathscr Q}^{(n)}_\phi-\hskip-2mm-\hskip-2mm\rightharpoonup {\mathscr Q}^{*}_\phi
\mbox{ and }
\hat {\mathscr Q}^{(n)}_\psi-\hskip-2mm-\hskip-2mm\rightharpoonup {\mathscr Q}^{*}_\psi
\mbox{ weakly * in }L^\infty((0,1)\times(0,m))
\mbox{ as } n\to\infty.
\end{gather}
Set
\begin{align*}
Q^*(\varphi,\psi)={\mathscr Q}^{*}\Big(\frac\varphi{\zeta_+^{(0)}},\psi\Big),\quad
(\varphi,\psi)\in(0,\zeta_+^{(0)})\times(0,m).
\end{align*}
Letting $n\to\infty$ in \eqref{t2itsupp-eq}--\eqref{t2itsupp-ubbc}
and using \eqref{t2-1-1}, \eqref{t2-1-3}, \eqref{t2-4} and \eqref{t2-5}, one can get
that $Q^*$ solves the following problem
\begin{align*}
&Q_{\varphi\varphi}^*
-(b(\tilde Q^{(0)})Q_{\psi}^{*})_\psi=0,
\quad&&(\varphi,\psi)\in(0,\zeta_{+}^{(0)})\times(0,m),
\\
&Q^{*}(0,\psi)=0,\quad&&\psi\in(0,m),
\\
&Q^{*}_\varphi(0,\psi)=0,\quad&&\psi\in(0,m),
\\
&Q^{*}_{\psi}(\varphi,0)=0,\quad&&\varphi\in(0,\zeta_{+}^{(0)}),
\\
&Q^{*}_{\psi}(\varphi,m)=-\frac{\Theta'_{+}(X^{(0)}_{+}(\varphi))
\cos\Theta_{+}(X^{(0)}_{+}(\varphi))}
{b(\tilde Q^{(0)}(\varphi,m)))A^{-1}_+(\tilde Q^{(0)}(\varphi,m))},
\quad&&\varphi\in(0,\zeta_{+}^{(0)}).
\end{align*}
It follows from the uniqueness in Proposition \ref{prop-b} that
$$
Q^*(\varphi,\psi)=Q^{(0)}(\varphi,\psi),\quad
(\varphi,\psi)\in(0,\zeta_{+}^{(0)})\times(0,m),
$$
which contradicts \eqref{t2-3} and \eqref{t2-4}.
Hence \eqref{t2-2} holds.
$\hfill\Box$\vskip 4mm

\subsection{Smooth sonic-supersonic flows}

In this subsection, we establish the existence of $C^{1,1}$ sonic-supersonic flows
under the assumption that $f_+\in C^{3}([0,l_+])$ satisfies \eqref{a-a-0} and \eqref{rr-1}.
The method is still a fixed point argument and the method of characteristics.
To this end, one should establish a priori $C^{1,1}$ estimates for the linearized problem
\eqref{b-1}--\eqref{b-9}.

\begin{lemma}
\label{rr-lemma1}
Assume that $(W,Z)$ is the unique weak solution to the problem
\eqref{hd-1}--\eqref{hd-6} in Lemma \ref{lemma-0}. Then,
\begin{align*}
%\label{rr-2}
|(W_\varphi,Z_\varphi)(\varphi,\psi)|\le\mu_{10}\varphi^{\lambda},
\quad(\varphi,\psi)\in(0,\zeta_{+})\times(0,m)
\end{align*}
with $\mu_{10}$ depending only on $\gamma$, $m$,
$\lambda$, $\delta_1$, $\delta_2$, $\delta_3$, $\mu_1$, $\mu_2$ and $\beta_2$.
\end{lemma}

\Proof
Set
$$
(U,V)(\varphi,\psi)=(W_\varphi,Z_\varphi)(\varphi,\psi),\quad
(\varphi,\psi)\in(0,\zeta_{+})\times(0,m).
$$
Then, $(U,V)$ is the unique weak solution to the following problem
\begin{align*}
&U_\varphi+b^{1/2}(\tilde Q)U_\psi
=\frac12b^{-1}(\tilde Q)p(\tilde Q)\tilde Q_\varphi U,
\quad&&(\varphi,\psi)\in(0,\zeta_{+})\times(0,m),
\\
&V_\varphi-b^{1/2}(\tilde Q)V_\psi
=\frac12b^{-1}(\tilde Q)p(\tilde Q)\tilde Q_\varphi V,
\quad&&(\varphi,\psi)\in(0,\zeta_{+})\times(0,m),
\\
&U(0,\psi)=0,\quad&&\psi\in(0,m),
\\
&V(0,\psi)=0,\quad&&\psi\in(0,m),
\\
&U(\varphi,0)+V(\varphi,0)=0,
\quad&&\varphi\in(0,\zeta_{+}),
\\
&U(\varphi,m)+V(\varphi,m)=\tilde h'(\varphi),
\quad&&\varphi\in(0,\zeta_{+}).
\end{align*}
Decompose $(U,V)$ as
$$
(U,V)(\varphi,\psi)=(U^{1},V^{1})(\varphi,\psi)+(U^{2},V^{2})(\varphi,\psi),\quad
(\varphi,\psi)\in(0,\zeta_{+})\times(0,m),
$$
where $(U^{1},V^{1})$ and $(U^{2},V^{2})$ solve the following problems
\begin{align*}
&U^{1}_\varphi+b^{1/2}(\tilde Q)U^{1}_\psi
=0,
\quad&&(\varphi,\psi)\in(0,\zeta_{+})\times(0,m),
\\
&V^{1}_\varphi-b^{1/2}(\tilde Q)V^{1}_\psi
=0,
\quad&&(\varphi,\psi)\in(0,\zeta_{+})\times(0,m),
\\
&U^{1}(0,\psi)=0,\quad&&\psi\in(0,m),
\\
&V^{1}(0,\psi)=0,\quad&&\psi\in(0,m),
\\
&U^{1}(\varphi,0)+V^{1}(\varphi,0)=0,
\quad&&\varphi\in(0,\zeta_{+}),
\\
&U^{1}(\varphi,m)+V^{1}(\varphi,m)=\tilde h'(\varphi),
\quad&&\varphi\in(0,\zeta_{+})
\end{align*}
and
\begin{align*}
&U^{2}_\varphi+b^{1/2}(\tilde Q)U^{2}_\psi
=\frac12b^{-1}(\tilde Q)p(\tilde Q)\tilde Q_\varphi(U^{1}+U^{2}),
\quad&&(\varphi,\psi)\in(0,\zeta_{+})\times(0,m),
\\
&V^{2}_\varphi-b^{1/2}(\tilde Q)V^{2}_\psi
=\frac12b^{-1}(\tilde Q)p(\tilde Q)\tilde Q_\varphi(V^{1}+V^{2}),
\quad&&(\varphi,\psi)\in(0,\zeta_{+})\times(0,m),
\\
&U^{2}(0,\psi)=0,\quad&&\psi\in(0,m),
\\
&V^{2}(0,\psi)=0,\quad&&\psi\in(0,m),
\\
&U^{2}(\varphi,0)+V^{2}(\varphi,0)=0,
\quad&&\varphi\in(0,\zeta_{+}),
\\
&U^{2}(\varphi,m)+V^{2}(\varphi,m)=0,
\quad&&\varphi\in(0,\zeta_{+}),
\end{align*}
respectively.
It follows from the same proof as for Lemma \ref{lemma-0} that
$$
|(U^{1},V^{1})(\varphi,\psi)|\le\frac12\mu_{10}\varphi^{\lambda},
\quad(\varphi,\psi)\in(0,\zeta_{+})\times(0,m),
$$
where $\mu_{10}$ depends only on $\gamma$, $m$,
$\lambda$, $\delta_1$, $\delta_2$, $\delta_3$, $\mu_1$, $\mu_2$ and $\beta_2$.
Then, Remark \ref{qqq9} gives
$$
|(U^{2},V^{2})(\varphi,\psi)|\le\frac12\mu_{10}\varphi^{\lambda},
\quad(\varphi,\psi)\in(0,\zeta_{+})\times(0,m).
$$
Therefore,
$$
|(U,V)(\varphi,\psi)|\le\mu_{10}\varphi^{\lambda},
\quad(\varphi,\psi)\in(0,\zeta_{+})\times(0,m).
$$
$\hfill\Box$\vskip 4mm

Now, we assume additionally that
\begin{align}
\label{rr-3}
|\tilde Q_{\varphi\varphi}(\varphi,\psi)|\le\beta_4\varphi^{\lambda},
\quad
|\tilde Q_{\varphi\psi}(\varphi,\psi)|\le\beta_4\varphi^{5\lambda/4+1/2},
\quad(\varphi,\psi)\in(0,\zeta_{+})\times(0,m)
\end{align}
with $\beta_4\ge1$ to be determined later.
Then,
\begin{align*}
|(\tilde W_\varphi,\tilde Z_\varphi)(\varphi,\psi)|\le\mu_{11}\beta_4\varphi^{\lambda},
\quad(\varphi,\psi)\in(0,\zeta_{+})\times(0,m),
\end{align*}
where $\mu_{11}$ depends only on $\gamma$, $m$,
$\lambda$, $\delta_1$, $\delta_2$, $\mu_1$, $\mu_2$ and $\beta_2$.

\begin{lemma}
\label{rr-lemma2}
Assume that \eqref{rr-3} holds, $0<\zeta_+\le\tau_3$, $\|(w-W^*,z-Z^*)\|_{\mathscr B}\le\mu_8\beta_2$
and
\begin{align*}
|((w-W^*)_\varphi,(z-Z^*)_\varphi)(\varphi,\psi)|\le N\varphi^{\lambda},
\quad(\varphi,\psi)\in(0,\zeta_{+})\times(0,m)
\end{align*}
with a constant $N\ge1$,
where $(W^*,Z^*)$ is the unique weak solution to the problem \eqref{hd-1}--\eqref{hd-6}
in Lemma \ref{lemma-0}.
Then for any $(\varphi,\psi)\in(0,\zeta_{+})\times(0,m)$,
\begin{gather}
\label{rr-lemma2888}
|((W-W^*)_\varphi,(Z-Z^*)_\varphi)(\varphi,\psi)|
\le\Big(\frac{2\mu_8\beta_2+\mu_4}{2\beta_1}\mu_{11}\beta_4\tau_3^{\lambda/4-1/2}
+\frac{\mu_5\beta_2}{\mu_5\beta_2+2\beta_1}N+\mu_{12}\Big)\varphi^{\lambda},
\end{gather}
where $(W,Z)$ is the unique weak solution to the problem \eqref{d-1}--\eqref{d-6}
in Proposition \ref{proposition-a}
and $\mu_{12}$ depends only on $\gamma$, $m$,
$\lambda$, $\delta_1$, $\delta_2$, $\delta_3$, $\mu_1$, $\mu_2$ and $\beta_2$.
\end{lemma}

\Proof
Let $({\mathscr W},{\mathscr Z})$
be the weak solution to the problem \eqref{scrd-1}--\eqref{scrd-6}
in the proof of Proposition \ref{proposition-a}.
Set
$$
(U,V)(\varphi,\psi)=({\mathscr W}_\varphi,{\mathscr Z}_\varphi)(\varphi,\psi),\quad
(\varphi,\psi)\in(0,\zeta_{+})\times(0,m).
$$
Then, $(U,V)$ solves the following problem
\begin{align*}
&U_\varphi+b^{1/2}(\tilde Q)U_\psi
=\frac14b^{-1}(\tilde Q)p(\tilde Q)(2\tilde Q_\varphi+\tilde W)
\big(U+F^{1}+F^{2}+F^{3}\big),
\quad&&(\varphi,\psi)\in(0,\zeta_{+})\times(0,m),
\\
&V_\varphi-b^{1/2}(\tilde Q)V_\psi
=\frac14b^{-1}(\tilde Q)p(\tilde Q)(2\tilde Q_\varphi-\tilde Z)
\big(V+G^{1}+G^{2}+G^{3}\big),
\quad&&(\varphi,\psi)\in(0,\zeta_{+})\times(0,m),
\\
&U(0,\psi)=0,\quad&&\psi\in(0,m),
\\
&V(0,\psi)=0,\quad&&\psi\in(0,m),
\\
&U(\varphi,0)+V(\varphi,0)=0,
\quad&&\varphi\in(0,\zeta_{+}),
\\
&U(\varphi,m)+V(\varphi,m)=0,
\quad&&\varphi\in(0,\zeta_{+}),
\end{align*}
where
\begin{align*}
F^{1}=&(2\tilde Q_\varphi+\tilde W)^{-1}\tilde W_\varphi
\big({\mathscr W}+(W^*+Z^*)+(z-Z^*)\big),
\\
F^{2}=&(2\tilde Q_\varphi+\tilde W)^{-1}\tilde W(z-Z^*)_\varphi,
\\
F^{3}=&(2\tilde Q_\varphi+\tilde W)^{-1}\tilde W\Big(b(\tilde Q)p^{-1}(\tilde Q)
\big(b^{-1}(\tilde Q)p(\tilde Q)\big)_\varphi-\frac12b^{-1}(\tilde Q)p(\tilde Q)\tilde Q_\varphi\Big)
\big({\mathscr W}+(W^*+Z^*)+(z-Z^*)\big)
\\
&\qquad+(2\tilde Q_\varphi+\tilde W)^{-1}\tilde W
(W^*+Z^*)_\varphi,
\\
G^{1}=&-(2\tilde Q_\varphi-\tilde Z)^{-1}\tilde Z_\varphi
\big((w-W^*)+(W^*+Z^*)+{\mathscr Z}\big),
\\
G^{2}=&-(2\tilde Q_\varphi-\tilde Z)^{-1}\tilde Z(w-W^*)_\varphi,
\\
G^{3}=&(2\tilde Q_\varphi-\tilde Z)^{-1}\tilde Z\Big(\frac12b^{-1}(\tilde Q)p(\tilde Q)\tilde Q_\varphi
-b(\tilde Q)p^{-1}(\tilde Q)\big(b^{-1}(\tilde Q)p(\tilde Q)\big)_\varphi\Big)
\big((w-W^*)+(W^*+Z^*)+{\mathscr Z}\big)
\\
&\qquad-(2\tilde Q_\varphi-\tilde Z)^{-1}\tilde Z(W^*+Z^*)_\varphi.
\end{align*}
Direct calculations show that
\begin{gather*}
|(F^{1},G^{1})(\varphi,\psi)|
\le\frac{2\mu_8\beta_2+\mu_4}{2\beta_1}\mu_{11}\beta_4\tau_3^{\lambda/4-1/2}\varphi^{\lambda},
\quad(\varphi,\psi)\in(0,\zeta_{+})\times(0,m),
\\
|(F^{2},G^{2})(\varphi,\psi)|\le\frac{\mu_5\beta_2}{\mu_5\beta_2+2\beta_1}N\varphi^{\lambda},
\quad(\varphi,\psi)\in(0,\zeta_{+})\times(0,m)
\end{gather*}
and
\begin{align*}
|(F^{3},G^{3})(\varphi,\psi)|\le\mu_{12}\varphi^{\lambda},
\quad(\varphi,\psi)\in(0,\zeta_{+})\times(0,m),
\end{align*}
where $\mu_{12}$ depends only on $\gamma$, $m$,
$\lambda$, $\delta_1$, $\delta_2$, $\delta_3$, $\mu_1$, $\mu_2$ and $\beta_2$.
Then, Remark \ref{qqq9} yields
\begin{align*}
&|(U,V)(\varphi,\psi)|\le\Big(\frac{2\mu_8\beta_2+\mu_4}{2\beta_1}\mu_{11}\beta_4\tau_3^{\lambda/4-1/2}
+\frac{\mu_5\beta_2}{\mu_5\beta_2+2\beta_1}N+\mu_{12}\Big)\varphi^{\lambda},
\quad(\varphi,\psi)\in(0,\zeta_{+})\times(0,m).
\end{align*}
$\hfill\Box$\vskip 4mm

\begin{remark}
Under the assumptions of Lemma \ref{rr-lemma2}, if, in addition,
\begin{align}
\label{rr-6}
\frac{2\mu_8\beta_2+\mu_4}{2\beta_1}\mu_{11}\beta_4\tau_3^{\lambda/4-1/2}
+\mu_{12}\le\frac{2\beta_1}{\mu_5\beta_2+2\beta_1}N,
\end{align}
then \eqref{rr-lemma2888} implies
\begin{align*}
|((W-W^*)_\varphi,(Z-Z^*)_\varphi)(\varphi,\psi)|\le N\varphi^{\lambda},
\quad(\varphi,\psi)\in(0,\zeta_{+})\times(0,m).
\end{align*}
Therefore, replacing ${\mathscr B}$ with
\begin{align*}
\tilde {\mathscr B}=&\Big\{(w,z)\in L^\infty((0,\zeta_{+})\times(0,m))
\times L^\infty((0,\zeta_{+})\times(0,m)):
\|(w-W^*,z-Z^*)\|_{\mathscr B}\le\mu_8\beta_2,
\\
&\qquad|((w-W^*)_\varphi,(z-Z^*)_\varphi)(\varphi,\psi)|\le N\varphi^{\lambda},
\quad(\varphi,\psi)\in(0,\zeta_{+})\times(0,m)\Big\}
\end{align*}
in the proof of Proposition \ref{prop-b}, we can show that
\end{remark}

\begin{proposition}
\label{rr-prop1}
Let the assumptions of Proposition \ref{prop-b} hold and
$(W,Z,Q)$ be the unique weak solution to the problem \eqref{b-1}--\eqref{b-9}.
Assume further that $0<\zeta_+\le\tau_3$, $N\ge1$
and \eqref{rr-6} holds.
Then $(W,Z)\in\tilde {\mathscr B}$ and
\begin{align}
\label{rr-prop1888}
|Q_{\varphi\varphi}(\varphi,\psi)|\le(N+\mu_{10})\varphi^{\lambda},
\quad
|Q_{\varphi\psi}(\varphi,\psi)|\le\mu_{13}(N+\mu_{10})\varphi^{5\lambda/4+1/2},
\quad(\varphi,\psi)\in(0,\zeta_{+})\times(0,m),
\end{align}
where $\mu_{13}$ depends only on $\gamma$, $m$,
$\lambda$, $\delta_1$, $\delta_2$, $\mu_1$, $\mu_2$ and $\beta_2$.
\end{proposition}

\begin{remark}
Under the assumptions of Proposition \ref{rr-prop1}, if, in addition,
\begin{align}
\label{rr-7}
\beta_4\ge\max\Big\{
(N+\mu_{10})\Big(\frac{\overline c\sqrt{\delta_2^2+1}}{c_*}\Big)^{\lambda+2},
\mu_{13}(N+\mu_{10})
\Big(\frac{\overline c\sqrt{\delta_2^2+1}}{c_*}\Big)^{5\lambda/4+3/2}\Big\},
\end{align}
then \eqref{rr-prop1888} leads to
\begin{align*}
&|Q_{\varphi\varphi}(\varphi,\psi)|\le\beta_4
\Big(\frac{c_*}{\overline c\sqrt{\delta_2^2+1}}\Big)^{\lambda+2}\varphi^{\lambda},
&&\quad(\varphi,\psi)\in(0,\zeta_{+})\times(0,m),
\\
&|Q_{\varphi\psi}(\varphi,\psi)|\le\beta_4
\Big(\frac{c_*}{\overline c\sqrt{\delta_2^2+1}}\Big)^{5\lambda/4+3/2}\varphi^{5\lambda/4+1/2},
&&\quad(\varphi,\psi)\in(0,\zeta_{+})\times(0,m).
\end{align*}
\end{remark}

We now choose
\begin{gather*}
\tau_3=\Big(\frac{\beta_1^2}{\mu_{11}(1+\mu_{13})(2\mu_8\beta_2+\mu_4)(\mu_5\beta_2+2\beta_1)}\Big)^{4/(\lambda-2)}
\Big(\frac{c_*}{\overline c\sqrt{\delta_2^2+1}}\Big)^{(5\lambda+6)/(\lambda-2)},
\\
\beta_4=(1+\mu_{13})\Big(2\mu_{10}+\frac{\mu_{12}(\mu_5\beta_2+2\beta_1)}{\beta_1}+1\Big)
\Big(\frac{\overline c\sqrt{\delta_2^2+1}}{c_*}\Big)^{5\lambda/4+3/2}
\end{gather*}
and
$$
N=\mu_{10}+\frac{\mu_{12}(\mu_5\beta_2+2\beta_1)}{\beta_1}+1.
$$
Then, \eqref{rr-6} and  \eqref{rr-7} hold.
Set
$$
\tilde \delta_0=\frac1{\overline c\sqrt{\delta_2^2+1}}\min\{\tau_1,\tau_2,\tau_3\}
$$
and
\begin{align*}
\tilde{\mathscr S}=\Big\{&({\mathscr Q}_{+},{\mathscr Q})\in
C([0,l_+])\times C^{1,1}([0,1]\times[0,m]):{\mathscr Q}_{+} \mbox{ satisfies }\eqref{qub+},
\mbox{ while } {\mathscr Q} \mbox{ satisfies }
\\
&\qquad
\mbox{ \eqref{c-1}--\eqref{c-2}
with \eqref{beta}, }|{\mathscr Q}_{\phi\phi}(\phi,\psi)|\le\beta_4
(c_*l_+)^{\lambda+2}\phi^{\lambda}\mbox{ and }
\\
&\qquad
|{\mathscr Q}_{\phi\psi}(\phi,\psi)|\le\beta_4
(c_*l_+)^{5\lambda/4+3/2}\phi^{5\lambda/4+1/2}\mbox{ for }
(\phi,\psi)\in(0,1)\times(0,m)\Big\}.
\end{align*}
Replacing $\delta_0$ and ${\mathscr S}$ with $\tilde\delta_0$ and $\tilde{\mathscr S}$
in the proof of Theorem \ref{theoremsuper}, one can prove that

\begin{theorem}
\label{rr-theorem}
Assume that $f_+\in C^{3}([0,l_+])$ satisfies \eqref{a-a-0} and \eqref{rr-1}.
If $0<l_+\le\tilde\delta_0$,
then the problem \eqref{supp-eq}--\eqref{supp-q} admits at least a solution
$Q\in C^{1,1}([0,\zeta_+]\times[0,m])$, which satisfies \eqref{theoremsuper-1} and
\begin{align}
\label{rr-theorem-1}
|Q_{\varphi\varphi}(\varphi,\psi)|\le\mu_{14}\varphi^{\lambda},
\quad
|Q_{\varphi\psi}(\varphi,\psi)|\le\mu_{14}\varphi^{5\lambda/4+1/2},
\quad
|Q_{\psi\psi}(\varphi,\psi)|\le\mu_{14}\varphi^{3\lambda/2+1},
\nonumber
\\
\quad(\varphi,\psi)\in(0,\zeta_{+})\times(0,m),
\end{align}
where $\tilde\delta_0$ and $\mu_{14}$ depend only on
$\gamma$, $m$,
$\lambda$, $\delta_1$, $\delta_2$ and $\delta_3$.
\end{theorem}

\subsection{Uniqueness of sonic-supersonic flows}

In this subsection, we study the uniqueness of sonic-supersonic flows.
Owing to the strong singularity of the equation at the sonic curve,
one can only prove the uniqueness of sonic-supersonic flows with
precise estimates as in the existence (Theorems \ref{theoremsuper} and \ref{rr-theorem}).
In Theorem \ref{rr-theorem}, we get a $C^{1,1}$ sonic-supersonic flow
under the assumption that $f_+\in C^{3}([0,l_+])$ satisfies \eqref{a-a-0} and \eqref{rr-1}.
This smooth sonic-supersonic flow will be shown to be unique if $f_+\in C^{4}([0,l_+])$ additionally.
Furthermore, it is also unique in the space of weak solutions given in Theorem \ref{theoremsuper}.
It is noted that the additional condition $f_+\in C^{4}([0,l_+])$ arises
from the nonlocal and implicit boundary condition \eqref{supp-ubbc}.

\begin{theorem}
\label{uu-theorem}
Assume that $f_+\in C^{4}([0,l_+])$ satisfies \eqref{a-a-0} and \eqref{rr-1}.
Then the problem \eqref{supp-eq}--\eqref{supp-q} admits at most one weak solution
$Q\in C^{0,1}([0,\zeta_+]\times[0,m])$ satisfying
\begin{gather}
\label{supun-1}
-M_2\varphi^{\lambda+2}\le Q(\varphi,\psi)\le-M_1\varphi^{\lambda+2},
\quad(\varphi,\psi)\in(0,\zeta_{+})\times(0,m),
\\
\label{supun-2}
-M_2\varphi^{\lambda+1}\le Q_\varphi(\varphi,\psi)
\le-M_1\varphi^{\lambda+1},
\quad
|Q_\psi(\varphi,\psi)|\le M_2\varphi^{3\lambda/2+1},
\quad(\varphi,\psi)\in(0,\zeta_{+})\times(0,m),
\end{gather}
where $M_1$ and $M_2\,(M_1\le M_2)$ are positive constants.
\end{theorem}

\Proof
Note that the problem \eqref{supp-eq}--\eqref{supp-q} is equivalent to 
the problem \eqref{phytran15}--\eqref{phytran20},
which has standard initial and boundary value conditions,
and \eqref{phytran15} is strictly hyperbolic away from the sonic curve.
So, it suffices to prove the uniqueness of weak solutions to
the problem \eqref{supp-eq}--\eqref{supp-q} for small $\zeta_+>0$.
Thanks to Theorem \ref{rr-theorem}, for sufficiently small $\zeta_+\in(0,1)$,
the problem \eqref{supp-eq}--\eqref{supp-q}
admits a solution $\tilde Q\in C^{1,1}([0,\zeta_+]\times[0,m])$ satisfying
\eqref{theoremsuper-1} and \eqref{rr-theorem-1}.
Therefore, one need only to prove that $\tilde Q$ is the unique weak
solution to the problem \eqref{supp-eq}--\eqref{supp-q}.
Assume that $\hat Q\in C^{0,1}([0,\zeta_+]\times[0,m])$
with \eqref{supun-1} and \eqref{supun-2} is a weak solution
to the problem \eqref{supp-eq}--\eqref{supp-q}.
Note that
$$
\mbox{div}\big(\hat Q_\varphi,-b(\hat Q)\hat Q_\psi\big)=
\mbox{div}\big(\hat Q_\psi,-\hat Q_\varphi\big)=0
$$
in $(0,\zeta_{+})\times(0,m)$ in the sense of distribution.
It follows from the theory of $L^\infty$ divergence-measure vector fields (\cite{pair}),
\eqref{supun-1} and \eqref{supun-2} that
$\hat Q_\varphi$ and $\hat Q_\psi$ have
$L^\infty$ trace on $[0,\zeta_+]\times\{\psi\}$ and $\{\varphi\}\times[0,m]$
for each $\psi\in[0,m]$ and $\varphi\in[0,\zeta_+]$.

Denote the functions given by \eqref{Xub+} corresponding to $\tilde Q$
and $\hat Q$ by $\tilde X_{+}$ and $\hat X_{+}$, respectively.
For convenience, we use $\mu>0$ and $M>0$ in the proof to denote
generic positive constants depending only on $\gamma$, $m$, $\lambda$,
$l_+$, $\delta_1$, $\delta_2$, $\delta_3$, $\|f_+^{(4)}\|_{L^\infty((0,l_+))}$, $M_1$ and $M_2$.
Set
$$
Q(\varphi,\psi)=\tilde Q(\varphi,\psi)-\hat Q(\varphi,\psi),
\quad(\varphi,\psi)\in[0,\zeta_{+}]\times[0,m]
$$
and
$$
X(\varphi)=\tilde X_{+}(\varphi)-\hat X_{+}(\varphi),
\quad(\varphi,\psi)\in[0,\zeta_{+}].
$$
Then $Q\in C^{0,1}([0,\zeta_+]\times[0,m])$ and $X\in C^{1}([0,\zeta_+])$ satisfy
$$
b(\tilde Q(\varphi,\psi))-b(\hat Q(\varphi,\psi))=h(\varphi,\psi)Q(\varphi,\psi),
\quad(\varphi,\psi)\in(0,\zeta_{+})\times(0,m),
$$
\begin{align}
\label{supun-6}
Q_{\psi}(\varphi,m)=e_1(\varphi)X(\varphi)+e_2(\varphi)Q(\varphi,m),
\quad\varphi\in(0,\zeta_{+})
\end{align}
and
\begin{align}
\label{supun-7}
X'(\varphi)=d_1(\varphi)X(\varphi)+d_2(\varphi)Q(\varphi,m),
\quad\varphi\in(0,\zeta_{+}),
\end{align}
where
$$
h(\varphi,\psi)=-\int_0^1K''(t\tilde Q(\varphi,\psi)+(1-t)\hat Q(\varphi,\psi))dt,
\quad(\varphi,\psi)\in(0,\zeta_{+})\times(0,m),
$$
$$
e_1(\varphi)=-\frac1{b(\tilde Q(\varphi,m)))A^{-1}_+(\tilde Q(\varphi,m))}
\int_0^1E_1(t\tilde X_{+}(\varphi)+(1-t)\hat X_{+}(\varphi))dt,
\quad\varphi\in(0,\zeta_{+}),
$$
$$
e_2(\varphi)=-\Theta'_{+}(\hat X_{+}(\varphi))\cos\Theta_{+}(\hat X_{+}(\varphi))
\int_0^1E_2(t \tilde Q(\varphi,m)+(1-t)\hat Q(\varphi,m))dt,
\quad\varphi\in(0,\zeta_{+}),
$$
$$
d_1(\varphi)=\frac1{A^{-1}_+(\tilde Q(\varphi,m))}
\int_0^1D_1(t\tilde X_{+}(\varphi)+(1-t)\hat X_{+}(\varphi))dt,
\quad\varphi\in(0,\zeta_{+}),
$$
$$
d_2(\varphi)=\cos\Theta_{+}(\hat X_{+}(\varphi))
\int_0^1D_2(t \tilde Q(\varphi,m)+(1-t)\hat Q(\varphi,m))dt,
\quad\varphi\in(0,\zeta_{+})
$$
with
$$
E_1(s)=\big(\Theta'_{+}(s)\cos\Theta_{+}(s)\big)',
\quad D_1(s)=\big(\cos\Theta_{+}(s)\big)',\quad0\le s\le l_+
$$
and
$$
E_2(s)=\Big(\frac1{b(s)A^{-1}_+(s)}\Big)',\quad
D_2(s)=\Big(\frac1{A^{-1}_+(s)}\Big)',\quad s<0.
$$
Making use of the conditions \eqref{a-a-0}, \eqref{rr-1} and $f_+\in C^{4}([0,l_+])$ on $f_+$
and the asymptotic behavior \eqref{theoremsuper-1} and \eqref{rr-theorem-1} on $\tilde Q$
and \eqref{supun-1} and \eqref{supun-2} on $\hat Q$ near $\varphi=0$,
one can show by direct calculations that for $(\varphi,\psi)\in(0,\zeta_{+})\times(0,m)$,
\begin{gather*}
|Q_\varphi(\varphi,\psi)|\le M\varphi^{\lambda+1},
\quad
|Q_\psi(\varphi,\psi)|\le M\varphi^{3\lambda/2+1},
\\
\mu\varphi^{-3\lambda/2-3}\le h(\varphi,\psi)\le M\varphi^{-3\lambda/2-3},
\quad
-M\varphi^{-3\lambda/2-4}\le h_\varphi(\varphi,\psi)\le-\mu\varphi^{-3\lambda/2-4},
\end{gather*}
and for $\varphi\in(0,\zeta_{+})$,
\begin{gather*}
|e_1(\varphi)|\le M\varphi^{3\lambda/2},
\quad |e_2(\varphi)|\le M\varphi^{\lambda/2-1},
\quad |e'_1(\varphi)|\le M\varphi^{\lambda/2+1},
\\
|e'_2(\varphi)|\le M\varphi^{\lambda/2-2},
\quad |d_1(\varphi)|\le M\varphi^{2\lambda+1},
\quad |d_2(\varphi)|\le M\varphi^{-\lambda/2-1},
\\
-M\varphi^{3\lambda/2+1}\le\tilde Q_\psi(\varphi,m)\le-\mu\varphi^{3\lambda/2+1},
\quad
|\tilde Q_{\varphi\psi}(\varphi,m)|\le M\varphi^{3\lambda/2}.
\end{gather*}
Here, the bound of $e'_1$ depends on $\|f_+^{(4)}\|_{L^\infty((0,l_+))}$
and it is noted that this is the only reason for the assumption that $f_+\in C^{4}([0,l_+])$.

First we estimate $X$.
It follows from \eqref{supun-7} and the regularity estimates of $d_1$ and $d_2$ that
\begin{align}
\label{supun-8}
|X(\phi)|\le M\int_0^\phi\varphi^{-\lambda/2-1}|Q(\varphi,m)|d\varphi,
\quad\phi\in(0,\zeta_{+}),
\end{align}
which implies
\begin{align}
\label{supun-9}
X^2(\phi)\le M\phi^2\int_0^\phi\varphi^{-\lambda-3}Q^2(\varphi,m)d\varphi,
\quad\phi\in(0,\zeta_{+})
\end{align}
and
\begin{align}
\label{supun-10}
|X'(\phi)|^2\le M\phi^{4\lambda+4}\int_0^\phi\varphi^{-\lambda-3}Q^2(\varphi,m)d\varphi
+M\phi^{-\lambda-2}Q^2(\phi,m),
\quad\phi\in(0,\zeta_{+}).
\end{align}
The H\"older inequality gives
\begin{align}
\label{supun-11}
\sup_{(0,m)}Q^2(\phi,\cdot)\le&\Big(\frac1m\int_0^m |Q(\phi,\psi)|d\psi+\int_0^m|Q_\psi(\phi,\psi)|d\psi\Big)^2
\nonumber
\\
\le&\Big(\frac1m\int_0^\phi\int_0^m|Q_\varphi(\varphi,\psi)|d\varphi d\psi+\int_0^m|Q_\psi(\phi,\psi)|d\psi\Big)^2
\nonumber
\\
\le&M\phi^{\lambda+3}\int_0^\phi\int_0^m\varphi^{-\lambda-2}Q_\varphi^2(\varphi,\psi)d\varphi d\psi
+M\int_0^mQ_\psi^2(\phi,\psi)d\psi,
\quad\phi\in(0,\zeta_{+}).
\end{align}
Then, one gets from \eqref{supun-9}--\eqref{supun-11} that
\begin{align}
\label{supun-12}
X^2(\phi)\le M\phi^3\int_0^\phi\int_0^m\varphi^{-\lambda-2}Q_\varphi^2(\varphi,\psi)d\varphi d\psi
+M\phi^{\lambda/2+2}\int_0^\phi\int_0^m\varphi^{-3\lambda/2-3}Q_\psi^2(\varphi,\psi)d\psi,
\quad\phi\in(0,\zeta_{+})
\end{align}
and
\begin{align}
\label{supun-13}
|X'(\phi)|^2\le& M\phi
\int_0^\phi\int_0^m\varphi^{-\lambda-2}Q_\varphi^2(\varphi,\psi)d\varphi d\psi
+M\phi^{9\lambda/2+4}\int_0^\phi\int_0^m\varphi^{-3\lambda/2-3}Q_\psi^2(\varphi,\psi)d\psi
\nonumber
\\
&\qquad+M\phi^{-\lambda-2}\int_0^mQ_\psi^2(\phi,\psi)d\psi,
\quad\phi\in(0,\zeta_{+}).
\end{align}

Next, we improve the estimate \eqref{supun-11} for $\|Q(\phi,\cdot)\|_{L^\infty((0,m))}$.
It follows from the definition of weak solutions that
\begin{align*}
\frac{d}{d\varphi}\int_0^m Q(\phi,\psi)d\psi=&
\int_0^\phi\Big(b(\tilde Q(\varphi,m))\pd{\tilde Q}{\psi}(\varphi,m)
-b(\hat Q(\varphi,m))\pd{\hat Q}{\psi}(\varphi,m)\Big)d\varphi
\\
=&\int_0^\phi(e_3(\varphi)X(\varphi)+e_4(\varphi)Q(\varphi,m))d\varphi,
\quad \phi\in(0,\zeta_{+}),
\end{align*}
where
\begin{gather*}
e_3(\varphi)=-\frac1{A^{-1}_+(\tilde Q(\varphi,m))}
\int_0^1E_1(t \tilde X_{+}(\varphi)+(1-t)\hat X_{+}(\varphi))dt,
\quad\varphi\in(0,\zeta_{+}),
\\
e_4(\varphi)=-\Theta'_{+}(\hat X_{+}(\varphi))\cos\Theta_{+}(\hat X_{+}(\varphi))
\int_0^1D_2(t \tilde Q(\varphi,m)+(1-t)\hat Q(\varphi,m))dt,
\quad\varphi\in(0,\zeta_{+})
\end{gather*}
with
$$
|e_3(\varphi)|\le M\varphi^{\lambda-1},
\quad |e_4(\varphi)|\le M\varphi^{\lambda/2-1},
\quad \varphi\in(0,\zeta_{+}).
$$
Therefore,
\begin{align*}
\Big|\int_0^m Q(\phi,\psi)d\psi\Big|\le&\phi\int_0^\phi|e_3(\varphi)X(\varphi)+e_4(\varphi)Q(\varphi,m)|d\varphi
\\
\le&M\phi^{\lambda+1/2}\Big(\int_0^\phi X^2(\varphi)d\varphi\Big)^{1/2}
+M\phi^{\lambda/2+1/2}\Big(\int_0^\phi Q^2(\varphi,m)d\varphi\Big)^{1/2},
\quad \phi\in(0,\zeta_{+}),
\end{align*}
which, together with \eqref{supun-12} and \eqref{supun-11}, implies
\begin{align*}
\Big(\int_0^m Q(\phi,\psi)d\psi\Big)^2
\le&M\phi^{2\lambda+5}\int_0^\phi\int_0^m\varphi^{-\lambda-2}Q_\varphi^2(\varphi,\psi)d\varphi d\psi
\\
&\qquad
+M\phi^{5\lambda/2+4}\int_0^\phi\int_0^m\varphi^{-3\lambda/2-3}Q_\psi^2(\varphi,\psi)d\psi,
\quad \phi\in(0,\zeta_{+}).
\end{align*}
Thus
\begin{align}
\label{supun-14}
\sup_{(0,m)}Q^2(\phi,\cdot)\le&\Big(\frac1m\Big|\int_0^m Q(\phi,\psi)d\psi\Big|+\int_0^m|Q_\psi(\phi,\psi)|d\psi\Big)^2
\nonumber
\\
\le&M\phi^{2\lambda+5}\int_0^\phi\int_0^m\varphi^{-\lambda-2}Q_\varphi^2(\varphi,\psi)d\varphi d\psi
+M\phi^{5\lambda/2+4}\int_0^\phi\int_0^m\varphi^{-3\lambda/2-3}Q_\psi^2(\varphi,\psi)d\psi
\nonumber
\\
&\qquad+M\int_0^mQ_\psi^2(\phi,\psi)d\psi,
\quad\phi\in(0,\zeta_{+}).
\end{align}
Moreover, it follows from \eqref{supun-8} and \eqref{supun-14} that
\begin{align}
\label{supun-888}
X^2(\phi)\le& M\phi^{\lambda/2+2}\int_0^\phi\varphi^{-3\lambda/2-3}Q^2(\varphi,m)d\varphi
\nonumber
\\
\le&M\phi^{\lambda+5}\int_0^\phi\int_0^m\varphi^{-\lambda-2}Q_\varphi^2(\varphi,\psi)d\varphi d\psi
+M\phi^{\lambda/2+2}\int_0^\phi\int_0^m\varphi^{-3\lambda/2-3}Q_\psi^2(\varphi,\psi)d\psi,
\quad\phi\in(0,\zeta_{+}),
\end{align}
which improves the estimate \eqref{supun-12}.

Now we are ready to prove the theorem by a weighted energy estimate.
Let $\phi\in(0,\zeta_{+}]$ be given.
Multiplying the equation of $\tilde Q$ by $\varphi^{-\lambda-1}\tilde Q_\varphi$
and $-\varphi^{-\lambda-1}\hat Q_\varphi$, respectively, and then integrating
over $(0,\phi)\times(0,m)$, one can get that
\begin{align}
\label{superunique-1}
0=&\int_0^\phi\int_0^m\varphi^{-\lambda-1}\tilde Q_{\varphi\varphi}
\tilde Q_\varphi d\varphi d\psi
-\int_0^\phi\int_0^m\varphi^{-\lambda-1}(b(\tilde Q)\tilde Q_{\psi})_\psi
\tilde Q_\varphi d\varphi d\psi
\nonumber
\\
=&\frac{1}{2}\int_0^m\varphi^{-\lambda-1}\tilde Q_\varphi^2d\psi\Big|_{\varphi=\phi}
+\frac{\lambda+1}{2}\int_0^\phi\int_0^m\varphi^{-\lambda-2}\tilde Q_\varphi^2d\varphi d\psi
-\int_0^\phi\varphi^{-\lambda-1}b(\tilde Q)
\tilde Q_\varphi\tilde Q_{\psi}d\varphi\Big|_{\psi=m}
\nonumber
\\
&\qquad
+\frac{1}2\int_0^m\varphi^{-\lambda-1}b(\tilde Q)\tilde Q^2_{\psi}d\psi\Big|_{\varphi=\phi}
-\frac12\int_0^\phi\int_0^m(\varphi^{-\lambda-1}b(\tilde Q))_{\varphi}
\tilde Q^2_{\psi}d\varphi d\psi
\end{align}
and
\begin{align}
\label{superunique-2}
0=&-\int_0^\phi\int_0^m\varphi^{-\lambda-1}\tilde Q_{\varphi\varphi}
\hat Q_\varphi d\varphi d\psi
+\int_0^\phi\int_0^m\varphi^{-\lambda-1}(b(\tilde Q)\tilde Q_\psi)_{\psi}
\hat Q_{\varphi}d\varphi d\psi.
\end{align}
It follows from the definition of weak solutions and a standard limit process that
\begin{align}
\label{superunique-3}
0=&\frac{1}{2}\int_0^m\varphi^{-\lambda-1}\hat Q_\varphi^2d\psi\Big|_{\varphi=\phi}
+\frac{\lambda+1}{2}\int_0^\phi\int_0^m\varphi^{-\lambda-2}\hat Q_\varphi^2d\varphi d\psi
-\int_0^\phi\varphi^{-\lambda-1}b(\hat Q)
\hat Q_\varphi\hat Q_{\psi}d\varphi\Big|_{\psi=m}
\nonumber
\\
&\qquad
+\frac{1}2\int_0^m\varphi^{-\lambda-1}b(\hat Q)\hat Q^2_{\psi}d\psi\Big|_{\varphi=\phi}
-\frac12\int_0^\phi\int_0^m(\varphi^{-\lambda-1}b(\hat Q))_{\varphi}
\hat Q^2_{\psi}d\varphi d\psi.
\end{align}
and
\begin{align}
\label{superunique-4}
0=&-\int_0^m \varphi^{-\lambda-1}\tilde Q_\varphi\hat Q_\varphi d\psi\Big|_{\varphi=\phi}
-(\lambda+1)\int_0^\phi\int_0^m\varphi^{-\lambda-2}
\tilde Q_\varphi\hat Q_\varphi d\varphi d\psi
+\int_0^\phi\int_0^m\varphi^{-\lambda-1}\tilde Q_{\varphi\varphi}
\hat Q_\varphi d\varphi d\psi
\nonumber
\\
&\qquad
+\int_0^\phi\varphi^{-\lambda-1}b(\hat Q)
\tilde Q_\varphi\hat Q_{\psi}d\varphi\Big|_{\psi=m}
-\int_0^\phi\int_0^m\varphi^{-\lambda-1}b(\hat Q)
\tilde Q_{\varphi\psi}\hat Q_{\psi}d\varphi d\psi.
\end{align}
A direct calculation and a standard limit process show that
\begin{align}
\label{superunique-5}
&\int_0^\phi\int_0^m\varphi^{-\lambda-1}(b(\tilde Q)\tilde Q_\psi)_{\psi}
\hat Q_{\varphi}d\varphi d\psi
\nonumber
\\
=&\int_0^\phi\varphi^{-\lambda-1}b(\tilde Q)\tilde Q_\psi
\hat Q_{\varphi}d\varphi\Big|_{\psi=m}
-\int_0^m \varphi^{-\lambda-1}b(\tilde Q)\tilde Q_\psi
\hat Q_{\psi} d\psi\Big|_{\varphi=\phi}
\nonumber
\\
&\qquad
+\int_0^\phi\int_0^m(\varphi^{-\lambda-1}b(\tilde Q))_\varphi
\tilde Q_\psi\hat Q_{\psi} d\varphi d\psi
+\int_0^\phi\int_0^m\varphi^{-\lambda-1}b(\tilde Q)
\tilde Q_{\varphi\psi}\hat Q_{\psi}d\varphi d\psi.
\end{align}
Integrating by parts leads to
\begin{align}
\label{superunique-6}
&\int_0^\phi\int_0^m\varphi^{-\lambda-1}hQ\tilde Q_{\varphi\psi}
\tilde Q_{\psi}d\varphi d\psi
\nonumber
\\
=&\frac12\int_0^m \varphi^{-\lambda-1}hQ\tilde Q^2_\psi d\psi\Big|_{\varphi=\phi}
-\frac12\int_0^\phi\int_0^m(\varphi^{-\lambda-1}hQ)_\varphi
\tilde Q^2_\psi d\varphi d\psi.
\end{align}
Summing up from \eqref{superunique-1} to \eqref{superunique-6} yields
\begin{align*}
0=I_1+I_2+J_1+J_2,
\end{align*}
where
\begin{align*}
I_1=&\frac{\lambda+1}{2}\int_0^\phi\int_0^m\varphi^{-\lambda-2}
\Big(\tilde Q_\varphi^2+\hat Q_\varphi^2-2\tilde Q_\varphi\hat Q_\varphi\Big)d\varphi d\psi
\\
=&\frac{\lambda+1}{2}\int_0^\phi\int_0^m\varphi^{-\lambda-2}Q_\varphi^2d\varphi d\psi,
\end{align*}
\begin{align*}
I_2=&-\frac12\int_0^\phi\int_0^m\Big((\varphi^{-\lambda-1}b(\tilde Q))_{\varphi}
\tilde Q^2_{\psi}+(\varphi^{-\lambda-1}b(\hat Q))_{\varphi}
\hat Q^2_{\psi}-2(\varphi^{-\lambda-1}b(\tilde Q))_\varphi
\tilde Q_\psi\hat Q_{\psi}\Big) d\varphi d\psi
\\
&\qquad
-\frac12\int_0^\phi\int_0^m(\varphi^{-\lambda-1}hQ)_\varphi
\tilde Q^2_\psi d\varphi d\psi
\\
&\qquad
-\int_0^\phi\int_0^m\varphi^{-\lambda-1}\Big(b(\hat Q)
\tilde Q_{\varphi\psi}\hat Q_{\psi}-b(\tilde Q)
\tilde Q_{\varphi\psi}\hat Q_{\psi}+hQ\tilde Q_{\varphi\psi}
\tilde Q_{\psi}\Big)d\varphi d\psi
\\
=&-\frac12\int_0^\phi\int_0^m(\varphi^{-\lambda-1}b(\tilde Q))_{\varphi}Q^2_{\psi}d\varphi d\psi
-\frac12\int_0^\phi\int_0^m\varphi^{-\lambda-1}h(\tilde Q_{\psi}+\hat Q_{\psi})
Q_\varphi Q_{\psi}d\varphi d\psi
\\
&\qquad
-\frac12\int_0^\phi\int_0^m(\varphi^{-\lambda-1}h)_{\varphi}(\tilde Q_{\psi}+\hat Q_{\psi})
Q Q_{\psi}d\varphi d\psi
-\int_0^\phi\int_0^m\varphi^{-\lambda-1}h\tilde Q_{\varphi\psi}QQ_{\psi}d\varphi d\psi,
\end{align*}
\begin{align*}
J_1=&\frac{1}{2}\int_0^m\varphi^{-\lambda-1}
\Big(\tilde Q_\varphi^2+\hat Q_\varphi^2-2\tilde Q_\varphi\hat Q_\varphi\Big)d\psi\Big|_{\varphi=\phi}
\\
&\qquad
+\frac{1}{2}\int_0^m\varphi^{-\lambda-1}
\Big(b(\tilde Q)\tilde Q^2_{\psi}-2b(\tilde Q)\tilde Q_{\psi}\hat Q_{\psi}
+b(\hat Q)\hat Q^2_{\psi}+hQ\tilde Q^2_\psi\Big)d\psi\Big|_{\varphi=\phi}
\\
=&\frac{1}{2}\int_0^m\varphi^{-\lambda-1}Q_\varphi^2d\psi\Big|_{\varphi=\phi}
+\frac{1}{2}\int_0^m\varphi^{-\lambda-1}b(\tilde Q)Q^2_{\psi}d\psi\Big|_{\varphi=\phi}
+\frac{1}{2}\int_0^m\varphi^{-\lambda-1}
h(\tilde Q_{\psi}+\hat Q_{\psi})QQ_\psi d\psi\Big|_{\varphi=\phi}
\end{align*}
and
\begin{align*}
J_2=&-\int_0^\phi\varphi^{-\lambda-1}
\Big(b(\tilde Q)\tilde Q_\varphi\tilde Q_{\psi}
-b(\tilde Q)\tilde Q_\psi\hat Q_{\varphi}
+b(\hat Q)\hat Q_\varphi\hat Q_{\psi}
-b(\hat Q)\tilde Q_\varphi\hat Q_{\psi}\Big)d\varphi\Big|_{\psi=m}
\\
=&-\int_0^\phi\varphi^{-\lambda-1}
\Big(b(\tilde Q)\tilde Q_\psi Q_{\varphi}
-b(\hat Q)\hat Q_\psi Q_\varphi\Big)d\varphi\Big|_{\psi=m}
\\
=&-\int_0^\phi\varphi^{-\lambda-1}
h\tilde Q_\psi QQ_{\varphi}d\varphi\Big|_{\psi=m}
-\int_0^\phi\varphi^{-\lambda-1}
b(\hat Q)Q_{\varphi}Q_\psi d\varphi\Big|_{\psi=m}
\\
=&-\frac12\varphi^{-\lambda-1}h\tilde Q_\psi Q^2\Big|_{(\varphi,\psi)=(\phi,m)}
+\frac12\int_0^\phi(\varphi^{-\lambda-1}h\tilde Q_\psi)_{\varphi}Q^2d\varphi\Big|_{\psi=m}
-\int_0^\phi\varphi^{-\lambda-1}
b(\hat Q)Q_{\varphi}Q_\psi d\varphi\Big|_{\psi=m}.
\end{align*}
Therefore,
\begin{align*}
&\frac{\lambda+1}{2}\int_0^\phi\int_0^m\varphi^{-\lambda-2}Q_\varphi^2d\varphi d\psi
-\frac12\int_0^\phi\int_0^m(\varphi^{-\lambda-1}b(\tilde Q))_{\varphi}Q^2_{\psi}d\varphi d\psi
\\
&\qquad
+\frac{1}{2}\int_0^m\varphi^{-\lambda-1}Q_\varphi^2d\psi\Big|_{\varphi=\phi}
+\frac{1}{2}\int_0^m\varphi^{-\lambda-1}b(\tilde Q)Q^2_{\psi}d\psi\Big|_{\varphi=\phi}
-\frac12\varphi^{-\lambda-1}h\tilde Q_\psi Q^2\Big|_{(\varphi,\psi)=(\phi,m)}
\\
=&\frac12\int_0^\phi\int_0^m\varphi^{-\lambda-1}h(\tilde Q_{\psi}+\hat Q_{\psi})
Q_\varphi Q_{\psi}d\varphi d\psi
+\frac12\int_0^\phi\int_0^m(\varphi^{-\lambda-1}h)_{\varphi}(\tilde Q_{\psi}+\hat Q_{\psi})
Q Q_{\psi}d\varphi d\psi
\\
&\qquad
+\int_0^\phi\int_0^m\varphi^{-\lambda-1}h\tilde Q_{\varphi\psi}QQ_{\psi}d\varphi d\psi
-\frac{1}{2}\int_0^m\varphi^{-\lambda-1}h(\tilde Q_{\psi}+\hat Q_{\psi})QQ_\psi d\psi\Big|_{\varphi=\phi}
\\
&\qquad
-\frac12\int_0^\phi(\varphi^{-\lambda-1}h\tilde Q_\psi)_{\varphi}Q^2d\varphi\Big|_{\psi=m}
+\int_0^\phi\varphi^{-\lambda-1}
b(\hat Q)Q_{\varphi}Q_\psi d\varphi\Big|_{\psi=m}.
\end{align*}
It follows from the regularity estimates of $\tilde Q$, $\hat Q$ and $h$,
the asymptotic behavior of $\tilde Q_\psi$ and $\tilde Q_{\varphi\psi}$ on the wall near $\varphi=0$
and the H\"older inequality that
\begin{align}
\label{superunique-7}
&\int_0^\phi\int_0^m\Big(\varphi^{-\lambda-2}Q_\varphi^2+\varphi^{-3\lambda/2-3}Q^2_{\psi}\Big)d\varphi d\psi
+\int_0^m\Big(\varphi^{-\lambda-1}Q_\varphi^2+\varphi^{-3\lambda/2-2}Q^2_{\psi}\Big)d\psi\Big|_{\varphi=\phi}
\nonumber
\\
\le&
M\int_0^\phi\int_0^m\varphi^{-\lambda-3}|Q_\varphi Q_{\psi}|d\varphi d\psi
+M\int_0^\phi\int_0^m\varphi^{-\lambda-4}|Q Q_{\psi}|d\varphi d\psi
+M\int_0^\phi\int_0^m\varphi^{-5\lambda/4-7/2}|Q Q_{\psi}|d\varphi d\psi
\nonumber
\\
&\qquad
+M\int_0^m\varphi^{-\lambda-3}|QQ_\psi| d\psi\Big|_{\varphi=\phi}
+M\int_0^\phi\varphi^{-\lambda-4}Q^2 d\varphi\Big|_{\psi=m}
+\int_0^\phi\varphi^{-\lambda-1}
b(\hat Q)Q_{\varphi}Q_\psi d\varphi\Big|_{\psi=m}
\nonumber
\\
\le&
M\int_0^\phi\int_0^m\varphi^{-3\lambda/4-5/2}Q_\varphi^2 d\varphi d\psi
+M\int_0^\phi\int_0^m\varphi^{-5\lambda/4-7/2}Q_{\psi}^2 d\varphi d\psi
+M\int_0^\phi\int_0^m\varphi^{-5\lambda/4-7/2}Q^2 d\varphi d\psi
\nonumber
\\
&\qquad
+M\int_0^m\varphi^{-\lambda-3}Q^2 d\psi\Big|_{\varphi=\phi}
+M\int_0^m\varphi^{-\lambda-3}Q_\psi^2 d\psi\Big|_{\varphi=\phi}
\nonumber
\\
&\qquad
+M\int_0^\phi\varphi^{-\lambda-4}Q^2 d\varphi\Big|_{\psi=m}
+\int_0^\phi\varphi^{-\lambda-1}b(\hat Q)Q_{\varphi}Q_\psi d\varphi\Big|_{\psi=m}.
\end{align}
To estimate the last term on the right side of \eqref{superunique-7},
one substitutes \eqref{supun-6} into this term,
then integrates by parts and uses the H\"older inequality
and the regularity estimates of $\hat Q$, $e_1$ and $e_2$ on the wall near $\varphi=0$
to get
\begin{align*}
&\int_0^\phi\varphi^{-\lambda-1}b(\hat Q)Q_{\varphi}Q_\psi d\varphi\Big|_{\psi=m}
\\
=&\int_0^\phi\varphi^{-\lambda-1}b(\hat Q(\varphi,m))(e_1(\varphi)X(\varphi)+e_2(\varphi)Q(\varphi,m))
Q_\varphi(\varphi,m)d\varphi
\\
=&\phi^{-\lambda-1}b(\hat Q(\phi,m))e_1(\phi)X(\phi)Q(\phi,m)
+\frac12\phi^{-\lambda-1}b(\hat Q(\phi,m))e_2(\phi)Q^2(\phi,m)
\\
&\qquad
-\int_0^\phi\big(\varphi^{-\lambda-1}b(\hat Q(\varphi,m))e_1(\varphi))'X(\varphi)Q(\varphi,m)d\varphi
-\int_0^\phi\varphi^{-\lambda-1}b(\hat Q(\varphi,m))e_1(\varphi)X'(\varphi)Q(\varphi,m)d\varphi
\\
&\qquad
-\frac12\int_0^\phi\big(\varphi^{-\lambda-1}b(\hat Q(\varphi,m))e_2(\varphi)\big)'Q^2(\varphi,m)d\varphi
\\
\le&M\phi^{-2}|X(\phi)Q(\phi,m)|+M\phi^{-\lambda-3}Q^2(\phi,m)
+M\int_0^\phi\varphi^{-\lambda-1}|X(\varphi)Q(\varphi,m)|d\varphi
\\
&\qquad+M\int_0^\phi\varphi^{-2}|X'(\varphi)Q(\varphi,m)|d\varphi
+M\int_0^\phi\varphi^{-\lambda-4}Q^2(\varphi,m)d\varphi
\\
\le&M\phi^{\lambda-1}X^2(\phi)+M\phi^{-\lambda-3}Q^2(\phi,m)
+M\int_0^\phi\varphi^{-\lambda/2-1}X^2(\varphi)d\varphi
+M\int_0^\phi\varphi^{-3\lambda/2-1}Q^2(\varphi,m)d\varphi
\\
&\qquad+M\int_0^\phi\varphi^{\lambda}|X'(\varphi)|^2d\varphi
+M\int_0^\phi\varphi^{-\lambda-4}Q^2(\varphi,m)d\varphi,
\end{align*}
which, together with \eqref{supun-13}--\eqref{supun-888}, yields
\begin{align}
\label{supun-16}
&\int_0^\phi\varphi^{-\lambda-1}b(\hat Q)Q_{\varphi}Q_\psi d\varphi\Big|_{\psi=m}
\nonumber
\\
\le&M(\phi^{\lambda+2}+\phi^{\lambda/2+5})\int_0^\phi\int_0^m\varphi^{-\lambda-2}Q_\varphi^2d\varphi d\psi
+M(\phi^{\lambda/2-1}+\phi^2)
\int_0^\phi\int_0^m\varphi^{-3\lambda/2-3}Q_\psi^2d\varphi d\psi
\nonumber
\\
&\qquad
+M\phi^{\lambda/2-1}\int_0^m\varphi^{-3\lambda/2-2} Q_\psi^2d\psi\Big|_{\varphi=\phi}.
\end{align}
Substituting \eqref{supun-16} into \eqref{superunique-7}
and using \eqref{supun-14}, we obtain
\begin{align*}
&\int_0^\phi\int_0^m\Big(\varphi^{-\lambda-2}Q_\varphi^2+\varphi^{-3\lambda/2-3}Q^2_{\psi}\Big)d\varphi d\psi
+\int_0^m\Big(\varphi^{-\lambda-1}Q_\varphi^2+\varphi^{-3\lambda/2-2}Q^2_{\psi}\Big)d\psi\Big|_{\varphi=\phi}
\\
\le&M\phi^{\lambda/4-1/2}
\int_0^\phi\int_0^m\varphi^{-\lambda-2}Q_\varphi^2 d\varphi d\psi
+M(\phi^{\lambda/4-1/2}+\phi^2)
\int_0^\phi\int_0^m\varphi^{-3\lambda/2-3}Q_\psi^2d\varphi d\psi
\\
&\qquad
+M\phi^{\lambda/2-1}\int_0^m\varphi^{-3\lambda/2-2} Q_\psi^2d\psi\Big|_{\varphi=\phi}.
\end{align*}
Owing to $\lambda>2$, one can conclude that
for sufficiently small $\phi\in(0,\zeta_{+}]$,
$$
Q_\varphi(\varphi,\psi)=0,\quad
Q_\psi(\varphi,\psi)=0,\quad(\varphi,\psi)\in(0,\phi)\times(0,m).
$$
Therefore,
$$
\tilde Q(\varphi,\psi)=\hat Q(\varphi,\psi),\quad(\varphi,\psi)\in(0,\phi)\times(0,m).
$$
$\hfill\Box$
\end{document}